Printed in the United Kingdom

# Approximation of Functions: Optimal Sampling and Complexity

David Krieg
*Faculty of Computer Science and Mathematics,*
*University of Passau, 94032 Passau, Germany*
*E-mail: david.krieg@uni-passau.de*

Mario Ullrich
*Department of Quantum Computing & Institute of Analysis,*
*Johannes Kepler University, 4040 Linz, Austria*
*E-mail: mario.ullrich@jku.at*

We consider approximation or recovery of functions based on a finite number of function evaluations. This is a well-studied problem in optimal recovery, machine learning, and numerical analysis in general, but many fundamental insights were obtained only recently. We discuss different aspects of the information-theoretic limit that appears because of the limited amount of data available, as well as algorithms and sampling strategies that come as close to it as possible.
We also discuss (optimal) sampling in a broader sense, allowing other types of measurements that may be nonlinear, adaptive and random, and present several relations between the different settings in the spirit of information-based complexity.
We hope that this article provides both, a basic introduction to the subject and a contemporary summary of the current state of research.

## 1. Introduction and overview

The numerical solution of an approximation problem can often be divided into two parts: data acquisition and reconstruction from this data. By sampling, we understand the (possibly random) process of data acquisition. The data might be acquired by means of physical measurements, by running computer simulations, or even via questionnaires. In any case, each such "measurement" comes at a certain cost and we can only acquire a limited finite amount of data. As such, the information about the problem is usually incomplete and this constitutes a fundamental limit of what can be achieved by any reconstruction based upon it.

We would like to understand whether this limitation on available measurements represents a serious drawback compared to approximations that use unrealistic-to-obtain or even complete information about the problem instance (like best approximations from a linear space, a manifold, etc.), but which constitute well-studied benchmark quantities in approximation theory and numerical analysis. This includes, on the one hand, the construction and analysis of sampling schemes and algorithms which come close to these benchmarks. On the other hand, it includes the study of inherent lower bounds on the amount of data necessary for specific model classes, i.e., the information complexity, and its dependence on the class of admissible measurements.

We try to give a systematic survey of this branch of modern computational mathematics which takes into account the admissible information. From our perspective, we operate at the intersection of the "classical" theories of numerical analysis, approximation theory and theoretical data science, and strive to establish connections between them.

Throughout this article, we assume that there is a *model class* $F$, typically functions $f\colon D \to \mathbb{R}$ on a set $D$, which describes the possible problem elements (or inputs) that we aim to approximate. The approximation does not have to be an element of the model class $F$, it can be in a larger space $Y$ (like $Y = L_2$), which we assume to be equipped with a (semi-)norm or metric in order to have a notion of error. We understand by an *algorithm* a mapping $A_n\colon F \to Y$ of the form

$$A_n = \Phi \circ N_n \quad \text{with} \quad N_n\colon F \to \mathbb{R}^n \qquad \text{and} \qquad \Phi\colon \mathbb{R}^n \to Y,$$

where $N_n$ is called the *measurement* or *information map* (aka *encoder*) and $\Phi$ is the *reconstruction map* (aka *decoder*). That is, we take $n$ real (or complex) measurements and use an, in principle arbitrary, mapping to reconstruct the target function $f$ using only these values and the a priori information that $f \in F$. Most of this article is devoted to *sampling algorithms* of the form

$$A_n(f) = \Phi\big(f(x_1), \ldots, f(x_n)\big)\,,$$

often with an explicit form of $\Phi$. Function values often serve as a *model* for acquired data, even if it is not known whether the data can be explained by an underlying function. In machine learning, e.g., it is typically assumed that samples are (noisy) evaluations of an unknown target function, which very much fits our setting. In general, function evaluations are usually among the more realistic classes of measurements and it is hence interesting to study the intrinsic power and limits of this type of algorithms.

We will discuss sampling algorithms in various settings and provide several theoretical error bounds, but also some (semi-)constructive sampling strategies for the points $x_1, \ldots, x_n \in D$. To describe the best possible approximation by sampling algorithms, it is often helpful to (re)state our results in terms of the *sampling numbers*

$$g_n(F, Y) := \inf_{\substack{x_1,\ldots,x_n \in D \\ \Phi\colon \mathbb{C}^n \to Y}} \sup_{f \in F} \Big\| f - \Phi\big(f(x_1), \ldots, f(x_n)\big) \Big\|_Y. \tag{1.1}$$

This quantity describes the *minimal worst-case error* achievable with sampling algorithms. It is the smallest possible number $\varepsilon \geq 0$ for which there exists a sampling algorithm that takes at most $n$ samples and whose $Y$-approximation error for any function from $F$ is at most $\varepsilon$. The most typical choices for $Y$ are the $L_p$-spaces or the space of bounded functions $B(D)$ with the sup-norm. See Remark 2.6 for a detailed comment on our motivation to study the *worst-case setting*.

The sampling numbers are well-studied for many (classical) classes of functions. Over the years, it has been understood that there is often no great loss of practical sampling algorithms compared to purely theoretical approximation benchmarks. In fact, for several of these model classes, the corresponding numerical and theoretical findings are probably a reason why they became important in the first place. In recent years, considerable progress has been made in explaining these phenomena from a broader perspective, i.e., not only by studying particular model classes, but rather by identifying structural properties of model classes that generally make it possible to find good sampling algorithms. Hopefully, the established general error bounds, algorithms and sampling strategies ultimately lead to new practicable approaches to the solution of 'difficult' problems. In any case, they provide insight into the capabilities and limitations of (optimal) sampling algorithms.

Here, we try to survey the state of the art in this area of research. We will focus on results that are general in the sense that they can be applied for a large variety of model classes. We do not want to assume any structure of the domain $D$ and we try to formulate minimal conditions on the model class $F$ required for the specific results. In particular, it is not the aim of this work to present a collection of results for specific model classes. While we will also discuss several examples along the way, this will mainly be for illustration and in order to discuss the power and weaknesses of the general theory.

The different sections of this article highlight different aspects of (optimal) approximation and sampling. We have tried to make these sections as self-contained as possible to allow selective reading. However, we also hope that it will prove suitable as (introductory) reading material in its entirety.

As a starter, we begin in Section 2 with a detailed discussion of the *easy example* of $L_p$-approximation of univariate Lipschitz functions. This classical example allows for a simple (graphical) explanation of the basic concepts, including the concepts of optimality and complexity as used in this article, and we will return to it repeatedly to illustrate new aspects.

In Section 3, we discuss to what extent optimal sampling algorithms can catch up with best approximation from a linear space. This section is solely based on particular *least squares algorithms*, and we show in Section 3.1 how error bounds for them relate to suitable *discretization inequalities*. We first focus on the existence of corresponding sampling points (and weights) and, in particular, present bounds by best approximation in the uniform norm (Section 3.2) and in $L_2$-norm (Section 3.3). A constructive approach to an implementation of these algorithms is discussed in Section 4.

In Section 5, we turn back to upper bounds for the sampling numbers by other theoretical approximation benchmarks. This time the emphasis is on the approximation on nonlinear spaces. First, and based on methods from *compressive sensing*, we compare the power of sampling algorithms in Section 5.1 to best sparse approximation from a given dictionary, which is still partially constructive. Sections 5.2 and 5.3 are concerned with bounds by the *(metric) entropy numbers* and the *Hilbert numbers*, respectively. The latter is a very small approximation benchmark that even lower bounds the error of best continuous approximation from a manifold. Unfortunately, the bounds from Sections 5.2 and 5.3 are not constructive and we do not know of an explicit algorithm that achieves these bounds.

By the methods of Sections 3 and 5, we obtain several *sharp* bounds on the sampling numbers, where by sharp we mean that there is some example that shows that a general bound cannot be (substantially) improved. One may summarize these results about the sampling numbers by six main results, which are collected in Appendix A.1, together with the (known) bounds for our frequently used and illustrative example of univariate Sobolev functions.

Some of the previously mentioned theoretical results are proven with the help of the *probabilistic method*. In fact, error bounds for sampling algorithms are often equivalent to or implied by certain properties of (possibly infinitely large) matrices, which can be handled by contemporary methods from random matrix theory. We discuss that in detail in Section 6, where we consider the setting that the sampling points are given as a realization of independent and identically distributed (i.i.d.) random variables. Still, we consider deterministic algorithms based on these sampling points and ask whether they can be 'good' in a worst-case sense, ideally with high probability. Again, we illustrate the setting with our easy example in Section 6.1, before discussing some general results in Section 6.2. In Section 6.3, we show how the general results can be improved for some important special classes, namely, for Sobolev spaces on general domains.

*Randomized sampling algorithms*, i.e., sampling algorithms that are allowed to use randomness tailored to the particular problem and for which new random sampling points can be drawn for each input, will be considered in Section 7. See again Section 7.1 for an illustration of the important concepts. Even more general bounds are available in this setting, where the demand of a guaranteed error bound is relaxed to an error bound with good probability or in expectation; see Section 7.2 where we discuss the *randomized sampling numbers*. Aspects of explicit constructions, particularly based on least squares or a *multilevel Monte Carlo method*, are discussed in Sections 7.3 and 7.4.

In Section 8, we shortly discuss issues related to *high-dimensional approximation*, meaning that we explicitly track the dependence of the minimal error (and complexity) on the number of input parameters. We will show in Section 8.1, again with an example close to our easy example of Lipschitz functions, that certain problems are practically impossible to solve in high dimensions. This phenomenon is usually called the *curse of dimensionality*. In Sections 8.2 and 8.3, we discuss how the general results apply in this context (leading to several sufficient conditions for tractability, and hence the absence of the curse of dimensionality) and present a collection of structural assumptions that allow for tractable problems, respectively.

In Section 9, we turn to a more general concept of sampling which is not focused solely on the acquisition of function values. We hence take a more general perspective on the problem of optimal sampling and introduce a framework that includes many types of approximation benchmarks and allows for their comparison. This includes a brief description of the theory of optimal algorithms (aka *information-based complexity*) and its most fundamental results in Section 9.1. In particular, we will introduce the concept of *adaptive algorithms* and discuss their advantage in comparison with non-adaptive, or even linear, algorithms in the context of (deterministic) linear measurements. We then make a detour and introduce *s-numbers* of operators (a generalization of singular values of matrices) in Section 9.2, which provides some powerful techniques for our purposes. Most notably, there exist bounds between all s-numbers. The relation to optimal algorithms and the implied bounds for minimal errors are discussed in Section 9.3.

In Section 10, we give a systematic comparison of minimal errors corresponding to algorithms based on different classes of admissible measurements. In Section 10.1, we compare (linear) sampling algorithms with general adaptive and/or randomized algorithms that have access to arbitrary linear measurements. In Section 10.2, we discuss the individual power of adaption and randomization for linear measurements in more detail. Finally, the class of continuous measurements and the power of adaption in this setting is outlined in Section 10.3.

**Notation.** For convenience, we collect some of the frequently used notation.

For two sequences, $(e_n)_{n\geq 0}$ and $(g_n)_{n\geq 0}$, we write $e_n \lesssim g_n$ or also $g_n \gtrsim e_n$ if $e_n \leq C\, g_n$ for a constant $C > 0$ and all $n \geq 2$, and $e_n \asymp g_n$ if $e_n \lesssim g_n$ and $g_n \lesssim e_n$. We use the notation $e_n \approx g_n$ to say that $\lim_{n\to\infty} e_n/g_n = 1$.

In what follows, $D$ may be any set. It will be the domain of the functions that we approximate. For example, it could be a bounded or unbounded subset of $\mathbb{R}^d$, but we restrict that only if needed.

For sets $D$ and $E$, we denote by $E^D$ the set of all functions $f\colon D \to E$. We use it most often for $E = \mathbb{C}$, in which case $E^D = \mathbb{C}^D$ is a $\mathbb{C}$-linear space. The space of bounded functions is denoted by $B(D) := \{f \in \mathbb{C}^D\colon \|f\|_\infty < \infty\}$ with the sup-norm $\|f\|_\infty := \sup_{x\in D} |f(x)|$, and we often abbreviate $B := B(D)$. For a measure space $(D, \mathcal{D}, \mu)$, we write $L_p = L_p(D, \mu)$, $1 \le p < \infty$ for the space of $p$-integrable functions on $D$, and $L_\infty = L_\infty(D, \mu)$ for the space of essentially bounded functions, both with the usual seminorms. We generally suppress the specification of the $\sigma$-algebra (since $\mu$ contains it as its domain anyways), and mostly consider elements of $L_p$ as functions, and not equivalence classes. Hence, we often speak of the '$L_p$-seminorm'. We write $X \hookrightarrow Y$ for two metric spaces $X, Y$ and say that $X$ *is embedded into* $Y$, if $X \subset Y$ and the identity (or inclusion) $\mathrm{id}\colon X \to Y$, $\mathrm{id}(f) = f$, is a continuous injection.

## 2. An easy example and some basic concepts

A main objective of this article is to study the optimality of sampling schemes and algorithms for the problem of approximating a function $f$. In order to discuss optimality, it is essential to fix a class $F$ of possible inputs $f \in F$. The class $F$ can be any class of functions that satisfy a given a priori knowledge. The emphasis is on general results say, for general subsets of Hilbert or Banach function spaces. However, to get a better idea of our concept(s) of optimality and complexity, we illustrate them by means of an easy example, namely, the problem of $L_p$-approximation of univariate Lipschitz functions. In this case, all definitions are easy to understand and the optimal algorithm is well known, and the proofs can almost be done by drawings. This allows for a simple explanation of the theoretical concepts, and is also helpful for highlighting some issues and misconceptions about optimality and randomness as used in this survey.

Of course, there are also other motivations for studying (toy) examples such as this one. Examples are needed to show that some general result is already *best possible*, they can disprove conjectures, and they can point us to the properties which are needed for a certain desirable statement. We will discuss more examples 'on the fly' in order to illustrate the power and limitations of the general theory of the upcoming sections.

The example we would like to discuss is the problem of $L_p$-approximation, $1 \le p \le \infty$, of functions from the class

$$F_{\mathrm{Lip}} := \Big\{ f\colon [0,1] \to \mathbb{R}\colon\ |f(x) - f(y)| \le \operatorname{dist}(x,y) \ \text{ for all } x, y \in [0,1] \Big\} \tag{2.1}$$

consisting of Lipschitz-continuous functions. To make the later computations as simple as possible, we consider these functions on the circle, aka 1-dimensional torus, and dist denotes the corresponding distance, i.e.,

$$\operatorname{dist}(x,y) = \min_{k \in \mathbb{Z}} |x + k - y| \quad \text{for} \quad x, y \in [0,1].$$

The results below would be essentially the same in the non-periodic case. As usual, we denote by $L_p := L_p([0,1], \lambda)$ the space of $p$-Lebesgue integrable functions on $[0,1]$, and by $L_p$-approximation we mean that the error of an approximation or algorithm is measured in the $L_p$-norm (or seminorm, since we usually talk about individual functions, not equivalence classes).

Though we always keep this example in mind, we introduce the concepts for general classes of functions $f\colon D \to \mathbb{C}$, defined on arbitrary sets $D$. So let $F \subset \mathbb{C}^D := \{f\colon D \to \mathbb{C}\}$ be such a class. We refer to $F$ as a **model class**. We want to find an approximation of some $f \in F$ solely based on function evaluations (and the assumption that $f \in F$). That is, we consider **measurement maps** (also called **encoders** or **information maps**) of the form

$$N_n(f) = \big(f(x_1), \dots, f(x_n)\big) \quad \text{with} \quad \mathcal{P}_n := \{x_1, \dots, x_n\} \subset D. \tag{2.2}$$

To turn this information into an approximation, we are looking for **reconstruction maps** (also called **decoders**)

$$\Phi\colon \mathbb{C}^n \to \mathbb{C}^D.$$

The composition of such a measurement and reconstruction map $A_n = \Phi \circ N_n$ is called a **sampling algorithm**. The approximation error is measured in a norm or seminorm $\|\cdot\|_Y$ of some normed space $Y$ (like $Y = L_p$), which of course only makes sense if we assume that $F \subset Y$ and that $\Phi$ maps to $Y$. The **worst-case error (w.c.e.)** of such an algorithm $A_n$ over $F$ is defined by

$$e(A_n, F, Y) := \sup_{f\in F} \|f - A_n(f)\|_Y . \tag{2.3}$$

The study of this quantity allows us to give reliable (guaranteed) a priori error bounds in the class $F$. See Remark 2.6 for a general comment on the worst-case setting and on the meaning of the class $F$. On the following pages, we will explain what we understand by optimal reconstruction maps and optimal sampling points, and illustrate those concepts for our easy example.

The error of a given algorithm should be compared to the best possible error for the given measurements, i.e., for a fixed set of sampling points $\mathcal{P}_n = \{x_1, \ldots, x_n\} \subset D$, we consider the **minimal (worst-case) error**

$$\mathrm{rad}(\mathcal{P}_n, F, Y) := \inf_{\Phi\colon \mathbb{C}^n\to Y}\ \sup_{f\in F} \left\|f - \Phi\big(f(x_1), \ldots, f(x_n)\big)\right\|_Y . \tag{2.4}$$

Hence, the infimum runs over all sampling algorithms $A_n = \Phi \circ N_n$ with the fixed measurement map $N_n$ from (2.2), and so this is the worst-case error (over $F$) of the best sampling algorithm based on the function values at $\mathcal{P}_n$. This is usually called the **radius of information** $N_n$, and sometimes $N_n$ is used in place of $\mathcal{P}_n$ in (2.4), esp. when we consider more general measurements, see Section 9.

Given a measurement map $N_n$ as in (2.2), i.e., a point set $\mathcal{P}_n$, a reconstruction map $\Phi^*\colon \mathbb{C}^n \to Y$ is called (worst-case) **optimal** for $N_n$ (over $F$) if we have

$$e(\Phi^* \circ N_n, F, Y) = \mathrm{rad}(\mathcal{P}_n, F, Y).$$

Usually, it is too hard to find an optimal reconstruction map and we are already happy if the above relation holds up to a constant independent of $n$, or up to a logarithmic factor in $n$. In this case, we call the reconstruction map almost optimal or near-optimal.

One may also refine the analysis by asking for optimality for every specific data $y = (y_1, \ldots, y_n) \in \mathbb{C}^n$. If the function $f$ is known only via the a priori knowledge $f \in F$ and the data $y_k = f(x_k)$ at known points $x_1, \ldots, x_n$, then

$$F^y := \big\{f \in F\colon f(x_k) = y_k,\ k = 1 \ldots, n\big\}$$

is the set of all the possible functions that $f$ could be equal to. Intuitively, the best we can do is to place the approximation $\Phi^*(y)$ in the 'middle' of the set $F^y$.

Formally, this means that the reconstruction is chosen as a *Chebyshev center* of $F^y$,

$$\Phi^*(y) \in \operatorname*{arg\,min}_{g\in Y} \sup_{f\in F^y} \|f-g\|_Y .$$

This means that $\Phi^*$ minimizes the worst-case error on all of the sets $F^y$, i.e.,

$$e(\Phi^* \circ N_n, F^y, Y) = \inf_{\Phi\colon \mathbb{C}^n\to Y} e(\Phi \circ N_n, F^y, Y) \qquad \text{for all } y \in \mathbb{C}^n.$$

Such a reconstruction map (or the corresponding algorithm) is called **instance optimal**; also the terms **central** and **strongly optimal** are used. Clearly, any instance optimal reconstruction map is also worst-case optimal, as one can see by taking the supremum over all $y$ in the last identity.

In practice, it is usually too difficult to find a Chebyshev center of $F^y$ and it is sometimes easier to choose $\Phi(y)$ as an arbitrary element of $F^y$. If $\Phi$ is chosen in such a manner, the reconstruction map is called **interpolatory**. The algorithm $A_n = \Phi \circ N_n$ then satisfies the two properties: $A_n(f) \in F$ and $A_n(f)(x_k) = y_k$ $(k \le n)$ for all $f \in F$, and also the algorithm $A_n$ is called interpolatory. Any interpolatory reconstruction map is instance optimal up to a factor of two, which follows from the triangle inequality. In particular, any interpolatory decoder $\Phi$ is also worst-case optimal up to a factor of two, i.e.,

$$e(\Phi \circ N_n, F, Y) \le 2 \cdot \mathrm{rad}(\mathcal{P}_n, F, Y).$$

If $F$ is the unit ball of some seminormed space $X$, a particular interpolatory reconstruction map is obtained if $\Phi(y)$ is defined as a minimizer (if it exists) of the seminorm $\|g\|_X$ under the constraints $g(x_k) = y_k$ for $k = 1, \dots, n$. This $\Phi$, and the corresponding $A_n$, are called **spline algorithm**, and $\Phi(y)$ is called a **spline**.

We now illustrate those concepts for our example. Fortunately, in this easy case, there is an explicit form of an instance optimal decoder $\Phi^* := \Phi^*_{\mathcal{P}_n,F_{\mathrm{Lip}}} : \mathbb{R}^n \to F_{\mathrm{Lip}}$ for any point set $\mathcal{P}_n$, which even outputs a function of $F_{\mathrm{Lip}}$. This can be seen by noting that, for every $y \in \mathbb{R}^n$, or more precisely, for every $y \in N_n(F_{\mathrm{Lip}})$, which should be understood as the function values of some $f \in F_{\mathrm{Lip}}$ at the (fixed) points $\mathcal{P}_n := \{x_1, \dots, x_n\} \subset [0,1]$, there are a largest and a smallest function that interpolate the data. These functions are the piecewise linear functions given by

$$h^+_{\mathcal{P}_n,y}(x) := \min_{i\le n} \big(y_i + \mathrm{dist}(x,x_i)\big) \quad \text{and} \quad h^-_{\mathcal{P}_n,y}(x) := \max_{i\le n} \big(y_i - \mathrm{dist}(x,x_i)\big),$$

see Figure 2.1.

In this example, the set of possible functions for given data $y$ at the point set $\mathcal{P}_n$ is given by

$$F^y_{\mathrm{Lip}} = \Big\{ f \in F_{\mathrm{Lip}} \colon h^-_{\mathcal{P}_n,y}(x) \le f \le h^+_{\mathcal{P}_n,y}(x) \Big\} .$$

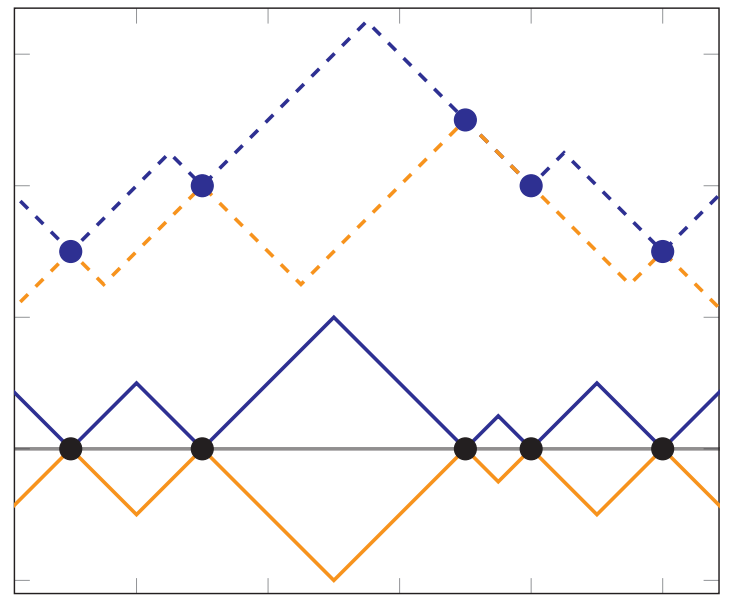

Blue dots: data points
Blue dashed line: $h^+_{\mathcal{P}_n,y}$
Orange dashed line: $h^-_{\mathcal{P}_n,y}$
Black dots: zero data points
Blue line: $h^+_{\mathcal{P}_n,0}$
Orange line: $h^-_{\mathcal{P}_n,0}$

Figure 2.1. The set $F^y_{\mathrm{Lip}}$ for some $\mathcal{P}_n$ and $y$, and for $y = 0$

Since there is no way to distinguish the functions from $F^y_{\mathrm{Lip}}$, it is natural to choose our reconstruction in the 'middle', i.e.,

$$\Phi^*(y) \;:=\; \frac{1}{2}\left(h^+_{\mathcal{P}_n,y} + h^-_{\mathcal{P}_n,y}\right).$$

And indeed, the function $\Phi^*(y)$ is a Chebyshev center of $F^y_{\mathrm{Lip}}$ within the space $Y = L_p$ for any $p$. Hence, for a fixed set of points $\mathcal{P}_n$, the algorithm

$$A^*_n(f) := \Phi^* \circ N_n(f) = \Phi^*\left(f(x_1), \ldots, f(x_n)\right) \tag{2.5}$$

is instance-optimal. Since $A^*_n(f) \in F_{\mathrm{Lip}}$ for every $f \in F_{\mathrm{Lip}}$, and since $A^*_n(f)(x_k) = f(x_k)$ for every $x_k \in \mathcal{P}_n$, the algorithm $A^*_n$ is also interpolatory. The function $A^*_n(f)$ is a piecewise linear function that interpolates the data, but usually with breakpoints different from the sampling points $x_k$, see Figure 2.2.

Another natural reconstruction is given by the piecewise linear interpolating function $\widetilde{\Phi}(y)$ with breakpoints equal to the sampling points. As we have $\widetilde{\Phi}(y) \in F^y_{\mathrm{Lip}}$ for any valid data $y \in N_n(F_{\mathrm{Lip}})$, also the algorithm $\widetilde{A}_n = \widetilde{\Phi} \circ N_n$ is interpolatory, and hence (instance and worst-case) optimal up to a factor of two; but it is not instance optimal in the sense of equality. In fact, the class $F_{\mathrm{Lip}}$ is the unit ball of the space of all Lipschitz-continuous functions with seminorm given by the Lipschitz-constant, and $\widetilde{A}_n(f)$ is the minimal-norm interpolant of the data, i.e., the *spline algorithm*. Both reconstructions $A^*_n(f)$ and $\widetilde{A}_n(f)$ are depicted in Figure 2.2.

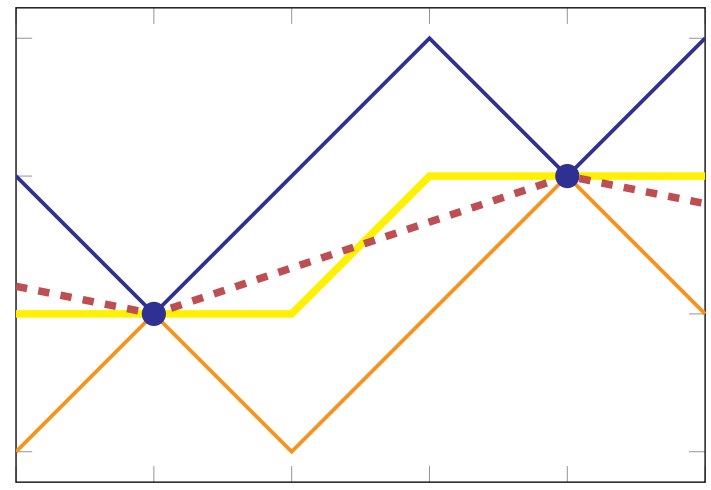

Blue dots: two data points
Blue line: $h^+_{\mathcal{P}_n,y}$
Orange line: $h^-_{\mathcal{P}_n,y}$
Yellow line: $A^*_n(f)$
Dotted line: $\widetilde{A}_n(f)$

Figure 2.2. Instance optimal reconstruction and interpolation for $F_{\mathrm{Lip}}$

Let us now consider the error of the optimal algorithm $A_n^*$. For any $f \in F_{\mathrm{Lip}}$ with data $y = N_n(f)$, we have

$$|f(x) - A_n^* f(x)| \leq \frac{1}{2}\left(h^+_{\mathcal{P}_n,y}(x) - h^-_{\mathcal{P}_n,y}(x)\right), \quad x \in [0,1],$$

which can be seen in Figure 2.2. By definition of the these functions, we have for all $i \leq n$ that

$$\frac{1}{2}\left(h^+_{\mathcal{P}_n,y}(x) - h^-_{\mathcal{P}_n,y}(x)\right) \leq \frac{1}{2}\Big[(y_i + \mathrm{dist}(x,x_i)) - (y_i - \mathrm{dist}(x,x_i)\Big]$$

and hence

$$|f(x) - A_n^* f(x)| \leq \mathrm{dist}(x, \mathcal{P}_n) := \min_{k \leq n} \mathrm{dist}(x, x_k). \tag{2.6}$$

It therefore holds that

$$e(A_n^*, F_{\mathrm{Lip}}, L_p) \leq \|\,\mathrm{dist}(\cdot, \mathcal{P}_n)\|_p. \tag{2.7}$$

For the particular function $f_0 := \mathrm{dist}(\cdot, \mathcal{P}_n)$, which indeed belongs to $F_{\mathrm{Lip}}$, we have $A_n^*(f_0) = 0$ and hence there is equality in (2.6) for all $x \in [0,1]$, and so also (2.7) is an equality. Such a function with 'large' error is also called a **fooling function**, since it 'fools' the algorithm $A_n^*$. In fact, since $N_n(f_0) = N_n(-f_0)$, no algorithm based on the sampling points $\mathcal{P}_n$ can distinguish between $f_0$ and $-f_0$. The function $f_0$ therefore 'fools' any such algorithm, and so its worst-case error is at least $\|f_0\|_p$. We hence obtain

$$\mathrm{rad}(\mathcal{P}_n, F_{\mathrm{Lip}}, L_p) = e(A_n^*, F_{\mathrm{Lip}}, L_p) = \|\,\mathrm{dist}(\cdot, \mathcal{P}_n)\|_p. \tag{2.8}$$

Here, we also encounter two more interesting facts. Fist, the 'worst' fooling function $f_0 = \mathrm{dist}(\cdot, \mathcal{P}_n)$ has data $y = 0$, and second, we have the formula

$$\mathrm{rad}(\mathcal{P}_n, F_{\mathrm{Lip}}, L_p) = \sup\left\{\|f\|_p : f \in F_{\mathrm{Lip}},\ f|_{\mathcal{P}_n} = 0\right\}, \tag{2.9}$$

which means that an optimal algorithm makes its largest error for $y = 0$. That is, the data $y = 0$ is the least informative data that we can obtain. Up to a factor two, this holds in much more generality for convex and symmetric classes $F$ in normed spaces $Y$, see, e.g., Traub, Wasilkowski and Woźniakowski (1988, Chapter 4). We come back to this in Section 9.

**Remark 2.1.** In more abstract terms, the quantity (2.9) is the radius in $L_p$ of the intersection of the set $F_{\mathrm{Lip}}$ with a linear space of codimension $n$. This geometric interpretation of minimal worst-case errors is not only curious, but can also be quite useful in the study of particular approximation problems, or to place these quantities in a broader picture, see Section 9 for details.

It remains to determine the best that can be done based on measurement maps of the form (2.2), i.e., we want to study optimal sampling point sets $\mathcal{P}_n$. For this, we define the **$n$-th minimal (worst-case) error** for approximation of functions from $F$ in $Y$ based on $n$ function evaluations by

$$g_n(F,Y) := \inf_{\mathcal{P}_n} \operatorname{rad}(\mathcal{P}_n, F, Y) = \inf_{A_n} e(A_n, F, Y), \tag{2.10}$$

or written out in detail,

$$g_n(F,Y) = \inf_{\substack{x_1,\dots,x_n \in D \\ \Phi\colon \mathbb{C}^n \to Y}} \sup_{f \in F} \left\| f - \Phi\left(f(x_1),\dots,f(x_n)\right)\right\|_Y.$$

This is called the **$n$-th sampling number**. The choice for the notation $g_n$ instead of the more natural $e_n$ will become clear when we discuss some other classical approximation benchmarks in the following sections (especially Section 9). The sampling numbers will be a main object of study in large parts of this work. Hence, let us make clear how we interpret this quantity:

- An upper bound $g_n < \varepsilon$ means that there is a sampling algorithm with an error at most $\varepsilon$ that uses at most $n$ function evaluations. This gives a hint of what is possible with a fixed sampling budget and sets a natural goal for finding such algorithms. A priori, this is not helpful for actual computations, as such an algorithm may not be known or implementable or stable. However, the proofs are often based on more or less explicit, or even easy, algorithms.
- A lower bound on $g_n$ shows that an algorithm with a smaller error cannot exist, no matter how we choose the sampling points or the reconstruction map.

We call $A_n^*$ an **optimal sampling algorithm** if

$$e(A_n^*, F, Y) = g_n(F,Y).$$

A set $\mathcal{P}_n^*$ of sampling points, or the corresponding measurement map $N_n^*$, is called **optimal** if it satisfies

$$\operatorname{rad}(\mathcal{P}_n^*, F, Y) = g_n(F,Y),$$

that is, if it is possible to build an (arbitrarily close to) optimal sampling algorithm $A_n^* = \Phi^* \circ N_n^*$ upon the sampling points $\mathcal{P}_n^*$ using some reconstruction map $\Phi^*$. Again, it is usually too difficult to find optimal sampling points and algorithms, respectively, and we are already happy if we find 'almost-optimal' sampling points where the previous relations hold up to a constant or logarithmic factor.

**Remark 2.2.** We ignore here (at first) the possibility to choose measurements *adaptively* based on previous measurements, which is an important subject in numerical analysis, and refer to the *general framework* introduced in Section 9. Just note that adaptive methods would not help for the easy example considered in this section (as well as many other examples), see Theorem 9.3.

For our example, we do not only know an optimal reconstruction for any point set $\mathcal{P}_n$, as discussed above, but also how to place the optimal points. It is probably not surprising, and maybe already clear from Figure 2.1, that equidistant nodes are optimal. In fact, since we are in a periodic setting (on the circle), any equidistant set will do, for instance, $\mathcal{P}_n^* := \{k/n \colon k = 1, \ldots, n\}$.

The following goes back at least to the seminal work of Sukharev (1978, 1979), see also (Novak 1988, Section 1.3.6 & 1.3.8).

**Proposition 2.3.** Let $F_{\mathrm{Lip}}$ from (2.1) and $\mathcal{P}_n^*$ consist of $n$ equidistant points on the circle. Then,

$$g_n(F_{\mathrm{Lip}}, L_p) \;=\; \mathrm{rad}(\mathcal{P}_n^*, F_{\mathrm{Lip}}, L_p) \;=\; \begin{cases} \dfrac{1}{2}\left(\dfrac{1}{1+p}\right)^{1/p} n^{-1} & \text{if} \quad 1 \le p < \infty, \\[2ex] \dfrac{1}{2}\, n^{-1} & \text{if} \quad p = \infty. \end{cases}$$

In particular, for any $1 \le p \le \infty$, we have

$$\frac{1}{4n} \le g_n(F_{\mathrm{Lip}}, L_p) \le \frac{1}{2n}.$$

*Proof.* Let $N_n$ be as in (2.2) with a general set $\mathcal{P}_n = \{x_1, \ldots, x_n\} \subset [0,1]$. By formula (2.8), we have

$$\mathrm{rad}(N_n) := \mathrm{rad}(N_n, F_{\mathrm{Lip}}, L_p) = \|\mathrm{dist}(\cdot, \mathcal{P}_n)\|_p \,.$$

Let us first consider the case $p = \infty$. Since the length of the union of the intervals $B_r(x) := [x-r, x+r]$ (considered on the circle) over $x \in \mathcal{P}_n$ is smaller than 1 for all $r < 1/(2n)$, there must be some $x \in [0,1]$ with $\mathrm{dist}(x, \mathcal{P}_n) \ge r$. Therefore,

$$\|\mathrm{dist}(\cdot, \mathcal{P}_n)\|_\infty \ge \frac{1}{2n}$$

for any $\mathcal{P}_n$, where equality holds when $\mathcal{P}_n = \mathcal{P}_n^*$.

For $p < \infty$, let $\lambda$ denote Lebesgue measure on $[0,1]$. Then,

$$\mathrm{rad}(N_n)^p = \int_0^1 \mathrm{dist}(x, \mathcal{P}_n)^p \; \mathrm{d}x = \int_0^\infty \lambda\Big(\{x \colon \mathrm{dist}(x, \mathcal{P}_n)^p \ge t\}\Big) \; \mathrm{d}t.$$

It holds that

$$\lambda\Big(\{x \colon \mathrm{dist}(x, \mathcal{P}_n)^p \ge t\}\Big) = 1 - \lambda\Big(\{x \colon \mathrm{dist}(x, \mathcal{P}_n)^p < t\}\Big) \ge 1 - 2n \cdot t^{1/p},$$

where equality holds if the sets $B_{t^{1/p}}(y)$ for $y \in \mathcal{P}_n$ are pairwise disjoint. Hence,

$$\begin{aligned} \mathrm{rad}(N_n)^p &\ge \int_0^{1/(2n)^p} \lambda\Big(\{x \colon \mathrm{dist}(x, \mathcal{P}_n)^p \ge t\}\Big) \; \mathrm{d}t \\ &\ge \frac{1}{(2n)^p} - 2n \int_0^{1/(2n)^p} t^{1/p} \, \mathrm{d}t = \frac{1}{(2n)^p} \, \frac{1}{1+p} \end{aligned}$$

with equality for $\mathcal{P}_n = \mathcal{P}_n^*$. This proves the statement. $\square$

Hence, for the classes $F_{\mathrm{Lip}}$ we are in a very favorable situation: We know optimal measurement maps and, for every measurement, we know optimal reconstruction maps. They are even "feasible" in the sense that they could be easily implemented on a computer.

The problem is, of course, that this is only a toy example. Usually, the classes of inputs are more complicated, e.g., because they are defined on more complex (high-dimensional) domains, or because they are not given by a (geometric) property that allows for easy determination of all elements that fit the data. In this case, the precise analysis of minimal errors and the construction of optimal algorithms is often out-of-reach. In fact, it is usually already quite challenging to prove the existence of good sampling points, i.e., to show good upper bounds for the sampling numbers $g_n$. In principle, this can be attempted as for the class $F_{\mathrm{Lip}}$ by obtaining a detailed understanding of the given model class and providing an analysis that is precisely tailored to the problem at hand. But model classes are numerous and complicated. This is why, in the following sections, we will rather aim for a general analysis that can be applied for 'any' model class, for universal proof strategies, and for comparison results of the sampling numbers with other approximation benchmarks that might be easier to handle.

Nonetheless, the example $F_{\mathrm{Lip}}$ is quite instructive and we will pick it up a few more times to illustrate new concepts. The first time is in Section 6, where we ask whether random points, once drawn and fixed, do allow for nearly-optimal algorithms. The second time is in Section 7, where we discuss the Monte Carlo methods for sampling recovery. Its third (and final) reappearance will be in Section 8, where the corresponding class of multivariate Lipschitz functions will be used to illustrate the curse of dimensionality.

To illustrate a few additional phenomena, and in order to have a larger variety of (counter-)examples for the upcoming general results, we now introduce a whole scale of model classes: the **Sobolev spaces** $W_q^s := W_q^s([0,1])$ with integrability index $q \in [1,\infty]$ and smoothness $s > 1/q$ (or $s \geq 1$ for $q = 1$). These are classes of continuous functions on the interval $[0,1]$. In the case $s \in \mathbb{N}$, and we mostly consider $s \in \mathbb{N}$, they are defined by

$$W_q^s := \left\{ f \in C^{s-1}([0,1]) \colon f^{(s-1)} \text{ abs. cont. and } f^{(s)} \in L_q \right\}, \tag{2.11}$$

where 'abs. cont.' means *absolutely continuous* and $f^{(k)}$ denotes the $k$-th order (weak) derivative. This space is equipped with the norm

$$\|f\|_{W_q^s}^q := \sum_{k=0}^{s} \|f^{(k)}\|_q^q. \tag{2.12}$$

We denote the unit ball of $W_q^s$ by $F_q^s$. Much is known about this simple, but rich, family of model classes. That is why they will serve as running examples for us to illustrate (problems with) the presented general results.

One can perform an analysis similar to the analysis described above for the class $F_q^s$. Optimal sampling points for these classes are again given by equidistant points (up to constants). More generally, and similar to the formula formula (2.8) for the class $F_{\mathrm{Lip}}$, the minimal error for any fixed point set is again described (up to constants) by the $L_\tau$-norm of the distance function for some $\tau \in [1,\infty]$. For details and references, we refer to the survey of Hinrichs, Krieg, Novak, Prochno and Ullrich (2020). In the way stated below, the formula for general point sets can be found in Krieg and Sonnleitner (2024).

**Proposition 2.4.** Let $1 \le p, q \le \infty$, $s > 1/q$ (or $s \ge 1$ for $q = 1$) and $F_q^s$ be the unit ball of the Sobolev space from (2.11). Then, for every set of points $\mathcal{P}_n \subset [0,1]$,

$$\mathrm{rad}(\mathcal{P}_n, F_q^s, L_p) \asymp \begin{cases} \|\,\mathrm{dist}(\cdot,\mathcal{P}_n)\|_\tau^s, & \text{if } \ p < q, \\ \|\,\mathrm{dist}(\cdot,\mathcal{P}_n)\|_\infty^{s+1/p-1/q}, & \text{otherwise}, \end{cases}$$

where $\tau := s \cdot (1/p - 1/q)^{-1}$. In particular, if $\mathcal{P}_n^*$ consists of $n$ equidistant points, then we have

$$g_n(F_q^s, L_p) \asymp \mathrm{rad}(\mathcal{P}_n^*, F_q^s, L_p) \asymp n^{-s+(1/q-1/p)_+}$$

with $(a)_+ := \max\{a, 0\}$.

That is, the minimal worst case error is of order $n^{-s+(1/q-1/p)_+}$. This seems to be known for a long time. Optimal upper bounds, based on polynomial spline approximation, are contained in Ciarlet (1978, Chapter 3), see also (Schumaker 1981, Corollary 6.26). The corresponding lower bounds are proven even for the larger class of randomized algorithms, as discussed in Section 7, by Mathé (1991, Theorem 4.3). We also refer to the works of Temlyakov (1993), Narcowich, Ward and Wendland (2005), Novak and Triebel (2006, Theorem 23 & Coro. 25) and Vybíral (2007), where similar results can be found in much more generality. In addition to the statement on the minimal worst case error, the proposition shows that for $L_p$-approximation with $p \ge q$, a set of sampling points is optimal if and only if the sup-norm of the distance function, which is more commonly known as the *fill distance* or the *covering radius* of the point set, is of order $\mathcal{O}(1/n)$. Geometrically speaking, this means that the point set can have no 'hole' larger than $\mathcal{O}(1/n)$. This changes in the case $p < q$, where optimality of the point set is not described properly by the fill distance. A fill distance $\mathcal{O}(1/n)$ is still sufficient in this case, but it is not necessary any more. A sampling set can be optimal even if is has a few holes larger than $\mathcal{O}(1/n)$, as long as they are not too large and not too many, so that the average size of the holes is still $\mathcal{O}(1/n)$. For example, this means that a typical set of random (i.i.d. uniformly distributed) sampling points will be optimal in the case $q < p$; such a set has most holes of size $\mathcal{O}(1/n)$ and only few of size $\Omega(\log n/n)$. All these considerations generalize to Sobolov spaces on general domains and manifolds, see Krieg and Sonnleitner (2025, 2024). We will pick up this topic later when we discuss the quality of i.i.d. samples in general, see Section 6.

**Remark 2.5 (Information-based complexity).** Many of the concepts that we presented in this section have been treated systematically in the monograph *Information-Based Complexity* by Traub *et al.* (1988) and Novak and Woźniakowski (2008). There, the reader can also find several general results on the possible optimality of non-adaptive sampling or linear reconstruction maps for a more general family of problems. We will (partially) discuss these matters in Section 9. The task of finding optimal reconstruction maps $\Phi$ for given data is often referred to as *optimal recovery*; this notion dates back to an influential survey of Micchelli and Rivlin (1978). We refer to the recent book of Foucart (2022) for a detailed treatment. As suggested by the title of this work, our focus is more on optimizing the way of acquiring the data, i.e., we optimize over the measurement maps $N_n$. This is more in the spirit of *information-based complexity (IBC)* than in the spirit of optimal recovery. Of course, the two areas are two sides of the same coin and it does not make a lot of sense to consider one without the other; we will certainly discuss the question of explicit reconstruction maps wherever possible. Other classical and a few recent works in these areas and on their intersection, which also highlight possible ambiguities, are given by Novak (1988), Traub and Woźniakowski (1992) and Woźniakowski (2019) in the context of IBC, and by Foucart (2021), Binev, Bonito, Cohen, Dahmen, DeVore and Petrova (2024*a*) and Binev, Bonito, DeVore and Petrova (2024*b*) with an emphasis on optimal recovery.

**Remark 2.6 (The worst-case setting).** Let us end with a general comment on the worst-case setting. In applications, there often is only one single function $f$ that has to be approximated, and not every function from some class $F$. On the other hand, the target function $f$ itself is usually unknown (in the sense that the mapping rule is unkown) and, prior to obtaining the data, there is only some structural a priori knowledge (like smoothness, monotonicity, convexity, permutation-invariance, etc.). For us, the two settings are one and the same. The a priori knowledge implicitly defines a function class $F$ of all the possibilities for the unknown function $f$. The worst-case error of the algorithm over $F$ is the best error bound that can be guaranteed under the given a priori knowledge.

Admittedly, there are also many applications where it is difficult to pin down the proper a priori assumptions about the target function and so the proper class $F$ is unclear. But even for such applications, we believe that a study of different (abstract and concrete) classes $F$ is beneficial. We want to gain insight how structural conditions of possible inputs affect the amount and type of measurements needed for reasonable approximations (as needed, e.g., for representations on computers). Such an analysis helps to identify beneficial properties and (in)evitable bottlenecks, and may consequently lead to new algorithms or a more careful modeling of the application. Of course, there are other interesting possibilities to model the uncertainty in the input $f$, mainly by introducing some form of randomness and considering average errors. In this work, however, we want to concentrate on a worst-case analysis, i.e., on guaranteed error bounds.

## 3. Least squares and best linear approximation

As discussed in Section 2, a sampling algorithm for the recovery of a function $f\colon D \to \mathbb{C}$ consists of two parts: the computation of $n$ samples $f(x_1),\ldots,f(x_n)$ and the application of a reconstruction map $\Phi\colon \mathbb{C}^n \to \mathbb{C}^D$. We first consider the situation where the *approximation space*, i.e., the image of the reconstruction map, is an $m$-dimensional linear space $V_m$ of complex-valued functions on $D$. Typical examples are spaces of polynomials, piecewise polynomials, trigonometric functions, etc. In this article, we mainly want to discuss results that hold for arbitrary linear spaces $V_m$. We therefore do not assume any structure of the domain $D$ or particular properties of the functions in $V_m$. While many (but not all) of the results in this section are well known for particular examples, they have been proven only recently and are therefore relatively unknown in the case of arbitrary spaces $V_m$.

The best approximation of some function $f$ with a function $g$ from $V_m$ and with error measured in the norm of some space $Y$, is given by an element of

$$\operatorname*{arg\,min}_{g\in V_m} \|f-g\|_Y$$

whenever this set is non-empty. If $f$ is only known via samples on a finite set $\mathcal{P}_n$, then the best approximation cannot be computed. It is natural to consider instead

$$\operatorname*{arg\,min}_{g\in V_m} \|f-g\|_* \tag{3.1}$$

with a seminorm $\|\cdot\|_*$ on the space of functions $\mathbb{C}^D$ that depends on the function $f$ solely via its function values $N_n(f)=(f(x))_{x\in\mathcal{P}_n}$ on the set $\mathcal{P}_n$. The hope is that the error of the algorithm (3.1) is close to the error of best approximation if the discrete norm $\|\cdot\|_*$ is close to the norm $\|\cdot\|_Y$ in some sense. Hence, the question arises naturally, how well the original error norm can be approximated by such discrete norms.

In addition, one wants that the discrete norm can be efficiently computed and/or that it is feasible to find a corresponding minimizer in (3.1). A particularly prominent choice in this respect is that $\|\cdot\|_*$ is the Euclidean norm of the function values, possibly with some weights. That is, for given points $\mathcal{P}_n := \{x_1,\ldots,x_n\}\subset D$ and corresponding weights $w:=(w_1,\ldots,w_n)\in\mathbb{R}_+^n$, we consider the seminorm

$$\|f\|_{\mathcal{P}_n^w} := \sqrt{\sum_{i=1}^n w_i\,|f(x_i)|^2}$$

on the space $\mathbb{C}^D$ of all functions $f\colon D\to\mathbb{C}$. For this choice of the seminorm $\|\cdot\|_*$, formula (3.1) defines the **(weighted) least squares algorithm** $A^w_{\mathcal{P}_n,V_m}\colon \mathbb{C}^D\to V_m$, given by

$$A^w_{\mathcal{P}_n,V_m}(f) \in \operatorname*{arg\,min}_{g\in V_m} \sum_{i=1}^n w_i\,|f(x_i)-g(x_i)|^2. \tag{3.2}$$

We usually suppress $V_m$ from the notation since we focus on varying $\mathcal{P}_n$ and $w$. It is easy to see that the expression on the right hand side indeed has a minimizer and the argmin is a nonempty set. In general, the minimizer need not be unique and we need to specify in (3.2) the chosen minimizer. However, we do not need to discuss this here since we will later only consider situations where the minimizer is unique, and hence write

$$A^w_{\mathcal{P}_n}(f) := A^w_{\mathcal{P}_n,V_m}(f) := \operatorname*{arg\,min}_{g\in V_m} \sum_{i=1}^n w_i \, |f(x_i) - g(x_i)|^2. \tag{3.3}$$

The algorithm is specified by $\mathcal{P}_n$, $V_m$ and $w$, as well as its domain and codomain as needed for any mapping. It thus gives the best approximation of $f$ from $V_m$ with respect to the seminorm $\|\cdot\|_{\mathcal{P}^w_n}$, which comes from the pseudo-inner product

$$\langle f, g\rangle_{\mathcal{P}^w_n} := \sum_{i=1}^n w_i \, f(x_i)\, \overline{g(x_i)}, \tag{3.4}$$

which is also defined on $\mathbb{C}^D$. The algorithm $A^w_{\mathcal{P}_n}$ can therefore also be viewed as the orthogonal projection from $\mathbb{C}^D$ onto $V_m$ with respect to this inner product, and so it is also clear that (under the uniqueness assumption) the least squares algorithm is a linear mapping and that it has the projection property (i.e., $A^w_{\mathcal{P}_n} f = f$ for $f \in V_m$).

At first sight, one might expect that the least squares method is limited to approximation in $L_2$-norm, i.e., for $Y = L_2$, since the norm $\|\cdot\|_{\mathcal{P}^w_n}$ mimics the $L_2$-norm. But the following considerations will show that quite useful error bounds can also be obtained for different norms if the points and weights are chosen such that suitable *discretization inequalities* are satisfied. This section contains

- some basics on least squares and discretization, and a first general error bound (Section 3.1),
- error bounds for (unweighted) least squares algorithms in terms of best uniform approximation on linear spaces (Section 3.2),
- error bounds for (weighted) least squares algorithms in terms of best $L_2$-approximation on linear spaces (Section 3.3).
- An algorithm for constructing points and weights for weighted least squares approximation will be discussed afterwards in Section 4.

As the least squares method is a linear algorithm, the analysis in this section will result in several upper bounds on the $n$**-th linear sampling numbers** of a model class $F$ in a semi-normed space $Y \subset \mathbb{C}^D$, defined by

$$g_n^{\mathrm{lin}}(F, Y) := \inf_{\substack{x_1,\ldots,x_n\in D\\ \varphi_1,\ldots,\varphi_n\in Y}} \sup_{f\in F} \Big\| f - \sum_{i=1}^n f(x_i)\, \varphi_i \Big\|_Y, \tag{3.5}$$

which is the worst-case error of the best linear algorithm that uses $n$ samples.

Depending on the model class $F$ and the error norm $\|\cdot\|_Y$, the restriction to linear recovery maps $\Phi$ can make a big difference or no difference at all. Results for nonlinear approximation will be discussed in Section 5. See also Section 9 for a treatment of more general types of algorithms.

### 3.1. *Least squares, discretization, and a first error bound*

Let us start with some (well-known) basics about least squares algorithms in general. If $\{b_1, \ldots, b_m\}$ is a basis of $V_m \subset \mathbb{C}^D$, then we can write each $g \in V_m$ in the form $g = \sum_{k=1}^m c_k b_k$ for some unique $c \in \mathbb{C}^m$. For such $g$, it holds

$$\|g\|^2_{\mathcal{P}_n^w} = c^* G c \tag{3.6}$$

with the Hermitian matrix

$$G := \big( \langle b_k, b_j \rangle_{\mathcal{P}_n^w} \big)_{j \le m, k \le m} \in \mathbb{C}^{m \times m}. \tag{3.7}$$

The matrix $G$ is called the **Gram matrix**. By definition, $g$ is a least squares approximation of $f \in \mathbb{C}^D$ if and only if $c \in \mathbb{C}^m$ is a minimizer of the quadratic functional

$$\Big\| f - \sum_{k=1}^m c_k b_k \Big\|^2_{\mathcal{P}_n^w} = \|f\|^2_{\mathcal{P}_n^w} - 2 \cdot \mathrm{Re}\Big( \sum_{k=1}^m \overline{c_k} \cdot \langle f, b_k \rangle_{\mathcal{P}_n^w} \Big) + c^* G c. \tag{3.8}$$

This expression has a unique minimizer for all $f \in \mathbb{C}^D$ if and only if the Gram matrix is positive definite. By (3.6), this is also equivalent to the fact that the seminorm $\|\cdot\|_{\mathcal{P}_n^w}$ is a proper norm on the space $V_m$, i.e., there is no function $g \in V_m \setminus \{0\}$ with $g(x_1) = \ldots = g(x_n) = 0$. It is clear that this is only possible for $n \ge m$. If this assumption is satisfied, the unique solution of (3.8) is given via the **normal equation**

$$G \cdot c = \gamma, \qquad \text{where} \quad \gamma := \big( \langle f, b_k \rangle_{\mathcal{P}_n^w} \big)_{k \le m}. \tag{3.9}$$

That is, the unique least squares estimator $A^w_{\mathcal{P}_n}(f) = \sum_{k=1}^m c_k b_k$ can be computed by solving this equation. Alternatively, if we introduce the matrix

$$M := \Big( w_i^{1/2} b_k(x_i) \Big)_{i \le n, k \le m} \in \mathbb{C}^{n \times m},$$

then we have $G = M^* M$ and $\gamma = M^* y$ with $y = (w_i^{1/2} f(x_i))_{i \le n}$, so the normal equation (in the case of uniqueness) is equivalent to

$$c = M^+ \cdot y, \qquad \text{where} \quad y = \Big( w_i^{1/2} f(x_i) \Big)_{i \le n}. \tag{3.10}$$

Here, $M^+ := (M^* M)^{-1} M^*$ is the **Moore–Penrose inverse** of $M$. This gives two ways to compute the least squares solution.

**Remark 3.1.** Computationally, it is a little faster to solve the normal equation (3.9) than to compute the Moore-Penrose inverse (3.10). We note that our choice of points and weights will often ensure that the normal equation is well-conditioned. However, since we will work in the worst-case setting, our sampling points will be good for every function $f \in F$ from a given model class $F$, and so the Moore-Penrose inverse, which does not depend on $f$, can be reused for all these functions. Computing the Moore-Penrose inverse can hence be regarded as offline cost, and then the online cost is just the time of the matrix vector multiplication (3.10), which is a lot faster than solving the normal equation. Computations are even faster in the presence of sparsity. We refer to Scott and Tůma (2025).

From (3.10), it is also apparent that the **stability** of the least squares method, that is, how small changes in the data affect the solution, is described by the operator norm of the matrix $M^+$. The operator norm can be considered with respect to different norms in $\mathbb{C}^n$ and $\mathbb{C}^m$. If it is considered with respect to the standard Euclidean norms, the operator norm equals the largest singular value of $M^+$, which is, in turn, given by the reciprocal of the smallest singular value of $G$, denoted by $\sigma_{\min}(G)$. The use of the Euclidean norm is meaningful, since $\|y\|_2 = \|f\|_{\mathcal{P}_n^w}$. However, regarding the norm of the solution, what we are actually interested in is the change of the solution function $g = \sum_k c_k b_k$, which shall be measured in the norm of $Y$. This means that $\sigma_{\min}(G)$ is the correct measure for stability only in the case $\|g\|_Y = \|c\|_2$, that is, if $Y$ is an inner product space and the basis $\{b_1, \ldots, b_m\}$ is orthonormal. In general, we need the operator norm with respect to the induced norm $\|c\|_* := \|g\|_Y$. That is, the stability of the least squares algorithm $A^w_{\mathcal{P}_n}$, in the sense of how changes in the function values of $f$ affect the solution $g = A^w_{\mathcal{P}_n}(f)$ in $Y$, is better described by the operator norm

$$\|M^+\|_* := \sup_{y \in \mathbb{C}^m \setminus \{0\}} \frac{\|M^+ y\|_*}{\|y\|_2} = \sup_{g \in V_m \setminus \{0\}} \frac{\|g\|_Y}{\|g\|_{\mathcal{P}_n}}.$$

In words, the stability constant $\|M^+\|_*$ is the smallest possible constant $K \geq 0$ in the (one-sided) **discretization inequality**

$$\|g\|_Y \leq K \cdot \|g\|_{\mathcal{P}_n^w}, \qquad g \in V_m. \tag{3.11}$$

If $\|\cdot\|_Y$ is a norm on $V_m$, under this condition, of course also $\|\cdot\|_{\mathcal{P}_n^w}$ is a norm on $V_m$ and the least squares estimator is uniquely defined. Let us summarize a bit in the following lemma.

**Lemma 3.2.** Let $V_m \subset \mathbb{C}^D$ be an $m$-dimensional linear space and let $\|\cdot\|_Y$ be a norm on $V_m$. If $\mathcal{P}_n = \{x_1, \ldots, x_n\} \subset D$ and $(w_1, \ldots, w_n) \in \mathbb{R}^n_+$ satisfy the discretization inequality (3.11), then the least squares algorithm $A^w_{\mathcal{P}_n}$ from (3.3) is a well-defined linear projection from $\mathbb{C}^D$ onto $V_m$ and satisfies

$$\left\| A^w_{\mathcal{P}_n}(f) - A^w_{\mathcal{P}_n}(g) \right\|_Y \leq K \cdot \|f - g\|_{\mathcal{P}_n^w}$$

for all $f, g \colon D \to \mathbb{C}$.

*Proof.* We already saw that in this case, the least squares algorithm is uniquely defined and an orthogonal projection w.r.t. the discrete inner product. This gives

$$\|A^w_{\mathcal{P}_n}(f) - A^w_{\mathcal{P}_n}(g)\|_Y \;=\; \|A^w_{\mathcal{P}_n}(f-g)\|_Y \;\le\; K\|A^w_{\mathcal{P}_n}(f-g)\|_{\mathcal{P}^w_n} \;\le\; K\|f-g\|_{\mathcal{P}^w_n},$$

where we used the orthogonal projection property in the last inequality and the discretization condition in the previous inequality.

□

It is well known that discretization inequalities, such as (3.11), do not only imply stability, but also lead to upper bounds on the error of the weighted least squares method. This principle has already been used by many authors, including Gröchenig (1999), Cohen, Davenport and Leviatan (2013), Migliorati, Nobile, von Schwerin and Tempone (2014), Temlyakov (2021), to name just a few examples. Also here, we want to make use of this principle. The foundation for the error bounds of this section is given by the following simple proposition.

**Proposition 3.3.** Let $D$ be a set, $Y \subset \mathbb{C}^D$ be a linear space with seminorm $\|\cdot\|_Y$, and $V_m$ be a finite-dimensional subspace of $Y$ on which $\|\cdot\|_Y$ is a norm. Moreover, let $\mathcal{P}_n := \{x_1, \dots, x_n\} \subset D$ and $w := (w_1, \dots, w_n) \in \mathbb{R}^n_+$ be such that

$$\|g\|_Y \;\le\; K \cdot \sqrt{\sum_{i=1}^n w_i\, |g(x_i)|^2} \qquad \text{for all} \quad g \in V_m \tag{3.12}$$

and some $K < \infty$. Then, for all $f \in Y$ and $g \in V_m$, we have

$$\big\|f - A^w_{\mathcal{P}_n}(f)\big\|_Y \;\le\; \big\|f-g\big\|_Y + K \cdot \sqrt{\sum_{i=1}^n w_i\, |f(x_i) - g(x_i)|^2}.$$

In other words, the error of the weighted least squares method is bounded by the error of best approximation in the norm $\|\cdot\|_* := \|\cdot\|_Y + K\|\cdot\|_{\mathcal{P}^w_n}$, that is,

$$\big\|f - A^w_{\mathcal{P}_n}(f)\big\|_Y \;\le\; \inf_{g \in V_m} \|f-g\|_*$$

for all $f\colon D \to \mathbb{C}$. Note that $g$ on the right-hand side is at our disposal, and we may choose a specific one, like a projection onto $V_m$. Also note that $K$ often depends on $n$ and $m$, respectively. The proof is immediate from the stability in Lemma 3.2.

*Proof.* By the triangle inequality

$$\|f - A^w_{\mathcal{P}_n}(f)\|_Y \;\le\; \|f-g\|_Y + \|g - A^w_{\mathcal{P}_n}(f)\|_Y,$$

and by Lemma 3.2,

$$\|g - A^w_{\mathcal{P}_n}(f)\|_Y \;=\; \|A^w_{\mathcal{P}_n}(g) - A^w_{\mathcal{P}_n}(f)\|_Y \;\le\; K\|f-g\|_{\mathcal{P}^w_n}.$$

□

In fact, it is not really essential for the error bound in Proposition 3.3 that the approximation space is a linear space, or that we consider the square norm of the function values in the definition of the least squares algorithm. We come back to this in Section 5.1, see Proposition 5.6, when we consider nonlinear approximation. For now, we stick to the linear setting.

A term like $\|f - g\|_Y$ in the error bound of Proposition 3.3 cannot be avoided as it comes from the fact that $A^w_{\mathcal{P}_n}$ maps to $V_m$. Therefore, we clearly have

$$\|f - A^w_{\mathcal{P}_n}(f)\|_Y \;\geq\; \inf_{g \in V_m} \|f - g\|_Y$$

for any $f \in Y$. The first hope for further bounding the right-hand side of the error bound in Proposition 3.3 may be that the second term on the right, i.e., $\|f - g\|_{\mathcal{P}^w_n}$, can also be bounded by $\|f - g\|_Y$ for a suitable choice of $\mathcal{P}_n$ and $w$. This can be done, e.g., in the case of uniform approximation as discussed in the next Section 3.2.

In general, however, it is out of reach to bound the $Y$-error of a (deterministic) sampling algorithm by best-approximation in $Y$. For example, in the important case $Y = L_2$, changing a function $f \in Y$ on the set $\mathcal{P}_n$ of sampling points does not change its norm in $Y$ and hence also not the error of best $L_2$-approximation. But the outcome of a corresponding sampling algorithm may change drastically. This is even true if $f$ is assumed to be continuous or smooth; the outcome of the algorithm can be manipulated drastically without a significant change of the error of best approximation in $L_2$. So the error of the least squares algorithm cannot possibly be bounded by the error of best approximation in $L_2$, at least not in this direct sense. It is one purpose of this work to discuss different remedies for this issue, and in particular, to study properties of suitable model classes $F$ such that the problem can be avoided. This will be discussed in the following sections. In particular, a direct bound for the $L_2$-error of the algorithm can be obtained by the error of best approximation in a larger norm, like the sup-norm (Section 3.2) or the norm in a reproducing kernel Hilbert space (Section 3.3); the latter also leads, indirectly, to bounds in terms of the error of best $L_2$-approximation.

Let us present a few remarks before discussing these approaches in detail.

**Remark 3.4 (Interpolatory algorithms).** In light of the considerations of Section 2, let us comment on the near-optimality of *interpolatory algorithms*. The least squares method described above will, in most cases, not be interpolatory, because we will mostly choose $n$ to be larger than $m$, i.e., the dimension of $V_m$. This approach has quite some advantages, as will become clear below. In fact, we can often work with a fixed $V_m$, which could also be helpful for applications of the results, because $V_m$ is usually part of the *architecture* of an algorithm and is specified (by the user) beforehand. However, already asking for exactly $n = m$ samples so that (3.12) is fulfilled would be a major restriction that would, in general, not allow for "good" bounds.

On the other hand, once we *know* good sampling points $\mathcal{P}_n$, i.e., a good measurement map for a model class $F$, it follows from earlier consideration that there exists a near-optimal interpolatory algorithm based on the same measurements. This interpolatory algorithm is different from the least squares algorithm and no longer maps to $V_m$. It is also unclear, in general, how to give an explicit formula for computing the corresponding reconstruction. An exception is the *representer theorem* for reproducing kernel Hilbert spaces with a known kernel, see Remark 3.28.

**Remark 3.5 (Hyperinterpolation).** Instead of the condition (3.12), the stronger condition

$$\|g\|_2 = \sqrt{\sum_{i=1}^{n} w_i |g(x_i)|^2} \qquad \text{for all} \quad g \in V_m \tag{3.13}$$

of exact discretization for equal or other positive weights $w_i$ often appears in the literature. The condition (3.13) is equivalent to the quadrature rule $Q_m f = \sum_{i=1}^{m} w_i f(x_i)$ being exact for all integrands from $V_m \cdot \overline{V}_m$. The most famous example is probably given by the space $V_m$ of polynomials of degree smaller than $m$ on an interval and the Gaussian quadrature points, in which case $n = m$ points suffice. In general, however, exact discretization requires $n$ to be much larger than $m$. Under the exact discretization condition, it is not hard to see that the (weighted) least squares method can alternatively be written explicitly as

$$A^{w}_{\mathcal{P}_n}(f) = \sum_{k=1}^{m} \langle f, b_k \rangle_{\mathcal{P}^w_n}\, b_k \tag{3.14}$$

where $\{b_k\}_{k=1}^{m}$ is any $L_2$-orthonormal basis of $V_m$. Such methods are often referred to as *hyperinterpolation*, a term which seems to have been coined by Sloan (1995). It refers to the fact that there will generally be too many points compared to the degrees of freedom for the method to be interpolating (though, of course, this is the case also for other methods such as least squares without exact discretization). Hyperinterpolation is studied, for instance, by Erdös and Turan (1937), Hua and Wang (1981), Sloan (1995), Reimer (2003), Dai (2006), Hansen, Atkinson and Chien (2009), see also the recent survey of An, Ran and Wu (2025). Without the exact discretization condition, the method (3.14) is different from the least squares method, but still of interest. It leads to a computationally faster algorithm but also to additional error terms compared to the least squares algorithm, see, e.g., An, Cai and Goda (2026, Theorem 3.3).

**Remark 3.6 (Noisy measurements).** Another advantage of the simple and linear method $A^w_{\mathcal{P}_n}$ is its robustness for *noisy measurements*. In fact, assume that we do not have access to exact measurements, but only to some noisy or corrupted values $y_i = f(x_i) + \delta_i$ for some unknown or even random $\delta := (\delta_1, \ldots, \delta_n) \in \mathbb{C}^n$. Then the output of the algorithm is equal to $\tilde{f} := A^w_{\mathcal{P}_n}(f_{\mathrm{cor}})$ with the corrupted function

$f_{\mathrm{cor}} = f + \sum_{i=1}^{n} \delta_i 1_{x_i}$. By Lemma 3.2,

$$\|A^w_{\mathcal{P}_n}(f_{\mathrm{cor}}) - A^w_{\mathcal{P}_n}(f)\|_Y \;\leq\; K \left\| \sum_{i=1}^{n} \delta_i 1_{x_i} \right\|_{\mathcal{P}_n^w}$$

and hence Proposition 3.3 gives

$$\left\| f - \tilde{f} \right\|_Y \;\leq\; \left\| f - g \right\|_Y + K \cdot \sqrt{\sum_{i=1}^{n} w_i \, |f(x_i) - g(x_i)|^2} + K \cdot \sqrt{\sum_{i=1}^{n} w_i \delta_i^2}$$

for all $g \in V_m$. With this in mind, the following results may be easily modified to take into account possible noise $\delta$, simply by adding $K$ times the weighted Euclidean norm of the noise in every error bound. We refer to Plaskota (1996) for a collection of results on noisy information from a complexity point of view.

### *3.2. Bounds by best uniform approximation*

In this section, we discuss how close sampling algorithms can get to the error of best uniform approximation from a given $m$-dimensional space $V_m$. For a function $f\colon D \to \mathbb{C}$, the error of **best uniform approximation** from $V_m$ is defined by

$$d(f, V_m)_\infty := \inf_{g \in V_m} \|f - g\|_\infty, \tag{3.15}$$

where

$$\|h\|_\infty \;:=\; \sup_{x \in D} |h(x)|, \quad \text{for} \quad h \in \mathbb{C}^D.$$

We also denote the **space of bounded functions** on a set $D$ by

$$B(D) := \left\{ f\colon D \to \mathbb{C}\colon\ \|f\|_\infty < \infty \right\},$$

and often just write $B := B(D)$. For a model class $F \subset B(D)$, the best (worst case over $F$) uniform approximation by any subspace of dimension $m$ is described by the **Kolmogorov widths** of $F$ in $B = B(D)$, given by

$$d_m(F, B) := \inf_{V_m} \sup_{f \in F} d(f, V_m)_\infty, \tag{3.16}$$

where the infimum is over all $m$-dimensional subspaces $V_m \subset B(D)$.

The sampling algorithms that achieve the best results in this context are, in fact, **unweighted least squares algorithms**. For a set of sampling points $\mathcal{P}_n := \{x_1, \ldots, x_n\} \subset D$ and a function $f\colon D \to \mathbb{C}$, the unweighted least squares algorithm is defined by

$$A_{\mathcal{P}_n}(f) \;:=\; \operatorname*{arg\,min}_{g \in V_m} \sum_{i=1}^{n} |f(x_i) - g(x_i)|^2. \tag{3.17}$$

In the terminology of the last section, this means that all the weights of the least squares algorithm are equal to one (or another constant positive value, since the algorithm does not change if we multiply all weights by the same positive constant). The algorithm is hence very simple and the challenge is 'only' to find good sampling points. Applying Proposition 3.3 to the unweighted least squares algorithm, we easily get the following error bound, which gives us a condition for good sampling points. Recall that we measure the error in $Y$, usually $Y = L_p$ or $Y = B = B(D)$.

**Proposition 3.7.** Let $D$ be a set, $Y \subset \mathbb{C}^D$ a linear space with seminorm $\|\cdot\|_Y$, $V_m$ be a finite-dimensional subspace of $Y$ on which $\|\cdot\|_Y$ is a norm, and $\mathcal{P}_n = \{x_1, \ldots, x_n\} \subset D$ be a finite multiset. Then, the conditions

$$\|g\|_Y \le C \cdot \|g\|_\infty \qquad \text{for all} \quad g \in Y,$$

and

$$\|g\|_Y \le K \cdot \sqrt{\frac{1}{n} \sum_{i=1}^{n} |g(x_i)|^2} \qquad \text{for all} \quad g \in V_m \tag{3.18}$$

for some $C, K < \infty$ imply that

$$\left\| f - A_{\mathcal{P}_n}(f) \right\|_Y \; \le \; (C + K) \cdot \inf_{g \in V_m} \|f - g\|_\infty$$

for all $f \in Y$. In particular, for $F \subset B(D)$, this implies that

$$g_n^{\mathrm{lin}}(F, Y) \le (C + K) \cdot d_m(F, B).$$

By a multiset, we mean that a point can appear more than once and is counted with the corresponding multiplicity. This is needed in special situations, e.g., if the domain $D$ is a finite set and $n > \#D$. Of course, multiple copies of the same point could be summarized into one element without changing the algorithm by using a weighted least squares algorithm (instead of an unweighted one) with weights proportional to the multiplicity of the point in the mulitset. But this only makes the formulas unnecessarily complicated, so we rather allow multisets.

*Proof.* First note that the algorithm (3.2) does not change if we multiply all weights by the same number. Hence, taking into account the conditions and Proposition 3.3 with $w_1 = \ldots = w_n = 1/n$, it remains to observe that

$$\sqrt{\frac{1}{n} \sum_{i=1}^{n} |f(x_i) - g(x_i)|^2} \le \|f - g\|_\infty,$$

which immediately gives the error bound. □

Hence, in order to obtain good error bounds from Proposition 3.7, we 'only' have to care about choosing points that satisfy the one-sided discretization bound (3.18) with a factor $K > 0$ as small as possible.

We first consider the case $Y = B(D)$ and $n = m$. In this case, a condition of the form (3.18) implies that the mapping $g \mapsto (g(x_i))_{i\le m}$ is a bijection from $V_m$ to $\mathbb{C}^m$. Hence, for any data $f(x_i)$, $i \le m$, there is a unique function $g \in V_m$ with $g(x_i) = f(x_i)$ for all $i \le m$, i.e., a unique interpolating function from $V_m$. By definition, the least squares estimator $A_{\mathcal{P}_n}(f)$ must be equal to this function.

The case $n = m$ has already been discussed by Novak (1988), see Remark 3.10. Here, the discretization bound relies on (a special case of) *Auerbach's lemma*, see e.g. Pietsch (1978, B.4.8). For convenience of the reader, and since we think the similarities to 'modern' methods are instructive, we present an elementary proof.

**Lemma 3.8 (Auerbach's lemma).** Let $V_n \subset B(D)$ be $n$-dimensional. Then, for every $\varepsilon > 0$, there exist $\varphi_1, \ldots, \varphi_n \in V_n$ and $x_1, \ldots, x_n \in D$ such that

$$\|\varphi_i\|_\infty = 1 + \varepsilon, \qquad \text{and} \qquad \varphi_i(x_j) = \delta_{ij}.$$

*Proof.* We consider a basis $v_1, \ldots, v_n$ of $V_n$ and set

$$d(t_1, \ldots, t_n) := \left|\det\big(v_i(t_j)\big)\right| \qquad \text{for} \qquad t_1, \ldots, t_n \in D.$$

There are $x_1, \ldots, x_n \in D$ such that $(1 + \varepsilon) \cdot d(x_1, \ldots, x_n) \ge \sup_{t_j} d(t_1, \ldots, t_t) > 0$. In particular, the matrix $(v_i(x_j))_{ij}$ is invertible and we can find $\varphi_1, \ldots, \varphi_n \in V_n$ such that

$$\sum_{j=1}^{n} v_i(x_j) \cdot \varphi_j = v_i \qquad \text{for} \qquad i = 1, \ldots, n.$$

This implies $\varphi_i(x_j) = \delta_{ij}$. Now, it follows from

$$\sum_{j=1}^{n} v_i(x_j) \cdot \varphi_j(t_k) = v_i(t_k)$$

that

$$\det\big(v_i(x_j)\big) \cdot \det\big(\varphi_j(t_k)\big) = \det\big(v_i(t_k)\big)\,.$$

By the maximizing properties of the $x_j$, we obtain

$$\left|\det\big(\varphi_j(t_k)\big)\right| \le 1 + \varepsilon \qquad \text{for all} \qquad t_1, \ldots, t_n \in D.$$

Now, for any $j = 1, \ldots, n$ and $x \in D$, we set $t_j := x$ and $t_k := x_k$ whenever $k \ne j$. We get

$$|\varphi_j(x)| = \left|\det\big(\varphi_j(t_k)\big)\right| \le 1 + \varepsilon.$$

Consequently, $\|\varphi_j\|_\infty \le 1 + \varepsilon$. □

In particular, the functions $\varphi_k$ from Auerbach's lemma form a basis of $V_n$ and we can write every $g \in V_n$ as

$$g = \sum_{k=1}^{n} g(x_k) \cdot \varphi_k.$$

This implies

$$\|g\|_\infty \le \sqrt{\sum_{k=1}^{n} \|\varphi_k\|_\infty^2} \cdot \sqrt{\sum_{i=1}^{n} |g(x_i)|^2}$$

and hence the *discretization* inequality (3.18) with $K = n + \varepsilon$ (with $\varepsilon > 0$ different from above, but still arbitrarily close to zero), i.e.,

$$\|g\|_\infty \le (n+\varepsilon) \cdot \sqrt{\frac{1}{n} \sum_{i=1}^{n} |g(x_i)|^2} \qquad \text{for all} \quad g \in V_n. \tag{3.19}$$

Together with Proposition 3.7 with $Y = B(D)$ and $C = 1$, this implies the following result of Novak (1988, 1.2.5).

**Theorem 3.9.** Let $D$ be a set and $\varepsilon > 0$. Then, for every $n$-dimensional subspace $V_n \subset B := B(D)$ there exist $\mathcal{P}_n = \{x_1, \dots, x_n\} \subset D$ such that $A_{\mathcal{P}_n}$ from (3.17) satisfies

$$\|f - A_{\mathcal{P}_n}(f)\|_\infty \le (n+1+\varepsilon) \cdot d(f, V_n)_\infty \qquad \text{for all} \quad f \in B.$$

In particular, for any $F \subset B$, we have

$$g_n^{\mathrm{lin}}(F, B) \le (n+1) \cdot d_n(F, B). \tag{3.20}$$

Novak (1988, 1.2.5) also showed that the constant $(n+1)$ in (3.20) is optimal in the sense that, for each $n$, there is a set $F$ with $g_n(F, B) = (n+1) \cdot d_n(F, B)$.

**Remark 3.10 (Discretization of the uniform norm).** The discretization inequality (3.19) implied by Auerbach's lemma immediately also implies the weaker discretization inequality

$$\|g\|_\infty \le (n+\varepsilon) \sup_{i=1,\dots,n} |g(x_i)| \qquad \text{for all} \quad g \in V_n. \tag{3.21}$$

Novak (1988, 1.2.5) proved Theorem 3.9 based on this weaker discretization condition. Recall that for $n = m$, the least squares algorithm $A_{\mathcal{P}_n}$ from (3.17) is *interpolatory*, and hence the same as used by Novak (1988).
We add that (sequences of) point sets $\mathcal{P}_n \subset D$ such that $\|f\|_\infty \le K_n \max_{x \in \mathcal{P}_n} |f(x)|$ for all $f \in V_n$ are well-studied, especially for $V_n$ being $n$-dimensional spaces of algebraic or trigonometric polynomials, respectively; see (Bos, Calvi, Levenberg, Sommariva and Vianello 2011, Kashin, Konyagin and Temlyakov 2023). If $K_n$ and $\#\mathcal{P}_n$ depend at most polynomially on $\dim(V_n) = n$, then the corresponding point sets $(\mathcal{P}_n)$ are sometimes called *(weakly) admissible meshes*. If $K_n \asymp 1$

they are called *norming sets* and, if additionally $\#\mathcal{P}_n \asymp n$, then they are called *optimal meshes*. We refer to Kroó (2019) and Dai and Prymak (2024) for the existence of optimal meshes for algebraic polynomials on arbitrary multivariate convex domains, see also Bos (2018) and Xu and Narayan (2023). For a discussion of corresponding least squares and Lagrange approximations, we refer to Calvi and Levenberg (2008) and Guo, Narayan and Zhou (2020). For other spaces $V_n$, there need not exist optimal meshes; there are even spaces $V_n$ where the relation $\log(\#\mathcal{P}_n) \asymp n$ is necessary (Kashin *et al.* 2023, Thm. 1.2) in order to obtain a norming set. This condition is also sufficient; for any space $V_n$, there exist norming sets $\mathcal{P}_n$ with $\log(\#\mathcal{P}_n) \asymp n$, see (Kashin *et al.* 2023, Thm. 1.3). It follows from the considerations below that $K_n \lesssim \sqrt{n}$ is the best bound that can be achieved for arbitrary spaces $V_n$ if we require $\#\mathcal{P}_n \asymp n$. This cannot be improved even if we allow $\#\mathcal{P}_n \asymp n^{42}$, see also Remark 3.19.

**Remark 3.11 (Interpolation with polynomials).** In the context of interpolation, especially with polynomials, the *Lebesgue constant* is a well-known characteristic for the quality of the approximation in the uniform norm. The Lebesgue constant is the norm of the *interpolation operator* as an operator on $B(D)$ or $C(D)$. The considerations above show: For any $n$-dimensional space $V_n$, we can find $n$ points such that their Lebesgue constant is at most of order $n$, and that it is not possible to improve this estimate in general. For particular spaces $V_n$, of course, there can be points with a much smaller Lebesgue constant. For example, for univariate polynomials of degree less than $n$ on the interval $[-1, 1]$, the Lebesgue constant of $n$ Chebychev points is bounded above $\frac{\pi}{2}\log(n) + 1$, and no other point set can be better up to a (small) additive constant. This can be traced back to the fundamental work of Bernstein (1931), see also Faber (1914), Erdős (1961) and Rivlin (1974). We also refer to Bernardi and Maday (1992) for error bounds in $L_2$ of interpolation with polynomials in Sobolev spaces.

Theorem 3.9 cannot be improved if we insist in $n = m$. However, if we allow some **oversampling**, i.e., $n > m$, then the factor $m$ (which equals $n$ above) can even be improved to $\sqrt{m}$ as we will see below. Moreover, we can then derive even better bounds in weaker error norms like $Y = L_p$ with $p < \infty$ (however, still in comparison with best approximation in $B(D)$). Recall that we consider the seminormed spaces $L_p = L_p(D, \mu)$, $1 \le p \le \infty$, where $\mu$ is any probability measure on $D$, as spaces of functions, and not equivalence classes, hence we talk about 'seminorms'. As before, we can choose $C = 1$ in Proposition 3.7.

The most interesting case in this respect might be $Y = L_2$ with norm $\|\cdot\|_2$. Here, we can employ the recent result of Chkifa, Dolbeault, Krieg and Ullrich (2026, Proposition 8), which states that, for any $V_m \subset L_2$ with $\dim(V_m) = m$ and any

$n > m$, there exist points $x_1, \ldots, x_n \in D$ such that

$$\left(1 - \sqrt{\frac{m}{n}}\right) \cdot \|g\|_2 \leq \sqrt{\frac{1}{n} \sum_{i=1}^{n} |g(x_i)|^2} \qquad \text{for all} \quad g \in V_m. \tag{3.22}$$

There is quite some history of such **unweighted discretizations** of $L_2$-norms, some of which we discuss in Remark 3.14 and in Section 4.1, where we present further discretization bounds and discuss the construction of the corresponding points. The inequality (3.22) can be proved by relatively elementary methods, along the lines of the proof of the main theorem of Batson, Spielman and Srivastava (2014). It leads to the following error bound, see Chkifa *et al.* (2026, Corollary 16).

**Theorem 3.12.** Let $(D, \mu)$ be a probability space. Then, for every $m$-dimensional subspace $V_m \subset L_2 = L_2(D, \mu)$ and every $n > m$, there exists $\mathcal{P}_n = \{x_1, \ldots, x_n\} \subset D$ such that $A_{\mathcal{P}_n}$ from (3.17) satisfies

$$\left\| f - A_{\mathcal{P}_n}(f) \right\|_2 \leq \left(1 + \frac{1}{1 - \sqrt{m/n}}\right) \cdot d(f, V_m)_\infty$$

for all $f \in L_2$. In particular, for $F \subset L_2$, we have

$$g_{2n}^{\mathrm{lin}}(F, L_2) \leq 5 \cdot d_n(F, B). \tag{3.23}$$

*Proof.* We apply Proposition 3.7 with $Y = L_2$, the given space $V_m$, and the points from (3.22), so that we have $C = 1$ and $K = (1 - \sqrt{m/n})^{-1}$. □

There are classes $F$ for which the inequality (3.23) is sharp up to constants, as shown by the following example.

**Example 3.13.** Let $F_q^s$ be the unit ball of the univariate Sobolev spaces, introduced in (2.11). It is known that

$$g_n(F_q^s, L_2) \asymp d_n(F_q^s, B) \asymp n^{-s+(1/q-1/2)_+}$$

for all $1 \leq q \leq \infty$ and $s \in \mathbb{N}$ (with $s > 1$ for $q = 1$, in order to have a compact embedding), see Appendix A.2. Hence, Theorem 3.12 even gives the optimal order for all these classes.

**Remark 3.14 (History).** The inequality (3.23) was first proved by Temlyakov (2021) for classes $F \subset C(D)$ of continuous functions for compact $D \subset \mathbb{R}^d$ and with larger constants. He used a similar approach, but relied on the weighted discretization inequality of Limonova and Temlyakov (2022), and hence used the weighted algorithm $A_{\mathcal{P}_n}^w$. Probably the most important difference is that the proof of the latter relies on Marcus, Spielman and Srivastava (2015), and is therefore *non-constructive*. In contrast, using the method of Batson *et al.* (2014) leads to a (semi-)*constructive* algorithm, as discussed in Section 4. A corresponding error bound together with an explicit algorithm for constructing the points for the unweighted least squares method was first given by Bartel, Schäfer and Ullrich (2023). The above version of the result is from Chkifa *et al.* (2026) and gives improved constants. Also note that a (semi-)constructive weighted version with similarly good constants is given by Dolbeault and Chkifa (2024).

We now turn back to approximation in the uniform norm. Another important feature of (3.22) is that it can be applied to any $\mu$. Hence, if we can find a probability measure $\mu^*$ on $D$ such that $\|g\|_\infty \le C \cdot \|g\|_{L_2(\mu^*)}$ for all $g \in V_m$, then we can apply (3.22) to this measure to get

$$\left(1 - \sqrt{\frac{m}{n}}\right) \cdot \|g\|_\infty \;\le\; C \cdot \sqrt{\frac{1}{n}\sum_{i=1}^{n} |g(x_i)|^2} \qquad \text{for all} \quad g \in V_m.$$

Fortunately, such a measure always exists with $C = \sqrt{m+\varepsilon}$ due to the following result of Kiefer and Wolfowitz (1960).

**Lemma 3.15.** Let $D$ be a set and $V_m$ be an $m$-dimensional subspace of $B(D)$. For every $\varepsilon > 0$, there exists a natural number $N \le m^2 + 1$, distinct points $x_1, \ldots, x_N \in D$ and positive weights $(\lambda_k)_{k=1}^N$ satisfying $\sum_{k=1}^N \lambda_k = 1$ such that, for all $f \in V_m$, we have

$$\|f\|_\infty \;\le\; \sqrt{m+\varepsilon} \cdot \left( \sum_{k=1}^{N} \lambda_k \, |f(x_k)|^2 \right)^{1/2}.$$

If $V_m \subset C(D)$ for some compact topological space $D$, then $\varepsilon = 0$ also works.

*Proof.* For a detailed proof, we refer to Krieg, Pozharska, Ullrich and Ullrich (2026*c*). Let us only highlight the ideas. Similarly as for Lemma 3.8, we choose points that maximize the determinant of some 'Gram matrix'. However, since we allow $N \ge m$ points and do not restrict to equal weights, we now maximize over a larger set and obtain a better bound. In fact, given a basis $b = (b_1, \ldots, b_m)^\top$ of $V_m$, as well as points $x_1, \ldots, x_N \in D$ and weights $\lambda_1, \ldots, \lambda_N > 0$, we define the *Gram matrices*

$$G_N^\lambda \;:=\; \Big( \langle b_i, b_j \rangle_{\mathcal{P}_N^\lambda} \Big)_{i,j \le m} \tag{3.24}$$

with $\langle \cdot, \cdot \rangle_{\mathcal{P}_N^\lambda}$ as in (3.4). It holds that

$$\sup_{f \in V_n \setminus \{0\}} \frac{\|f\|_\infty^2}{\|f\|_{\mathcal{P}_N^\lambda}^2} \;=\; \sup_{x \in D} \; b(x)^* (G_N^\lambda)^{-1} b(x).$$

To bound the right hand side by $m + \varepsilon$ for some points and weights, we choose a positive-definite matrix $G \in \mathbb{C}^{m\times m}$ such that

$$\det(G) \;=\; \sup\left\{\det(G_N^\lambda)\colon x_1,\dots,x_N \in D,\; \sum_{k=1}^N \lambda_k = 1,\; \lambda_k > 0,\; N \in \mathbb{N}\right\} \;>\; 0.$$

We then define the polynomial

$$p_x(\alpha) := \det\big((1-\alpha)\cdot G + \alpha\cdot b(x)b(x)^*\big)\cdot \det(G)^{-1},$$

so that $p_x'(0) = \mathrm{Tr}\left(G^{-1}(b(x)b(x)^* - G)\right) = b(x)^* G^{-1} b(x) - n$, where Tr denotes the trace and we used its invariance under circular shifts. By the maximizing property of $G$, we have $p_x'(0) \le 0$, and thus also $b(x)^* G^{-1} b(x) \le n$ for all $x \in D$.

This inequality remains true up to the given tolerance $\varepsilon$ if the matrix $G$ is replaced by a sufficiently good approximation of the form (3.24). By Caratheodory's theorem we obtain that we can take $N \le \dim(\mathrm{span}(V_m\cdot \overline{V}_m)) + 1 \le m^2+1$.

□

With this, we obtain the required discretization inequality (3.18) with $Y = B(D)$, i.e., the inequality

$$\left(1-\sqrt{\frac{m}{n}}\right)\cdot \|g\|_\infty \;\le\; \sqrt{\frac{m}{n}\sum_{i=1}^n |g(x_i)|^2} \qquad \text{for all} \quad g \in V_m, \tag{3.25}$$

see again Chkifa *et al.* (2026, Proposition 8). (The $\varepsilon$ from above can be removed by replacing (3.22) with a slightly stronger inequality.)

Although the method to generate such points is again constructive to some extent, finding the measure $\mu^*$ is usually problematic. This problem is closely related to so-called **D-optimal designs** and is discussed in detail, e.g., by Kiefer (1959, 1974) and Huan, Jagalur and Marzouk (2024). See also Bartel, Kämmerer, Pozharska, Schäfer and Ullrich (2025) for a treatment that is close (in language) to this work.

The inequality (3.25), together with Proposition 3.7 with $Y = B(D)$, imply the following result which improves upon Theorem 3.9 due to the *oversampling*. With worse constants, Theorem 3.16 first appeared in Krieg *et al.* (2026*c*); see Chkifa *et al.* (2026, Corollary 16) for the version as below.

**Theorem 3.16.** Let $D$ be a set. Then, for every $m$-dimensional subspace $V_m \subset B = B(D)$ and $n > m$, there exist $\mathcal{P}_n = \{x_1,\dots,x_n\} \subset D$ such that $A_{\mathcal{P}_n}$ from (3.17) satisfies

$$\|f - A_{\mathcal{P}_n}(f)\|_\infty \;\le\; \left(1 + \frac{\sqrt{m}}{1-\sqrt{m/n}}\right)\cdot d(f, V_m)_\infty$$

for all $f \in B$. In particular, for $F \subset B(D)$ and $n \ge 3$, we have

$$g_{2n}^{\mathrm{lin}}(F,B) \;\le\; 4\sqrt{n}\cdot d_n(F,B). \tag{3.26}$$

The case of $L_p$-approximation now follows immediately from the above consid-erations. First, note that $\|g\|_p \le \|g\|_2$ for $p \le 2$ and any probability space $(D, \mu)$. For $2 < p < \infty$, we simply use that, for point sets $\mathcal{P}_n^1$ and $\mathcal{P}_n^2$ satisfying (3.22) and (3.25), respectively, we can take its union $\mathcal{P}_n := \mathcal{P}_n^1 \cup \mathcal{P}_n^2 = \{x_1, \ldots, x_{2n}\}$ so that

$$\left(1 - \sqrt{\frac{m}{n}}\right) \cdot \|g\|_p \;\le\; m^{\frac{1}{2}-\frac{1}{p}} \cdot \sqrt{\frac{1}{n} \sum_{i=1}^{2n} |g(x_i)|^2} \tag{3.27}$$

for all $g \in V_m$. We obtain the following result of Krieg *et al.* (2026c, Theorem 20) with improved constants of Chkifa *et al.* (2026).

**Theorem 3.17.** Let $D$ be a set, $(D, \mu)$ be a probability space and $1 \le p \le \infty$. Then, for every $m$-dimensional subspace $V_m \subset L_2 \cap B(D)$ and every $n > m$, there exists $\mathcal{P}_{2n} = \{x_1, \ldots, x_{2n}\} \subset D$ such that $A_{\mathcal{P}_{2n}}$ from (3.17) satisfies

$$\|f - A_{\mathcal{P}_{2n}}(f)\|_p \;\le\; \left(1 + \frac{\sqrt{2}}{1 - \sqrt{m/n}} \cdot m^{(1/2-1/p)_+}\right) \cdot d(f, V_m)_\infty$$

for all $f \in L_2 \cap B$. In particular, for $F \subset L_2 \cap B$, we have

$$g_{4n}^{\mathrm{lin}}(F, L_p) \;\le\; 6 \cdot n^{(1/2-1/p)_+} \cdot d_n(F, B) \tag{3.28}$$

with $(a)_+ := \max\{0, a\}$.

Also the results for $L_p$-approximation are sharp in the sense that the factor $n^{(1/2-1/p)_+}$ can in general not be replaced by a lower order term.

**Example 3.18.** For $p \ge 2$, consider the unit ball $F_q^s$ of the univariate Sobolev spaces from (2.11), with $q = 1$, for which

$$g_n(F_1^s, L_p) \asymp n^{-s+1-1/p} \quad \text{and} \quad d_n(F_1^s, B) \asymp n^{-s+1/2}$$

for $s \in \mathbb{N}$ with $s > 1$, see Appendix A.2. For $p \le 2$, we consider $F_\infty^s$, $s \in \mathbb{N}$, in which case we have

$$g_n(F_\infty^s, L_p) \asymp d_n(F_\infty^s, B) \asymp n^{-s}.$$

Hence, (3.28) is sharp for all $1 \le p \le \infty$.

**Remark 3.19 (Superlinear oversampling).** As we have seen, the factors in the various error bounds and discretization inequalities could be reduced a lot if instead of $m = \dim(V_m)$ samples, we use $n = 2m$ samples. For example, the factor $(n + 1)$ from (3.20) is reduced to $4\sqrt{n}$ in (3.26). It is a natural question whether the bounds can be improved further by allowing a larger amount of samples, for example, $42m$ points or even $m^{42}$ points. However, in general, the order $\sqrt{m}$ cannot be improved with any polynomial amount of samples. A corresponding bound would lead to a contradiction for the Sobolev classes $F_1^s$ with fractional smoothness $s > 1$ close to one; see Krieg *et al.* (2026c, Section 5.3) for details.

A similar phenomenon occurs for $L_p$-approximation. So the factor $\sqrt{m}$ can only be improved with a superpolynomial amount of samples. And indeed, Kashin *et al.* (2023, Theorem 5.2) show that it can be reduced to a constant if we use exponentially many samples, see also Remark 3.10.

### *3.3. Bounds by best $L_2$-approximation*

The results of the last subsection are quite general, and *local* in the sense that we are free to choose the subspace $V_m$, and are therefore in a sense *easy-to-use*. However, the error bounds for $L_p$-norms by best approximation in the sup-norm are a somewhat odd comparison. Although they lead to quite sharp results for some classes, they do not for others. In addition, for a given model class $F$, the good subspaces $V_m$ for approximation in the sup-norm are often unknown, while it is usually much easier to find good subspaces for $L_2$-approximation. And even if they were known, uniform approximation is often "much harder" than approximation in $L_2$, especially when it comes to spaces of smooth functions in high dimensions, see Section 8.

Ideally, we would like to compare the $L_2$-error of least squares (or another method) directly with the best approximation in $L_2$, but this is in general not possible, as noted in the discussion after Proposition 3.3. However, and this will be the main scheme of this section, one can derive something close to that if there exists a whole sequence of "good" subspaces. This alternative line of research has been initiated by Krieg and Ullrich (2021*a*), with inspiration taken from Cohen and Migliorati (2017) and Hinrichs, Krieg, Novak, Prochno and Ullrich (2021*b*), which will be discussed in Section 7 and Remark 6.15, respectively. We will say more on the history of the presented results below in Remark 3.27.

The necessity of working with a whole sequence of approximation spaces makes the analysis more complicated. In order to provide motivation for the reader to deal with the necessary theory and technicalities, we first present a corollary of our analysis, which we nonetheless call a theorem and which we hope to be easily accessible. Roughly speaking, it says that, for any nested sequence of subspaces $(V_m)_{m\in\mathbb{N}}$ of $L_2$, which are 'good enough' for approximating the target function $f$, we find weighted least squares algorithms whose $L_2$-errors are of the same magnitude as the error of **best $L_2$-approximation** on the spaces $V_m$, which is defined by

$$d(f,V_m)_2 := \inf_{g\in V_m} \|f-g\|_2. \tag{3.29}$$

By a nested sequence of subspaces of $L_2$, we mean that $V_m$ is an $m$-dimensional space of square-integrable functions and $V_m \subset V_{m+1}$ for all $m \in \mathbb{N}$. Such a sequence gives rise to a unique orthonormal system $\{b_k\}_{k\in\mathbb{N}_0}$ in $L_2$ such that $V_m = \operatorname{span}\{b_k \colon 0 \le k < m\}$, and of course, vice versa, any such orthonormal sequence defines a nested sequence of subspaces. We say that $f \in L_2$ has a **pointwise convergent Fourier series** if

$$f(x) = \sum_{k\in\mathbb{N}_0} \langle f, b_k\rangle_2 \, b_k(x) \quad \text{for all} \quad x \in D. \tag{3.30}$$

The following is a special case of Proposition 3.30. It follows from results of Chkifa *et al.* (2026), and is an improvement of (Dolbeault, Krieg and Ullrich 2023) in the sense of better constants. See also Remarks 3.23 and 3.27 for more comments on the history. We will discuss in Section 4 when and how the points and weights of the weighted least squares algorithms can be constructed.

**Theorem 3.20.** Let $(D,\mu)$ be a measure space, let $(V_m)_{m\in\mathbb{N}}$ be a nested sequence of subspaces of $L_2(D,\mu)$ with $\dim(V_m) = m$, and let $\varepsilon > 0$. Then, for each $m \in \mathbb{N}$, there exist $4m$ points $\mathcal{P} \subset D$ and weights $w \in \mathbb{R}_+^{4m}$ such that the weighted least squares algorithm $A_m := A^w_{\mathcal{P},V_{2m}}$ from (3.2) satisfies

$$\|f - A_m f\|_2^2 \le C_\varepsilon \cdot m^{-1-\varepsilon} \sum_{k>m} k^\varepsilon \, d(f, V_k)_2^2$$

for any $f \in L_2$ with pointwise convergent Fourier series (3.30). Here, $C_\varepsilon > 0$ is a constant that solely depends on $\varepsilon$. In particular, for any $\alpha > 1/2 + \varepsilon/2$, we have

$$d(f, V_m)_2 \lesssim m^{-\alpha} \quad \Longrightarrow \quad \|f - A_m f\|_2 \lesssim m^{-\alpha}.$$

In addition, if $\mu$ is a probability measure and $V_m$ contains the constant functions, the weights can be chosen such that $\|w\|_1 \le C_\varepsilon$.

We point out that the points and weights, and hence the algorithm $A_m$ in Theorem 3.20, do not depend on the rate $\alpha$, only on the chosen sequence of approximation spaces and the error measure $\mu$, as well as the parameter $\varepsilon > 0$. The latter can be chosen freely and may be assumed fix (e.g., $\varepsilon = 0.01$). This means that the algorithm is a **universal algorithm** in the sense that its worst case error has the optimal rate of convergence for a whole family of function classes. In order to illustrate this property for an easy example, let us consider $L_2$-approximation on an interval $D$ with respect to the Lebesgue measure and the algorithm $A_m$ obtained from choosing $(V_m)_{m\in\mathbb{N}}$ as the spaces of polynomials of degree less than $m$; then the worst case error of $A_m$ has the optimal rate of convergence $m^{-k}$ for the classes $C^k(D)$ of $k$-times continuously differentiable functions for all $k \in \mathbb{N}$.

**Remark 3.21 (Noise).** By Remark 3.6, the algorithm from Theorem 3.20 responds well to noise in the function values if the noise is measured in the weighted Euclidean norm. The weights are the same as the weights of the algorithm and hence not at our disposal. This is one reason why a bounded sum of weights is of interest; it ensures that stability and robustness also holds in terms of ‘pointwise inaccuracies’, i.e., if the noise is measured in the (more natural) maximum norm.

The property will also be important for the actual computation of the sampling points, since it will allow us to truncate the infinite basis $\{b_k\}$.

We now discuss how to derive results like Theorem 3.20. The main difference to Section 3.2 is to bound the right-hand side of the error bound of Proposition 3.3, i.e., the discrete norm of $f-g$, from above by a different norm than the supremum norm. In fact, the norm $\|\cdot\|_H$ of a Hilbert space $H$ is used. So in addition to the one-sided discretization (3.12) on $V_m$, i.e., in addition to

$$\|g\|_Y \;\le\; K\sqrt{\sum_{i=1}^n w_i\,|g(x_i)|^2} \qquad \text{for all} \quad g\in V_m,$$

we want a reverse *discretization* inequality of the form

$$\sqrt{\sum_{i=1}^n w_i\,|h(x_i)|^2} \;\le\; \gamma\cdot\|h\|_H \tag{3.31}$$

for the functions of the form $h=f-g$, where $f\in F$ is the target function and $g\in V_m$ is some good approximation of $f$. The class of functions where (3.31) should hold is typically an infinite-dimensional object and it turns out that (3.31) is often the difficult part. A difference of (3.31) compared to bounding the discrete norm by the supremum norm is that the Hilbert space norm $\|h\|_H$ will be larger than the sup-norm, but in turn, this will allow for a much smaller factor $\gamma\in o(m^{-1/2})$.

The key result in the analysis is a result on $L_2$-approximation in reproducing kernel Hilbert spaces $H$ and for particular approximation spaces $V_m$. The results for general model classes $F$, general approximation spaces $(V_m)_{m\in\mathbb{N}}$, and more general error norms $Y\neq L_2$ will be developed thereafter by (comparably basic) embeddings and estimates.

### 3.3.1. *$L_2$-approximation on Hilbert spaces*

We start by discussing some basics of reproducing kernel Hilbert spaces which are embedded in some $L_2$-space. For a comprehensive treatment, we refer to Aronszajn (1950) and Steinwart and Scovel (2012).

Let $D$ be an arbitrary set, let $\mu$ be a measure on $D$ and let $L_2=L_2(D,\mu)$. A **reproducing kernel Hilbert space (RKHS)** $H$ on $D$ is a Hilbert space of complex-valued functions on $D$ such that function evaluation

$$\delta_x\colon H\to\mathbb{C},\qquad f\mapsto f(x)$$

at each point $x\in D$ is a continuous linear functional. By the Riesz representation theorem, this means that there exist functions $k_x\in H$ for all $x\in D$ such that

$$f(x)=\langle f,k_x\rangle_H \qquad \text{for all} \quad f\in H \quad \text{and} \quad x\in D. \tag{3.32}$$

We write

$$k_H(x, y) := \langle k_x, k_y \rangle_H, \qquad x, y \in D,$$

for the **reproducing kernel** of $H$. In order to be able to discuss $L_2$-approximation on $H$, we will of course assume that $H \subset L_2$, and moreover that we have an *injective* embedding

$$I\colon H \hookrightarrow L_2. \tag{3.33}$$

Here, $L_2$ is meant as the Hilbert space of equivalence classes of functions, i.e., injectivity means that, whenever two functions from $H$ are equal $\mu$-a.e., then they are also pointwise equal. In other words, $H$ is only allowed to contain one representer of each $\mu$-equivalence class. We are used to the injectivity assumption being included in the word 'embedding' and hence, in the following, this assumption is included when we say that $H$ is embedded into $L_2$.

For the following analysis, we will additionally assume that the the kernel $k_H$ has **finite trace**, that is,

$$\int_D k_H(x, x)\, \mathrm{d}\mu(x) < \infty. \tag{3.34}$$

This condition implies that the embedding (3.33) is a Hilbert-Schmidt operator (and in particular compact), see for instance Steinwart and Scovel (2012, Lemma 2.3). The converse is also true (i.e., if the embedding is Hilbert-Schmidt, then the kernel has finite trace), which will become clear soon.

Indeed, let us not yet assume (3.34) and only assume that $I$ is injective and compact. Then the operator $T := I^*I$ is a compact, injective and positive definite operator on the Hilbert space $H$. By the spectral theorem, this implies that there exists a countable or finite orthonormal basis $\{e_k\}_k$ of $H$ and corresponding positive eigenvalues $\lambda_k$, $0 \le k < \dim(H)$, i.e., $Te_k = \lambda_k e_k$, such that $(\lambda_k)_k$ is decreasing and tends to zero. A simple calculation gives

$$\langle e_j, e_k \rangle_2 = \langle Ie_j, Ie_k \rangle_2 = \langle e_j, I^*Ie_k \rangle_H = \lambda_k \delta_{k,j},$$

i.e., that the system $\{e_k\}_k$ is also orthogonal in $L_2$. The numbers $\sigma_k := \|e_k\|_2 = \lambda_k^{1/2}$ are called the **singular values** of the embedding $T$. Instead of $\{e_k\}_k$, we prefer to work with the $L_2$-normalized system, i.e., $b_k := \sigma_k^{-1} e_k$, which we call the **singular functions** of the embedding. For any function $f \in H$, we have

$$f = \sum_k \langle f, e_k \rangle_H\, e_k = \sum_k \langle f, b_k \rangle_2\, b_k$$

with convergence in $H$, and hence also pointwise and in $L_2$. In particular, the norm in $H$ can be written as

$$\|f\|_H^2 = \sum_k \left| \frac{\langle f, b_k \rangle_2}{\sigma_k} \right|^2.$$

By Steinwart and Scovel (2012, Theorem 3.1), in the situation of a compact

embedding (3.33), we have the pointwise **Mercer decomposition** of the kernel,

$$k_H(x,y) = \sum_k \sigma_k^2 \cdot b_k(x)\,\overline{b_k(y)} \qquad \text{for all} \quad x,y \in D. \tag{3.35}$$

From this, we also see that the kernel has finite trace if and only if the singular values are square-summable; the latter is a possible definition for an operator to be a Hilbert Schmidt operator.

Having finished with our little introduction to RKHSs, we now consider the *weighted least squares algorithm* $A^w_{\mathcal{P}_n}$ from (3.3), for the $m$-dimensional space

$$V_m \;:=\; \operatorname{span}\{b_k \colon k < m\} \tag{3.36}$$

that is spanned by the basis functions corresponding to the $m$ largest $\sigma_k$, where $m \le \dim(H)$. For $f \in L_2$, the best $L_2$-approximation from the space $V_m$ is given by the $L_2$-orthogonal projection

$$P_m f \;:=\; \sum_{k<m} \langle f, b_k\rangle_2\; b_k.$$

Due to the special form of the spaces $V_m$, the function $P_m f$ is also the $H$-orthogonal projection of $f \in H$ onto $V_m$, so it is also the best approximation of $f$ from $V_m$ in the $H$-norm.

Analogously to Proposition 3.3, we get the following conditional error bound, with a small improvement by using orthogonality instead of the triangle inequality. We denote with $V_m^\perp$ the orthogonal complement of $V_m$ in $H$. (For the particular $V_m$, it does not matter whether orthogonality is considered in $H$ or $L_2$.)

**Proposition 3.22.** Let $(D,\mu)$ be a measure space and let $H$ be a RKHS that is (injectively) embedded into $L_2(D,\mu)$ and whose kernel has a finite trace (3.34). As above, let $(\sigma_m)_{m\in\mathbb{N}_0}$ denote the singular values of the embedding and $(V_m)_{m\in\mathbb{N}}$ the corresponding eigenspaces. For some $m,n \in \mathbb{N}$, let $\mathcal{P}_n := \{x_1,\dots,x_n\} \subset D$ and $w := (w_1,\dots,w_n) \in \mathbb{R}^n_+$ be such that

$$\|g\|_2 \;\le\; K \sqrt{\sum_{i=1}^n w_i\, |g(x_i)|^2} \qquad \text{for all} \quad g \in V_m, \tag{3.37}$$

and

$$\sqrt{\sum_{i=1}^n w_i\, |h(x_i)|^2} \;\le\; \gamma\, \|h\|_H \qquad \text{for all} \quad h \in V_m^\perp. \tag{3.38}$$

Then, for all $f \in H$, the weighted least squares algorithm $A^w_{\mathcal{P}_n} := A^w_{\mathcal{P}_n,V_m}$ from (3.3) satisfies

$$\begin{aligned} \left\| f - A^w_{\mathcal{P}_n}(f) \right\|_2^2 \;&\le\; \|f - P_m f\|_2^2 \;+\; (K\gamma)^2\, \|f - P_m f\|_H^2 \\ &\le\; \left(\sigma_m^2 + K^2\gamma^2\right) \cdot \inf_{g\in V_m} \|f-g\|_H^2 \end{aligned}$$

*Proof.* Since $P_m f$ and $A^w_{\mathcal{P}_n}(f)$ are in $V_m$ and $f - P_m f$ is orthogonal to $V_m$, we have

$$\|f - A^w_{\mathcal{P}_n} f\|_2^2 = \|f - P_m f\|_2^2 + \|P_m f - A^w_{\mathcal{P}_n} f\|_2^2.$$

By Lemma 3.2, the second summand is bounded as

$$\|P_m f - A^w_{\mathcal{P}_n} f\|_2^2 = \|A^w_{\mathcal{P}_n} f - A^w_{\mathcal{P}_n} P_m f\|_2^2 \le K^2 \sum_{i=1}^n w_i \, |(f - P_m f)(x_i)|^2 ,$$

so that (3.38) implies the first bound. The second bound is obtained since it holds for all $h \in V_m^\perp$ that

$$\|h\|_2 \le \sigma_m \|h\|_H.$$

□

So, the main task is to find points and weights such that the two discretization conditions (3.37) and (3.38) are satisfied with factors $K, \gamma > 0$ as small as possible. The latest result in this direction is from Chkifa *et al.* (2026, Corollary 6). See Section 4.1 for a detailed discussion. It implies that, for any $m \in \mathbb{N}$ (with $m \le \dim(H)$) and $n \ge m$, there exist points $x_1, \ldots, x_n \in D$ and weights $w_1, \ldots, w_n > 0$ such that (3.37) and (3.38) hold simultaneously with

$$K = \left(1 - \sqrt{\frac{m-1}{n}}\right)^{-1} \quad \text{and} \quad \gamma = \sigma_m + \sqrt{\frac{1}{n} \sum_{k \ge m} \sigma_k^2}. \tag{3.39}$$

In particular, if we choose the number $n$ of points linear in the dimension $m$ of the space $V_m$, like $n = 2m$, we see that

$$K \asymp 1 \quad \text{and} \quad \gamma \asymp \sigma_m + \sqrt{\frac{1}{m} \sum_{k \ge m} \sigma_k^2}. \tag{3.40}$$

In fact, under the additional assumption that $\mu$ is a probability measure and $1 \in V_m$, we can additionally ensure a bounded sum of the weights, namely

$$\sum_{i=1}^n w_i \le \left(1 + \sqrt{\frac{1}{n\sigma_m^2} \sum_{k \ge m} \sigma_k^2}\right)^2 = \frac{\gamma^2}{\sigma_m^2}, \tag{3.41}$$

see Chkifa *et al.* (2026, Proposition 7) and note that under the assumption $1 \in V_m$ all functions in $V_m^\perp$ have zero integral.

**Remark 3.23 (History).** A similar discretization inequality has already been obtained and used (implicitly) by Krieg and Ullrich (2021*a*), where the authors used probabilistic arguments to show that (3.40) is achieved with high probability by $n \asymp m \log(m)$ random points (taken i.i.d. with respect to a distribution tailored especially to this problem). We will discuss that in more detail in Section 6.2.

The logarithmic oversampling has been resolved by Dolbeault *et al.* (2023), at the expense of an enormous constant, using a (non-constructive) approach based on work of Marcus *et al.* (2015) together with a technique of Nitzan, Olevskii and Ulanovskii (2016). Even before that, another paper by Nagel, Schäfer and Ullrich (2022) had removed the logarithmic oversampling factor at the expense of an additional logarithmic factor in the discretization bound. These authors already used the results of Batson *et al.* (2014) in their original (symmetric) form, see Proposition 4.1 below. The final result (3.39) from Chkifa *et al.* (2026) also uses the framework of Batson *et al.* (2014), but in a non-symmetric (and in a certain sense dimension-independent) form, see Theorem 4.2. In contrast to the results based on Marcus *et al.* (2015), the results based on Batson *et al.* (2014) are, to some extent, constructive. This is discussed further in Section 4.

Together with Proposition 3.22, the discretization result (3.39) gives the following error bound, see Chkifa *et al.* (2026, Corollaries 10 & 11) which may be regarded as the key result of this section.

**Theorem 3.24.** Let $(D, \mu)$ be a measure space and let $H$ be a RKHS that is (injectively) embedded into $L_2(D, \mu)$ and whose kernel has a finite trace (3.34). Moreover, let $V_m$ as in (3.36), i.e., the eigenspace corresponding to the $m$ largest singular values $\sigma_k$. Then, for all $m \in \mathbb{N}$ and $n \geq m$, there exist $\mathcal{P}_n = \{x_1, \ldots, x_n\} \subset D$ and $w := (w_1, \ldots, w_n) \in \mathbb{R}_+^n$ such that

$$\|f - A^w_{\mathcal{P}_n}(f)\|_2 \leq \frac{\sqrt{2}}{1-r}\left(\sigma_m + \sqrt{\frac{1}{n}\sum_{k\geq m}\sigma_k^2}\right) \cdot \inf_{g\in V_m} \|f - g\|_H, \tag{3.42}$$

for any function $f \in H$, where $r = \sqrt{(m-1)/n}$.

If $\mu$ is a probability measure and $1 \in V_m$, we can additionally ensure a bounded sum of the weights as in (3.41).

*Proof.* By (3.39), we can apply Proposition 3.22 with

$$K = \left(1 - \sqrt{\frac{m-1}{n}}\right)^{-1} \quad \text{and} \quad \gamma = \sigma_m + \sqrt{\frac{1}{n}\sum_{k\geq m}\sigma_k^2},$$

together with the rough estimate $\sigma_m^2 \leq K^2\gamma^2$. □

Theorem 3.24 implies bounds on the sampling numbers of model classes $F = B_H$, which are the unit ball of a RKHS $H$. In fact, using $g = 0$ above gives an upper bound that depends on $f$ only via $\|f\|_H$. We explicitly state the corresponding bounds for $n = 2m$ and $n = m$, where we obtain

$$g_{2m}^{\mathrm{lin}}(B_H, L_2) \le 5 \left( \sigma_m + \sqrt{\frac{1}{2m} \sum_{k \ge m} \sigma_k^2} \right) \tag{3.43}$$

and

$$g_{m}^{\mathrm{lin}}(B_H, L_2) \le 3m \left( \sigma_m + \sqrt{\frac{1}{m} \sum_{k \ge m} \sigma_k^2} \right) \tag{3.44}$$

in terms of the singular numbers $(\sigma_m)_{m\in\mathbb{N}_0}$ of the embedding $H \hookrightarrow L_2$, respectively. The second bound (which corresponds to interpolatory algorithms) may be better than (3.43) if the $\sigma_m$ decay very fast, like for Hilbert spaces with a Gaussian kernel. In this case, the oversampling factor two has a worse effect than the extra factor $m$ in the error bound. In the case of polynomial decay, (3.43) is sharp in an asymptotic sense, as can be seen again by specific examples.

In fact, the upper bound (3.43) cannot be improved apart from the precise value of the constants 2 and 5. Indeed, for any monotone and square-summable sequence $(\sigma_m)_{m\in\mathbb{N}_0}$, there exists a RKHS $H$ with corresponding singular numbers and for which

$$g_{m/8}^{\mathrm{lin}}(B_H, L_2) \ge \frac{1}{8} \left( \sigma_m + \sqrt{\frac{1}{m} \sum_{k \ge m} \sigma_k^2} \right) \tag{3.45}$$

holds for all sufficiently large $m$, namely, for all $m \in \mathbb{N}$ for which the expression in brackets is smaller than $\sigma_0$; see Hinrichs, Krieg, Novak and Vybíral (2022, Theorem 2) and Krieg and Vybíral (2023, Section 3.3). The lower bound is fulfilled for Sobolev spaces (of generalized smoothness) on an interval.

The following example of a tensor product problem is of special value for us. It was our original motivation for developing the theory of this section and hence can be seen as the 'historical origin' of the corresponding results. It shall also serve as an advertisement for Theorem 3.24 since for this example, the generic theorem gives better upper bounds than any other known bounds for sampling algorithms, even those that have been tailored specifically to tensor product problems.

**Example 3.25 (Tensor product spaces).** Let $H$ be a RKHS on $D$ that is compactly embedded into $L_2(D, \mu)$ and let $F$ be its unit ball. We consider $L_2$-approximation on the unit ball $F_d$ of the $d$-fold tensor product $H_d$ of $H$, which is a RKHS on the domain $D^d$. We assume that $g_n(F, L_2) \lesssim n^{-\alpha}$ for some $\alpha > 0$. The famous Smolyak algorithm, see Smolyak (1963), gives the estimate

$$g_n^{\mathrm{lin}}(F_d, L_2) \lesssim n^{-\alpha} \log^{(\alpha+1)(d-1)} n. \tag{3.46}$$

Examples of such tensor product spaces are the **Sobolev spaces of (dominating) mixed smoothness** $\alpha > 1/2$, denoted by $H_d = \mathbf{W}_2^\alpha$, which are obtained if $D = [0, 1]$ and $H = W_2^\alpha$ is the Sobolev space from (2.11), see Dũng, Temlyakov and Ullrich (2018). For these spaces, it is known that the error bound (3.46) for the Smolyak algorithm can be improved, see Sickel and Ullrich (2007); the exponent of the logarithm can be reduced to $(\alpha+1/2)(d-1)$. It is also shown that no further improvements are possible for Smolyak's algorithm, which led to the conjecture that this rate is optimal, see (Dũng *et al.* 2018, Conjecture 5.26) and *Outstanding Open Problem 1.4* therein.

Due to Babenko (1960) and Mityagin (1962), it is known that $d_n(F_d, L_2) \asymp n^{-\alpha} \log^{\alpha(d-1)} n$. In fact, this holds for arbitrary tensor products under the assumption that $d_n(F, L_2) \asymp n^{-\alpha}$, see e.g. Krieg (2018, Theorem 1). With Theorem 3.24, particularly (3.43), we now obtain

$$g_n^{\mathrm{lin}}(F_d, L_2) \lesssim n^{-\alpha} \log^{\alpha(d-1)} n \quad \text{if } \alpha > 1/2, \tag{3.47}$$

which was first proved by Dolbeault *et al.* (2023). This bound is asymptotically optimal, since $g_n^{\mathrm{lin}} \geq d_n$.

By now, we have an algorithm that can construct sampling points that achieve the optimal rate of convergence (and hence beats the sparse grid construction) in polynomial time, see Section 4. However, we still do not know any explicit expression for such sampling points and it seems that there should be a "simpler" geometric or algebraic construction for this classical problem.

**Example 3.26 (Exponential convergence).** In the case of exponential error rates, (3.44) gives better results than (3.43). Indeed, whenever the singular values $(\sigma_m)_{m\in\mathbb{N}_0}$ of the embedding $H \hookrightarrow L_2$ satisfy

$$\sigma_m \lesssim \omega^m \tag{3.48}$$

for some $\omega \in (0,1)$, then (3.43) implies that

$$g_m^{\mathrm{lin}}(B_H, L_2(\mu)) \lesssim m \cdot \omega^m, \tag{3.49}$$

i.e., the sampling numbers decay exponentially with the same base as the singular numbers. Such exponential rates occur, for example, for Hilbert spaces with a Gaussian kernel or an Hermite kernel, see, e.g., Buhmann (2003), Wendland (2005), Rasmussen and Williams (2006), Rieger and Zwicknagl (2010), Fasshauer, Hickernell and Woźniakowski (2012), Irrgeher, Kritzer, Pillichshammer and Woźniakowski (2016). The easiest example is probably given by the Gaussian kernel on the real line,

$$k_H(x,y) = \exp\left(-\frac{\gamma^2}{2}(x-y)^2\right), \qquad x, y \in \mathbb{R},$$

with shape parameter $\gamma > 0$ and approximation in $L_2(\mu)$ where $\mu$ the Gaussian probability measure with Lebesgue density $\pi^{-1/2} e^{-x^2}$. In this case, the singular values of the embedding $H \hookrightarrow L_2(\mu)$ are known precisely; it holds that

$$\sigma_m = \sqrt{1-\omega^2} \cdot \omega^m, \qquad \text{where} \quad \omega = \frac{\gamma}{\sqrt{1+\gamma^2+\sqrt{1+2\gamma^2}}},$$

see, e.g., Rasmussen and Williams (2006, Section 4.3.1)

**Remark 3.27 (History II).** To the best of our knowledge, the first general upper bound for the sampling numbers of Hilbert classes in $L_2$ in terms of the singular values has been obtained by Wasilkowski and Woźniakowski (2001). There it was proved for all $n, m \in \mathbb{N}_0$ that

$$g_n^{\mathrm{lin}}(B_H, L_2) \le \sigma_m + \sqrt{\frac{m}{n} \sum_{k=0}^{\infty} \sigma_k^2}. \tag{3.50}$$

Note that the singular values are equal to the 'minimal worst-case error using arbitrary linear information', which is the way the result is formulated in the aforementioned paper. This result has been improved by Kuo, Wasilkowski and Woźniakowski (2008), who used multilevel algorithms to show that a polynomial decay of the singular values of order $\alpha > 1/2$ implies a polynomial decay of the sampling numbers of order at least $\frac{2\alpha}{2\alpha+1} \cdot \alpha$. The subsequent improvements are based on Proposition 3.22 and their history was already outlined in Remark 3.23.

All these upper bounds require the square-summability of the singular numbers, that is, the finiteness of the trace (3.34). In this case, non-trivial lower bounds can be proven based on the technique of decomposable kernels from Novak and Woźniakowski (2001), see Novak and Woźniakowski (2010, 2012) for further results, or based on a variant of the Schur product theorem from Vybíral (2020), see Krieg and Vybíral (2023) for a survey of this technique; all these lower bounds have in common that they even hold for the (simpler) problem of numerical integration. Without the finite trace assumption, it has been shown by Hinrichs, Novak and Vybíral (2008) that there cannot be any general upper bound on the sampling numbers in terms of the singular numbers.

**Remark 3.28 (Kernel interpolation).** If we fix the sampling points $x_1, \ldots, x_n$, it is known that the smallest possible worst case error is achieved by the *spline algorithm*

$$S_n(f) := \underset{g \in H:\, g(x_i)=f(x_i)}{\operatorname{argmin}} \|g\|_H,$$

that is,

$$\inf_{\varphi_1,\ldots,\varphi_n \in L_2} \sup_{f \in B_H} \Big\| f - \sum_{i=1}^n f(x_i)\, \varphi_i \Big\|_{L_2} = \sup_{f \in B_H} \Big\| f - S_n(f) \Big\|_{L_2},$$

see e.g. Traub and Woźniakowski (1980, Theorem 5.1). The function $S_n(f)$ is also known as the minimal norm interpolant and, by the famous *representer theorem*, can be expressed as a linear combination of the kernel functions $k(\cdot, x_i)$, see for example Wainwright (2019, Proposition 12.32). Namely, we have $S_n(f) = \sum_{j=1}^n c_j k(\cdot, x_j)$ where the coefficients solve the linear system

$$\sum_{j=1}^{n} c_j \, k(x_i, x_j) \; = \; f(x_i) \quad \text{for all} \quad i \le n.$$

Hence, in the context of RKHSs, the spline algorithm is also referred to as *kernel interpolation*. Therefore, our upper bounds are true not only for the least squares algorithm, but also for the kernel-based approximation $S_n(f)$. Both types of algorithms are common in the context of learning, see for example the seminal paper of Cucker and Smale (2002).

### *3.3.2. The non-Hilbert setting*

We now discuss how Theorem 3.24 can be used in a non-Hilbert setting and for approximation in other norms than $L_2$. We start by dropping the Hilbertian assumption on $F$, which was first done by Krieg and Ullrich (2021*b*). The key idea, which was properly worked out in (Dolbeault *et al.* 2023, Section 6), is to construct an 'artificial' Hilbert space that acts as a go-between.

**Lemma 3.29.** Let $(b_k)_{k\in\mathbb{N}_0}$ be an orthonormal system in $L_2(D, \mu)$ and let $(\sigma_k)_{k\in\mathbb{N}_0}$ be a non-negative, monotone, square-summable sequence. Then, there exists a set $D_0 \subset D$ with $\mu(D \setminus D_0) = 0$ such that the set $H$ of functions $f \in L_2(D_0, \mu)$ with pointwise convergent Fourier series and

$$\|f\|_H^2 \; := \; \sum_{k=0}^{\infty} \left| \frac{\langle f, b_k \rangle_2}{\sigma_k} \right|^2 \; < \; \infty$$

is a RKHS on $D_0$. The embedding $H \hookrightarrow L_2(D_0, \mu)$ is injective, the families $(b_k)_{k\in\mathbb{N}_0}$ and $(\sigma_k)_{k\in\mathbb{N}_0}$ are its singular functions and singular values and the kernel of $H$ has finite trace. Moreover, for any function $f \in L_2(D_0, \mu)$ with pointwise convergent Fourier series and $m \in \mathbb{N}$, it holds that

$$\|f - P_{2m} f\|_H \; \le \; \sqrt{2 \sum_{k=m+1}^{\infty} \frac{d(f, V_k)_2^2}{k \cdot \sigma_{4k}^2}}\,, \tag{3.51}$$

where, as above, $V_k = \operatorname{span}\{b_j \colon 0 \le j < k\}$ and $P_{2m}$ is the $L_2$-orthogonal projection onto $V_{2m}$. In particular, $f \in H$ whenever the above expression is finite.

*Proof.* By assumption, the function

$$k(x, x) \; := \; \sum_{k\in\mathbb{N}_0} \sigma_k^2 \cdot |b_k(x)|^2$$

has a finite integral. Hence, there exists a set $D_0$ with $\mu(D \setminus D_0) = 0$ such that

$$k(x,x) < \infty \quad \text{for all} \quad x \in D_0.$$

Under this condition, it is well-known that $H$, as defined above, is a RKHS with (pointwise convergent) reproducing kernel

$$k(x,y) := \sum_{k\in\mathbb{N}_0} \sigma_k^2 b_k(x)\overline{b_k(y)},$$

see, for instance, Steinwart and Scovel (2012, Lemma 2.6). The square-summability implies the finite trace and the pointwise convergence of the Fourier series implies the injectivity of the embedding. So it only remains to prove the estimate (3.51). For this, we compute

$$\begin{aligned}
\|f - P_{2m}f\|_H^2 &= \sum_{k=2m}^{\infty} \left|\frac{\langle f, b_k\rangle_2}{\sigma_k}\right|^2 = \sum_{\ell=1}^{\infty} \sum_{k=2^\ell m}^{2^{\ell+1}m-1} \left|\frac{\langle f, b_k\rangle_2}{\sigma_k}\right|^2 \\
&\leq \sum_{\ell=1}^{\infty} \frac{1}{\sigma_{2^{\ell+1}m-1}^2} \sum_{k=2^\ell m}^{2^{\ell+1}m-1} |\langle f, b_k\rangle_2|^2 \leq \sum_{\ell=1}^{\infty} \frac{d(f, V_{2^\ell m})_2^2}{\sigma_{2^{\ell+1}m-1}^2} \\
&\leq \sum_{\ell=1}^{\infty} \frac{1}{\sigma_{2^{\ell+1}m-1}^2 2^{\ell-1}m} \sum_{k=2^{\ell-1}m+1}^{2^\ell m} d(f, V_k)_2^2 \\
&\leq \sum_{\ell=1}^{\infty} \sum_{k=2^{\ell-1}m+1}^{2^\ell m} \frac{2d(f, V_k)_2^2}{k\sigma_{4k}^2} \leq \sum_{k=m+1}^{\infty} \frac{2d(f, V_k)_2^2}{k\sigma_{4k}^2}.
\end{aligned}$$

The finiteness of this expression implies, by definition of $H$, that $f - P_{2m}f \in H$, and since $P_{2m}f \in H$, hence also $f \in H$.

□

The introductory Theorem 3.20 is now proved by applying Theorem 3.24 to a Hilbert space as constructed in Lemma 3.29. We write down a more general (and more ugly) form of the result with more freedom in the choice of parameters.

**Proposition 3.30.** Let $(D,\mu)$ be a measure space, let $(V_m)_{m\in\mathbb{N}}$ be a nested sequence of subspaces of $L_2(D,\mu)$, and let $(\sigma_k)_{k\in\mathbb{N}_0}$ be a non-negative, monotone, square-summable sequence. Then, for each $m \in \mathbb{N}$, there exist $4m$ points $\mathcal{P} \subset D$ and weights $w \in \mathbb{R}_+^{4n}$ such that the weighted least squares algorithm $A_m = A_{\mathcal{P},V_{2m}}^w$ satisfies

$$\|f - A_m f\|_2 \leq 13\sqrt{\frac{1}{m}\sum_{\ell>m}\sigma_\ell^2 \cdot \sum_{k>m} \frac{d(f,V_k)_2^2}{k\cdot\sigma_{4k}^2}}$$

for any function $f \in L_2$ with pointwise convergent Fourier series (3.30). If $\mu$ is a probability measure and $1 \in V_{2m}$, we can additionally ensure that

$$\sum_{i=1}^{4m} w_i \;\leq\; \left(1 + \sqrt{\frac{1}{n\sigma_{2m}^2} \sum_{k \geq 2m} \sigma_k^2}\right)^2 .$$

*Proof.* Let $H$ be the RKHS on $D_0 \subset D$ from Lemma 3.29 with the basis induced by the subspaces $(V_m)_{m\in\mathbb{N}}$. Then, by Theorem 3.24 there exist $4m$ points $\mathcal{P} \subset D_0$ and weights $w \in \mathbb{R}_+^{4n}$ such that the weighted least squares algorithm $A_{4m} = A^w_{\mathcal{P},V_{2m}}$ satisfies for any $f \in H$,

$$\|f - A_m f\|_2 \;\leq\; 5\left(\sigma_{2m} + \sqrt{\frac{1}{2m} \sum_{k \geq 2m} \sigma_k^2}\right) \cdot \|f - P_{2m} f\|_H .$$

As in (3.41), the weights can be chosen to be bounded under the additional assumptions. We additionally simplify by

$$\sigma_{2m}^2 \;\leq\; \frac{1}{m} \sum_{k=m+1}^{2m} \sigma_k^2 \;\leq\; \frac{1}{m} \sum_{k>m} \sigma_k^2$$

to get

$$\|f - A_m f\|_2 \;\leq\; 9\sqrt{\frac{1}{m} \sum_{k>m} \sigma_k^2} \cdot \|f - P_{2m} f\|_H . \tag{3.52}$$

Now, if $f \in L_2$ has a pointwise convergent Fourier series, by Lemma 3.29, we have

$$\|f - P_{2m} f\|_H \;\leq\; \sqrt{2 \sum_{k=m+1}^{\infty} \frac{d(f,V_k)_2^2}{k\sigma_{4k}^2}} . \tag{3.53}$$

If this expression is infinite, the upper bound in the corollary is infinite as well and trivially true. It the expression is finite, we have $f \in H$ and so (3.52) and (3.53) can be combined to obtain the stated estimate, i.e.,

$$\|f - A_m f\|_2 \;\leq\; 9\sqrt{2}\sqrt{\frac{1}{m} \sum_{\ell>m} \sigma_\ell^2}\sqrt{\sum_{k>m} \frac{d(f,V_k)_2^2}{k\sigma_{4k}^2}} .$$

To be precise, the statement is obtained if all the functions are considered on $D_0$ instead of $D$. But none of the expressions in the aformentioned inequality (neither the output of the algorithm nor the $L_2$-norms and errors) change if we consider functions on $D$ instead of $D_0$.

This proves Proposition 3.30. The simplified form of the introductory Theorem 3.20 is obtained by considering $\sigma_k = k^{-1/2-\varepsilon/2}$ for $k \in \mathbb{N}$ (and $\sigma_0 = 1$). □

The sequence $(\sigma_k)_{k\in\mathbb{N}_0}$ in Proposition 3.30 is a sequence of parameters that can be chosen freely. How should we choose it? If nothing is known about the error of best approximation, it seems best to choose the sequence with a slow decay, so that it is just barely two-summable, like $\sigma_k = k^{-0.51}$, since then the upper bound is finite for as many functions $f$ as possible. As we have seen in Theorem 3.20, a faster decay of the error of best approximation is detected automatically by the algorithm and reflected in the error rate.

On the other hand, we are often interested in the situation that the target function comes from a given model class $F$ (like a smoothness class) and then there is an 'optimal' choice both for the sequence of approximation spaces $V_m$ and for the parameters $\sigma_k$, which we describe in the following. We remark already that this 'optimal' choice is usually not practical, however, it leads to the best known theoretical bounds for the sampling numbers.

We need the following minimal assumption on the model class $F$.

**Assumption A.** Let $D$ be a set and $\mu$ be a measure on $D$, and let $F$ be a class of complex-valued square-integrable functions on $D$. We assume that there exists a topology on $F$ such that:

- (A.1) $F$ is separable.
- (A.2) For every $x \in D$, function evaluation $f \mapsto f(x)$ is a continuous functional on $F$.
- (A.2) The identity mapping $F \to L_2(D, \mu)$, $f \mapsto f$, is continuous.

Although this assumption may look a little peculiar at first, one quickly realizes that it is a very weak assumption. For example, it is satisfied whenever

- $D$ is a compact topological space,
- $\mu$ is a finite Borel measure,
- $F$ is a compact subset of $C(D)$.

However, continuity of the functions, finiteness of the measure, or compactness of the domain is not really required, this is just an example. Another example would be any countable subset of $L_2(D, \mu)$; here the assumptions are satisfied if $F$ is equipped with the discrete metric.

For such classes $F$, the right choice of the approximation spaces and the parameter sequence leads to the following bound of the sampling numbers, which was obtained by Krieg, Pozharska, Ullrich and Ullrich (2025*b*); see also Dolbeault *et al.* (2023) for an earlier and only slightly less general version. The bound is in terms of the **Kolmogorov widths** of a class $F$ in $L_2$ which are defined by

$$d_n(F, L_2) := \inf_{V_m} \sup_{f \in F} d(f, V_m)_2, \tag{3.54}$$

where the infimum is over all $m$-dimensional subspaces $V_m \subset L_2$.

**Corollary 3.31.** There is an absolute constant $c \in \mathbb{N}$ such that, for all $D$, $\mu$ and $F$ that satisfy Assumption A, and for all $m \in \mathbb{N}$, we have

$$g_{cm}^{\mathrm{lin}}(F, L_2) \leq \frac{1}{\sqrt{m}} \sum_{k>m} \frac{d_k(F, L_2)}{\sqrt{k}}.$$

*Proof.* We only sketch the proof. For the details, we refer to (Krieg *et al.* 2025*b*, Section 2). We want to apply Proposition 3.30 for a good choice of the approximation spaces and the parameter sequence. It is natural to choose the sequence of approximation spaces such that $V_m^*$ attains or almost attains the infimum in (3.54). Such a sequence is not necessarily nested, of course, so we cannot apply Theorem 3.20 directly. However, we can simply replace $V_m^*$ by the space

$$V_m := \operatorname{span}\Bigg( \bigcup_{\ell < \lfloor \log_2(m) \rfloor} V_{2^\ell}^* \Bigg).$$

Then the sequence $(V_m)_{m\in\mathbb{N}}$ is nested and it is easy to see that it is still almost optimal in the sense that

$$\sup_{f\in F} d(f, V_{4m})_2 \leq (1+\varepsilon)\cdot d_m(F, L_2)$$

for some small $\varepsilon > 0$ and all $m$. In order to minimize the error bound, we now choose a decreasing parameter sequence $(\sigma_k)_{k\in\mathbb{N}_0}$ with

$$\sigma_{16k} := \sqrt{\frac{d_k(F, L_2)}{\sqrt{k}}}.$$

This sequence can be assumed to be square-summable; otherwise the error bound is infinite and hence trivial. For this parameter sequence, the error bound from Proposition 3.30 is such that both infinite sums are (essentially) equal; the error bound takes the form

$$\|f - A_m f\|_2 \leq \frac{c_2}{\sqrt{m}} \sum_{k > c_1 m} \frac{d_k(F, L_2)}{\sqrt{k}} \tag{3.55}$$

with some universal constants $0 < c_1 < 1 < c_2$. These constants can be swallowed by replacing the index $m$ with $cm$ for appropriate $c \in \mathbb{N}$.

There is one additional issue: We need to ensure that all functions from the class $F$ have a pointwise convergent Fourier series for the error bound to be valid for all $f \in F$. This is where Assumption A comes in. Namely, we fix some countable dense subset $F_0 \subset F$, which exists by the separability assumption. Assuming the finiteness of the aforementioned error bound, the Rademacher Menchov theorem can be used to conclude that pointwise convergence holds on some set $D_* \subset D$ with $\mu(D \setminus D_*) = 0$, simultaneously for every function $f \in F_0$. Thus, (3.55) is proved for all functions $f \in F_0$ (first only on the domain $D_*$, but then clearly also on the full domain $D$). Finally, the two continuity assumptions from Assumption A can be used to see that the error bound, in fact, holds on all $f \in F$. Indeed, the assumptions imply that the left hand side of (3.55) depends continuously on $f$ and so we can write the left hand side as a limit over functions from $F_0$.

□

In particular, Corollary 3.31 implies that for sufficient decay of the Kolmogorov numbers, the sequence of sampling numbers has the same decay. Indeed, if $d_n(F, L_2) \lesssim n^{-\alpha} \log^{\beta} n$ for some $\alpha > 1/2$ and $\beta \in \mathbb{R}$, then the corollary gives $g_n^{\mathrm{lin}}(F, L_2) \lesssim n^{-\alpha} \log^{\beta} n$.

**Remark 3.32 (Bound for Hilbert spaces).** Let us compare the bound from Corollary 3.31 with the Hilbert case. If $F$ is the unit ball of a RKHS as above, then we have $d_n(F, L_2) = \sigma_n$, see, e.g., Pietsch (1987) or Section 9.2. So in the Hilbert case we have the formula

$$g_{4m}^{\mathrm{lin}}(F, L_2) \;\leq\; 9 \sqrt{\frac{1}{m} \sum_{k>m} d_k(F, L_2)^2} \tag{3.56}$$

for all $m \in \mathbb{N}$. This follows from the bound (3.43), which can be simplified by using the estimate $\sigma_{2m}^2 \leq \frac{1}{m} \sum_{k=m+1}^{2m} \sigma_k^2$, which holds for all $m \in \mathbb{N}$ due to the monotonicity of the singular numbers. Up to constants, the bound (3.56) was established by Dolbeault *et al.* (2023). Inserting various examples, the reader will realize that the difference between the two bounds is rather insignificant, i.e., the result for general classes is almost as strong as the result for Hilbert spaces. Still, there is *some* difference. For example, there exist non-Hilbert classes $F$ that satisfy Assumption A with $d_n(F, L_2) \asymp n^{-1/2} \log^{-1} n$ for which it holds that $g_n(F, L_2) \asymp 1$, that is, the sampling numbers have no decay, see Dolbeault *et al.* (2023, Example 30). Such a thing cannot happen in the RKHS case; if $F$ were a Hilbert class, the aforementioned bound would imply $g_n(F, L_2) \lesssim n^{-1/2} \log^{-1/2} n$ in this case.

We now discuss how the above analysis can also be used to obtain error bounds in other norms $Y \neq L_2$. This is possible if $Y$ satisfies the following assumptions.

**Assumption B.** Let $(D, \mu)$ be a measure space and $L_2 = L_2(D, \mu)$ and let $Y$ be a seminormed space of complex-valued functions on $D$. Assume that

(B.1) $Y \cap L_2$ is complete w.r.t. the seminorm $\|\cdot\|_* := \|\cdot\|_Y + \|\cdot\|_2$.
(B.2) If $f, g \in Y$ and $f = g$ holds $\mu$-a.e., then $\|f\|_Y = \|g\|_Y$.

Clearly, the assumptions are satisfied for $Y = L_p(D, \mu)$ with $1 \leq p \leq \infty$, which is our main interest. Other possible choices for $Y$ can be the space of continuous functions $C(D)$ or the energy space $H^1(D)$ on compact domains with the Lebesgue measure.

The general idea, which was already used by Pozharska and Ullrich (2022) for $Y = B(D)$, is to lift the results from $L_2$-approximation to $Y$-approximation, using the quantity

$$\Lambda_m \;:=\; \Lambda(V_m, Y) \;:=\; \sup_{f \in V_m,\, f \neq 0} \frac{\|f\|_Y}{\|f\|_2}. \tag{3.57}$$

Depending on the choice of $V_m$, $G$ and $L_2(\mu)$, the $\Lambda_m$ correspond to the inverse of the Christoffel function, or the best constant in Bernstein or Nikol'skii inequalities.

Their appearance is no surprise, as it is an important quantity in approximation theory, especially for $G = B(D)$ and $V_m$ a space of polynomials, see e.g. Freud (1969), Nevai (1986) or Gröchenig (2020). It already turned out to be instrumental for discretization and sampling, also in more general contexts. We refer to the early work of Kuo, Wasilkowski and Woźniakowski (2009*a*), and to the more recent surveys of Kashin, Kosov, Limonova and Temlyakov (2022) and Temlyakov (2018*b*).

**Remark 3.33** ($Y = L_p$)**.** If $\mu$ is a finite measure and $V_m$ is an $m$-dimensional space of bounded and measurable functions, then the numbers $\Lambda(V_m, B)$ can be expressed in terms of an arbitrary $L_2$-orthonormal basis $\{b_k\}_{k=0}^{m-1}$ of $V_m$. We have

$$\Lambda(V_m, B) = \Big\| \sum_{k=0}^{m-1} |b_k|^2 \Big\|_\infty^{1/2}. \tag{3.58}$$

Moreover, we always have $\Lambda(V_n, B) \in \Omega(m^{1/2})$ since

$$\Lambda(V_m, B(D)) \geq \left( \frac{1}{\mu(D)} \int_D \sum_{k=0}^{m-1} |b_k|^2 \,\mathrm{d}\mu \right)^{1/2} = \sqrt{\frac{m}{\mu(D)}}.$$

The optimal behavior $\Lambda(V_m, B) \in \Theta(m^{1/2})$ occurs for spaces spanned by trigonometric monomials, Chebyshev polynomials and spherical harmonics if $\mu$ is their corresponding orthogonality measure. In addition, it is met for certain Walsh and (Haar) wavelet spaces. In other cases, $\Lambda(V_m, B)$ might be larger. For $2 \leq p \leq \infty$, an easy interpolation shows

$$\Lambda(V_m, L_p) \leq \Lambda(V_m, B)^{1-2/p} = \Big\| \sum_{k=0}^{m-1} |b_k|^2 \Big\|_\infty^{1/2-1/p}.$$

Results for $Y$-approximation are obtained from the results for $L_2$-approximation by the simple estimate

$$\begin{aligned} \|f - Af\|_Y &\leq \|f - P_m f\|_Y + \|P_m f - Af\|_Y \\ &\leq \|f - P_m f\|_Y + \Lambda_m \cdot \|P_m f - Af\|_2 \\ &\leq \|f - P_m f\|_Y + \Lambda_m \cdot \|f - Af\|_2, \end{aligned} \tag{3.59}$$

which holds for any $f \in Y \cap L_2$ and any $Af \in V_m$, as well as the following lemma.

**Lemma 3.34.** Let $(D, \mu)$ be a measure space and let $Y \subset \mathbb{C}^D$ be a seminormed space that satisfies Assumption B. Let $(V_m)_{m\in\mathbb{N}}$ be a nested sequence of subspaces of $Y \cap L_2(D, \mu)$. For any $m \in \mathbb{N}$ and $f \in Y \cap L_2$, we have

$$\|f - P_m f\|_Y \leq 2 \sum_{k > \lfloor m/4 \rfloor} \frac{\Lambda_{4k}}{k} \cdot d(f, V_k)_2$$

with $\Lambda_m$ from (3.57).

*Proof.* Note that there is nothing to prove if the right hand side of the inequality is infinite so we may assume that it is finite. We first show that $P_m f$ converges in $Y \cap L_2$. We define $I_\ell = \{k \in \mathbb{N} \colon 2^\ell < k \le 2^{\ell+1}\}$ for $\ell \in \mathbb{N}_0$. Let $m \in I_s$ and $n \in I_t$ with $s \le t$. Then

$$\begin{aligned} \|P_n f - P_m f\|_Y \;&\le\; \sum_{\ell=s}^{t} \Big\| \sum_{k \in I_\ell \colon m<k\le n} \langle f, b_k\rangle_2\, b_k \Big\|_Y \\ &\le \sum_{\ell=s}^{t} \Big\| \sum_{k\in I_\ell} \langle f, b_k\rangle_2\, b_k \Big\|_2 \cdot \Lambda_{2^{\ell+1}} \;\le\; \sum_{\ell\ge s} d(f, V_{2^\ell})\Lambda_{2^{\ell+1}}. \end{aligned}$$

Further we use the fact that the sequence $(d(f, V_k))_{k=1}^\infty$ is decreasing, and therefore

$$d(f, V_{2^\ell}) \le \frac{1}{2^{\ell-1}} \sum_{j\in I_{\ell-1}} d(f, V_j),$$

and hence, since the sequence $(\Lambda_k)_{k=1}^\infty$ is increasing,

$$\begin{aligned} \|P_n f - P_m f\|_Y \;&\le\; 2 \sum_{\ell\ge s} \sum_{j\in I_{\ell-1}} \frac{d(f, V_j)\Lambda_{2^{\ell+1}}}{2^\ell} \\ &\le 2 \sum_{\ell\ge s} \sum_{j\in I_{\ell-1}} \frac{d(f, V_j)\Lambda_{4j}}{j} \;\le\; 2 \sum_{j>\lfloor m/4\rfloor} \frac{d(f, V_j)\Lambda_{4j}}{j}, \end{aligned}$$

which is convergent due to the assumption. This means that $(P_m f)$ is a Cauchy sequence in $Y$. Since we also know that $P_m f$ converges to $f$ in $L_2$, we obtain that $(P_m f)$ is a Cauchy sequence in $Y \cap L_2$ and converges to some $g \in Y \cap L_2$ by Assumption B. In particular, $f = g$ almost everywhere. Thus, using again Assumption B, we obtain

$$\|f - P_m f\|_Y \;=\; \|g - P_m f\|_Y \;\le\; 2 \sum_{j>\lfloor m/4\rfloor} \frac{d(f, V_j)\Lambda_{4j}}{j},$$

which yields the statement.

□

Together, inequality (3.59) and Lemma 3.34 take the form

$$\|f - Af\|_Y \;\le\; 2 \sum_{k>\lfloor m/4\rfloor} \frac{\Lambda_{4k}}{k} \cdot d(f, V_k)_2 \;+\; \Lambda_m \cdot \|f - Af\|_2 \qquad (3.60)$$

for $f \in Y \cap L_2$ and $Af \in V_m$. It only remains to plug in your favorite $L_2$-approximation result. Various results obtained in this way are discussed in detail by Krieg *et al.* (2025*b*). It turns out that this seemingly crude method gives optimal results for surprisingly many examples. Let us present a reformulation in terms of the sampling numbers similar to Corollary 3.31, which also proves the result in combination with (3.60), see also (3.55).

**Theorem 3.35.** There are absolute constants $c, C \in \mathbb{N}$ such that, for all $D$, $\mu$ and $F$ that satisfy Assumption A, all seminormed spaces $Y \subset \mathbb{C}^D$ that satisfy Assumption B, all nested sequences $(V_m)_{m\in\mathbb{N}}$ of subspaces of $Y \cap L_2(D,\mu)$, and all $m \in \mathbb{N}$, we have

$$g_{cm}^{\mathrm{lin}}(F,Y) \;\leq\; C \cdot \sum_{k>m} \frac{d(F,V_k)_2}{\sqrt{k}} \cdot \max\left\{\frac{\Lambda_{4k}}{\sqrt{k}}, \frac{\Lambda_{4m}}{\sqrt{m}}\right\}$$

with $\Lambda_m = \Lambda(V_m, Y)$ from (3.57).

Let us finally present a simple and (hopefully) easy-to-use result in the spirit of the introductory Theorem 3.20.

**Corollary 3.36.** Let $(D,\mu)$ be a measure space and let $Y \subset \mathbb{C}^D$ be a seminormed space that satisfies Assumption B. Let $(V_m)_{m\in\mathbb{N}}$ be a nested sequence of subspaces of $Y \cap L_2(D,\mu)$ such that

$$\Lambda(V_m, Y) = \sup_{f\in V_m,\, f\neq 0} \frac{\|f\|_Y}{\|f\|_2} \;\lesssim\; m^{\beta}$$

for some $\beta \geq 0$, and let $\varepsilon > 0$. Then the weighted least squares algorithm $A_m$ from Theorem 3.20 satisfies

$$\|f - A_m f\|_Y \;\lesssim\; m^{-\alpha+\beta}$$

for any function $f \in Y \cap L_2$ with pointwise convergent Fourier series (3.30) and with $d(f, V_n)_2 \lesssim n^{-\alpha}$ for some $\alpha > \max\{\beta, 1/2 + \varepsilon/2\}$.

We point out that this means that the algorithm from Theorem 3.20 (see also Theorem 4.4 below) can be universal in two ways: One and the same algorithm is good for any approximation rate $\alpha$ (where 'approximation rate' can often be replaced by 'smoothness') and at the same time potentially also for various error norms $Y$.

As above in Example 3.25, the Sobolev spaces of mixed smoothness are of particular interest to us, as they show some interesting behavior and are subject to various open problems.

**Example 3.37 (Mixed smoothness).** For the spaces $\mathbf{W}_2^{\alpha}$ of mixed smoothness from Example 3.25, the least squares algorithm from Theorem 3.20 has the optimal rate of convergence for $L_p$-approximation for all $\alpha > 1/2$ and $1 \leq p \leq \infty$. We refer to Krieg *et al.* (2025*b*) for details. The spaces $\mathbf{W}_q^{\alpha}$ can also be considered for integrability parameters $q \neq 2$, in which case the discussion of optimal algorithms becomes more diverse. In some parameter regions, better bounds are obtained by the Smolyak algorithm, while in other regions, we do not know a matching lower

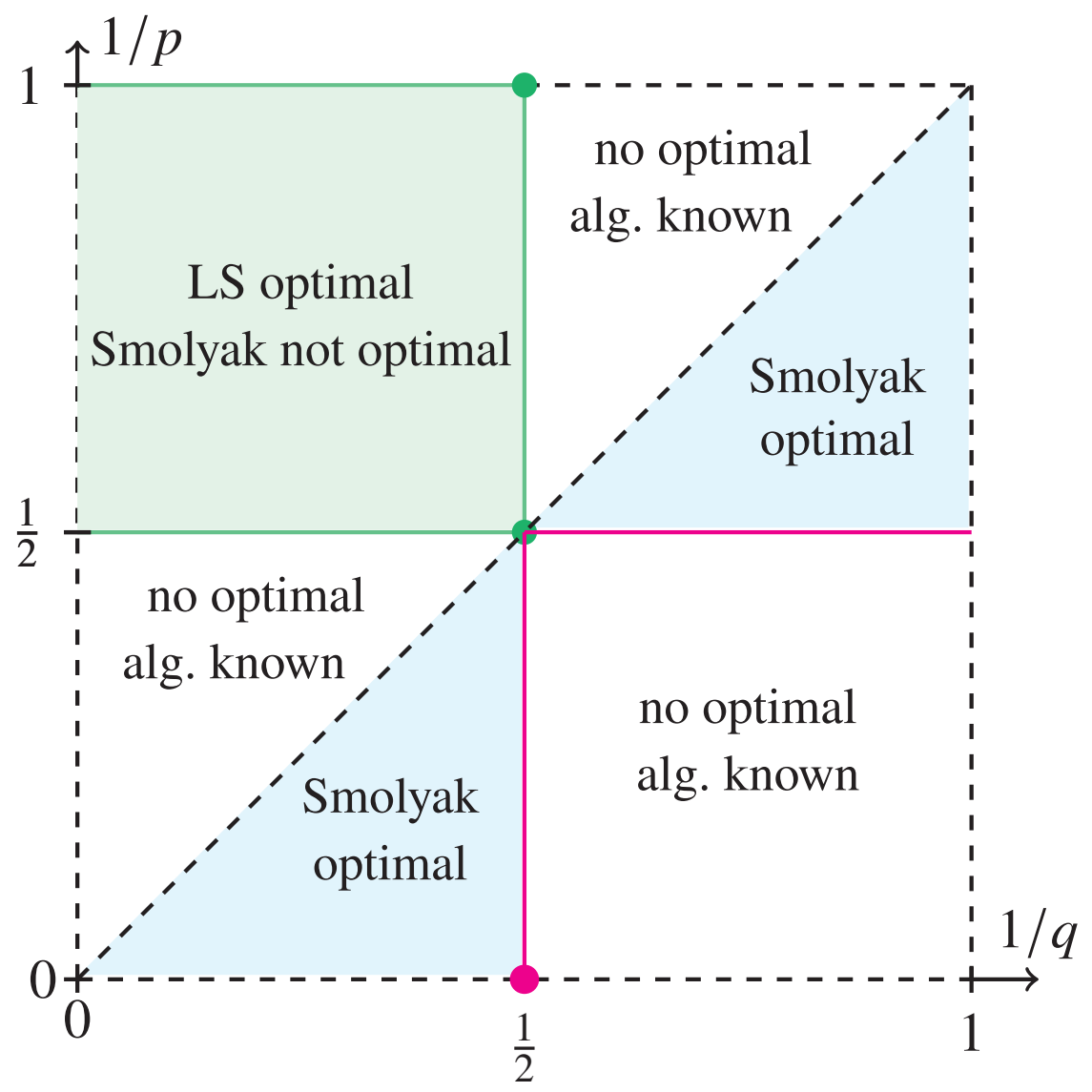

- Green area: Least squares is optimal, Smolyak is not optimal.
- Blue area: Smolyak is optimal, known LS bounds are not optimal.
- Magenta lines: Both algorithms are optimal.

Figure 3.1. Linear sampling algorithms for $L_p$-approximation on $\mathbf{W}_q^\alpha$.

bound and, possibly, none of the two algorithms gives optimal error bounds. The state of the art for linear sampling algorithms is summarized by Figure 3.1, see again Krieg *et al.* (2025*b*) for details and further remarks. We also refer to Temlyakov and Ullrich (2021, 2022) for the case of *small smoothness* $1/2 < \alpha < 1/p$.

## 4. Construction of sampling points

We now discuss whether and how one can actually obtain the points and weights behind the theorems of Section 3 numerically.

### 4.1. *Discretization*

The error bounds for least squares algorithms from Section 3 are all based on (different) discretization inequalities. We now present the algorithm that can be used to construct the corresponding points and weights. At the same time, we discuss a bit the origins of these results.

Discretization of (integral and uniform) norms has a long and rich history. We only mention results that hold for general linear spaces and are therefore relevant for the general theory from Section 3. We refer to the comprehensive treatments of Temlyakov (2018*b*), Gröchenig (2020), Kashin *et al.* (2022), and Limonova, Malykhin and Temlyakov (2024), and the references therein for more details.

A big part of the literature (until recently) is concerned with *two-sided discretization* which means that some continuous norm, usually an $L_2$-norm on a finite-dimensional space, is bounded above and below by the same discrete norm. For the $L_2$-norm, this problem was essentially solved by Batson *et al.* (2014), even in a constructive way. These authors had a different application in mind (namely, graph sparsification) and hence used a different language and worked under some inessential additional assumptions. Explicitly (up to the constants), a result of the following form was given by Dai, Prymak, Temlyakov and Tikhonov (2019, Theorem 2.13), with an additional assumption $V_m \subset L_4$ later removed by Dai, Prymak, Shadrin, Temlyakov and Tikhonov (2021, Theorem 6.4) and Limonova and Temlyakov (2022, Theorem 1.1). See Chkifa *et al.* (2026) for the version stated here with a direct proof along the lines of Batson *et al.* (2014).

**Proposition 4.1.** Let $(D, \mu)$ be a measure space, $V_m \subset L_2(D, \mu)$ an $m$-dimensional space of functions. Then, for any $n \geq m$, there exist points $x_1, \ldots, x_n \in D$ and weights $w_1, \ldots, w_n > 0$ such that we have

$$\left(1 - \sqrt{\frac{m-1}{n}}\right) \|f\|_2 \;\leq\; \sqrt{\sum_{i=1}^{n} w_i |f(x_i)|^2} \;\leq\; \left(1 + \sqrt{\frac{m-1}{n}}\right) \|f\|_2$$

for all $f \in V_m$.

In our context, especially for the results of Section 3.3, the following modifications of this result are particularly important:

- The spaces on which the lower and upper discretization bounds hold should be allowed to be different spaces.
- We want the upper bound to hold on possibly infinite-dimensional spaces and hence the factor $m$ on the right should be replaced with something dimension-independent.

- We want to discretize norms different from $L_2$.
- In Section 3.2, we want the weights to be equal, as this allows for the use of unweighted least squares methods.

All the required discretization bounds are implied by the following result of Chkifa *et al.* (2026), see Corollary 6 and Proposition 7 there. It is based on the points and weights constructed by Algorithm 1 that will be described below. The proof, as well as the algorithm, use the framework developed by Batson *et al.* (2014), but generalizes the latter in the aforementioned regards. Note that the modification for different upper and lower discretization spaces (in the finite-dimensional setting) has already been undertaken in the context of matrix reconstruction by Boutsidis, Drineas and Magdon-Ismail (2014).

**Theorem 4.2.** Let $(D, \mu)$ be a measure space, $V_m \subset L_2(D, \mu)$ an $m$-dimensional space of functions, and $H_m$ a reproducing kernel Hilbert space that is (injectively) embedded into $L_2(D, \mu)$ with

$$A \;:=\; \sup_{f \in H_m} \frac{\|f\|_2}{\|f\|_{H_m}}$$

and whose kernel $k_m$ has a finite trace

$$K \;:=\; \int_D k_m(x, x)\, d\mu(x).$$

Then, for any $n \geq m$, there exist points $x_1, \ldots, x_n \in D$ and weights $w_1, \ldots, w_n > 0$ such that we have

$$\left(1 - \sqrt{\frac{m-1}{n}}\right) \|g\|_2 \;\leq\; \sqrt{\sum_{i=1}^{n} w_i |g(x_i)|^2} \qquad \text{for all } g \in V_m$$

and

$$\sqrt{\sum_{i=1}^{n} w_i |f(x_i)|^2} \;\leq\; \left(A + \sqrt{\frac{K}{n}}\right) \|f\|_{H_m} \qquad \text{for all } f \in H_m.$$

These bounds are satisfied for any set of points and weights returned by Algorithm 1, where $a$ is an $L_2$-orthonormal basis of $V_m$, $b$ is an $H_m$-orthonormal basis of $H_m$, and the oracle is arbitrary.

Moreover, if $\mu$ is a probability measure and $\int_D f\, d\mu = 0$ for all $f \in H_m$, then we can add the constant function $A$ as first element of the basis $b$ in Algorithm 1, and additionally obtain

$$\sum_{i=1}^{n} w_i \;\leq\; \left(1 + \sqrt{\frac{K}{n \cdot A^2}}\right)^2.$$

**Algorithm 1** Construction of points and weights.

1: **Input:**

$$m \geq 2, \quad n \geq m, \quad N \in \mathbb{N} \cup \{\infty\},$$

$$a \in L_2(D,\mu)^m \text{ orthonormal}, \quad b \in L_2(D,\mu)^N.$$

We also assume we have an **oracle** that suggests sampling points.

2: **Set up:**

$$J = \int bb^* \, d\mu \in \mathbb{C}^{N \times N},$$

$$r = \sqrt{\frac{m-1}{n}}, \quad s = \sqrt{\frac{\mathrm{Tr}(J)/\lambda_{\max}(J) - 1}{n}}, \quad \delta = \frac{1-r}{n}, \quad \zeta = \frac{1+s}{n},$$

$$G = \frac{\delta m}{r} \cdot I_m, \quad \Gamma = -\frac{\zeta \cdot \mathrm{Tr}(J)}{s} \cdot I_N, \quad \ell = u = 0.$$

3: **for** $i = 1, \ldots, n$ **do**

4:

$$L = (G - \ell I)^{-1}, \quad U = (uJ - \Gamma)^{-1}.$$

$$L_\delta = (G - (\ell + \delta) I)^{-1}, \quad U_\zeta = ((u+\zeta)J - \Gamma)^{-1}.$$

$$V = \frac{L_\delta^2}{\mathrm{Tr}(L_\delta - L)} - L_\delta, \quad W = \frac{U_\zeta J U_\zeta}{\mathrm{Tr}(JU - JU_\zeta)} + U_\zeta.$$

If $\mathrm{Tr}(J)/\lambda_{\max}(J) < 1 + 1/n$, we instead put $W = n/\mathrm{Tr}(J) \cdot I_N$.

5: Suggest a sampling point $x_i \in D$ using the oracle.

6: **if**

$$a(x_i)^* V a(x_i) \;\geq\; b(x_i)^* W b(x_i),$$

**then**

7: Choose $w_i > 0$ such that $w_i^{-1}$ is between those two values and update

$$G \leftarrow G + w_i a(x_i) a(x_i)^*, \quad \Gamma \leftarrow \Gamma + w_i b(x_i) b(x_i)^*,$$

$$\ell \leftarrow \ell + \delta, \qquad u \leftarrow u + \zeta.$$

8: **else**

9: go to line 5.

10: **end if**

11: **end for**

12: **Output:** Points $x_1, \ldots, x_n \in D$ and weights $w_1, \ldots, w_n > 0$ such that

$$\lambda_{\min}\Big( \sum w_i a(x_i) a(x_i)^* \Big) \geq (1-r)^2, \quad \lambda_{\max}\Big( \sum w_i b(x_i) b(x_i)^* \Big) \leq (1+s)^2 \, \|J\|.$$

Algorithm 1 chooses the points and weights one by one. At each step, an oracle suggests a new point which is then tested based on the spectral properties of corresponding Gram matrices. In order for the algorithm to be numerically feasible, two conditions are apparent:

1. The dimension of $b$ has to be finite (and not too large).
2. An oracle has to be available which suggests sampling points that get accepted by the algorithm at a decent rate.

Both issues are discussed in more detail in Section 4.2, when we discuss the application of Theorem 4.2 to least squares approximation. But let us already mention that, in case of an infinite basis, optimal recovery results can often still be obtained by applying Algorithm 1 to a truncated basis $b$. Regarding the second issue, the following can be said, see again Chkifa *et al.* (2026) for details:

- In each step, the set of points that would get accepted by the algorithm has a positive measure w.r.t. $\mu$.
- If we suggest the points randomly w.r.t. the Christoffel density associated to $V_m$, they will be accepted with a decent probability. Precisely, in the case $n \geq 2m$ and if we choose $w_i$ minimal at each step, the acceptance probability is at least $0.06/m^2$. We can even get an acceptance probability of order $1/m$ with insubstantial changes in the discretization bounds by slightly changing the parameter $\delta$, see Chkifa *et al.* (2026, Section 5.1) for details.
- We can also choose the oracle based on intuition or based on a large set of good sampling points of which we want to find a smaller subset. Whenever the algorithm terminates, the point set will satisfy the bounds from Theorem 4.2.

Several interesting discretization inequalities can be derived from Theorem 4.2 via different choices of $H_m$, which we now discuss shortly.

First, we note that Proposition 4.1 is a special case of Theorem 4.2, which is obtained by choosing $H_m = V_m$, which is a reproducing kernel Hilbert space when equipped with the $L_2$-norm, and in which case we have $A = 1$ and $K = m$.

Second, if $\mu$ is a probability measure, then one-sided discretization bounds with equal weights can be obtained by inserting for $H_m$ just a one-dimensional space of constant functions. Here, we have to look at the construction algorithm, which simplifies in this case. Some of the matrices become just real numbers, and if we choose the weights $w_i$ maximal in each step, they will all be equal. This eventually gives a discretization inequality of the form

$$\left(1 - \sqrt{\frac{m}{n}}\right) \cdot \|g\|_2 \;\leq\; \sqrt{\frac{1}{n}\sum_{i=1}^{n} |g(x_i)|^2} \qquad \text{for all } g \in V_m, \tag{4.1}$$

see Chkifa *et al.* (2026, Proposition 8) for details. Recall that this kind of one-sided unweighted discretization was the basis for the results from Section 3.2. Let us mention again that, up to worse constants, the result (4.1) can already be found in the paper of Bartel *et al.* (2023) together with a construction algorithm.

A third interesting application, and the basis of the results from Section 3.3, is the following. If we are given a RKHS $H$ that is embedded into $L_2(D,\mu)$ with finite trace, i.e., with square-summable singular values $(\sigma_k)_{k\in\mathbb{N}_0}$, an interesting choice is to take $V_m$ as the space spanned by the first $m$ singular functions and $H_m$ as the orthogonal complement of $V_m$. In this case, we have $A=\sigma_m$ and $K=\sum_{k\geq m}\sigma_k^2$. This means that we obtain the bounds

$$\left(1-\sqrt{\frac{m-1}{n}}\right)\|g\|_2 \leq \sqrt{\sum_{i=1}^n w_i|g(x_i)|^2} \qquad \text{for all } g\in V_m$$

and

$$\sqrt{\sum_{i=1}^n w_i|f(x_i)|^2} \leq \left(\sigma_m+\sqrt{\frac{1}{n}\sum_{k\geq m}\sigma_k^2}\right)\|f\|_H \quad \text{for all } f\in V_m^\perp.$$

*4.2. Weighted least squares*

As noted already by Bartel *et al.* (2023) and Dolbeault and Chkifa (2024), and discussed in the previous section as well as in Chkifa *et al.* (2026), the construction from Batson *et al.* (2014) can be used more or less directly to obtain (weighted or unweighted) least squares algorithms whose error in the $L_2$-norm is bounded in terms of the error of best uniform approximation on a given finite dimensional space $V_m$, that is, an error bound as in Theorem 3.12. In case a probability measure $\tilde{\mu}$ is known for which $\|\cdot\|_\infty \lesssim \sqrt{m}\|\cdot\|_{L_2(\tilde{\mu})}$ is satisfied on $V_m$, the same construction can be done twice (once for $\mu$ and once for $\tilde{\mu}$) to obtain sampling points that satisfy also the $L_p$-error bounds from Theorem 3.17.

We now discuss whether and how one can also obtain the points and weights behind the theorems from Section 3.3, where we bounded the error of weighted least squares algorithms in terms of the error of best $L_2$-approximation. We do this in the setting of the introductory Theorem 3.20 and Corollary 3.36, i.e., we assume that we are given a nested sequence $(V_m)_{m\in\mathbb{N}}$ of subspaces of $L_2$ with a corresponding orthonormal basis $(b_k)_{k\in\mathbb{N}_0}$. This is already the first important assumption: We need to know this basis. It is a very interesting open problem to generalize the construction (or invent a new one) to the situation where the basis is not orthonormal, or where the functions are not even linearly independent.

Then, in theory, the bounds from the aforementioned theorems are obtained by the points and weights returned by Algorithm 1 for the choice $a=(b_k)_{k<m}$ and $b=(k^{-1/2-\varepsilon/2}b_k)_{k\geq m}$, since according to Theorem 4.2, these points and weights satisfy the conditions of Proposition 3.22. Recall that Algorithm 1 chooses the points one by one following two basic steps:

1 A sampling point is suggested.

2 The point is tested against a certain criterion and accepted or rejected accordingly.

The first step can be done in numerous ways. One way is to use points from a large finite set (like a fine grid or a lattice) that is known, or at least conjectured, to contain some good set of sampling points. Of course, this is only practical if this 'candidate set' is not extremely large, since if we are unlucky, we have to test all the candidate points in the second step. However, we know of several examples where the size of a good candidate set should be not too bad, e.g., a regular grid on $[0,1]^2$ with $m^2$ points, a rank-1 lattice on $[0,1]^d$ of size $m^2$, or a sparse grid on $[0,1]^d$ of a quasi-linear size in $m$. We note that, even if this candidate set can be based on intuition and conjectures only, the second step is 'water-tight', and so we will either end up with a point set that is guaranteed to satisfy the desired error bound or otherwise the algorithm will tell us that all the suggested points are not acceptable.

Another method for the first step is to suggest the sampling point randomly. From the general theory, see again Chkifa *et al.* (2026), we know that an acceptable point exists and further that, if we suggest the point randomly according to the **Christoffel density**

$$\varrho_m(x) := \frac{1}{m} \sum_{k=0}^{m-1} |b_k(x)|^2, \tag{4.2}$$

it will be accepted with a probability of order $1/m^2$, or even $1/m$ for a slightly changed parameter $\delta$. This makes the first step feasible whenever we are able to generate points from the Christoffel density, leading to an expected number of $\mathcal{O}(mn)$ of random points that need to be tested in total.

**Remark 4.3.** The Christoffel density, which is clearly a probability density with respect to $\mu$, is an important quantity in various contexts and discussed in more detail in Section 7.3 when we talk about the randomized setting. Note that the sup-norm of the Christoffel density already appeared above in the context of lifting results for $L_2$-approximation to approximation results in the sup-norm, see (3.58).

The second step of the construction, in theory, involves computations with the whole basis sequence. We can deal with this issue by a truncation of the basis. The trick is to consider the function values of the target function as noisy samples of the truncated Fourier series and to exploit the robustness of the least squares method against noise, see Remark 3.6. In order to make this work, we need two assumptions: First, we assume that the measure $\mu$ is finite and that $V_m$ contains the constant functions. This will make it possible to have a bounded sum of weights, leading to the desired robustness against noise. For simplicity we take a normalized measure, i.e., a probability measure. Second, we want a nice pointwise behavior of the basis functions. Again, we give this assumption in terms of the quantity

$$\Lambda(V_m, B(D))^2 \;=\; \sup_{f\in V_m\setminus\{0\}} \frac{\|f\|_\infty^2}{\|f\|_2^2} \;=\; \Big\| \sum_{k=0}^{m-1} |b_k|^2 \Big\|_\infty .$$

Namely, we assume

$$\Lambda(V_m, B(D)) \lesssim m^\theta \tag{4.3}$$

for some (known) $\theta \geq 1/2$. Recall that $\theta = 1/2$ is best possible and that this will happen for any bounded orthonormal system but also for systems that are bounded on average like the Haar system, see Remark 3.33. However, let us also point out already that a larger $\theta > 1/2$ will not result in a worse error bound for $L_2$-approximation, only in stronger requirements on the 'smoothness' of $f$ and in a higher computational effort to find the points and weights.

Fortunately, this concludes the list of the required assumptions to obtain a practical construction:

- (C1) $\mu$ is a probability measure and approximation is with respect to an orthonormal basis of $L_2(D,\mu)$, containing the constant function.
- (C2) The basis satisfies a condition of the form (4.3).
- (C3) We can either generate samples from the Christoffel density or we know a good and reasonably-sized candidate set.

Note that these assumptions are not entirely independent and often coincide. In particular, under the assumption (4.3) the Christoffel sampling can be done by rejection sampling from the underlying probability measure $\mu$ with an acceptance probability of order $m^{-(2\theta-1)}$.

Under these conditions, we have the following constructive version of Theorem 3.20 and Corollary 3.36, which is again from Chkifa *et al.* (2026), at least in the case $p = 2$. This time, we also add the possibility of extra logarithmic factors to allow for a finer analysis of particular examples (like the function spaces of mixed smoothness from Example 3.25).

**Theorem 4.4.** Let $(D,\mu)$ be a probability space and let $\{b_k\}_{k\in\mathbb{N}_0}$ with $b_0 = 1$ be an orthonormal system in $L_2(D,\mu)$ such that (4.3) holds. Choose some $\alpha_0 > \theta$ and $t \in (\theta, \alpha_0)$. For $m \in \mathbb{N}$, consider the weighted least squares algorithm $A_m = A^w_{\mathcal{P}_n}$ with $n = 2m$ points and weights from Algorithm 1 for the bases

$$a = (b_0, \ldots, b_{m-1}) \quad \text{and} \quad b = (m^{-t} b_0, k^{-t} b_k)_{m\leq k<M},$$

with truncation index $M \geq m^{\alpha_0/(\alpha_0-\theta)}$. Then, the following holds for any function $f \in L_2(D,\mu)\cap B(D)$ with pointwise convergent Fourier series. If

$$d(f, V_m)_2 \;\lesssim\; m^{-\alpha} \log^\beta m \tag{4.4}$$

for some $\alpha \ge \alpha_0$ and $\beta \ge 0$, then

$$\|f - A_m f\|_p \lesssim m^{-\alpha+\theta(1-2/p)_+} \cdot \log^\beta m \tag{4.5}$$

for all $1 \le p \le \infty$. The hidden constant in (4.5) only depends on the constant in (4.3) and (4.4), and the numbers $\alpha, \beta, \theta$ and $p$.

Note that the polynomial rate in (4.5) reduces to $\alpha$ in the case $p \le 2$, resulting in optimal approximation rates for $L_2$-approximation for all sufficiently large $\alpha$ (i.e., for all sufficiently 'smooth' functions $f$). Whether the error bounds are good also for $L_p$-approximation with $p > 2$ depends much on the basis $(b_k)_{k\in\mathbb{N}_0}$ and the corresponding parameter $\theta$. Assuming that we can sample from the underlying measure $\mu$ and that the maximal time to compute a function value of $b_k$ is polynomial in $k \in \mathbb{N}$, the time to compute the points and weights from Theorem 4.4 is polynomial in $m$, see Chkifa *et al.* (2026, Section 5) for a fast implementation of Algorithm 1 and an analysis of the runtime.

*Proof.* It suffices to prove the result in the case $p = 2$. For $p \neq 2$, the error bound is obtained from the case $p = 2$ via the lifting result (3.60) together with the observation from Remark 3.33 that

$$\Lambda(V_m, L_p) \le \Lambda(V_m, B(D))^{(1-2/p)_+} \lesssim m^{\theta(1-2/p)_+}.$$

So let $p = 2$. We split the error into three parts:

$$\|f - A_m f\|_2 \le \|f - P_M f\|_2 + \|P_M f - A_m P_M f\|_2 + \|A_m P_M f - A_m f\|_2$$

and show that each term is $\mathcal{O}(m^{-\alpha} \log^\beta m)$.

For the first term, this is obvious. The second term of the above estimate is the 'main part' of the error bound; it is the error of the least squares algorithm for functions from $V_M$, the universe for which $A_m$ was constructed. Theorem 3.20 gives us

$$\|P_M f - A_m P_M f\|_2 \lesssim m^{-\alpha} \log^\beta m,$$

as desired. The third term regards the stability of the least squares method. Here, we also use that $V_m$ contains the constant functions, so that, according to Theorem 3.20, we have $\sum w_i \lesssim 1$. By Lemma 3.2, we obtain

$$\|A_m P_M f - A_m f\|_2 \le K \|f - P_M f\|_{\mathcal{P}_n^w} \lesssim \|f - P_M f\|_\infty.$$

Recall from (3.40) that $K$ is an absolute constant. Finally, the truncation error in the sup-norm can be bounded via Lemma 3.34, which gives us

$$\|f - P_M f\|_\infty \lesssim M^{-\alpha+\theta} \log^\beta M.$$

The truncation parameter $M$ was chosen so that the right hand side of the last inequality is in $\mathcal{O}(m^{-\alpha} \log^\beta m)$, so the proof is complete.

□

Let us finish the section with some numerical examples for the sampling points and weights from Theorem 4.4; more can be found in (Chkifa *et al.* 2026).

Figures 4.1 and 4.2 show the sampling points for $D = [0,1]^2$ with the Lebesgue measure that are obtained by applying the construction from Theorem 4.4. In all cases $(b_k)_{k\in\mathbb{N}_0}$ is a certain ordering of the Fourier basis $(e^{2\pi i(rx+jy)})_{r,j\in\mathbb{Z}}$, in which case we have $\theta = 1/2$. We choose $\alpha_0 = 3/2$ and $t = 1$, which means that the constructed points will give optimal approximation rates if the rate of best $L_2$-approximation in this basis is at least $3/2$.

Specifically, we consider

$$a = \left( e^{2\pi i(rx+jy)} \colon\ r,j \in \mathbb{Z} \ \text{ with } \ \sigma_{r,j} \le R \right)$$

and

$$b = \left( \frac{1}{\sigma_{r,j}} \cdot e^{2\pi i(rx+jy)} \colon\ r,j \in \mathbb{Z} \ \text{ with } \ R < \sigma_{r,j} \le R' \ \text{ or } \ (r,j) = (0,0) \right)$$

with $R' = 1000$, for $\sigma_{r,j} = 2\pi(r^2 + j^2)$ (isotropic smoothness; Figure 4.1) and $R = 80$, as well as $\sigma_{r,j} = \max\{1, 2\pi|r|\} \cdot \max\{1, 2\pi|j|\}$ (mixed smoothness; Figure 4.2) and $R = 50$, respectively. The definition of $\sigma_{0,0}$ is modified, and plays the role of the factor $m^t$ in Theorem 4.4. Note that the condition $M \ge m^{\alpha_0/(\alpha_0-\theta)}$ is satisfied for this choice of $R$ and $R'$.

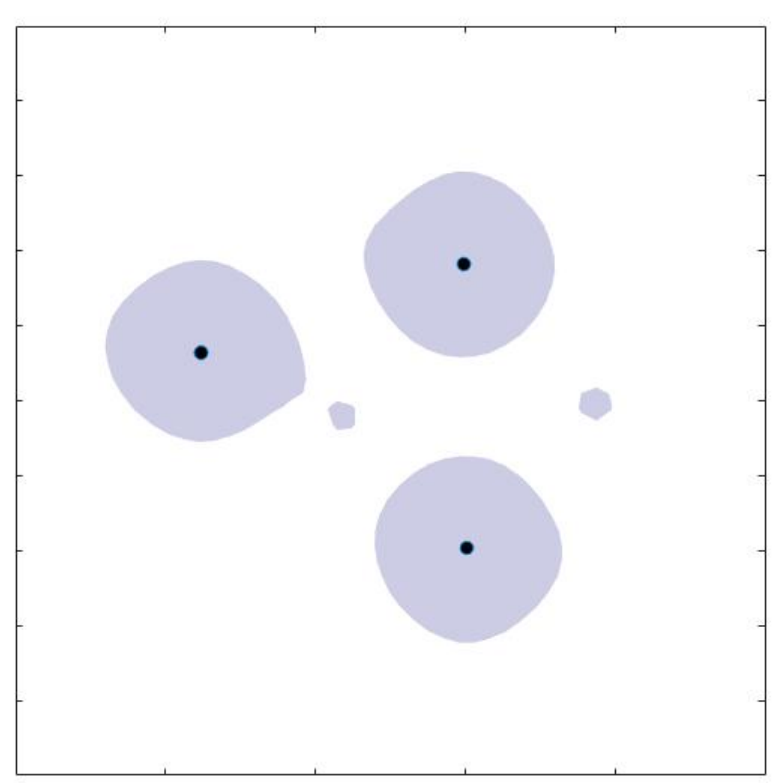

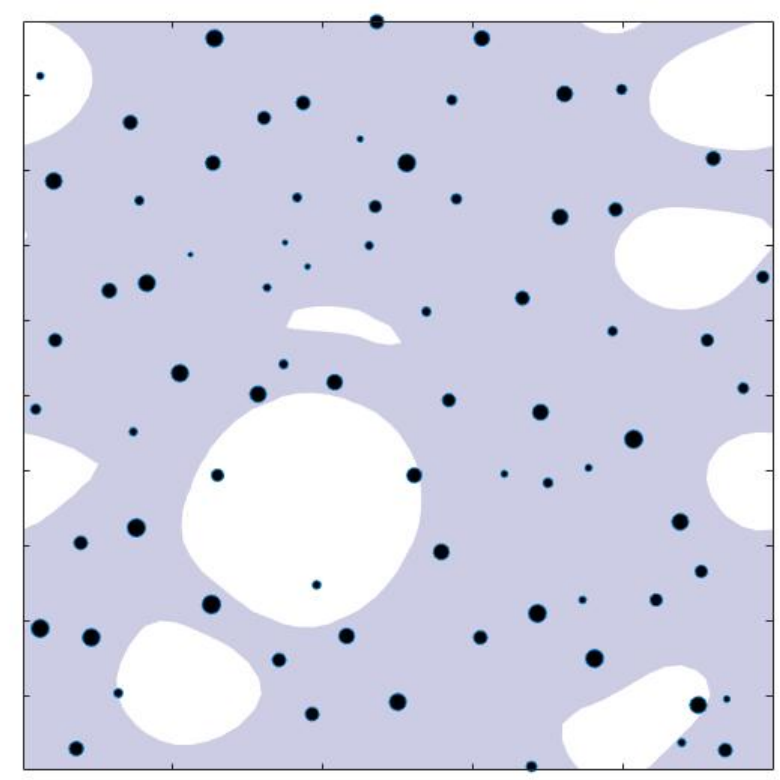

Figure 4.1. Points generated according to Theorem 4.4 for the Fourier basis and isotropic smoothness. Colored area: points that would be rejected in the next step.

Note that we choose the weights in the definition of the basis $b$ of order $k^{-t} = 1/k$ (in the mixed case up to a logarithmic factor, which does not change the result for $\alpha \ge 3/2$), but not precisely equal as in Theorem 4.4; it seems more natural for this example to choose them constant for all frequencies with the same distance to the origin. This way, the point set will give optimal rates for Sobolev classes of isotropic smoothness $s \geq 3$ and mixed smoothness $s \geq 3/2$, respectively. In both cases, the points in the algorithm are suggested randomly and the weight in Step 7 is chosen minimal. We note that the Christoffel density is just the uniform distribution in this case, and that the true acceptance rate of the suggested points was a lot higher than the theoretical lower bound $0.06/m^2$ in our tests; usually we only had to suggest two to ten points. In all pictures, the size of the points is proportional to the size of the corresponding weight.

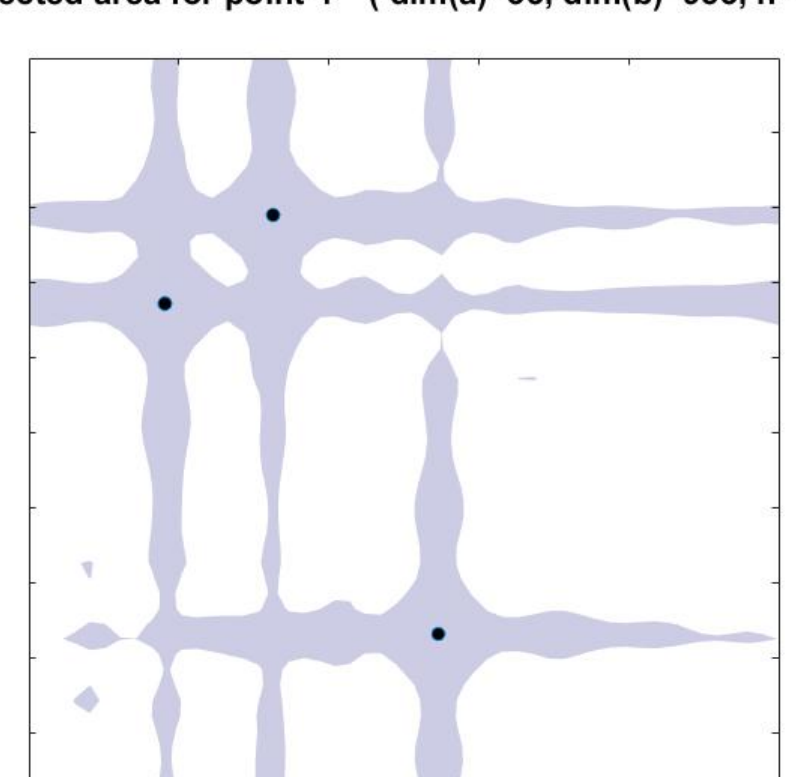

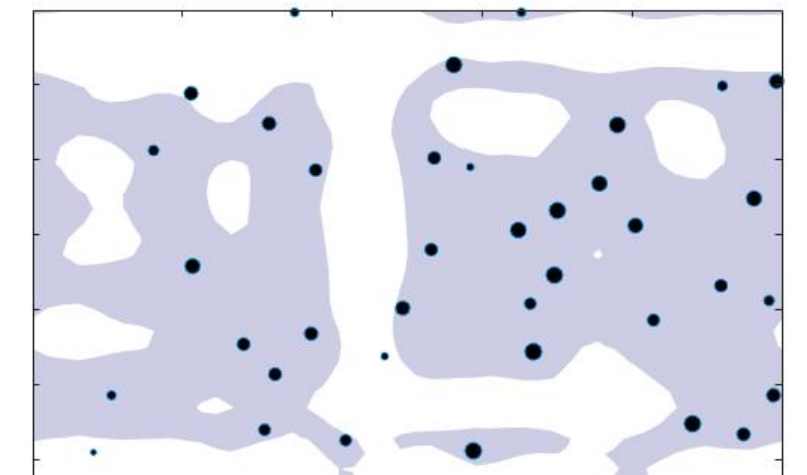

Figure 4.2. Points generated according to Theorem 4.4 for the Fourier basis and mixed smoothness. Colored area: points that would be rejected in the next step.

In both cases, we first plot the first three points constructed by the algorithm, together with the area of all points that would get accepted by the algorithm in the next step. The dark area is the area of all points that would get rejected. It is interesting to see that this area of “blocked” points takes roughly the form of balls around the existing points in the isotropic case, and roughly the form of hyperbolic crosses in the mixed case.

As a second picture in each case, we plot the full point set together with the area of points that would have been accepted or rejected for the last point. This time, the white areas (i.e., the areas of acceptable points) resemble a little the “large empty balls” in the isotropic case and the “large empty hyperbolic crosses” in the mixed case.

## 5. Optimal sampling and nonlinear approximation

In the previous sections, we discussed sampling algorithms for approximation with the goal of getting close to the error of best approximation from a linear space. We now discuss bounds by benchmarks from nonlinear approximation which may be much smaller, as they also allow for more general *approximation procedures*.

The results presented here are a rather loose collection of different bounds for sampling algorithms and sampling numbers. Recall that the $n$-th sampling number of a model class $F$ in a semi-normed space $Y \subset \mathbb{C}^D$ is the worst-case error of the best algorithm that uses $n$ samples, defined by

$$g_n(F, Y) := \inf_{\substack{x_1,\dots,x_n \in D \\ \Phi\colon \mathbb{C}^n \to Y}} \sup_{f \in F} \Big\| f - \Phi(f(x_1), \dots, f(x_n)) \Big\|_Y .$$

In this section, it will often be important that we allow nonlinear reconstruction mappings $\Phi$ and the linear sampling numbers $g_n^{\mathrm{lin}}(F, Y)$ from Section 3.1 will be of less importance. So far, we are far from obtaining the same level of understanding as in the linear setting and do not really know what to expect from such general comparisons. Although there are some very *convincing examples* and applications, the theory presented in the following certainly comes with a smaller variety of them. In this section, we consider upper bounds in terms of

- the error of best sparse approximation (Section 5.1),
- metric entropy (Section 5.2), and
- several other approximation benchmarks, including the error of best continuous approximation from a manifold (Section 5.3).

Another natural benchmark in this context are the *Gelfand widths*. Upper bounds on the sampling numbers in terms of Gelfand widths are particularly interesting from an information-theoretic point of view since this leads to a comparison of the power of different types of sampling: function evaluations versus arbitrary linear measurements. Several such bounds are obtained as corollaries of the bounds presented in this section. This will be discussed later in Sections 9 and 10, see especially Section 10.1.

### *5.1. Bounds by best sparse approximation*

Instead of an approximation from a linear space of dimension $m$, we now aim for an approximation by a $m$-sparse function in a linear space $V_N$ of much larger dimension $N \gg m$. With $m$-sparse functions, we mean functions whose expansions in a fixed basis have the property that all but $m$ coefficients are zero (or close to zero); we will give precise definitions soon. Hence, we aim for an approximation in a union of $m$-dimensional linear spaces; a nonlinear space that does not allow for the direct use of the techniques from Section 3. We ask the following questions: Are there sampling algorithms with an error close to the error of best approximation by an $m$-sparse function? If yes, how many sampling points do we need for this?

Undoubtedly, there is a large variety of results about sparse approximation in the literature, especially in the area of **compressive sensing**. We refer to the book by Foucart and Rauhut (2013), the seminal work of Candès, Romberg and Tao (2006) and Donoho (2006), and early key results by Kashin (1977) and Garnaev and Gluskin (1984). Though the most classical setting here is to recover high-dimensional vectors from few linear measurements (and not function values), several of the results also deal with the recovery of functions from function values, which is the setting that we are interested in.

Here we present a bound which is based on a result about the square root Lasso algorithm by Petersen and Jung (2022), an algorithm introduced by Belloni, Chernozhukov and Wang (2011). The result we present is essentially from Moeller, Neumayer, Pozharska, Sommerfeld and Ullrich (2025*a*) and Moeller, Stasyuk and Ullrich (2025*b*). Comparable bounds are available also for the (weak) orthogonal matching pursuit algorithm which is discussed by Dai and Temlyakov (2024) (see esp. Corollary 5.1 there) and also in the paper by Moeller *et al.* (2025*a*). In fact, interesting results can be obtained already by considering unweighted least squares algorithms on a nonlinear space. We discuss this shortly towards the end of the section, where we also comment a bit more on the history, see Remark 5.5.

We now rigorously introduce the setting. Let $(D, \mu)$ be a probability space and for some (typically huge) number $N \in \mathbb{N}$, assume that $\mathcal{B}_N := \{b_1, \ldots, b_N\}$ is an orthonormal system in $L_2(D, \mu)$ that has the property

$$\max_{k \le N} \|b_k\|_\infty \;\le\; \Theta$$

for some absolute constant $\Theta \ge 1$. We say that $\mathcal{B}_N$ is a **bounded orthonormal system (BOS)** and call $\Theta$ the BOS-constant. Examples are families of Fourier or the Walsh basis functions on $D = [0,1]^d$ with the Lebesgue measure, in which case we have $\Theta = 1$. For an integrable function $f\colon D \to \mathbb{C}$, define

$$\hat{f}(k) \;:=\; \int f(x)\,\overline{b_k(x)}\,d\mu(x), \qquad k \le N. \tag{5.1}$$

By $V_N$ we denote the linear span of $\mathcal{B}_N$. On the space $V_N$, a norm is defined by

$$|||f||| \;:=\; \sum_{k=1}^{N} |\hat{f}(k)|.$$

By assumption, we have that $\|f\|_\infty \le \Theta \cdot |||f|||$ for $f \in V_N$, but a small norm $|||f|||$ additionally implies some sort of sparsity, as will be discussed below.

For $m \le N$, the set of $m$**-sparse** functions from $V_N$ is

$$V_N^{(m)} \;:=\; \Big\{ f \in V_N \,:\, \#\{k \le N \colon \hat{f}(k) \ne 0\} \le m \Big\}. \tag{5.2}$$

The error of **best $m$-term approximation** of $f \in V_N$ is defined by

$$\sigma_m(f) := \sigma_m(f, V_N) := \inf_{g \in V_N^{(m)}} |||f - g||| = \min_{\substack{I \subset \mathbb{N} \\ \#I = m}} \sum_{k \in \mathbb{N} \setminus I} |\hat{f}(k)|.$$

Of course, best $m$-term approximation can also be considered with respect to other norms than the 'Wiener algebra norm' $|||\cdot|||$. We use this norm not because it is the most interesting, but rather because it is the norm in terms of which we can prove the best results for sampling recovery. Similarly good results in terms of best approximation in the sup-norm (which is a smaller quantity, but which also leads to additional factors in the error bounds) are obtained by Moeller *et al.* (2025*a*).

Assume that we are given sampling points $\mathcal{P}_n = \{x_1, \ldots, x_n\}$ in $D$ and a parameter $\lambda \geq 0$. For any function $f \colon D \to \mathbb{C}$, the **square root Lasso algorithm** is defined as

$$R_{\mathcal{P}_n,\lambda}(f) \in \operatorname*{arg\,min}_{g \in V_N} \left( |||g||| + \lambda \sqrt{\sum_{i=1}^{n} |g(x_i) - f(x_i)|^2} \right). \tag{5.3}$$

A reason to use the square root Lasso algorithm (and not, for instance, constrained $\ell_1$-minimization) is that the parameter $\lambda$, and hence the algorithm, can be chosen independent of the noise level, see e.g. Adcock, Brugiapaglia and Webster (2022, Section 6.6) for a discussion. In our setting, this means that the algorithm requires no a priori bound on the approximability of the target function $f$ by functions from the space $V_N$ in the uniform norm.

The following error bound is similar to Proposition 2.5 and Theorem 2.8 from Moeller *et al.* (2025*a*). It is implied by Theorem 3.1 of Petersen and Jung (2022). The proof is a technical translation of this result and shifted to the end of the section.

**Proposition 5.1.** Let $(D, \mu)$ be a probability space and let $V_N$ be an $N$-dimensional space that is spanned by a bounded orthonormal system in $L_2(D, \mu)$ and $m \leq N/2$. Assume that the point set $\mathcal{P}_n = \{x_1, \ldots, x_n\}$ is such that

$$\frac{3}{4} \|g\|_2^2 \leq \frac{1}{n} \sum_{i=1}^{n} |g(x_i)|^2 \leq \frac{5}{4} \|g\|_2^2 \qquad \text{for all } g \in V_N^{(2m)}. \tag{5.4}$$

Then the square root Lasso algorithm (5.3) with $\lambda = 2\sqrt{m/n}$ satisfies

$$\|f - R_{\mathcal{P}_n,\lambda}(f)\|_p \leq C\, m^{(1/2-1/p)_+} \inf_{g \in V_N} \left( \frac{\sigma_m(g)}{\sqrt{m}} + \|f - g\|_\infty \right) \tag{5.5}$$

for any function $f \in L_2(D, \mu) \cap B(D)$ and all $1 \leq p \leq \infty$, where $C > 0$ is a constant that only depends on the BOS-constant $\Theta$.

Before we proceed to discuss the existence of points $\mathcal{P}_n$ as in (5.4), we mention two immediate implications of the bound (5.5). The first is obtained by considering only sparse functions $g \in V_N^{(m)}$, for which $\sigma_m(g) = 0$, so that we get

$$\|f - R_{\mathcal{P}_n,\lambda}(f)\|_p \;\le\; C\, m^{(1/2-1/p)_+} \inf_{g \in V_N^{(m)}} \|f - g\|_\infty. \tag{5.6}$$

That is, the error of the algorithm is bounded by the error of best uniform approximation by an $m$-sparse function from $V_N$. We point to the similarity of this error bound to the error bound from Theorem 3.17. Particularly, for $N = m$, the bound (5.6) is in terms of the error of best uniform approximation from $V_m$ and essentially equal to the bound of Theorem 3.17. The possible gain happens, of course, in the case $N > m$, where the best approximation from an $m$-dimensional space is replaced with a smaller quantity: the best $m$-term approximation from a space of much larger dimension $N$. On the other hand, also the discretization condition (5.4) is stronger than the discretization condition of Theorem 3.17 and hence it is no surprise that (5.4) will require a larger sample. However, the required sampling size is only larger by logarithmic factors in $m$ and $N$, as we discuss below, which makes it possible to deal with spaces of dimension $N \gg m$. The second implication is obtained by putting $g = P_N f$, which shows

$$\|f - R_{\mathcal{P}_n,\lambda}(f)\|_p \;\le\; C\, m^{(1/2-1/p)_+} \left( \frac{1}{\sqrt{m}} \sum_{k=m+1}^{N} |\hat{f}(k)| + \|f - P_N f\|_\infty \right). \tag{5.7}$$

Since $N$ is usually huge compared to $m$, for $p < \infty$, this already gives “good” bounds under the very mild assumption of a bounded sum of Fourier coefficients; a bound which is especially useful in tractability studies, see Section 8.

We turn to the discussion of the property (5.4) of the point set $\mathcal{P}_n$. In recent papers, this property is discussed under the name **universal discretization**, referring to the fact that the point set $\mathcal{P}_n$ provides discretization of a continuous norm universally for a whole family of $2m$-dimensional linear spaces. The term universal discretization was coined by Temlyakov (2018*a*). The universal discretization condition (5.4) is equivalent to the **restricted isometry property (RIP)** of order $2m$ of the ‘Fourier matrix’

$$A \;:=\; \frac{1}{\sqrt{n}} \big( b_j(x_i) \big)_{j \le N, i \le n}\,. \tag{5.8}$$

Namely, it is equivalent to

$$\frac{3}{4}\, \|c\|_2^2 \;\le\; \|Ac\|_2^2 \;\le\; \frac{5}{4}\, \|c\|_2^2$$

for all vectors $c \in \mathbb{C}^N$ with at most $2m$ non-zero entries; of course, the constants $3/4$ and $5/4$ chosen here are quite arbitrary. The restricted isometry property has been first introduced by Candès and Tao (2005).

Finding measurement matrices (5.8) that satisfy the restricted isometry property is a challenging problem. Results are usually proven by probabilistic arguments.

Already Kashin (1977), who considered measurement matrices $A \in \mathbb{C}^{n\times N}$ different from (5.8), noticed that very good results are obtained with high probability by choosing the matrix randomly. There also exist deterministic constructions, see for example DeVore (2007), but the bounds obtained in this way are much worse. Regarding the Fourier matrices of the form (5.8), it was proved by Candès and Tao (2006) that $\mathcal{O}(m \log^6 N)$ random points satisfy the RIP with high probability. The upper bound on the number of required samples has been subsequently improved by Rudelson and Vershynin (2008), Cheraghchi, Guruswami and Velingker (2013), Bourgain (2014), Chkifa, Dexter, Tran and Webster (2018), Haviv and Regev (2017), and Brugiapaglia, Dirksen, Jung and Rauhut (2021), resulting in a bound of order $\mathcal{O}(m \log^2 m \log N)$. The following is implied by Theorem 2.3 of Brugiapaglia *et al.* (2021), see also Adcock *et al.* (2022, Section 6.4.3).

**Proposition 5.2.** Let $(D,\mu)$ be a probability space and let $V_N$ be an $N$-dimensional space that is spanned by a bounded orthonormal system in $L_2(D,\mu)$. For $m \in \mathbb{N}$, there exist

$$n \;\le\; C\, m \log^2 m \log N$$

points $\mathcal{P}_n = \{x_1,\dots,x_n\}$ such that (5.4) holds, i.e.,

$$\frac{3}{4}\,\|g\|_2^2 \;\le\; \frac{1}{n}\sum_{i=1}^{n} |g(x_i)|^2 \;\le\; \frac{5}{4}\,\|g\|_2^2 \qquad \text{for all } g \in V_N^{(2m)}. \tag{5.9}$$

Here, $C > 0$ is a constant that only depends on the BOS-constant $\Theta$.

In fact, the stated amount of i.i.d. random points with distribution $\mu$ satisfy the statement from Proposition 5.2 with high probability. In this sense, the result is semi-constructive. However, it is not fully constructive: No way is known of constructing a point set with an even remotely similar size that satisfies (5.4) with absolute certainty. Also for the randomly drawn points, it is not computationally feasible to verify whether a particular realization satisfies the required RIP.

For the error of the square root Lasso algorithm, Propositions 5.1 and 5.2 imply the following; note the similarity to (Moeller *et al.* 2025*a*, Theorem 2.8) and (Moeller *et al.* 2025*b*, Theorem 3.7).

**Theorem 5.3.** Let $(D,\mu)$ be a probability space and let $V_N$ be an $N$-dimensional space that is spanned by a bounded orthonormal system in $L_2(D,\mu)$. For $m \le N/2$, there exist

$$n \;\le\; C\, m \log^2 m \log N$$

points $\mathcal{P}_n = \{x_1,\dots,x_n\}$ such that, for $\lambda = 2\sqrt{m/n}$, it holds

$$\|f - R_{\mathcal{P}_n,\lambda}(f)\|_p \;\le\; C\, m^{(1/2-1/p)_+} \inf_{g\in V_N} \left( \frac{\sigma_m(g)}{\sqrt{m}} + \|f-g\|_\infty \right)$$

for any function $f \in L_2(D,\mu) \cap B(D)$ and all $1 \le p \le \infty$ and with a constant

$C > 0$ that only depends on the BOS-constant $\Theta$. In particular,

$$\|f - R_{\mathcal{P}_n,\lambda}(f)\|_p \;\leq\; C\, m^{(1/2-1/p)_+} \inf_{g\in V_N^{(m)}} \|f-g\|_\infty. \tag{5.10}$$

Hence, it is possible to get close to the error of best $m$-sparse uniform approximation in an $N$-dimensional space with a sampling algorithm that uses $\mathcal{O}(m\log^2 m \log N)$ samples. However, recall that this only holds for spaces $V_N$ that are spanned by a bounded orthonormal basis. In this sense, the results of this section are weaker compared to those of the previous ones, which often work for general (linear) subspaces. Progress in this regard has been made via the concept of weighted sparsity, see Candès, Wakin and Boyd (2008) or the recent comprehensive treatments of Adcock *et al.* (2022) and Adcock, Brugiapaglia, Dexter and Moraga (2024) with a focus on algebraic polynomials.

Let us discuss an easy example.

**Example 5.4.** Let $\{b_k\colon k\in\mathbb{Z}^d\}$ be the trigonometric system, i.e.,

$$b_k(x) = e^{2\pi i\langle k,x\rangle}, \quad x\in[0,1]^d,$$

which is an orthonormal basis of $L_2([0,1]^d)$ with the Lebesgue measure. Then, (5.1) defines the usual Fourier coefficients of $f\in L_2([0,1])$,

$$\hat{f}(k) = \int_{[0,1]^d} f(x)\cdot e^{-2\pi i\langle k,x\rangle}\,dx, \quad k\in\mathbb{Z}^d. \tag{5.11}$$

Assume that the model class $F$ is such that all functions have a nice summability of their Fourier coefficients,

$$F \;\subseteq\; \left\{ f\in C([0,1]^d)\colon \sum_{k\in\mathbb{Z}^d} \frac{|\hat{f}(k)|}{w_k} \leq 1 \right\} \tag{5.12}$$

for some Fourier weights $(w_k)_{k\in\mathbb{Z}^d}$ with $w_k\to 0$ for $\|k\|_2\to\infty$. Note that this condition implies that any $f\in F$ is pointwise equal to its Fourier series. In such a case, Theorem 5.3 gives for all $m,d\in\mathbb{N}$ that

$$g_{Cdm\log^3 m}(F,L_2) \;\leq\; \frac{w_m^*}{\sqrt{m}} + w_{m^d}^*, \tag{5.13}$$

where $(w_n^*)_{n\in\mathbb{N}}$ shall be the monotonically decreasing rearrangement of the Fourier weights, and $C$ is a universal constant. Indeed, we use the theorem for $N=m^d$, $V_N$ being the space of trigonometric polynomials with frequencies corresponding to the $N$ largest Fourier weights, and for $g=P_N f$ being the orthogonal projection onto $V_N$. Then it is clear that

$$\|f-P_N f\|_\infty \;\leq \sum_{k\colon w_k\leq w_N^*} \left|\hat{f}(k)\right| \;\leq\; w_N^*$$

and

$$\sigma_m(P_N f) \leq \sum_{k\colon w_k \leq w_m^*} \left|\hat{f}(k)\right| \leq w_m^*,$$

which gives the result. Classes of the form (5.12) might be called weighted Wiener spaces or also Fourier-Barron spaces and sampling recovery on such classes is considered in much more detail by Kolomoitsev, Lomako and Tikhonov (2023), Jahn, Ullrich and Voigtlaender (2023) and Moeller *et al.* (2025*a*,*b*).
The bound (5.13) seems particularly useful in high dimensions since for a normalized problem with $w_1^* = 1$, the first summand can be bounded by $m^{-1/2}$, independently of the dimension, and the number of samples is only logarithmic in the index of the second summand $w_N$; this means that it is no problem if an error demand $w_N^* \leq \varepsilon$ requires $N$ to be (super-)exponentially large in the dimension. The consequences of such a bound for the tractability of several function approximation problems in high dimensions are studied in Krieg (2024), see also Section 8. It is also quite interesting to note that neither the optimal error rates from Jahn *et al.* (2023) and Moeller *et al.* (2025*a*) nor the tractability results from Krieg (2024) can be achieved with any linear algorithm. Hence, for classes of the form (5.12), in general, nonlinear sampling algorithms are better than linear ones (like the least squares method considered in the last section).

**Remark 5.5 (History).** The applicability of the compressive sensing theory of Candès *et al.* (2006) and Donoho (2006) to sampling recovery went largely unnoticed for quite some time. This is probably due to the fact that classical compressive sensing results are concerned with a finite-dimensional setting, while in sampling recovery, it is usually an undesirable assumption that the target function belongs to some finite-dimensional space $V_N$. However, since those results often allow for noisy measurements, it is usually possible to interpret the function values of the target function $f \notin V_N$ as noisy measurements of some $g \in V_N$ and thus use the corresponding results.

This simple yet somewhat hidden path, used, for example, already in Rauhut and Ward (2016), has recently been explored in a number of papers. It seems that a general connection of sampling numbers to compressive sensing was first noticed by Jahn *et al.* (2023), where the authors used a result by Rauhut and Ward (2012) about *basis pursuit denoising* to prove new results on the order of convergence for sampling recovery on mixed smoothness classes with a low integrability index.

We note that the error bound (5.10) in terms of the error of best $m$-sparse approximation can also be proven with the unweighted least squares algorithm instead of the square-root Lasso algorithm. For this, it is enough to combine Proposition 5.2 with the following simple proposition on the error of least squares (and other) methods on nonlinear sets. Namely, for an arbitrary set of functions

$\Sigma \subset \mathbb{C}^D$, we may consider the algorithm

$$A^*_\Sigma f \;\in\; \operatorname*{arg\,min}_{g\in\Sigma} \|f-g\|_*,$$

provided that a minimizer exists, with a seminorm $\|\cdot\|_*$ on $\mathbb{C}^D$, like a weighted $\ell_2$-norm of function values as in (3.2). If a minimizer does not exist, it is easy to modify the following result using an almost-minimizer.

**Proposition 5.6.** Let $D$ be a set, $Y \subset \mathbb{C}^D$ be a semi-normed space, $\Sigma \subset Y$ and $\|\cdot\|_*$ be another seminorm such that

$$\|g\|_Y \;\le\; K\cdot\|g\|_* \qquad \text{for all} \quad g\in\Sigma-\Sigma \tag{5.14}$$

and some $K<\infty$. Then, for all $f\in Y$ and $g\in\Sigma$, we have

$$\left\|f-A^*_\Sigma f\right\|_Y \;\le\; \left\|f-g\right\|_Y + 2K\cdot\|f-g\|_*.$$

The proof is similar to the linear case from Proposition 3.3.

*Proof.* By the triangle inequality and since $g-A^*_\Sigma f\in\Sigma-\Sigma$, it holds

$$\begin{aligned}\left\|f-A^*_\Sigma f\right\|_Y &\le \|f-g\|_Y+\|g-A^*_\Sigma f\|_Y\\ &\le \|f-g\|_Y+K\|g-A^*_\Sigma f\|_*\\ &\le \|f-g\|_Y+K\big(\|g-f\|_*+\|f-A^*_\Sigma f\|_*\big).\end{aligned}$$

By the definition of $A^*_\Sigma f$, it holds that

$$\|f-A^*_\Sigma f\|_* \;\le\; \|f-g\|_*,$$

so that the claim is proved.

□

If we use this proposition for $\|\cdot\|_*$ being the unweighted $\ell_2$ norm of function values at the sampling points from Proposition 5.2, the approximation set $\Sigma=V_N^{(m)}$, and the error norm $Y=L_p$, the algorithm $A^*_\Sigma$ is the unweighted least squares algorithm on $V_N^{(m)}$ and the error bound (5.10) is obtained with the square root Lasso algorithm $R_{\mathcal{P}_n,\lambda}(f)$ replaced by the least squares algorithm $A^*_\Sigma$. However, this bound seems to be more of theoretical value, since it is generally hard to find such a minimizer. We refer to Cohen, Dolbeault, Mula and Somacal (2023) for more on *nonlinear least squares* and its relation to compressive sensing.

We conclude with the proof of Proposition 5.1.

*Proof.* In terms of the Fourier matrix (5.8), the assumption means that $A$ has the restricted isometry property of order $2m$ with a constant $1/4$. By (Foucart and Rauhut 2013, Theorem 6.13), this implies the $\ell_2$-robust null space property of order $m$ as required by (Petersen and Jung 2022, Theorem 3.1) with a constant $\tau\le 3/2$ and for the Euclidean norm $\|\cdot\|_2$, and therefore also the $\ell_1$-robust null space property for the norm $\sqrt{m}\|\cdot\|_1$ with the same value of $\tau$. Therefore, for $y \in \mathbb{C}^N$, we consider

$$c^{\#} \in \operatorname*{arg\,min}_{z\in\mathbb{C}^N}\Big(\|z\|_1 + 2\sqrt{m}\,\|Az-y\|_2\Big)$$

and define

$$\sigma_m(c)_1 := \inf_{\substack{I\subset\{1,\dots,N\}\\ \#I\le m}} \sum_{k\in\{1,\dots,N\}\setminus I} |c_k|.$$

We obtain for any $c \in \mathbb{C}^N$ that

$$\|c^{\#}-c\|_2 \le C\left(\frac{\sigma_m(c)_1}{\sqrt{m}} + \|Ac-y\|_2\right)$$

and

$$\|c^{\#}-c\|_1 \le C\sqrt{m}\left(\frac{\sigma_m(c)_1}{\sqrt{m}} + \|Ac-y\|_2\right)$$

with a universal constant $C > 0$. If we put $y = \frac{1}{\sqrt{n}}(f(x_i))_{i\le n}$, then $c^{\#}$ is precisely the vector of coefficients of $R_{\mathcal{P}_n,\lambda}(f)$. If we further let $g = \sum_{k=1}^N c_k b_k$, then $\sigma_m(c)_1 = \sigma_m(g)$ and

$$\|Ac-y\|_2 = \sqrt{\frac{1}{n}\sum_{i=1}^{n}|g(x_i)-f(x_i)|^2} \le \|f-g\|_\infty$$

and so

$$\|R_{\mathcal{P}_n,\lambda}(f)-g\|_2 = \|c^{\#}-c\|_2 \le C\left(\frac{\sigma_m(g)}{\sqrt{m}} + \|f-g\|_\infty\right)$$

and

$$\|R_{\mathcal{P}_n,\lambda}(f)-g\|_\infty \le \Theta\|c^{\#}-c\|_1 \le C\sqrt{\Theta^2 m}\left(\frac{\sigma_m(g)}{\sqrt{m}} + \|f-g\|_\infty\right).$$

Interpolation gives

$$\|R_{\mathcal{P}_n,\lambda}(f)-g\|_p \le C\,(\Theta^2 m)^{(1/2-1/p)_+}\left(\frac{\sigma_m(g)}{\sqrt{m}} + \|f-g\|_\infty\right).$$

From the triangle inequality and $\|\cdot\|_p \le \|\cdot\|_\infty$, it follows that

$$\|f-R_{\mathcal{P}_n,\lambda}(f)\|_p \le (C+1)\,(\Theta^2 m)^{(1/2-1/p)_+}\left(\frac{\sigma_m(g)}{\sqrt{m}} + \|f-g\|_\infty\right) \tag{5.15}$$

for any function $g \in V_N$.

□

*5.2. Bounds by entropy numbers*

We now turn to some non-constructive upper bounds on the (linear) sampling numbers in terms of other *approximation benchmarks*. In this section we use the **(metric) entropy numbers** of $F$ in $Y$,

$$\varepsilon_n(F,Y) := \inf_{y_1,\ldots,y_{2^n}\in Y} \sup_{f\in F} \min_{i=1,\ldots,2^n} \|f - y_i\|_Y \tag{5.16}$$

with $Y = B := B(D)$. In other words, $\varepsilon_n(F,Y)$ represents the smallest $\varepsilon$ such that $F$ can be covered by $2^n$ balls in $Y$ of radius $\varepsilon$. One may interpret this as the minimal worst-case error over all binary encoders $N_n(f) \in \{0,1\}^n$ (where $N_n(f)$ can encode the index $i$ of a closest center $y_i$), and they are therefore often used as ultimate benchmark for compression, and sometimes also for approximation by (general) nonlinear methods, see for example (Cohen, DeVore, Petrova and Wojtaszczyk 2022). The entropy numbers encode certain properties of the sets $F \subset Y$: For example, $\varepsilon_n$ is bounded (or a null sequence) if and only if $F$ is bounded (or totally bounded). The speed of decay of $\varepsilon_n$ can therefore be seen as kind of *degree of compactness*. Clearly, the entropy numbers are the inverse quantity to the (dyadic) *covering numbers*, i.e., the minimal number $n$ such that $F$ can be covered by $2^n$ balls of a given radius.

**Remark 5.7 (Notation).** Note that we differ from the standard notation of entropy numbers by using the symbol $\varepsilon$ instead of $e$ and by an index shift, i.e., $\varepsilon_n = e_{n+1}$. We use $\varepsilon$ because $e$ is already taken for the error of algorithms and we use the index shift to be consistent with the other error quantities in this manuscript, where $n$ is always the number of measurements and we start counting with $n = 0$.

We will comment on the general theory of widths and s-numbers in Section 9.2. But for comparison with the previous sections, let us already mention **Carl's inequality** due to Carl (1981*b*), which says that, for $F = B_X$ being the unit ball of a normed space, we have

$$\sup_{k\le n} (k+1)^s \cdot \varepsilon_k(F,Y) \le c_s \cdot \sup_{k\le n} (k+1)^s \cdot d_k(F,Y), \tag{5.17}$$

where $c_s < \infty$ only depends on $s > 0$. This bound also holds for unit balls of $p$-Banach spaces, see Hinrichs, Kolleck and Vybíral (2016), which are linear spaces $X$ which are complete with respect to a quasinorm that, instead of the triangle inequality, only satisfies $\|f+g\|_X^p \le \|f\|_X^p + \|g\|_X^p$ for all $f, g \in X$. This means that the following upper bounds by entropy can be translated into upper bounds by Kolmogorov widths.

The first bound we want to discuss is the following result by Kosov and Temlyakov (2025, Remark 2.2). Note that we simplify the bound for convenience and that the authors also present a variant for $s = 1/2$ with an additional logarithmic (in $n$) factor.

**Theorem 5.8.** Let $D \subset \mathbb{R}^d$ be compact, $F \subset C(D)$ be convex and symmetric, and $(D,\mu)$ be a probability space. Then, for all $1 \le p < \infty$ and $s \neq 1/2$, we have

$$g_n(F, L_p) \le C \cdot n^{-s'/p} \left( \sup_{k \le n} (k+1)^s \cdot \varepsilon_k(F, B) \right)^{1/p} \tag{5.18}$$

with $s' := \min\{s, 1/2\}$ and $C > 0$ only depending on $s \neq 1/2$ and $p$.

Theorem 5.8 is based on a result about discretization of the $L_p$-norm on the class $F$. In fact, it is used that

$$g_n(F, L_p) \le 2 \cdot \inf_{x_1,\dots,x_n \in D} \sup_{f \in F} \left| \|f\|_p^p - \frac{1}{n} \sum_{i=1}^n |f(x_i)|^p \right|^{1/p},$$

and the right hand side is bounded as in (5.18), see (Kosov and Temlyakov 2025, Theorem 1.3 & Remark 2.2). It is achieved with high probability by i.i.d. sampling points with distribution $\mu$ and the heart of the proof is an application of Dudley's entropy bound. We also refer to the results of Kosov and Tikhonov (2025) on discretization and approximation in (Orlicz) norms different from $L_p$.

The theorem is best in situations where the entropy numbers decay slowly. It cannot give a bound smaller than $n^{-1/(2p)}$, no matter how fast the $\varepsilon_n$ decay, and it does not give any rate for $p = \infty$. However, it shows a rate of convergence of the sampling numbers for any $p < \infty$ under very weak assumptions on the entropy numbers. Especially for $L_1$-approximation, the rate of decay of the sampling numbers is at least as the rate of the entropy numbers in $B(D)$ in the regime $s \in (0, 1/2)$.

In particular, Theorem 5.8 can improve Theorem 3.17 in cases where $d_n(F, B(D))$ decays slowly, at least for convex and symmetric $F$ and for possibly nonlinear sampling algorithms. To illustrate this, let us add a corollary. Combining Carl's inequality with Theorem 5.8, we obtain the following.

**Corollary 5.9.** Let $D \subset \mathbb{R}^d$ be compact, $F \subset C(D)$ be compact, convex and symmetric, and $(D,\mu)$ be a probability space. Then, for all $1 \le p < \infty$ and $0 < s \neq \frac{1}{2}$, we have

$$d_n(F, B) \asymp n^{-s} \quad \Longrightarrow \quad g_n(F, L_p) \lesssim \begin{cases} n^{-s/p}, & \text{for } s < \frac{1}{2}, \\ n^{-1/(2p)}, & \text{for } s > \frac{1}{2}. \end{cases}$$

Note that Theorem 3.17 only gives a meaningful upper bound for $s > \frac{1}{2} - \frac{1}{p}$, while Corollary 5.9 also works for smaller $s$. In the case $p = 1$ and $s < \frac{1}{2}$, Corollary 5.9 and Theorem 3.17 give the same result. However, the entropy numbers might be smaller than the Kolmogorov numbers, and hence Corollary 5.9 may be weaker than Theorem 5.8. However, for classical examples, we did not find a case where this matters when we consider the resulting bounds on the sampling numbers.

**Example 5.10.** Let us illustrate the bound of Theorem 5.8 for an easy example. For $0 < s \leq 1$, we consider the class of univariate Hölder-continuous functions

$$F^s := \left\{ f\colon [0,1] \to [-1,1] : \sup_{x \neq y} \frac{|f(x) - f(y)|}{|x-y|^s} \leq 1 \right\}.$$

It is known and easy to prove (by considering appropriate piecewise constant functions) that $\varepsilon_n(F^s, B([0,1])) \asymp n^{-s}$. Theorem 5.8 hence gives

$$g_n(F^s, L_p) \lesssim \begin{cases} n^{-s/p}, & \text{for } s < \frac{1}{2}, \\ n^{-1/(2p)}, & \text{for } s > \frac{1}{2}. \end{cases}$$

On the other hand, a simple piecewise linear interpolation at equidistant nodes (for the upper bound) together with considering *fooling functions* of the form $f(x) = \min_{i \leq n} |x - x_i|^s$ (for the lower bound) gives

$$g_n(F^s, L_p) \asymp n^{-s}.$$

So Theorem 5.8 gives a sharp result if and only if $s < 1/2$ and $p = 1$.

**Example 5.11.** An interesting example is given by $L_p$-approximation for the family $F = W_p^s$ of (univariate) Sobolev classes with the same integrability parameter. Here, we discuss Sobolev spaces with fractional smoothness, which are defined using *Fourier methods* or *differences*, and usually require the restriction $1 < p < \infty$. We omit the details and refer to Novak and Triebel (2006) for precise definitions and more results about sampling numbers of Sobolev classes. It is known that

$$\varepsilon_n(W_p^s, B) \asymp n^{-s} \quad \text{and} \quad d_n(W_p^s, B) \asymp n^{-s+(1/p-1/2)_+}$$

for $s > 1/p$, where $B := B([0,1])$, see (Pietsch 2007, 6.7.8.13 & 6.7.8.15) or Appendix A.2. Considering the regime $p > 2$ and $\frac{1}{p} < s < \frac{1}{2}$ of 'small smoothness', Theorem 3.17 and Theorem 5.8 give

$$g_n^{\text{lin}}(W_p^s, L_p) \lesssim n^{-s+1/2-1/p} \quad \text{versus} \quad g_n(W_p^s, L_p) \lesssim n^{-s/p},$$

respectively. The second bound is smaller for $s < \frac{1}{2} - \frac{1}{2(p-1)}$, the first is smaller for $s > \frac{1}{2} - \frac{1}{2(p-1)}$. On the other hand, for this basic example, it is known that

$$g_n(W_p^s, L_p) \asymp g_n^{\text{lin}}(W_p^s, L_p) \asymp n^{-s}$$

for all $s > 1/p$ and $p$, see Proposition 2.4. So none of the above-mentioned bounds is optimal for $p > 2$. So far, we only obtain optimal bounds in this scale for $p = 2$ from Theorem 3.17. It seems an important and challenging open problem to extend the general analysis of sampling numbers in some way that gives optimal bounds for the whole family of these basic examples.

There is another bound on the sampling numbers in terms of entropy for uni-form approximation that also reflects fast decay of $\varepsilon_n$. However, and similar to Theorems 3.9 and 3.16, it comes with an additional factor that makes it only useful if the decay is fast enough. The following result of Ullrich (2026) is based on a comparatively elementary argument.

**Theorem 5.12.** Let $D$ be a set and $F \subset B := B(D)$ be convex and bounded. Then, for every $n \in \mathbb{N}_0$, we have

$$g_n(F, B) \;\leq\; 2\,(n+1) \cdot \varepsilon_n(F, B).$$

Moreover, if $F = B_X$ is the unit ball of a real $p$-Banach space $X$ with $0 < p \leq 1$, then

$$g_n(F, B) \;\leq\; (n+1)^{1/p} \cdot \varepsilon_n(F, B).$$

We also add that Ullrich (2026) considered real-valued functions, and thereby saved the factor 2. For comparison with Theorem 5.8, note that $g_n(F, L_p) \leq g_n(F, B)$ whenever $(D, \mu)$ is probability space.

The proof of this result is almost copy-and-paste from (Pietsch 1978, Theorem 12.3.2), it only employs the additional observation that the proof there only requires functionals from a norming set. This is part of the following lemma.

**Lemma 5.13.** Let $D$ be a set, $F \subset B(D)$ be bounded. Then, for every $\delta > 0$ and $n \in \mathbb{N}_0$, we can find $f_0, h_0, \ldots, f_n, h_n \in F$ and $x_0, \ldots, x_n \in D$ such that, for $k = 0, \ldots, n$, we have

$$f_k(x_j) = h_k(x_j) \quad \text{for } j < k, \qquad \text{and} \qquad |f_k(x_k) - h_k(x_k)| \geq g_k(F, B) - \delta.$$

If $F$ consists of real-valued functions, then $g_k(F, B)$ can be replaced by $2 \cdot g_k(F, B)$.

*Proof.* The proof is by induction. Let $k \in \{0, \ldots, n\}$ and assume that $f_j, h_j$ and $x_j$ for $j < k$ have already been found. Define

$$M_k \;:=\; \left\{ (f, h) \in F \times F \colon\, f(x_j) = h(x_j) \quad \text{for } \; 0 \leq j < k \right\}.$$

Then it holds that

$$g_k(F, B) \leq \sup_{(f,h) \in M_k} \left\| f - h \right\|_\infty.$$

This error bound is achieved by any *interpolatory algorithm* $A_k$ of the form $A_k(f) = \Phi(f(x_1), \ldots, f(x_k))$ with the property

$$A_k(f) \;\in\; \{h \in F \colon\, h(x_i) = f(x_i) \text{ for all } 0 \leq j < k\}$$

for all $f \in F$. Hence, we can choose $f_k, h_k \in F$ with $(f_k, h_k) \in M_k$ and

$$\|f_k - h_k\|_\infty \;>\; g_k(F, B) - \delta.$$

This gives $x_k \in D$ with

$$|f_k(x_k) - h_k(x_k)| \;>\; g_k(F, B) - \delta.$$

This finishes the proof. The improvements for real-valued $F$ are due to the existence of a central algorithm, as has been explained, e.g., by Ullrich (2026).

□

We can now prove Theorem 5.12.

*Proof.* First note that for every $\rho < g_n(F,B)$, by Lemma 5.13, there exist $f_0, g_0, \dots, f_n, g_n \in F$ and $x_0, \dots, x_n \in D$ such that $f_k(x_j) = g_k(x_j)$ for $j < k$ and $|f_k(x_k) - g_k(x_k)| > \rho$ for $k = 0, \dots, n$.

For $\xi = (\xi_0, \dots, \xi_n) \in \{0,1\}^{n+1}$, we define the functions

$$h_\xi := \frac{1}{n+1} \sum_{i=0}^{n} (\xi_i f_i + (1-\xi_i) g_i) .$$

By convexity, $H_n := \{h_\xi : \xi \in \{0,1\}^{n+1}\} \subset F$ with $\# H_n = 2^{n+1}$.

For $\xi, \xi' \in H_n$, $\xi \neq \xi'$, let $k := \min\{j : \xi_j \neq \xi'_j\}$. Then,

$$\begin{aligned} \left| h_\xi(x_k) - h_{\xi'}(x_k) \right| &= \frac{1}{n+1} \left| \sum_{i=k}^{n} (\xi_i - \xi'_i) f_i(x_k) + (\xi'_i - \xi_i) g_i(x_k) \right| \\ &= \frac{1}{n+1} \left| \xi_k - \xi'_k \right| |f_k(x_k) - g_k(x_k)| \geq \frac{\rho}{n+1}. \end{aligned}$$

This implies $\left\| h_\xi - h_{\xi'} \right\|_\infty \geq \frac{\rho}{n+1}$ for all $\xi \neq \xi' \in \{0,1\}^{n+1}$.

Since a union of $2^n$ balls can only cover $F$ if one of the balls contains more than one function from $H_n$, we must have $\rho \leq 2(n+1)\varepsilon_n(F,B)$ for all $\rho < g_n(F,B)$, which yields the result.

The proof of the second part works almost exactly as above with the only difference that the $(n+1)$ in the definition of $h_\xi$ must be replaced by $(n+1)^{1/p}$ to ensure $h_\xi \in B_X$. Since we assume $X$ to be real, we also avoid the factor 2.

□

### 5.3. *Bounds by other benchmarks*

There is a big variety of approximation benchmarks in the literature. There are Kolmogorov widths, best $m$-term widths and entropy numbers, which we all encountered before, but also approximation numbers, Gelfand widths, manifold widths, Bernstein numbers, etc. We now discuss a bound on the sampling numbers by another approximation benchmark, which is of particular theoretical interest.

Indeed, many of the aforementioned (and other) benchmark quantities fall into the family of s-numbers, as introduced by Pietsch (1974), and further discussed in Section 9.2; others do not fall directly into this family but are still bounded from above and below by certain s-numbers. It is hence quite desirable to find upper bounds for the sampling numbers in terms of the smallest possible s-numbers. This way, we get upper bounds in terms of *any* of the aforementioned benchmarks.

The smallest s-numbers are the *Hilbert numbers* $h_n(F,Y)$ defined below. Hence, the main goal in this section is to bound the sampling numbers by the Hilbert numbers. But before we even come to the Hilbert numbers, we define a variant of the **Grothendieck numbers** of $F$ in $B(D)$, namely, for a convex and symmetric subset $F \subset B(D)$, we put

$$\Gamma_n(F) \;:=\; \sup\left\{ \det\Big( \big(f_i(x_j)\big)_{i,j=1}^{n+1} \Big)^{\frac{1}{n+1}} \;:\; f_k \in F,\; x_k \in D \right\}. \tag{5.19}$$

This definition differs from the *original* Grothendieck numbers, which are defined using the supremum over arbitrary continuous linear functionals $L_j(f_i)$, $L_j \in B(D)'$, instead of only function evaluations $f_i(x_j)$. This is a crucial concept for the analysis of (finite-dimensional) Banach spaces, see the monograph of Pietsch (2007, 6.1.11.4 & 6.3.13.10). In the above form, it was also used by Novak (1995*a*, page 127). For general subsets $F \subset B(D)$, we simply put

$$\Gamma_n(F) \;:=\; \Gamma_n\left(\mathrm{co}(F)\right),$$

where $\mathrm{co}(F)$ is the convex and symmetric hull of $F$. Note that, if $F$ is convex, then we have the simple relation $\mathrm{co}(F) = \frac{1}{2}(F-F)$. We obtain the following simple lemma.

**Lemma 5.14.** Let $D$ be a set and $F \subset B := B(D)$ be bounded. Then, for all $n \in \mathbb{N}_0$, we have

$$g_n(F,B) \;\le\; 2\cdot \Gamma_n(F).$$

We note that the factor 2 is not needed for real-valued functions.

*Proof.* First note that for every $\rho_k < g_{k-1}(F,B)$, $k=1,\dots,n+1$, by Lemma 5.13, there exist $f_1, g_1, \dots, f_{n+1}, g_{n+1} \in F$ and $x_1,\dots,x_{n+1} \in D$ such that $f_k(x_j) = g_k(x_j)$ for $j<k$ and $|f_k(x_k) - g_k(x_k)| > \rho_k$. We define $p_k := \frac{f_k - g_k}{2} \in \mathrm{co}(F)$ for $k=1,\dots,n+1$. By construction, we have $p_k(x_j)=0$ for $j<k$, i.e., $\big(p_k(x_j)\big)_{i,j}$ is a triangular matrix, and $2|p_k(x_k)| > \rho_k$ for $k=1,\dots,n+1$. The determinant of a triangular matrix is the product of its diagonal entries, and therefore

$$\prod_{k=1}^{n+1} \frac{\rho_k}{2} \le \det\Big( \big(p_i(x_j)\big)_{i,j=1}^{n+1} \Big) \le \Gamma_n(F)^{n+1}.$$

Taking the supremum over all $\rho_k < g_{k-1}(F,B)$, and $g_n \le g_k$ for $k \le n$, proves the claim.

□

The Grothendieck numbers (of $F$ in $B(D)$) can be bounded by various quantities, including the Hilbert numbers. The **Hilbert numbers** of $F \subset Y \subset \mathbb{C}^D$ are defined

by

$$h_n(F,Y) := \sup\Big\{\sigma_n\big(RA\big)\colon R\in\mathcal{L}(Y,\ell_2) \text{ with } \|R\|\le 1, \text{ and } A\in\mathcal{L}(\ell_2,Y) \text{ with } A(B_{\ell_2})+x\subset F \text{ for some } x\in F\Big\}, \tag{5.20}$$

where $\sigma_n$, $n\in\mathbb{N}_0$, is the $(n+1)$-st largest singular value of the operator $RA\in\mathcal{L}(\ell_2,\ell_2)$. In particular, $\sigma_0(RA)=\|RA\|$. (In fact, it would be enough in the above supremum to consider operators on $\ell_2^{n+1}$, i.e., matrices $RA\in\mathbb{C}^{(n+1)\times(n+1)}$. We omit the details.)

The Hilbert numbers do not look very handy at first. One may interpret them as *minimal errors for the hardest Hilbert-sub-problem*, and as such they allow us to use quite some *heavy machinery*. (Note that numerical analysis, and functional analysis in general, are usually much easier over Hilbert spaces.)

These numbers have been introduced by Bauhardt (1977) as the smallest *s-number* (of an operator) in the sense of Pietsch (1974). We refer to Section 9.2 for a detailed discussion of this subject. Here, we use a variant that has been introduced by Krieg, Novak and Ullrich (2025*a*) for the purpose of treating approximation over more general sets $F$ (the original definition corresponds to $F$ being the unit ball of some Banach space). The $h_n$ are not only the smallest s-numbers, but also a lower bound for other approximation quantities. In particular, let us mention the bounds

$$h_n(F,B)\le\min\big\{d_n(F,B),\ \delta_n(F,B),\ 2\cdot\varepsilon_n(F,B)\big\} \tag{5.21}$$

with $B:=B(D)$, whenever $F$ is the unit ball of a (real) normed space, with the Kolmogorov widths $d_n$ from (3.16), the entropy numbers $\varepsilon_n$ from (5.16), and the manifold widths $\delta_n$, see Section 9.2 (Proposition 9.12 and Remark 9.17). The **manifold widths** of $F\subset Y$ are defined by

$$\delta_n(F,Y) := \inf_{\substack{N_n\in C(F,\mathbb{R}^n)\\ \Phi\in C(\mathbb{R}^n,Y)}}\ \sup_{f\in F}\big\|f-\Phi\left(N_n(f)\right)\big\|, \tag{5.22}$$

where $C(X,Y)$ denotes the class of continuous mappings from $X$ to $Y$. (Here, a topology on $F$ is required. If $F$ is the unit ball of a normed space, we use the norm topology.) These represent the worst-case error that can be achieved with a continuous encoder and decoder. The manifold widths have been introduced by DeVore, Kyriazis, Leviatan and Tikhomirov (1993), though they are equivalent up to constants to a concept already introduced by Alexandroff (1956). We refer to Stesin (1975) for an early treatment and DeVore, Howard and Micchelli (1989) for more on nonlinear widths. The following result of Krieg *et al.* (2025*a*, Theorem 1.2) hence also implies bounds by the above quantities. However, this bound is purely theoretical, we do not know explicitly any algorithm that achieves the bound.

**Theorem 5.15.** Let $D$ be a set and $F\subset B:=B(D)$ be convex and bounded.

Then, for all $n \in \mathbb{N}_0$, we have

$$g_n(F,B) \ \le\ 2\,(n+1)^{3/2} \cdot \left( \prod_{k=0}^{n} h_k(F,B) \right)^{\frac{1}{n+1}}.$$

If $F$ is the unit ball of a Banach space[1], then we can replace the factor $(n+1)^{3/2}$ with $(n+1)$.

As above, the factor 2 is not needed for real-valued functions.

The geometric mean appearing on the right hand side of Theorem 5.15 is often of the same order of magnitude as the number $h_n(F,B)$, for instance, if the Hilbert numbers are of polynomial decay. In particular, if $F$ is the unit ball of a Banach space, we can use $n^\alpha \le e^\alpha\,(n!)^{\alpha/n}$ for $\alpha \ge 0$, to obtain an inequality

$$g_n(F,B) \ \le\ 2\,e^\alpha\, n^{-\alpha+1} \ \cdot \sup_{k\le n}\,(k+1)^\alpha \cdot s_k(F,B), \tag{5.23}$$

for all $\alpha > 1$ and $s_k \in \{h_k, 2\cdot\varepsilon_k, \delta_k, d_k\}$.

*Proof.* By Lemma 5.14, it remains to bound $\Gamma_n$ in terms of the $h_n$. For this, let $\varepsilon > 0$, $p_k := \frac{f_k - g_k}{2} \in \mathrm{co}(F) = \frac{F-F}{2}$ with $f_1, g_1, \ldots, f_{n+1}, g_{n+1} \in F$, and $x_1, \ldots, x_{n+1} \in D$ be given such that

$$\det\Big( \big(p_i(x_j)\big)_{i,j=1}^{n+1} \Big)^{1/(n+1)} \ge \Gamma_n(F) - \varepsilon.$$

Moreover, let $p_0 := \frac{1}{n+1}\sum_{i\le n+1} \frac{f_i+g_i}{2}$. Defining the operators

$$A(\xi) \ :=\ \frac{1}{n+1}\sum_{i=1}^{n+1} \xi_i p_i \in Y, \quad \xi = (\xi_i)_{i\le n+1} \in \ell_2^{n+1}, \tag{5.24}$$

and

$$R(f) \ :=\ \frac{1}{\sqrt{n+1}} \big(f(x_i)\big)_{i\le n+1} \in \ell_2^{n+1}, \qquad f \in Y,$$

we observe that $\|R\colon B(D) \to \ell_2^{n+1}\| \le 1$ and, for all $\xi \in [-1,1]^{n+1}$, due to convexity, it holds

$$A(\xi) + p_0 \ =\ \frac{1}{n+1}\sum_{i=1}^{n+1} \left( \frac{1+\xi_i}{2} f_i + \frac{1-\xi_i}{2} g_i \right) \in F$$

and $p_0 \in F$. In particular, $A(B_{\ell_2^{n+1}}) + p_0 \subset F$. Since we have the matrix represent-

[1] Theorems 1.2 and 6.4 in (Krieg *et al.* 2025*a*) additionally state that the factor $(n+1)^{3/2}$ can be replaced by $\sqrt{n+1}$ if $F$ is the unit ball of a Hilbert space, but the claim that the proof of Theorem 3.3 transfers to this setting is not correct. It is hence open whether the statement about the improved factor in the Hilbert case is correct.

ation

$$RA = (n+1)^{-3/2}\big(p_i(x_j)\big)_{i,j=1}^{n+1},$$

we obtain for $k = 0, \dots, n$,

$$\sigma_k\Big(\big(p_i(x_j)\big)_{i,j=1}^{n+1}\Big) = (n+1)^{3/2}\,\sigma_k\big(RA\big) \le (n+1)^{3/2}\,h_k(F,B).$$

The determinant is the product of the singular values and therefore

$$\Gamma_n(F) - \varepsilon \le \det\Big(\big(p_i(x_j)\big)_{i,j=1}^{n+1}\Big)^{1/(n+1)} \le (n+1)^{3/2}\left(\prod_{k=0}^{n} h_k(F,B)\right)^{1/(n+1)}.$$

This proves the claim with $\varepsilon \to 0$. If $F$ is the unit ball of a Banach space, then we can use $A(\xi) := \frac{1}{\sqrt{n+1}}\sum_{i=1}^{n+1}\xi_i p_i$ and $p_0 = 0$, to have $A(B_{\ell_2^{n+1}}) \subset F$. The remaining proof is as above with $(n+1)^{3/2}$ replaced by $(n+1)$.

□

**Remark 5.16.** In Sections 3 and 5, we gave direct upper bounds for the sampling numbers $g_n(F,B)$ in terms of the numbers $d_n(F,B)$, $\varepsilon_n(F,B)$, and $h_n(F,B)$. Any of the three quantities can be bounded in some way by any of the other two, see Pietsch (2007, 6.2.3.14 & 6.2.4.5). So, indirectly, we obtain various additional upper bounds for the sampling numbers in terms of each of these three quantities. The bounds obtained in this way are worse than the direct bounds, so we spare the reader the details. A message of this is, however, that there is no strict hierarchy between the three approaches.

## 6. Almost-optimality of i.i.d. sampling

As we have seen, and as the experienced reader will know, it is often very difficult and not known how to obtain optimal sampling points. This is especially the case for high-dimensional model classes. In this section, we discuss how good random sampling points are compared with optimally chosen sampling points. That is, we ask whether random points, once drawn and fixed, do allow for nearly-optimal algorithms. Ideally, with high probability and based on a feasible distribution of points. In particular, it turns out that i.i.d. sampling points often suffice, sometimes even uniformly distributed in the domain. Since using randomness in numerical software is typically not considered an issue, this gives a *semi-constructive* approach to optimal sampling. In some cases, there is even a feasible way to check whether the random sampling points, once they are drawn, satisfy the desired (near-)optimality. So, in this case, the random approach may even be fully constructive.

Apart from the semi-constructive way of obtaining optimal sampling points, there is a second motivation to study the quality of random sampling points. Indeed, in many applications, the data is already given and the sampling points cannot be chosen freely. We think that random sampling points, and especially i.i.d. samples, are a good and popular model of the way this given data has been obtained (for instance, it is a standard assumption in the area of machine learning). The study of the quality of such i.i.d. samples will tell us how good 'typical' sampling points are compared to optimal sampling points and whether or not one should expect good results when working with the given data, or whether it is worth trying to acquire new data at carefully chosen positions. Note that we consider random data $y_i = f(x_i)$ with random $x_i$ and without additional additive noise. But recall that many of our results can be extended to noisy data, see Remark 3.6.

It is quite common in (numerical) analysis to prove existence of an object via the **probabilistic method**. That is, instead of deterministic construction, a random construction is provided that satisfies the desired properties with positive or even high probability. This implies the existence of a deterministic construction. We have seen such an example, for instance, when we discussed the restricted isometry property for random sampling points in Section 5.1. Therefore, many results about the quality of random sampling points are (partially hidden) already contained in the literature.

Personally, we started our work about the quality of random information in Hinrichs *et al.* (2020) and Hinrichs *et al.* (2021*b*). A more recent survey is given by Sonnleitner and Ullrich (2023). In this section, we will proceed as follows:

- We illustrate the question for the introductory example of Lipschitz functions in Section 6.1.
- We discussed in Section 3.3 how $\mathcal{O}(m)$ *optimal samples* can result in an error close to the error of best $L_2$-approximation from an $m$-dimensional space. In Section 6.2, we discuss how many i.i.d. samples with respect to a given distribution are required to achieve such error bounds. In particular, we identify sampling distributions for which $\mathcal{O}(m \log m)$ i.i.d. samples suffice.

- In general, the extra logarithmic factor cannot be avoided when using i.i.d. sampling points. In Section 6.3, we identify some situations where $\mathcal{O}(m)$ i.i.d. samples suffice, namely, $L_p$-approximation in certain Sobolev spaces.

**Remark 6.1 (Randomized algorithms).** Let us point out that we do *not* discuss randomized algorithms in this section. The sampling points, and therefore the corresponding sampling algorithm, are only randomly drawn once, and should then be "good" for all functions from the given model class. In contrast, when considering randomized algorithms, we allow to draw a new random sample for each input. If the random sampling points are only to be used for a single function, the algorithmic procedure is the same. In this case, there is no need to demand that the random sample is good for all functions from the model class *at once* and it is enough to study the weaker Monte Carlo error criterion as discussed in Section 7. For some general motivation as to why we consider the worst case setting, we also refer to Remark 2.6. Also in the Monte Carlo setting, it is of interest to compare i.i.d. samples with optimal random samples. We refer to Remark 7.4 for a little more information and references on this research question. More general randomized algorithms are discussed in Sections 9 and 10.2.

### *6.1. The easy example revisited*

We discuss again the introductory example of approximating univariate periodic Lipschitz functions from the class $F_{\mathrm{Lip}}$ as defined in (2.1). This time, however, we assume that the sampling points cannot be chosen optimally, but instead are given as i.i.d. random points. In this example, it is natural to consider points $\mathcal{P}_n = \{x_1, \ldots, x_n\} \subset [0,1]$ that are independent and uniformly distributed (i.i.d.) in $D = [0,1]$. This example is from Hinrichs *et al.* (2020).

First, we recall from Section 2 that the minimal worst-case error that can be achieved based on the sampling points $\mathcal{P}_n$ is given by the radius of information

$$\mathrm{rad}(\mathcal{P}_n, F, L_p) := \inf_{\Phi\colon \mathbb{R}^n \to L_p} \sup_{f \in F} \big\| f - \Phi\big(f(x_1), \ldots, f(x_n)\big) \big\|_p .$$

In order to characterize the (typical) quality of the i.i.d. random sampling points, we hence study the **expected radius of random information**.

The model class $F_{\mathrm{Lip}}$ is very convenient for studying the power of random information because of the precise formula for the radius of information for every point set $\mathcal{P}_n$ in terms of the distance function from (2.8), i.e.,

$$\mathrm{rad}(\mathcal{P}_n, F_{\mathrm{Lip}}, L_p) = \| \operatorname{dist}(\cdot, \mathcal{P}_n) \|_p .$$

This makes it possible to compute precisely the $p^{\mathrm{th}}$ moment of the minimal worst-case error $\mathrm{rad}(\mathcal{P}_n, F_{\mathrm{Lip}}, L_p)$ and the expected value of $\mathrm{rad}(\mathcal{P}_n, F_{\mathrm{Lip}}, L_\infty)$.

**Proposition 6.2.** Let $F_{\mathrm{Lip}}$ from (2.1), and let $\mathcal{P}_n$ consist of independent and uniformly distributed points $x_1, \ldots, x_n \in [0,1]$. Then,

$$\sqrt[p]{\mathbb{E}\left[\mathrm{rad}(\mathcal{P}_n, F_{\mathrm{Lip}}, L_p)^p\right]} \;=\; \frac{1}{2}\sqrt[p]{\frac{n!}{(p+1)\cdots(p+n)}} \;\approx\; \frac{1}{2n}\sqrt[p]{\Gamma(p+1)}$$

for all $1 \le p < \infty$, and for $p = \infty$,

$$\mathbb{E}\left[\mathrm{rad}(\mathcal{P}_n, F_{\mathrm{Lip}}, L_\infty)\right] \;=\; \frac{1}{2n}\sum_{i=1}^{n}\frac{1}{i} \;\approx\; \frac{1}{2n}\log n.$$

*Proof.* We first consider the case $p < \infty$. Let $\mathcal{P}_n = \{x_1, \ldots, x_n\}$ and write $\mathrm{rad}(\mathcal{P}_n) := \mathrm{rad}(\mathcal{P}_n, F_{\mathrm{Lip}}, L_p)$. Recall that

$$\mathrm{rad}(\mathcal{P}_n)^p = \int_{[0,1]} \mathrm{dist}(x, \mathcal{P}_n)^p \;\mathrm{d}x$$

and note $\mathrm{dist}(\cdot,\cdot)\colon [0,1] \times [0,1]^n \to \mathbb{R}$ is continuous and thus Borel measurable. Using Tonelli's theorem, we obtain

$$\mathbb{E}\left[\mathrm{rad}(\mathcal{P}_n)^p\right] = \int_{[0,1]} \mathbb{E}\left[\,\mathrm{dist}(x, \mathcal{P}_n)^p\right] \;\mathrm{d}x.$$

We will show that the integrand of the latter integral is constant. Let us fix $x \in [0,1]$ and note that $\mathrm{dist}(x, \mathcal{P}_n) \in [0, 1/2]$. For any $t \in [0, 1/2]$, we have

$$\mathrm{dist}(x, \mathcal{P}_n) \ge t \quad\Longleftrightarrow\quad \forall i \in \{1, \ldots, n\} : x_i \notin [x - t, x + t].$$

Thus

$$\mathbb{P}[\mathrm{dist}(x, \mathcal{P}_n) \ge t] = \left(1 - 2t\right)^n.$$

The substitution $s = 1 - 2t^{1/p}$ and integration by parts yield

$$\begin{aligned}
\mathbb{E}\left[\,\mathrm{dist}(x, \mathcal{P}_n)^p\right] &= \int_0^{2^{-p}} \mathbb{P}[\mathrm{dist}(x, \mathcal{P}_n)^p \ge t] \;\mathrm{d}t = \int_0^{2^{-p}} \left(1 - 2t^{1/p}\right)^n \;\mathrm{d}t \\
&= \frac{p}{2^p}\int_0^1 s^n (1-s)^{p-1} \;\mathrm{d}s = \frac{1}{2^p}\frac{n!}{(p+1)\cdots(p+n)},
\end{aligned}$$

which implies the statement.

In the case $p = \infty$, we have from (2.8) that

$$\mathrm{rad}(\mathcal{P}_n) \;=\; \|\,\mathrm{dist}(\cdot, \mathcal{P}_n)\,\|_\infty,$$

so the minimal worst-case error for a given point set $\mathcal{P}_n$ is equal to half of the largest distance between neighboring points. That is, we have to study the expected value of the largest distance between $n$ random points on the circle, or equivalently, the longest piece of a stick of unit length that we break at $n - 1$ random places (the first random point on the circle turns the circle into an interval). This is of course an old and well studied problem, see, e.g., (Holst 1980). The longest piece of the broken stick has an expected length

$$\frac{1}{n}\sum_{i=1}^{n}\frac{1}{i} \approx \frac{\log n}{n}.$$

□

By recalling the formulas for the error achieved by optimal sampling points, as given in Proposition 2.3, we arrive at the following conclusion.

**Corollary 6.3.** Let $F_{\mathrm{Lip}}$ from (2.1) and $x_1, x_2, \ldots$ be independent and uniformly distributed points on $[0,1]$. Consider the random $n$-point set $\mathcal{P}_n = \{x_1, \ldots, x_n\}$. Then,

$$\mathbb{E}\left[\mathrm{rad}(\mathcal{P}_n, F_{\mathrm{Lip}}, L_p)\right] \asymp \begin{cases} g_n(F_{\mathrm{Lip}}, L_p) & \text{if } p < \infty, \\ g_{n/\log n}(F_{\mathrm{Lip}}, L_p) & \text{if } p = \infty. \end{cases}$$

That is, in expectation, i.i.d. uniformly distributed points are optimal up to a constant for $L_p$-approximation with $p < \infty$ and optimal up to a logarithmic factor for $L_\infty$-approximation.

*Proof.* By Proposition 2.3 and Proposition 6.2, we have

$$g_n(F_{\mathrm{Lip}}, L_p) := \inf_{\mathcal{P}_n} \mathrm{rad}(\mathcal{P}_n, F_{\mathrm{Lip}}, L_\infty) \asymp \frac{1}{n} \asymp \mathbb{E}\left[\mathrm{rad}(\mathcal{P}_{n\log n}, F_{\mathrm{Lip}}, L_\infty)\right].$$

In the case $p < \infty$, we write $\mathrm{rad}(\mathcal{P}_n) := \mathrm{rad}(\mathcal{P}_n, F_{\mathrm{Lip}}, L_p)$ and use Jensen's inequality to get

$$\frac{1}{n} \lesssim \inf_{\mathcal{P}_n} \mathrm{rad}(\mathcal{P}_n) \lesssim \mathbb{E}\left[\mathrm{rad}(\mathcal{P}_n)\right] \lesssim \sqrt[p]{\mathbb{E}\left[\mathrm{rad}(\mathcal{P}_n)^p\right]} \lesssim \frac{1}{n}.$$

□

The loss in the case $p = \infty$ above is related to the so-called *coupon collector's problem*, and can also be interpreted in the following way: If we divide the domain into $m$ equal "boxes" and take $n$ i.i.d. uniformly distributed points in the domain, we must have $n$ of order $m \log m$ in order to have a good chance of hitting all the boxes; and we must hit all the boxes in order to have a a good approximation in the uniform norm. This is different for $L_p$-approximation: Here it is no problem if we miss a few boxes as long as we hit most of them.

The interpretation in terms of the coupon collector problem is helpful because the latter also applies to various other examples, leading to the necessity of a logarithmic oversampling when i.i.d. sampling points have to be used. For instance, it also appears for uniform approximation of multivariate $C^k$-functions as discussed by Bauer, Devroye, Kohler, Krzyżak and Walk (2017), and $L_p$-approximation on Sobolev spaces $W_q^s$ on domains and manifolds with $q \leq p$ as discussed by Krieg and Sonnleitner (2025, 2024). I.i.d. samples for $L_p$-approximation on Sobolev spaces will be discussed further in Section 6.3. First, however, we want to give some results for general model classes $F$ similar to the results from Section 3.

### *6.2. Least squares and i.i.d. sampling*

We discuss the quality of random sampling points that are chosen i.i.d. with respect to a given distribution. We start by looking at the problem of $L_2$-approximation in a *reproducing kernel Hilbert space*. As in Section 3.3, results in this setting can then be used to obtain error bounds also for non-Hilbert classes and in error norms other than $L_2$.

We shortly recall some basics of the Hilbert setting, see Section 3.3.1 for more details. We assume that we are given a reproducing kernel Hilbert space $H$ of complex-valued functions on some non-empty set $D$, whose kernel we denote by $k_H$. We discuss approximation in $L_2 = L_2(D, \mu)$, where $\mu$ is a measure on $D$, and assume that we have an injective embedding $H \hookrightarrow L_2$. For our analysis, we again need that the embedding has a finite trace, i.e., that the kernel is integrable on the diagonal (3.34). In this situation, the embedding $H \hookrightarrow L_2$ is also compact and, by the singular value decomposition, there exists a unique square-summable decreasing sequence $(\sigma_k)_{k\in\mathbb{N}_0}$ (the singular values) and a countable or finite orthogonal basis $(b_k)_k$ of $H$ which is orthonormal in $L_2$, and which gives the Mercer decomposition (3.35) of the kernel.

We assume that we are given i.i.d. samples with respect to a $\mu$-density $\varrho$, which, in general, we cannot choose by ourselves. However, it is helpful to first look at the case where we can choose the sampling density. In this situation, we know an optimal choice of the density, which is given by

$$\varrho_m^*(x) := \frac{1}{2}\left(\frac{1}{m}\sum_{k<m}|b_k(x)|^2 + \frac{1}{\sum\limits_{k\geq m}\sigma_k^2}\sum_{k\geq m}\sigma_k^2|b_k(x)|^2\right) \tag{6.1}$$

for given $m \in \mathbb{N}$. By 'optimal', we mean that using $\mathcal{O}(m\log m)$ i.i.d. samples with respect to this distribution gives an error bound of the optimal order $\sigma_m$ with high probability (as we shall see below), whereas, in general, no sampling distribution can result in an error of this order with $o(m\log m)$ i.i.d. samples (as we shall see in Section 6.3). This optimal density was introduced in Krieg and Ullrich (2021*a*) as part of the probabilistic method to prove the *existence* of $\mathcal{O}(m\log m)$ samples that give an error of order $\sigma_m$, see Ullrich (2020), Krieg and Ullrich (2021*b*), Kämmerer, Ullrich and Volkmer (2021), Dolbeault *et al.* (2023) for further uses of this density. We add that the density $\rho_m^*$ sometimes reduces to a 'simple' one. For example, it is just constant if $(b_k)_k$ are trigonometric monomials on the unit cube or spherical harmonics on a manifold, but it is more complicated (and often intractable) in other settings, see Examples 6.7–6.9 below.

We attempt to give an intuition as to why this density is good for least squares approximation on $V_m = \operatorname{span}\{b_0, \dots, b_{m-1}\}$ of functions from $H$. This seems easier if we instead look at this problem as a problem of recovering functions $g = P_m f \in V_m$ (the orthogonal projection of $f$) from noisy samples (since we only have function values of $f$ instead of $g$). This gives the following intuition:

- The first part of the sum (6.1) is the norm of function evaluation at $x$ on $(V_m, \|\cdot\|_2)$. It therefore ensures that we take 'enough' of those samples that would give the important information for recovering $g$ with exact information.
- The second part of the sum (6.1) is the norm of function evaluation at $x$ on $(V_m^{\perp}, \|\cdot\|_H)$. If this norm is large, it means that function values of $f \in H$ are only badly approximated by function values of the orthogonal projection $g \in V_m$. Hence, in such areas, we have large 'noise' and need more samples in order to counterbalance it.

For general densities $\varrho$, we can say something about the quality of i.i.d. samples under the condition

$$M \;:=\; \inf\left\{c\colon\, \rho_m^*(x) \le c \cdot \rho(x) \quad \text{for all } \; x \in D\right\} \;<\; \infty. \tag{6.2}$$

Roughly speaking, we show that $\mathcal{O}(Mm \log(Mm))$ i.i.d. samples with density $\varrho$ give an error of the optimal order $\sigma_m$. This is foremost a theoretical insight since it is hard to determine the characteristic $M$ for concrete examples; it is usually already very difficult to find the singular value decomposition $(\sigma_k, b_k)_k$ for a given RKHS $H$. But let us already point out that in the non-Hilbert setting discussed afterwards, we can choose the basis $(b_k)$ and the numbers $(\sigma_k)$ freely, at least to some extent, and hence also obtain control over the numbers $M$.

**Proposition 6.4.** Let $(D, \mu)$ be a measure space and let $H$ be a RKHS that is (injectively) embedded into $L_2(D, \mu)$ with square-summable singular values $(\sigma_k)_{k\in\mathbb{N}_0}$. Let $\mathcal{P}_n = \{x_1, \dots, x_n\}$ be a set of $n$ i.i.d. random points with $\mu$-density $\varrho$ and

$$\frac{n}{\ln(n)} \;\ge\; 128\,(2+t)\,M \cdot \max\left\{m\,,\; \frac{\sum_{k\ge m} \sigma_k^2}{\sigma_m^2}\right\} \tag{6.3}$$

with $M$ as in (6.2). Then the weighted least squares estimator (3.2) with points $\mathcal{P}_n$ and weights $w_i = \varrho(x_i)^{-1}$ on the space $V_m = \operatorname{span}\{b_0, \dots, b_{m-1}\}$ satisfies with probability at least $1 - 8n^{-t}$ that

$$\|f - A_{\mathcal{P}_n}^{w} f\|_2 \;\le\; 2\,\sigma_m \inf_{g\in V_m} \|f - g\|_H \qquad \text{for all } f \in H.$$

We note that the second expression of the maximum (6.3) is often of order $m$, for instance, if $\sigma_m \asymp m^{-\alpha} \log^{\beta} m$ for some $\alpha > 1/2$. Proposition 6.4 generalizes the result of Krieg and Ullrich (2021*a*) and the proof follows the same lines. It combines the conditional error bound from Proposition 3.22 with a concentration result of Mendelson and Pajor (2006) for the average of (infinite-dimensional) i.i.d. random matrices of rank one. The precise formulation of the following lemma is from Sonnleitner and Ullrich (2023, Lemma 3.3).

**Lemma 6.5.** Let $N \geq 3$, $R > 0$ and $y_1, \ldots, y_N$ be i.i.d. random sequences from $\ell_2$ satisfying $\|y_i\|_2^2 \leq R^2$ almost surely and $\|E\| \leq 1$, where $E = \mathbb{E}(y_i y_i^*)$ and $\|\cdot\|$ denotes the spectral norm. Then

$$\mathbb{P}\left(\left\|\frac{1}{N}\sum_{i=1}^N y_i y_i^* - E\right\| > \frac{1}{2}\right) \leq \frac{4}{N^c},$$

whenever

$$\frac{N}{\ln(N)} \geq 64\,(2+c)\,R^2.$$

We now prove Proposition 6.4.

*Proof.* First note that $\varrho(x_i) > 0$ almost surely for all $x_i$. We consider normalized weights $w_i = (n\varrho(x_i))^{-1}$, which does not change the algorithm. We apply Lemma 6.5 twice. First we apply it to the random vectors

$$y_i := \varrho(x_i)^{-1/2}(b_0(x_i), \ldots, b_{m-1}(x_i))$$

which satisfy

$$\mathbb{E}(y_i y_i^*) = I_m, \qquad \|y_i\|_2^2 \leq 2\,Mm \quad \text{a.s.}$$

For the computation of the expected value, we used $M < \infty$, since this condition implies that $b_k = 0$ on the set $\{\rho = 0\}$ for all $k \in \mathbb{N}_0$. We put

$$G := \frac{1}{n}\sum_{i=1}^n y_i y_i^*.$$

Lemma 6.5 implies that the smallest singular value of $G$ is bounded below by $1/2$ with probability at least $1 - 4n^{-t}$. This translates into the discretization condition (3.37) of Proposition 3.22 with $K = \sqrt{2}$ since for $g = \sum c_k b_k \in V_m$, we have

$$\sum_{i=1}^n w_i |g(x_i)|^2 = c^* G c \qquad \text{and} \qquad \|g\|_2^2 = \|c\|_2^2.$$

Next, we apply the concentration inequality to the random vectors

$$y_i := (\sigma_m^2 \varrho(x_i))^{-1/2}(\sigma_m b_m(x_i), \sigma_{m+1} b_{m+1}(x_i), \ldots)$$

which satisfy

$$\mathbb{E}(y_i y_i^*) = \sigma_m^{-2} \cdot \operatorname{diag}(\sigma_m^2, \sigma_{m+1}^2, \ldots), \qquad \|y_i\|_2^2 \leq 2\,M\,\sigma_m^{-2} \sum_{k \geq m} \sigma_k^2 \quad \text{a.s.}$$

We put

$$\Gamma := \frac{1}{n}\sum_{i=1}^{n} y_i y_i^*.$$

Lemma 6.5 implies that $\|\Gamma\| \le 3/2$ with probability at least $1-4n^{-t}$. This translates into the discretization condition (3.38) of Proposition 3.22 with $\gamma = \sqrt{3/2}\,\sigma_m$ since for every $h = \sum \sigma_k c_k b_k \in V_m^{\perp}$, we have

$$\sum_{i=1}^{n} w_i |h(x_i)|^2 = \sigma_m^2 \cdot c^* \Gamma c \qquad \text{and} \qquad \|h\|_H^2 = \|c\|_2^2.$$

So with probability at least $1 - 8n^{-t}$, both discretization conditions are satisfied and Proposition 3.22 can be applied.

□

Let us now transfer this result to the non-Hilbert setting from Section 3.3. This will allow us to choose a 'convenient' basis for the approximation for which the quantity $M$ from (6.2), and hence the number of required i.i.d. samples, can be controlled.

Let $(D,\mu)$ be a measure space. We consider approximation based on i.i.d. samples with respect to a given $\mu$-density $\varrho$. The measure $\mu$ itself does not have to be finite, but if $\mu$ is itself a probability measure, e.g., the normalized Lebesgue measure of some compact domain, then a particular interesting example seems to be $\varrho = 1$, in which case the approximation is performed based on i.i.d. points w.r.t. the measure $\mu$.

The approximation is performed with respect to a nested sequence of subspaces $(V_m)_{m\in\mathbb{N}}$ of $L_2(D,\mu)$ with corresponding orthonormal system $(b_k)_{k\in\mathbb{N}_0}$, that is, $V_m = \operatorname{span}\{b_k \colon 0 \le k < m\}$ for $m \in \mathbb{N}$. We characterize the quality of the i.i.d. samples in terms of the quantity

$$\mathcal{R}_m := \inf\{c \colon \rho_m(x) \le c \cdot \rho(x) \quad \text{for all } \ x \in D\} \tag{6.4}$$

with the probability density

$$\rho_m(x) := \frac{1}{m}\sum_{k=0}^{m-1} |b_k(x)|^2 = \frac{1}{m}\sup_{f\in V_m\setminus\{0\}} \frac{|f(x)|^2}{\|f\|_2^2}.$$

The density $\rho_m$ is called the *Christoffel density* associated to $V_m$ and is also used and discussed in Sections 4.2 and 7.3. Note that, in the case of a probability measure $\mu$ and $\varrho = 1$, this quantity equals the quantity $\frac{1}{m}\Lambda(V_m, B(D))^2$ from Section 3.3 and Section 4, but we choose a lighter and normalized notation here. It is clear that $\mathcal{R}_m \ge 1$. Moreover, $\mathcal{R}_m \in \mathcal{O}(1)$ occurs for various examples, see Remark 3.33. In the following, we assume that there are some $C > 0$ and $r \ge 0$ such that

$$\mathcal{R}_m \le C m^r \tag{6.5}$$

for all $m \in \mathbb{N}$. We want to show that $\mathcal{O}(m^{r+1}\log m)$ i.i.d. random samples can give an approximation of the same order as best approximation on $V_m$. Hence, for a given sequence $(V_m)_{m\in\mathbb{N}}$ and $\alpha > 0$, we define the approximation classes

$$E_\alpha := \big\{ f \in L_2 \colon f \text{ has a pointwise convergent Fourier series (3.30)} \\ \text{and } d(f, V_m)_2 \lesssim m^{-\alpha} \big\} \tag{6.6}$$

with

$$\|f\|_{E_\alpha} := \sup_{m\in\mathbb{N}} m^\alpha \cdot d(f, V_m)_2.$$

The following is implied by Proposition 6.4.

**Theorem 6.6.** There are universal constants $B, b > 0$ such that the following holds. Let $(D,\mu)$ be a measure space, $\rho$ a probability density and $(V_m)_{m\in\mathbb{N}}$ be a nested sequence of subspaces of $L_2(D,\mu)$ with $\dim V_m = m$ such that (6.5) holds for some $C > 0$ and $r \geq 0$, and let $n, m, t \geq 2$ with

$$\frac{n}{\ln(n)} \geq b\, C\, t\, m^{r+1}.$$

Then, with probability at least $1 - n^{-t}$, the weighted least squares estimator $A^w_{\mathcal{P}_n}$ from (3.3) with $n$ i.i.d. points $\mathcal{P}_n = \{x_1, \dots, x_n\}$ distributed according to $\varrho \cdot d\mu$ and weights $\varrho(x_i)^{-1}$ satisfies

$$\|f - A^w_{\mathcal{P}_n} f\|_2 \leq B \cdot m^{-\alpha} \cdot \|f\|_{E_\alpha}$$

for all $f \in E_\alpha$ and all $\alpha \geq 0.51 + r/2$.

As the reader may have guessed, the number 0.51 in Theorem 6.6 can be replaced by any other number $1/2 + \varepsilon$ with $\varepsilon > 0$ (in which case the constants depend on $\varepsilon$). We refer to Bartel (2023) and Cohen and Migliorati (2017) for related results (in a randomized setting) that also use $\mathcal{R}_m$ as a reference quantity.

*Proof.* This is obtained when Proposition 6.4 is applied to the Hilbert space from Lemma 3.29 with singular values $\sigma_k^2 = k^{-(r+1.01)}$. In this case, the quantity $M$ from (6.2) satisfies

$$M \leq c\, C\, m^r,$$

with a universal constant $c$ and $C$ from (6.5), which can be verified by a dyadic decomposition of the sum in (6.2) as done several times in Section 3.3. □

We stress that the probabilistic statement of Theorem 6.6 holds uniformly for all functions from $E_\alpha$ (and all $\alpha$). Statements of the form "For all $f \in F$, the following error bound holds with high probability" are weaker and belong to the realm of randomized algorithms (with the Monte Carlo or probabilistic error criterion) as discussed in Section 7. Theorem 6.6 gives a statement of the form "With high probability, the following error bound holds for all $f \in F$". Only the (subtle) change of the order of quantifiers makes it possible not to draw a new random sample for each approximation.

In the following examples, we always consider approximation with respect to a probability measure $\mu$, and sampling with respect to the distribution $\mu$, i.e., the density $\varrho = 1$.

**Example 6.7 (Trigonometric polynomials).** Assume that $\mu$ is the uniform distribution on $[0,1]^d$ and that each $V_m = \mathcal{T}_{\Lambda_m}$ is a space of trigonometric polynomials on $[0,1]^d$ with frequencies from a set $\Lambda_m \subset \mathbb{Z}^d$ of cardinality $m$, with $\Lambda_m \subset \Lambda_{m+1}$. Then, (6.5) holds for $C = 1$ and $r = 0$, and so Theorem 6.6 gives that $\mathcal{O}(m \log m)$ random samples give the optimal degree of approximation with respect to trigonometric polynomials with frequencies from $\Lambda_m$, where the hidden constant is independent of the dimension.

**Example 6.8 (Spherical harmonics).** Assume that $\mu$ is the uniform distribution on the $d$-dimensional sphere and the spaces $(V_m)_{m\in\mathbb{N}_0}$ are given by the spherical harmonics. More generally, we may consider the normalized volume measure on a compact Riemannian manifold and approximation with respect to the eigenspaces of the Laplace-Beltrami operator (ordered according to the size of the eigenvalues). Then, (6.5) again holds for $r = 0$, see Mhaskar (2010, Lemma 5.2), and so Theorem 6.6 implies that $\mathcal{O}(m \log m)$ random samples from the uniform distribution give the optimal degree of approximation with respect to the spherical harmonics. This time, however, we are not aware of the correct dimension-dependence of $C$.

**Example 6.9 (Algebraic polynomials).** Assume that $\mu$ is the uniform distribution on $[-1,1]^d$ and that each $V_m = \mathbb{P}_{\Lambda_m}$ is a space of algebraic polynomials on $[-1,1]^d$ with exponents from a lower set $\Lambda_m \subset \mathbb{N}_0^d$ of cardinality $m$. By a lower set, we mean that a component-wise inequality $\nu' \le \nu$ and $\nu \in \Lambda_m$ imply $\nu' \in \Lambda_m$. In this case, it is known that (6.5) holds with $C = 1$ and $r = 1$, see Chkifa, Cohen, Migliorati, Nobile and Tempone (2015, Lemma 3.1). Hence, Theorem 6.6 gives that $\mathcal{O}(m^2 \log m)$ random samples (with a dimension-independent constant) ensure the optimal degree of approximation with respect to polynomials from the spaces $V_m$. If instead the error is measured with respect to the tensorized Chebyshev measure (and hence also the samples are taken i.i.d. with respect to the Chebyshev measure), then it is known that (6.5) holds with $C = 2^d$ and $r = 0$. This means that $\mathcal{O}(m \log m)$ samples suffice for the optimal degree of approximation in this error norm, but that the hidden constant can have a very bad dimension-dependence. It is also known that (6.5) holds with $C = 1$ and $r = \ln 3/\ln 2 - 1$, see Chkifa *et al.* (2015, Lemma 3.3), which means that we can instead take $\mathcal{O}(m^{\ln 3/\ln 2} \log m)$ samples with a dimension-independent constant.

Let us finally note that the ultimate comparison would probably be the comparison of the expected radius of information (i.e., the expected quality of the random sampling points) with the sampling numbers (i.e., the quality of optimal sampling points), as we did for the example $F = F_{\mathrm{Lip}}$. However, this was only possible since we considered a particular and easy example, and we do not know to what extent this is possible for general classes $F$. Theorem 6.6 'only' shows that $\mathcal{O}(n \log n)$ random samples are as good as $n$ optimal samples in many cases. It is clearly of interest to extend the result on the expected radius of information to more general settings. But it seems at least as interesting to understand when and why the logarithmic oversampling is indeed needed. This is the subject of the next section, where we will see that the oversampling is indeed needed for many examples (e.g., for $L_2$-approximation on RKHSs and approximation classes $E_\alpha$). But it is not needed in other examples.

### *6.3. Do we need the logarithmic oversampling?*

In Section 6.2, we have seen many examples where $\mathcal{O}(n \log n)$ i.i.d. sampling points with high probability give an approximation as good as $n$ optimal sampling points. We now discuss whether the logarithmic oversampling can be avoided. As we do not know a satisfactory theory that answers this question for general (abstract) model classes $F$, we discuss this question at least for a general family of examples. As we shall see, the answer to whether or not already $\mathcal{O}(n)$ points are optimal depends very much on the 'geometry' of the problem.

The example is the problem of $L_p$-approximation on Sobolev spaces $W_q^s(D)$ on general bounded domains $D$ by means of i.i.d. and uniformly distributed sampling points. (By *domain* we mean any open and nonempty set.) Similar to the univariate case, that was already discussed shortly in Section 2, for integer smoothness $s \in \mathbb{N}$ with $s > d/q$, the spaces $W_q^s(D)$ are Banach spaces of continuous real-valued functions on $D$ defined by

$$W_q^s(D) := \left\{ f \in C(D) \colon f^{(k)} \in L_q \text{ for all } k \in \mathbb{N}_0^d \text{ with } |k| \leq s \right\}, \tag{6.7}$$

where $f^{(k)}$ denotes the $k$-th weak partial derivative with total order $|k| = k_1 + \ldots + k_d$. We consider $W_q^s(D)$ with the seminorm

$$|f|_{W_q^s} := \left( \sum_{|k|=s} \int_D |f^{(k)}(x)|^q \, dx \right)^{1/q} \tag{6.8}$$

and norm

$$\|f\|_{W_q^s} := \|f\|_\infty + |f|_{W_q^s}$$

and denote the unit ball of $W_q^s := W_q^s(D)$ by $F_q^s := F_q^s(D)$.

For this example, the decay of the sampling numbers $g_n(F_q^s, L_p)$, i.e., the error that can be achieved with optimally chosen sampling points, is well known, at least for Lipschitz domains $D$. We refer to Narcowich *et al.* (2005) and Novak and Triebel (2006) for the corresponding results. On the other hand, as shown by Krieg and Sonnleitner (2024) and detailed below, i.i.d. uniformly distributed sampling points give precisely the same order of convergence if and only if $p < q$, while a logarithmic oversampling is required in the case $p \geq q$. The same is true for Sobolev spaces on connected and compact Riemannian manifolds, see Krieg and Sonnleitner (2025).

Behind these results is a geometric characterization of the radius of information of arbitrary point sets $\mathcal{P}_n = \{x_1, \ldots, x_n\} \subset D$. Recall that the radius of information is the minimal worst case error that can be achieved with the given set of sampling points, i.e.,

$$\mathrm{rad}(\mathcal{P}_n, F_q^s, L_p) := \inf_{\Phi\colon \mathbb{R}^n \to L_p} \sup_{f \in F_q^s} \left\| f - \Phi\big(f(x_1), \ldots, f(x_n)\big) \right\|_{L_p}.$$

The characterization is in terms of the distance function

$$\mathrm{dist}(\,\cdot\,, \mathcal{P}_n)\colon D \to \mathbb{R}, \qquad \mathrm{dist}(x, \mathcal{P}_n) := \min_{k \leq n} \mathrm{dist}(x, x_k).$$

We already have seen a similar characterization in the introductory example of univariate Lipschitz functions, see (2.8), which corresponds to the case $s = 1$, $q = \infty$ and $d = 1$. It generalizes as follows, see Krieg and Sonnleitner (2024, Theorem 0.1) and Krieg and Sonnleitner (2025, Theorem 1).

**Theorem 6.10.** Let $D$ be a bounded convex domain in $\mathbb{R}^d$ or a connected compact Riemannian manifold of dimension $d \in \mathbb{N}$. Let $1 \leq p < q \leq \infty$ and $s \in \mathbb{N}$ with $s > d/q$. For any $n \in \mathbb{N}$ and $\mathcal{P}_n = \{x_1, \ldots, x_n\} \subset D$, it holds that

$$\mathrm{rad}\left(\mathcal{P}_n, F_q^s, L_p\right) \asymp \left\| \mathrm{dist}(\cdot, \mathcal{P}_n) \right\|_{L_\gamma}^s, \tag{6.9}$$

where $\gamma := s \cdot (1/p - 1/q)^{-1}$ and the implied constants are independent of $n$ and $\mathcal{P}_n$. In the case $p \geq q$, it holds that

$$\mathrm{rad}\left(\mathcal{P}_n, F_q^s, L_p\right) \asymp \left\| \mathrm{dist}(\cdot, \mathcal{P}_n) \right\|_\infty^{s - d(1/q - 1/p)}. \tag{6.10}$$

The quantity $\| \mathrm{dist}(\cdot, \mathcal{P}_n) \|_\infty$, which characterizes the error of optimal algorithms in the case $p \geq q$, is the radius of the largest 'empty' ball (i.e., the largest ball that does not contain a sampling point) and commonly known as the **covering radius**, the *mesh norm* or the *fill distance*. It is quite common in the literature to use the covering radius to bound errors of sampling-based algorithms, see, e.g., Wu and Schaback (1993), Mhaskar (2010), Narcowich *et al.* (2005), Arcangéli, López de Silanes and Torrens (2007), Brauchart, Dick, Saff, Sloan, Wang and Womersley (2015), Ehler, Graef and Oates (2019). The novelty of the result from Krieg and Sonnleitner (2024) is that, in the case $p < q$, the covering radius can be replaced with the smaller quantity $\| \mathrm{dist}(\cdot, \mathcal{P}_n) \|_{L_\gamma}$, that is, with the radius of an 'average' empty ball. This quantity is sometimes called the *distortion* of the point set; it also plays an important role in the theory of the quantization of measures. The replacement of the covering radius with the distortion is important for the study of i.i.d. sampling (and maybe even more important for other irregular sampling point sets), since the distortion of random points is significantly smaller than their covering radius. In particular, the following is known; we refer to Krieg and Sonnleitner (2025, 2024) for a proof and further references.

**Proposition 6.11.** Let $D$ be a bounded convex domain in $\mathbb{R}^d$ or a connected compact Riemannian manifold of dimension $d \in \mathbb{N}$. Let $x_1, x_2, \ldots$ be independent and uniformly distributed random vectors on $D$ and let $0 < \alpha < \infty$. Consider the random $n$-point set $\mathcal{P}_n = \{x_1, \ldots, x_n\}$. Then

$$\Big(\mathbb{E}\, \| \operatorname{dist}(\cdot, \mathcal{P}_n) \|_{L_\gamma}^\alpha \Big)^{1/\alpha} \asymp \begin{cases} n^{-1/d} & \text{if } 0 < \gamma < \infty, \\ \Big(\frac{n}{\log n}\Big)^{-1/d} & \text{if } \gamma = \infty, \end{cases}$$

whereas

$$\inf_{\mathcal{P}_n} \| \operatorname{dist}(\cdot, \mathcal{P}_n) \|_{L_\gamma} \asymp n^{-1/d} \qquad \text{for all } 0 < \gamma \leq \infty.$$

This immediately gives the following corollary, see also Corollary 6.3.

**Corollary 6.12.** Let $D$ be a bounded convex domain in $\mathbb{R}^d$ or a connected compact Riemannian manifold of dimension $d \in \mathbb{N}$. Let $1 \leq p, q \leq \infty$ and $s \in \mathbb{N}$ with $s > d/q$. Let $x_1, x_2, \ldots$ be independent and uniformly distributed on $D$. Consider the random $n$-point set $\mathcal{P}_n = \{x_1, \ldots, x_n\}$. Then

$$\mathbb{E}\left[ \operatorname{rad}(\mathcal{P}_n, F_q^s, L_p) \right] \asymp \begin{cases} g_n(F_q^s, L_p) & \text{if } q > p, \\ g_{n/\log n}(F_q^s, L_p) & \text{if } q \leq p. \end{cases}$$

That is, in expectation, i.i.d. uniformly distributed points are asymptotically as good as optimal points for $L_p$-approximation on $F_q^s$ if and only if $p < q$.

The Hilbert setting that we discussed in Proposition 6.4 of Section 6.2 corresponds to the case $p = q = 2$ (i.e., $\gamma = \infty$), in which case the error is given by the covering radius. Hence, $\Omega(n \log n)$ i.i.d. random points are required. This is even true for any probability distribution on $D$. Indeed, it follows from classical results on the coupon collector problem that for any distribution, at least $\Omega(n \log n)$ i.i.d. random samples are required to hit each of $n$ disjoint balls in $D$, and hence to achieve a covering radius of the same order as with $n$ optimal points. Hence, in the general setting of Proposition 6.4, the logarithmic oversampling cannot be avoided when function values are sampled at i.i.d. locations.

Unfortunately, the algorithm that was used by Krieg and Sonnleitner (2024) to prove the upper bounds is not very practical since it uses a decomposition of the domain based on the sampling points that is hard to compute. We now present a slightly different algorithm which seems to be be easy to implement and which also achieves the upper bounds from Theorem 6.10 (and hence also the bounds from Corollary 6.12).

In principle, this algorithm should work for general domains $D$ if they are sufficiently regular (e.g., convex). However, we avoid many technicalities by considering only $D = (0,1)^d$. The algorithm consists of two main steps.

The first step is to split the domain into subdomains that are well-suited to the set of sampling points $\mathcal{P}_n$. This may be regarded as precomputation since the same domain decomposition can be used for every $f \in W_q^s(D)$. The decomposition is described by Algorithm 2. For a simple description, we use a certain terminology: For $k \in \mathbb{N}$, by the $k^d$ subcubes of a cube, we refer to the closed cubes of equal volume and pairwise disjoint interior that are each contained in the (closure of the) cube $Q$ and whose union equals $Q$. The parameter $\ell$ in the algorithm is a constant that we have to choose based on the smoothness $s$, for instance, $\ell = s + 1$.

**Algorithm 2** Cube splitting algorithm

1: **Input:** point set $\mathcal{P}_n = \{x_1, \ldots, x_n\} \subset (0,1)^d$, $\ell \in \mathbb{N}$
2: **Set up:** Pool_1 = $\{(0,1)^d\}$, Pool_2 = $\emptyset$
3: **while** Pool_1 $\neq \emptyset$ **do**
4: **for** $Q \in$ Pool_1 **do**
5: **if** all $2^d$ subcubes of $Q$ satisfy that each of their $(2\ell)^d$ subcubes contains a point of $\mathcal{P}_n$
**then**
6: replace $Q$ in Pool_1 by the $2^d$ subcubes
7: **else**
8: move $Q$ from Pool_1 to Pool_2
9: **end if**
10: **end for**
11: **end while**
12: **Output:** Pool_2 = $\{Q_1, \ldots, Q_m\}$

The algorithm returns a collection of cubes $\mathcal{Q} = \{Q_1, \ldots, Q_m\}$ that cover $(0,1)^d$ and which have pairwise disjoint interior (but do not need to have equal volume). It can happen that $\mathcal{Q} = \{(0,1)^d\}$. In this case, there must be a cube of constant side-length $1/(4\ell)$ inside $(0,1)^d$ that does not contain any point of $\mathcal{P}_n$ and so the right hand sides of the bounds in Theorem 6.10 are in $\Omega(1)$. This means that we can simply choose the zero function as an approximation and the upper bound that is claimed by the theorem is satisfied for the reconstruction map $\Phi = 0$.

In the following, we hence assume that $\mathcal{Q} \neq \{(0,1)^d\}$. In this case, by construction of $\mathcal{Q}$, each cube $Q_i \in \mathcal{Q}$ contains at least $(2\ell)^d$ points from $\mathcal{P}_n$, one in each of its $(2\ell)^d$ subcubes. In particular $m \lesssim n$. Due to the geometric progression of the side-lengths, it is also clear that the number of cubes that have to be tested in line 5 in total is also $\mathcal{O}(n)$. Given this collection of cubes, we choose a point set $\mathcal{P}_n^{(j)}$ for each cube $Q_j$ that consists of $\ell^d$ 'well-spread' points by splitting $Q_j$ into only $\ell^d$ subcubes and taking one point with each coordinate in the left half of the corresponding interval. This point set has covering radius (in the supremum norm) at most $\mathrm{length}(Q_j)/\ell$ and separation distance at least $\mathrm{length}(Q_j)/(2\ell)$, where $\mathrm{length}(Q)$ denotes the side-length of a cube $Q$.

The second step is the actual approximation of a function $f \in W_q^s((0,1)^d)$, which we denote by $A_n(f)$. On each cube $Q_j$, we approximate the function $f$ by a polynomial $\pi_j(f) \in \mathbb{P}_s$ of total degree at most $s$. We use an approximation that only depends on the function values of $f$ at $\mathcal{P}_n^{(j)}$ and satisfies

$$\|f - \pi_j(f)\|_\infty \le C \inf_{\pi \in \mathbb{P}_s} \|f - \pi\|_\infty \tag{6.11}$$

for all continuous functions $f\colon Q_j \to \mathbb{R}$ and some constant $C > 0$ that only depends on $s$ and $d$. This could be a unweighted least squares algorithm or an moving least squares algorithm. We do not fix the method and only note that such methods exist, which follows, for instance, from Wendland (2005, Theorem 4.7). The approximation $A_n(f)$ shall simply be the piecewise polynomial given by

$$A_n f(x) = \pi_j(f)(x), \qquad \text{for } x \in Q_j. \tag{6.12}$$

Since the attentive reader of Section 3.2 might be a little irritated by (6.11) at first sight, we remark that the factor "$\sqrt{m}$" that appeared in similar upper bounds of Section 3.2, see Theorem 3.16, is now just a constant since also the dimension $m$ of the approximation space $\mathbb{P}_s$ is constant.

For the proof of the optimality of this algorithm, the following lemma about optimal polynomial approximation on cubes is essential. It is a consequence of the Bramble-Hilbert lemma or generalized Poincaré inequality, see Brenner and Scott (2008, Lemma 4.3.8) or Maz'ya (1985, Lemma 1.1.11), as well as Krieg and Sonnleitner (2024, Lemma 1.9) for this precise formulation.

**Lemma 6.13.** For any $1 \le q \le \infty$ and $s \in \mathbb{N}$ with $s > d/q$, there is a constant $c > 0$ such that the following holds. For any cube $Q \subset \mathbb{R}^d$ with side-length $\mathrm{length}(Q) \le 1$ and any $f \in W_q^s(Q)$, there is a polynomial $\pi$ of degree at most $s$ such that

$$\sup_{x \in Q} |(f - \pi)(x)| \le c\ \mathrm{length}(Q)^{s-d/q}\ |f|_{W_q^s(Q)}.$$

We now prove that the algorithm (6.12) gives the optimal error bound from Theorem 6.10 and Corollary 6.12 in the case $D = (0,1)^d$.

*Proof.* In the following, "$\lesssim$" means that the upper bound holds up to a constant factor that can depend only on $d$ and $s$ (possibly via $\ell$), but not on the point set $\mathcal{P}_n$. The case $\mathcal{Q} = \{(0,1)^d\}$, in which case we defined $A_n f = 0$, has already been discussed, so let us consider $\mathcal{Q} \neq \{(0,1)^d\}$. In this case, (6.11) and Lemma 6.13

imply

$$\sup_{x\in Q_j} |(f - \pi_j(f))(x)| \lesssim \text{length}(Q_j)^{s-d/q}|f|_{W^s_q(Q_j)}. \tag{6.13}$$

For $p < \infty$, this gives

$$\begin{aligned}\|f - A_n f\|^p_{L_p((0,1)^d)} &\lesssim \sum_{j=1}^m \int_{Q_j} \text{length}(Q_j)^{sp-dp/q}|f|^p_{W^s_q(Q_j)}\,dx \\ &= \sum_{j=1}^m \text{length}(Q_j)^{sp-dp/q+d}|f|^p_{W^s_q(Q_j)}.\end{aligned} \tag{6.14}$$

We start with the (more interesting) case $p < q$. We apply Hölder's inequality with conjugate indices $(1 - p/q)^{-1}$ and $q/p$ and continue the estimate by

$$\leq \left(\sum_{j=1}^m \text{length}(Q_j)^{(sp-dp/q+d)(1-p/q)^{-1}}\right)^{1-p/q} \left(\sum_{j=1}^m |f|^q_{W^s_q(Q_j)}\right)^{p/q}.$$

Since $(1/p - 1/q) = s/\gamma$ and $(sp - dp/q + d)(1 - p/q)^{-1} = \gamma + d$,

$$\|f - A_n f\|_{L_p((0,1)^d)} \lesssim \left(\sum_{j=1}^m \text{length}(Q_j)^{\gamma+d}\right)^{s/\gamma} |f|_{W^s_q((0,1)^d)}.$$

Now we observe that each $Q_j$ contains a subcube of side-length length$(Q_j)/(4\ell)$ that does not contain any sampling point (otherwise, it would have been split into further pieces by Algorithm 2). This means that

$$\int_{Q_j} \text{dist}(x, \mathcal{P}_n)^\gamma\,dx \gtrsim \text{length}(Q_j)^{\gamma+d}$$

and therefore

$$\sum_{j=1}^m \text{length}(Q_j)^{\gamma+d} \lesssim \big\|\text{dist}(\,\cdot\,, \mathcal{P}_n)\big\|^\gamma_{L_\gamma((0,1)^d)}$$

and the proof in the case $p < q$ is complete.

In the case $\infty > p \geq q$, we continue from (6.14) via

$$\text{length}(Q_j) \lesssim \|\text{dist}(\cdot, \mathcal{P}_n)\|_\infty \tag{6.15}$$

to get

$$\|f - A_n f\|_{L_p((0,1)^d)} \lesssim \|\text{dist}(\cdot, \mathcal{P}_n)\|_\infty^{s-d/q+d/p} \left(\sum_{j=1}^m |f|^p_{W^s_q(Q_j)}\right)^{1/p},$$

so that it only remains to use $\|\cdot\|_{\ell_p} \leq \|\cdot\|_{\ell_q}$ to obtain the stated inequality. In the remaining case $p = \infty$, the bound follows directly from (6.13) via (6.15) and

$$|f|_{W_q^s(Q_j)} \leq |f|_{W_q^s((0,1)^d)}.$$

□

**Remark 6.14 (High dimensions).** The above algorithm is rather useless if the dimension $d$ is large, as the number $m$ of considered cubes is exponential in the dimension. On the other hand, the approximation problem on isotropic Sobolev classes suffers from the curse of dimensionality, see Section 8, so there cannot be *any* feasible algorithm in high dimensions for the model class $F_q^s(D)$.

**Remark 6.15 (Gaussian information).** There is a different type of random measurements where the logarithmic oversampling in the Hilbert case is not required, namely, Gaussian information. As in Proposition 6.4, let $H$ be a RKHS that is embedded into $L_2(D, \mu)$ with square-summable singular values $(\sigma_k)_k$ and corresponding $L_2$-orthonormal basis $(b_k)_k$. Instead of function evaluations, we consider information of the form

$$\ell_i(f) = \sum_k g_{ik} \langle f, b_k \rangle_2, \qquad 1 \leq i \leq n,$$

with independent standard Gaussian variables $g_{ik}$, which is well-defined almost surely for all $f \in H$. It is proven by Hinrichs *et al.* (2021*b*) that the least squares algorithm $A_n$ on $V_m = \operatorname{span}\{b_0, \ldots, b_{m-1}\}$ based on $n$ i.i.d. Gaussian measurements satisfies with high probability that

$$\|f - A_n f\|_2 \leq 2\sigma_m, \qquad \text{for all } f \in H$$

with the condition (6.3) replaced by

$$n \geq C \max\left\{4m\,,\ \frac{\sum_{k\geq m} \sigma_k^2}{\sigma_m^2}\right\}. \tag{6.16}$$

Hence, for singular values with polynomial decay, in which case the maximum in (6.16) is of order $m$, we see that $\mathcal{O}(m)$ pieces of i.i.d. Gaussian information are as good as $m$ pieces of optimal information. In particular, for $L_2$-approximation on $W_2^s(D)$, i.i.d. Gaussian information, if available, gives better results than i.i.d. sampling of function values.

Gaussian information is maybe most natural in the case that $D = \{1, \ldots, N\}$ with the counting measure and when $f \in \mathbb{R}^N$, so that the information $\ell_i(f)$, $i \leq n$, up to rescaling, is given by $n$ coordinates of $f$ in random directions which are i.i.d. and uniformly distributed on the unit sphere. We refer to Hinrichs, Prochno and Sonnleitner (2023) for a continuation of the study of the power of Gaussian information.

## 7. Randomized sampling algorithms

Randomized algorithms, also known as *Monte Carlo methods*, are a large class of algorithms which are allowed to use randomly chosen measurements. Even though we used randomness before, this is different from the algorithms discussed so far, since previously, we first fixed a set of sampling points (sometimes deterministically, but sometimes also randomly) and used the same set of sampling points for every function from the input class $F$. Now we use different random sampling points for each $f \in F$. That is, a **randomized algorithm** or **Monte Carlo method** $M_n$ for sampling recovery is a mapping of the form

$$M_n \colon F \times \Omega \to Y, \quad M_n(f,\omega) = \Phi^{\omega}\big(f(x_1(\omega)), \ldots, f(x_n(\omega))\big), \tag{7.1}$$

where the sampling points are random variables $x_i \colon \Omega \to D$ on some probability space $(\Omega, \mathbb{P})$ and also the recovery map $\Phi^{\omega} \colon \mathbb{C}^n \to Y$ can be random.

Another crucial difference to our previous approach is the way we measure the error of the Monte Carlo method. Assuming the measurability of the mapping $\omega \mapsto \|f - M_n(f,\omega)\|_Y$ for each $f \in F$, the **Monte Carlo error** of $M_n$ can be defined as

$$e^{\mathrm{rms}}(M_n, F, Y) \ := \ \sup_{f \in F} \sqrt{\mathbb{E}\Big[\|f - M_n(f)\|_Y^2\Big]}.$$

The Monte Carlo error is also called the **worst-case root-mean-square error**. That is, in order to have a small Monte Carlo error, the goal is for each individual $f \in F$ to have a good chance of a small error. This differs from Section 6, where we chose random sampling points once to be used for all $f \in F$ and wanted a good chance to have a small error for every $f \in F$ simultaneously.

The downside of randomized algorithms compared to deterministic algorithms is rather clear: Monte Carlo methods do not allow for reliable error guarantees, in the sense that error bounds only hold with certain probability, and that a realization of a randomized algorithm may have small error for some $f \in F$, but not for all at once. But on the upside, due to the larger class of algorithms and the weaker error criterion, we may hope to obtain better error bounds for Monte Carlo methods than for deterministic algorithms, as well as error bounds for larger input classes $F$. And indeed, we will see below that this is often the case.

In this section, we aim to provide a glance into the Monte Carlo approach for function recovery and to highlight some differences to the deterministic setting. The literature on sampling recovery using randomized algorithms is vast, and we do not even try to accomplish a comprehensive overview in this single section. We refer to Novak (1988), Novak and Woźniakowski (2012), Heinrich (1994), Kunsch (2017) for other general results and surveys on randomized approximation, with an emphasis on complexity as treated in this work. A more extensive survey on least squares methods based on random sampling, as discussed below, was given by Adcock (2025). We also recommend comparing with some machine learning literature, such as Christmann and Steinwart (2008), Belkin (2021) and Bach (2024), for a closely related perspective.

We proceed as follows. As in the deterministic setting, we first discuss in Section 7.1 a simple example to illustrate some general ideas and phenomena. Turning to general classes $F$ in Section 7.2, we start again by discussing the minimal error that can be achieved with arbitrary reconstruction schemes based on $n$ random samples. That is, we analyze the benchmark quantity

$$g_n^{\mathrm{rms}}(F,Y) \;:=\; \inf_{M_n} e^{\mathrm{rms}}(M_n,F,Y).$$

The infimum is taken over all randomized algorithms of the form (7.1) for which $\omega \mapsto \|f - M_n(f,\omega)\|_Y$ is measurable for all $f \in F$ (so that the expected value in the definition of the Monte Carlo error exists). We fully concentrate on the case of $L_2$-approximation, $Y = L_2$, since this is the so far best understood case.

Afterwards, we present a selection of explicit sampling strategies and algorithms and discuss how close they get to the theoretical benchmark. We start in Section 7.3 with the discussion of weighted least squares algorithms based on different kinds of sampling designs: i.i.d. sampling, subsampling, greedy sampling, and volume sampling. The weighted least squares algorithms can be used to obtain the theoretically optimal error bounds from Section 7.2. In Section 7.4, we also discuss the multilevel Monte Carlo approach based on hierarchical Christoffel sampling, which is a little less general than the least squares method but can be of advantage from a computational point of view.

### 7.1. *The easy example revisited*

Let us revisit the introductory example of Lipschitz functions from Section 2, i.e., the problem of $L_p$-approximation on the class

$$F_{\mathrm{Lip}} := \Big\{ f\colon\, [0,1] \to \mathbb{C}\colon\;\; |f(x) - f(y)| \le \mathrm{dist}(x,y) \;\text{ for all } x,y \in [0,1] \Big\}.$$

As discussed, using Monte Carlo methods, we may hope to achieve a smaller error than with deterministic algorithms. For this example, however, the hope is not justified. Namely, the following holds.

**Proposition 7.1.** For all $n \in \mathbb{N}$, we have

$$g_n^{\mathrm{rms}}(F_{\mathrm{Lip}}, L_1) \;\ge\; \frac{1}{16n}. \tag{7.2}$$

In particular, for all $1 \le p \le \infty$, it holds that

$$g_n(F_{\mathrm{Lip}}, L_p) \;\asymp\; g_n^{\mathrm{rms}}(F_{\mathrm{Lip}}, L_p) \;\asymp\; \frac{1}{n}. \tag{7.3}$$

Hence, the Monte Carlo approach is not superior to deterministic approximation

(for which equidistant points are optimal) for the $L_p$-approximation of Lipschitz functions. Similar things will happen for other problems, too. In fact, there are even some general lower bounds for Monte Carlo methods that show that the *maximal gain* over deterministic algorithms is limited, see Section 10. In any case, finding the right lower bound for particular model classes is often the hard part. A general lower bound for the error of Monte Carlo methods in terms of Bernstein numbers is given by Kunsch (2016), see Section 9.3.

But how does one prove lower bounds for randomized algorithms? For deterministic algorithms, we used the following technique: Find a pair of *fooling functions* $f, g \in F$ with information $f(x_i) = g(x_i)$, $i \le n$, and for which $\|f - g\|_Y$ is large. Then the algorithm cannot distinguish $f$ and $g$ and makes a worst-case error at least $\frac{1}{2}\|f - g\|_Y$. For randomized algorithms, this does not work: If the sampling points $x_1, \ldots, x_n$ invoke enough randomness, for any pair of functions $f \neq g$, we have a (usually very high) chance of being able to distinguish them by the obtained function values. So randomized algorithms are not so easily fooled and a different proof technique is needed.

In general, lower bounds for randomized algorithms are proven by **Bakhvalov's proof technique** (also called **Yao's principle**) due to seminal work by Bakhvalov (1959) and Yao (1977). The worst-case error of a randomized algorithm is compared with the average-case error of deterministic algorithms. For a probability measure $\nu$ on the input class $F$, the *average-case error* of a deterministic measurable algorithm $A_n \colon F \to Y$ is defined as

$$e^{\mathrm{avg}}(A_n, \nu, Y) \;:=\; \sqrt{\int_F \|f - A_n(f)\|_Y^2 \,\mathrm{d}\nu(f)},$$

and the $n$-th minimal average-case error is

$$g_n^{\mathrm{avg}}(\nu, Y) \;:=\; \inf_{A_n} e^{\mathrm{avg}}(A_n, \nu),$$

where the infimum runs over all deterministic measurable algorithms that use at most $n$ sampling points. Now, for every randomized algorithm $M_n \colon F \times \Omega \to Y$, it holds that

$$\begin{aligned} e^{\mathrm{rms}}(M_n, F)^2 \;&=\; \sup_{f \in F} \int_\Omega \|f - M_n(f, \omega)\|_Y^2 \,\mathrm{d}\mathbb{P}(\omega) \\ &\ge\; \int_F \int_\Omega \|f - M_n(f, \omega)\|_Y^2 \,\mathrm{d}\mathbb{P}(\omega)\,\mathrm{d}\nu(f) \\ &=\; \int_\Omega \int_F \|f - M_n(f, \omega)\|_Y^2 \,\mathrm{d}\nu(f)\,\mathrm{d}\mathbb{P}(\omega) \\ &\ge\; \inf_{\omega \in \Omega} \; e^{\mathrm{avg}}(M_n(\,\cdot\,, \omega), \nu)^2. \end{aligned}$$

Here $M_n(\,\cdot\,, \omega)$ is a deterministic algorithm that uses at most $n$ samples. Here we assume that the integrand $\|f - M_n(f,\omega)\|_Y$ is a measurable function on $F \times \Omega$, which is the case, in particular, for all finitely supported measures $\nu$. It follows that

$$g_n^{\mathrm{rms}}(F,Y) \;\geq\; g_n^{\mathrm{avg}}(\nu,Y) \tag{7.4}$$

for any such measure $\nu$, i.e., the Monte Carlo error of any randomized algorithm is larger than the average error of the best deterministic algorithm.

The 'art' of applying this proof technique mainly consists in finding good measures $\nu$ on the function class $F$ for which the minimal average-case error becomes large. The following lemma is similar to several propositions by Novak (1988, Section 2.2), see also Mathé (1991) and Heinrich (1994).

**Lemma 7.2.** Let $F$ be a class of integrable functions on $D$ and let $f_1, \dots, f_{2n}$ be functions with mutually disjoint support such that $\sum_{j=1}^{2n} \delta_j f_j \in F$ for all signs $\delta_j \in \{-1,1\}$ and $\|f_j\|_1 \geq \varepsilon$ for all $j \leq 2n$. Then

$$g_n^{\mathrm{rms}}(F,L_1) \;\geq\; n\varepsilon.$$

*Proof.* The purpose of this proof is to illustrate ideas and we will hence be a little informal. We consider the average-case error with respect to the uniform distribution $\nu$ on the (finite) set

$$F_0 \;:=\; \left\{ \sum_{j=1}^{2n} \delta_j f_j \colon \delta_j \in \{-1,1\} \right\}.$$

A deterministic algorithm that computes $n$ function values at (possibly adaptively chosen) sampling points $x_1, \dots, x_n$ can, at best, identify $n$ of the signs $\delta_j$ of the unknown function $f \in F_0$ and has no information about the other signs. So it is intuitively clear that the average-case error is minimized by the algorithm

$$A_n^*(f) \;=\; \sum_{j \in M(f)} \delta_j f_j,$$

where $M(f)$ is the set of all indices $j$ for which $\{x \in D \colon f_j(x) \neq 0\}$ contains a sampling point. That is, $A_n^*$ is the best algorithm that uses the sample $\{x_1, \dots, x_n\}$. For this average-optimal algorithm, we have for any $f \in F_0$ that

$$\|f - A_n^*(f)\|_1 \;\geq\; n\varepsilon$$

and so also the average-case error is larger than $n\varepsilon$. Therefore, any deterministic algorithm has an average-case error of at least $n\varepsilon$ and the relation (7.4) gives

$$g_n^{\mathrm{rms}}(F,L_1) \;\geq\; g_n^{\mathrm{avg}}(\nu,L_1) \;\geq\; n\varepsilon.$$

□

Now the proof of Proposition 7.1 is easy:

*Proof.* We partition the domain into $2n$ intervals of length $(2n)^{-1}$ and apply

Lemma 7.2 for the hat functions $f_1, \dots, f_{2n}$, where $f_j$ is supported in the $j$-th interval, linear on each half of the interval, and equal to $(4n)^{-1}$ at the midpoint. Then all signed sums of the functions $f_j$ are in the class $F_{\mathrm{Lip}}$ and $\|f_j\|_1 = (16n^2)^{-1}$ for all $j \le 2n$.

This gives the lower bound for $g_n^{\mathrm{rms}}(F_{\mathrm{Lip}}, L_1)$. For the second statement, we invoke

$$g_n^{\mathrm{rms}}(F_{\mathrm{Lip}}, L_1) \;\le\; g_n^{\mathrm{rms}}(F_{\mathrm{Lip}}, L_p) \;\le\; g_n(F_{\mathrm{Lip}}, L_p)$$

and the upper bounds from the deterministic case from Proposition 2.3.

□

There are examples different from $F_{\mathrm{Lip}}$ where the Monte Carlo error of randomized algorithms can hugely improve upon the error of any deterministic algorithm. This will become apparent also from the general results discussed in the next subsection. But let us already mention an example. Namely, for $s > 0$, let us consider the class

$$F_2^s \;:=\; \left\{ f \in C([0,1])\colon\, \sum_{k\in\mathbb{Z}} \left(1 + |k|^{2s}\right) |\hat{f}(k)|^2 \le 1 \right\}.$$

With $\hat{f}(k)$ we denote the $k$-th Fourier coefficient of $f$. The class $F_2^s$ is, up to the continuity assumption and periodicity, the unit ball of a Sobolev space $W_2^s$ from (2.11). But note that this definition makes sense also for $s \notin \mathbb{N}$. In the case $s > 1/2$, the set is a compact subset of $C([0,1])$ and we have

$$g_n(F_2^s, L_2) \;\asymp\; g_n^{\mathrm{rms}}(F_2^s, L_2) \;\asymp\; n^{-s},$$

so randomization is not helpful. In the case $s \le 1/2$, however,

$$g_n(F_2^s, L_2) \;\asymp\; 1 \quad \text{and} \quad g_n^{\mathrm{rms}}(F_2^s, L_2) \;\asymp\; n^{-s},$$

so randomized algorithms still have a small error, while deterministic algorithms are useless. This is discussed (in much more generality) by Heinrich (2008). The upper bound for randomized algorithms also follows from the general bounds in Section 7.2.

## 7.2. *How good are optimal random samples?*

In this section, we want to discuss the benchmark quantity

$$g_n^{\mathrm{rms}}(F, L_2) \;:=\; \inf_{\Phi, \mathcal{P}_n} \; \sup_{f\in F} \sqrt{\mathbb{E}\Big[ \|f - \Phi(f(x_1), \dots, f(x_n)\|_2^2 \Big]}.$$

Here, $L_2 = L_2(D, \mu)$ for an arbitrary measure space $(D, \mu)$ and $F \subset \mathbb{C}^D$ is a class of square integrable functions. The infimum is taken over all random choices of sampling points $\mathcal{P}_n = \{x_i \colon i \le n\} \subset D$ and mappings $\Phi\colon \mathbb{C}^n \to L_2$ for which the expected value exists. The numbers $g_n^{\mathrm{rms}}(F, L_2)$ are called the **randomized** **sampling numbers** and represent the minimal Monte Carlo error in the $L_2$-norm that can be achieved with $n$ samples.

The randomized sampling numbers are the probabilistic counterpart of the sampling numbers for deterministic algorithms, which we studied in Section 3. Since every deterministic algorithm can be viewed as a randomized algorithm with constant random variables, it is clear that we have

$$g_n^{\mathrm{rms}}(F, L_2) \leq g_n(F, L_2). \tag{7.5}$$

In the deterministic setting, we were able to compare the sampling numbers with the error of best approximation from a linear space, represented by the Kolmogorov numbers $d_n(F, L_2)$, see (3.54). In particular, we obtained for all $\alpha > 1/2$ and $\beta \in \mathbb{R}$ that

$$d_n(F, L_2) \lesssim n^{-\alpha} \log^{\beta}(n) \quad \Longrightarrow \quad g_n(F, L_2) \lesssim n^{-\alpha} \log^{\beta}(n) \tag{7.6}$$

under some mild assumptions on the input class $F$, see Section 3.3 for details.

In the randomized setting, due to the relation (7.5), we may hope to obtain better bounds: Can we remove the condition $\alpha > 1/2$ on the polynomial decay? Can we even remove the assumption on the tail of the Kolmogorov numbers and obtain bounds for each individual $n$? And can we remove the "mild assumptions" on the class $F$? In all three cases, the answer is yes. We start with the end of the story. The following is an immediate consequence of a recent breakthrough result of Dolbeault and Cohen (2022*a*).

**Theorem 7.3.** There are universal constants $b \in \mathbb{N}$ and $C > 0$, such that for every class $F$ of square-integrable functions and all $n \in \mathbb{N}$ it holds that

$$g_{bn}^{\mathrm{rms}}(F, L_2) \leq C \cdot d_n(F, L_2).$$

In fact, a stronger version of the result is true for any $n$-dimensional subspace $V_n$ of $L_2$. The result is proved by a weighted least squares algorithm. We discuss this in more detail in Section 7.3. For the moment, let us only discuss the randomized sampling numbers.

We first discuss the sharpness of Theorem 7.3. For this, note that we do not know a general lower bound that is specific to pointwise sampling. We can only employ the more general bounds on algorithms that have access to arbitrary linear measurement maps, see Section 9, which may be much smaller. However, it turns out that if $F$ is the unit ball of a reproducing kernel Hilbert space that is compactly embedded into $L_2$, we have (almost) a matching lower bound. In fact, by a result of Novak (1992), we have that

$$g_n^{\mathrm{rms-meas}}(F, L_2) \geq \frac{1}{\sqrt{2}} \cdot d_{2n}(F, L_2), \tag{7.7}$$

where $g_n^{\mathrm{rms-meas}}$ is the minimal Monte Carlo error that can be achieved with meas-urable algorithms $M_n\colon F\times\Omega\to L_2$. Recall that, in the definition of $g_n^{\mathrm{rms}}$, we only assumed measurability in the second argument and therefore considered a larger class of randomized algorithms. This ensured that the deterministic algorithms are a subclass of the randomized algorithms. Since the upper bound in Theorem 7.3 can be achieved with a measurable algorithm, we may conclude that the error of measurable randomized sampling algorithms is characterized (up to universal constants) by the Kolmogorov numbers for Hilbert classes $F$.

For more general sets $F$, we only have a lower bound in terms of the *Bernstein widths*, which is due to Kunsch (2016) and presented in Theorem 9.23 in Section 9.3 (applied for $S=\mathrm{id}_Y$ and $Y=L_2$). This lower bound gives the sharpness of the bound from Theorem 7.3 for many other examples, namely, for all classes $F$ where the Bernstein widths and the Kolmogorov numbers are of a similar size. In general, however, the Bernstein widths can be much smaller than the Kolmogorov numbers, and so there is a gap between the general upper and lower bounds (which can partially be explained by the superiority of adaptive sampling in the randomized setting, see Section 10.2). It is an interesting open problem whether a general upper bound as in Theorem 7.3 is possible with the Kolmogorov widths replaced by smaller approximation benchmarks, ideally in terms of Bernstein widths.

We now discuss the history of Theorem 7.3. Indeed, the result has several forerunners. It was proven by Wasilkowski and Woźniakowski (2007) that the randomized sampling numbers must always have the same polynomial decay as the Kolmogorov numbers, namely,

$$
\begin{aligned}
&d_n(F,L_2) \lesssim n^{-\alpha}\log^{\beta}(n)\\
&\qquad\Longrightarrow\quad g_n^{\mathrm{rms}}(F,L_2) \lesssim n^{-\alpha}\log^{\beta}(n)(\log\log n)^{\alpha+1/2}
\end{aligned}
\tag{7.8}
$$

for all $\alpha\ge 0$ and $\beta\in\mathbb{R}$. The next result is by Cohen and Migliorati (2017) and gives that

$$
g^{\mathrm{rms}}_{bn\log n}(F,L_2) \le C\, d_n(F,L_2),
\tag{7.9}
$$

with universal constants $b,C>0$, which implies a slightly worse asymptotic bound, but holds for every individual $n$ and hence makes a decay assumption on the Kolmogorov widths obsolete. Going back to asymptotic analysis, Krieg (2019) removed excessive logarithmic factors, resulting in a statement of the form

$$
\begin{aligned}
&d_n(F,L_2) \lesssim n^{-\alpha}\log^{\beta}(n)\\
&\qquad\Longrightarrow\quad g_n^{\mathrm{rms}}(F,L_2) \lesssim n^{-\alpha}\log^{\beta}(n)
\end{aligned}
\tag{7.10}
$$

again for all $\alpha\ge 0$ and $\beta\in\mathbb{R}$. Finally, Dolbeault and Cohen (2022*a*) could remove the excessive logarithmic factor in (7.9), resulting in Theorem 7.3, which also implies (7.10) but does not require a decay assumption on the Kolmogorov widths. A further improvement of Theorem 7.3 is due to Dolbeault and Chkifa (2024), who, apart from advances regarding practical implementation, improved the constants in Theorem 7.3 and obtained results also for small oversampling factors $b\ge 1$.

Namely, it is shown that, for any $m \in \mathbb{N}$ and $n \geq m$, it holds that

$$g_n^{\mathrm{rms}}(F, L_2) \leq \left(1 + \frac{1}{(1 - 1/\sqrt{r})^2}\right) \cdot d_m(F, L_2), \tag{7.11}$$

where $r = n/(m-1)$. The algorithms used for (7.8) and (7.10) are multilevel Monte Carlo algorithms, the other results are proven with randomized weighted least squares algorithms. We will discuss these algorithms in the following sections.

Let us comment again on the difference to the deterministic results from Section 3.3, particularly Corollary 3.31, which bounds the deterministic sampling numbers by the Kolmogorov widths in $L_2$. A big difference of Theorem 7.3 in comparison with Corollary 3.31 is that the tail of the sequence of Kolmogorov numbers does not play any role for optimal Monte Carlo algorithms. This leads to good bounds also in the case of a decay with polynomial rate $\alpha \leq 1/2$, where no general bound for the deterministic sampling numbers can exist, see Hinrichs *et al.* (2008) and Hinrichs *et al.* (2022, Theorem 6). But it can also make a big difference even in the case of polynomial decay $\alpha > 1/2$. Indeed, the disappearance of the tail implies that some high-dimensional approximation problems that suffer from the curse of dimensionality if only deterministic algorithms are used, become tractable in the randomized setting, see Section 8.

In the following sections, we turn to the discussion of algorithms. We start with least squares algorithms in Section 7.3 since these algorithms are used to prove Theorem 7.3. In Section 7.4, we shortly discuss multilevel Monte Carlo algorithms, which might have advantages from a computational point of view.

**Remark 7.4 (Predescribed sampling measures).** The emphasis of this section (or this survey in general) is on optimal sampling strategies. We therefore do not discuss here (i.i.d.) sampling with respect to a general given (and non-optimal) measure as we did in Section 6. However, let us add that this is a very important subject, since, in *real-world applications*, the sampling measure may not be at our disposal. Especially in a machine learning context, one may even have different measures to handle, sometimes called *source* and *target distribution*. This area of research is called *domain (or covariate shift) adaptation*, see the influential works of Shimodaira (2000), Ben-David, Blitzer, Crammer and Pereira (2006) and Ganin and Lempitsky (2015). One prominent approach to this problem is to assume the availability of values of the Radon-Nikodym derivative of the measures, and take their size into account. We refer to Bartel (2023) and Krieg, Novak and Sonnleitner (2022) for work close to our setting (and language), as well as to Gizewski, Mayer, Moser, Nguyen, Pereverzyev, Pereverzyev, Shepeleva and Zellinger (2022), Nguyen, Zellinger and Pereverzyev (2024), Zellinger, Kindermann and Pereverzyev (2023) for other recent contributions. Practical aspects and references in this direction are given by Dinu, Holzleitner, Beck, Nguyen, Huber, Eghbal-Zadeh, Moser, Pereverzyev, Hochreiter and Zellinger (2023).

### 7.3. *Randomized least squares algorithms*

We work in a very similar setting to Section 3.1, and only shortly recall the basic definitions. Let $(D,\mu)$ be a measure space and let $V_m$ be an $m$-dimensional linear space of complex-valued square-integrable functions on $D$. Given points $\mathcal{P}_n := \{x_1,\dots,x_n\} \in D$ and corresponding weights $w := (w_1,\dots,w_n) \in \mathbb{R}_+^n$, we consider the *weighted least squares algorithm*

$$A_{\mathcal{P}_n}^w(f) \;=\; \operatorname*{arg\,min}_{g\in V_m} \sum_{i=1}^n w_i \,|f(x_i) - g(x_i)|^2 \tag{7.12}$$

for functions $f\colon D \to \mathbb{C}$, see (3.3). That is, $A_{\mathcal{P}_n}^w(f)$ is the best approximation of $f$ from $V_m$ with respect to the seminorm

$$\|f\|_{\mathcal{P}_n^w} \;:=\; \sqrt{\sum_{i=1}^n w_i |f(x_i)|^2}.$$

In the deterministic setting, we always need two properties of the points and weights to ensure a small error of the least squares algorithm $A_{\mathcal{P}_n}^w$ on the input class $F \subset \mathbb{C}^D$. The first property was a discretization condition of the form

$$\|g\|_2 \;\le\; K\|g\|_{\mathcal{P}_n^w} \qquad \text{for all } g \in V_m \tag{7.13}$$

with some constant $K > 0$, see (3.11). This condition ensures that $A_{\mathcal{P}_n}^w$ is well defined but it also ensures the *stability* of the reconstruction map, see Lemma 3.2 and Remark 3.6, which remains important in the randomized setting. For easier reference, we refer to a condition like (7.13) as *reconstructional stability*.

The second condition takes different forms, depending on the desired error bound. In any case, it was needed to ensure that

$$\|f - g\|_{\mathcal{P}_n^w} \;\le\; \varepsilon \qquad \text{for all } f \in F \tag{7.14}$$

for some small number $\varepsilon > 0$, and some good approximation $g = g(f) \in V_m$ of $f$, like its orthogonal projection in $L_2$ or its best approximation in the uniform norm. This condition can be interpreted in the way that changing from $f \in F$ to $g(f) \in V_m$ (or vice versa) only leads to small observation errors. For easier reference, we refer to a condition like (7.14) as *observational stability*.

Together, the two conditions imply

$$\|f - A_{\mathcal{P}_n}^w(f)\|_2 \;\le\; \|f-g\|_2 \,+\, K\varepsilon \qquad \text{for all } f \in F,$$

which is a simple computation:

$$\begin{aligned}\|f - A_{\mathcal{P}_n}^w(f)\|_2 \;&\le\; \|f-g\|_2 + \|g - A_{\mathcal{P}_n}^w(f)\|_2 \;=\; \|f-g\|_2 + K\|A_{\mathcal{P}_n}^w(g-f)\|_{\mathcal{P}_n^w}\\ &\le\; \|f-g\|_2 + K\|g-f\|_{\mathcal{P}_n^w} \;\le\; \|f-g\|_2 + K\varepsilon.\end{aligned}$$

In the randomized setting, the condition (7.14) can be replaced by a probabilistic upper bound or an upper bound for the expected value. A uniform upper bound that holds simultaneously for all $f \in F$ is not required any more. For example, if we replace (7.14) with

$$\mathbb{E}\,\|f-g\|^2_{\mathcal{P}_n^w} \;\leq\; \varepsilon^2 \qquad \text{for all } f \in F \tag{7.15}$$

we get a corresponding bound

$$\mathbb{E}\,\|f-A^w_{\mathcal{P}_n}(f)\|_2 \;\leq\; \|f-g\|_2 \,+\, K\varepsilon \qquad \text{for all } f \in F.$$

The relaxed observational stability (7.15) makes it possible to consider very large input classes $F$, which is a big advantage compared to the deterministic setting. For example, if we choose the points i.i.d. with some $\mu$-density $\varrho$ and the weights $w_i = (n\varrho(x_i))^{-1}$, we get

$$\mathbb{E}\,\|f\|^2_{\mathcal{P}_n^w} \;=\; \frac{1}{n}\sum_{i=1}^{n} \mathbb{E}\left(\frac{|f(x_i)|^2}{\varrho(x_i)}\right) \;=\; \|f\|_2^2 \tag{7.16}$$

for any square-integrable function $f \in L_2$. This leads to a bound of the form (7.15) for any $f \in L_2$ for which $\|f-g\|_2^2 \leq \varepsilon$, and hence any function that has a good approximation $g$. The observation (7.16) suggests that there might be no harm in imposing a condition on observational stability for all $f \in L_2$, namely

$$\mathbb{E}\,\|f\|^2_{\mathcal{P}_n^w} \;\leq\; \beta^2\,\|f\|_2^2 \qquad \text{for all } f \in L_2, \tag{7.17}$$

in order to obtain error bounds for any $f \in L_2$. This is done in the following proposition of Dolbeault and Cohen (2022*a*).

**Proposition 7.5.** Let $V_m$ be a finite-dimensional subspace of $L_2$. Consider random points $\mathcal{P}_n := \{x_1, \ldots, x_n\} \in D$ and weights $w := (w_1, \ldots, w_n) \in \mathbb{R}^n_+$. Assume that the reconstructional stability (7.13) is satisfied almost surely and that the observational stability (7.17) holds. Then the least squares algorithm (7.12) satisfies

$$\left(\mathbb{E}\left\|f-A^w_{\mathcal{P}_n}(f)\right\|_2^2\right)^{1/2} \;\leq\; C\cdot \inf_{g\in V_m}\|f-g\|_2 \tag{7.18}$$

for all $f \in L_2$, where $C = \sqrt{1+K^2\beta^2}$.

*Proof.* Let $P$ be the orthogonal projection onto (i.e., the best approximation from) $V_m$ in $L_2$. Recall that $A^w_{\mathcal{P}_n}$ is linear with $A^w_{\mathcal{P}_n}(g) = g$ for $g \in V_m$. By orthogonality and the reconstructional stability, we have

$$\begin{aligned}
\|f-A^w_{\mathcal{P}_n}(f)\|_2^2 &= \|f-Pf\|_2^2 + \|Pf-A^w_{\mathcal{P}_n}(f)\|_2^2 \\
&= \|f-Pf\|_2^2 + \|A^w_{\mathcal{P}_n}(f-Pf)\|_2^2 \\
&\leq \|f-Pf\|_2^2 + K^2\|A^w_{\mathcal{P}_n}(f-Pf)\|^2_{\mathcal{P}_n^w} \\
&\leq \|f-Pf\|_2^2 + K^2\|f-Pf\|^2_{\mathcal{P}_n^w},
\end{aligned}$$

almost surely. It remains to take the expected value and use the observational stability,

$$\mathbb{E}\,\|f - Pf\|^2_{\mathcal{P}_n^w} \;\leq\; \beta^2\,\|f - Pf\|_2^2.$$

□

It remains to find random choices of points and weights that satisfy the observational stability (7.17) and reconstructional stability (7.13). Clearly, the number $n$ of points should be as small as possible, ideally linear in the dimension $m$ of the space $V_m$.

If the measure $\mu$ that we use to quantify the error is a probability measure, a first natural choice is to choose the points i.i.d. according to $\mu$ and consider weights $w_i = 1/n$, that is, an unweighted least squares algorithm. Unweighted randomized least squares algorithm are studied, for instance, by Cohen *et al.* (2013) and Migliorati *et al.* (2014). Depending on the space $V_m$, they can lead to almost optimal sampling size, but in general, a sampling size $n$ much larger than the dimension $m$ is needed.

A first general way to find $n \asymp m \log m$ such points is from Cohen and Migliorati (2017). It is based on **Christoffel sampling**, that is, on i.i.d. random draws with respect to the **Christoffel density**

$$\varrho_m(x) \;:=\; \frac{1}{m}\sum_{k=1}^{m} |b_k(x)|^2,$$

where $(b_1, \ldots, b_m)$ is an arbitrary $L_2$-orthonormal basis of $V_m$. This is an important and well-studied object, especially in the context of orthogonal polynomials, and we refer to Szegö (1975) and Nevai (1986) for more on this. Recall that we already used this probability density in Section 4.2. In machine learning, where often $D$ is a large finite set and the basis $(b_1, \ldots, b_m)$ is hence a matrix, this kind of sampling is also known as **leverage score sampling**. For this, we recommend the early contribution by Chatterjee and Hadi (1986), as well as the more recent work of Ma, Mahoney and Yu (2015) and references therein.

To be precise, the desired properties are achieved by *conditional* Christoffel sampling as follows. Let $\mathcal{Z}_n = \{x_1, \ldots, x_n\}$ be a random point set consisting of $n = 10\, m \log(4m)$ independent random points distributed according to the Christoffel density. We take a random draw of $\mathcal{Z}_n$ and accept the obtained point set if the reconstructional stability (7.13) is satisfied with the corresponding weights $w_i = (n\varrho_m(x_i))^{-1}$, namely, if the Gram matrix $G_{\mathcal{Z}_n}$ from (3.7) satisfies

$$\lambda_{\min}(G_{\mathcal{Z}_n}) \;\geq\; \frac{1}{2}.$$

If this is not the case, we draw another independent realization of $\mathcal{Z}_n$, until the procedure terminates. A concentration inequality for the random matrix $G_{\mathcal{Z}_n}$, for instance (Tropp 2012, Theorem 1.1), implies that

$$\mathbb{P}\left(\lambda_{\min}(G_{\mathcal{Z}_n}) \geq \frac{1}{2}\right) \geq \frac{1}{2},$$

which means that, on average, it will take two random draws of $\mathcal{Z}_n$ until the procedure terminates. Let us call the random set of $n$ points that is obtained in this way by $\mathcal{Y}_n$. Then, of course, $\lambda_{\min}(G_{\mathcal{Y}_n}) \geq \frac{1}{2}$ is satisfied almost surely. Moreover, the random point set $\mathcal{Y}_n$ satisfies the following. Recall that this result is due to Cohen and Migliorati (2017).

**Theorem 7.6.** Let $V_m$ be an $m$-dimensional subspace of $L_2$. The set $\mathcal{Y}_n = \{x_1, \ldots, x_n\}$ from above with $n = 10\, m \log(4m)$ random points and weights $w_i := (n\rho_m(x_i))^{-1}$ satisfies

$$\begin{aligned} \|g\|_2^2 &\leq 2\|g\|_{\mathcal{Y}_n^w}^2 \qquad \text{a.s. for all } g \in V_m, \\ \mathbb{E}\,\|f\|_{\mathcal{Y}_n^w}^2 &\leq 2\|f\|_2^2 \qquad \text{for all } f \in L_2. \end{aligned}$$

Consequently, the corresponding least squares algorithm (7.12) satisfies

$$\mathbb{E}\,\|f - A_{\mathcal{Y}_n}^w f\|_2^2 \leq 5 \cdot \inf_{g \in V_m} \|f - g\|_2^2 \qquad \text{for all } f \in L_2.$$

**Remark 7.7 (Alternatives to Christoffel sampling).** Generating samples from the Christoffel density, as needed for the algorithm of Theorem 7.6, can be problematic, especially for irregular domains $D \subset \mathbb{R}^d$, where an analytic expression for an orthonormal basis of $V_m$ is not available even for "simple" spaces such as spaces of algebraic polynomials. Alternative sampling strategies for the case that Christoffel sampling is not possible are discussed, for instance, by Adcock and Cardenas (2020), Migliorati (2021) and Dolbeault and Cohen (2022*b*).

Starting from Theorem 7.6, it is possible to obtain the existence of a random construction with a linear amount of samples based on the recent solution of the Kadison-Singer conjecture by Marcus *et al.* (2015), or more precisely, an implied result about frames by Nitzan, Olevskii and Ulanovskii (2013). Indeed, the aforementioned result implies that we can partition the random set $\mathcal{Y}_n$ from above into sets $\mathcal{X}_1, \ldots, \mathcal{X}_r$, each with a cardinality at most $bm$ with some universal constant $b$, such that the reconstructional stability

$$\|g\|_2 \leq K\, \|g\|_{\mathcal{X}_i^w} \qquad \text{for all } g \in V_m,$$

is satisfied for any of those sets $\mathcal{X}_i$ with weight $w(x) = (|\mathcal{X}_i|\varrho_m(x))^{-1}$ for each $x \in \mathcal{X}_i$ and a universal constant $K$. By randomly drawing a subset $\mathcal{P}_{bm}$ from the collection of those sets $\mathcal{X}_i$, each with probability $|\mathcal{X}_i|/|\mathcal{Y}_n|$, we can also preserve the observational stability

$$\mathbb{E}\,\|f\|_{\mathcal{P}_{bm}^w}^2 \leq 2\|f\|_2^2 \qquad \text{for all } f \in L_2.$$

This gives the following breakthrough result of Dolbeault and Cohen (2022*a*).

**Theorem 7.8.** For any $m$-dimensional space $V_m \subset L_2$, there exists a random sample $\mathcal{P}_n = \{x_1, \dots, x_n\} \subset D$ consisting of $n \leq bm$ points and random weights $w_1, \dots, w_n > 0$ such that the least squares estimator (7.12) satisfies

$$\mathbb{E}\, \|f - A^w_{\mathcal{P}_n}(f)\|_2^2 \leq C \cdot \inf_{g \in V_m} \|f - g\|_2^2 \qquad \text{for all } f \in L_2.$$

Here, $b, C > 0$ are universal constants.

The downside of Theorem 7.8 is that it does not come with a construction of the points and weights as also the result of Marcus *et al.* (2015) is not constructive. A major progress is hence due to Dolbeault and Chkifa (2024), where a constructive way of finding a linear amount of samples is derived from the framework of Batson *et al.* (2014). This approach has the additional advantage of significantly improved constants $b$ and $C$, also providing results for very small oversampling factors $b \geq 1$, including the interpolation case $b = 1$. Namely, for every $n \geq m$, a set of $n$ random points and weights are constructed such that the following holds. This is Theorem 1.1 from Dolbeault and Chkifa (2024).

**Theorem 7.9.** For any $m$-dimensional space $V_m \subset L_2$ and any $n \geq m$, the least squares algorithm (7.12) with points $\mathcal{P}_n := \{x_1, \dots, x_n\} \in D$ and weights $w := (w_1, \dots, w_n) \in \mathbb{R}^n_+$ sampled according to Algorithm 3 satisfies

$$\mathbb{E}\, \|f - A^w_{\mathcal{P}_n}(f)\|_2^2 \leq \left(1 + \frac{1}{(1 - 1/\sqrt{r})^2}\right) \cdot \inf_{g \in V_m} \|f - g\|_2^2 \qquad \text{for all } f \in L_2,$$

where $r = n/(m-1)$.

For instance, an oversampling factor $r = 4$ in Theorem 7.9 leads to best approximation up to a factor

$$\mathbb{E}\, \|f - A^w_{\mathcal{P}_n}(f)\|_2^2 \leq 5 \cdot \inf_{g \in V_m} \|f - g\|_2^2 \qquad \text{for all } f \in L_2,$$

and hence the same bound as in Theorem 7.6 with now $4(m-1)$ instead of $10\, m \log(4m)$ samples. For the reader's convenience, we also copy the algorithm by Dolbeault and Chkifa (2024), namely Algorithm 2 with $\kappa = 0$, that gives the points and weights from Theorem 7.9. As the weights are chosen one by one, the sampling design is a kind of *greedy* sampling strategy.

Other successful random sampling techniques for weighted least squares approximation are based on determinantal point processes (DPPs), referred to as *determinantal sampling* and *volume sampling*. This is a well-known sampling technique in machine learning and statistics, see, for instance, the work of Kulesza and Taskar (2012) and Lavancier, Møller and Rubak (2014). Here, the samples are

**Algorithm 3** Construction of points and weights (Dolbeault and Chkifa 2024)

**Input:** probability measure $\mu$, $b(x) = (b_1(x), \dots, b_m(x))^\top$, and $n \in \mathbb{N}$.

$\boldsymbol{A}_0 = \boldsymbol{0} \in \mathbb{R}^{m\times m}, \quad \ell_0 = -m, \quad r = n/(m-1)$

**for** $i = 1, \dots, n$ **do**

Let $\boldsymbol{Y}_i = (\boldsymbol{A}_{i-1} - \ell_{i-1}\boldsymbol{I})^{-1}$

Update $\ell_i = \ell_{i-1} + 1/\sqrt{r}$

Let $\boldsymbol{Z}_i = (\boldsymbol{A}_{i-1} - \ell_i\boldsymbol{I})^{-1}$

Let $\boldsymbol{W}_i = \left(\mathrm{Tr}(\boldsymbol{Z}_i) - \mathrm{Tr}(\boldsymbol{Y}_i)\right)^{-1}\boldsymbol{Z}_i^2 - \boldsymbol{Z}_i$ and $\varrho_i : x \mapsto b(x)^*\boldsymbol{W}_i b(x)$

Define $R_i(x) = \varrho_i(x)\mathbf{1}_{\varrho_i(x)\geq 0}$

Draw $x_i$ from probability density $\dfrac{R_i(x)}{\Gamma_i}d\mu(x)$, where $\Gamma_i = \displaystyle\int_X R_i(x)d\mu(x)$

Let $w_i = 1/R_i(x_i) = 1/\varrho_i(x_i)$

Update $\boldsymbol{A}_i = \boldsymbol{A}_{i-1} + w_i b(x_i)b(x_i)^*$

**end for**

**Output:** Sample $\{x_1, \dots, x_n\}$, weights $\{w_1, \dots, w_n\}$.

drawn from a joint distribution of the form

$$(x_1, \dots, x_n) \sim \frac{1}{Z}\det\left(K(x_i, x_j)\right)_{i,j\leq n} d\mu^{\otimes n}(x_1, \dots, x_n),$$

where $K\colon D \times D \to \mathbb{C}$ is a reproducing kernel on $D$ and $Z > 0$ a normalization constant. For instance, the $L_2$-kernel of the approximation space $V_m$ is used, i.e.,

$$K(x, y) = \sum_{k=1}^{m} b_k(x)\overline{b_k(y)}.$$

Drawing $n = m$ samples from the latter distribution leads to a sampling design which is similar to the i.i.d. Christoffel sampling from Cohen and Migliorati (2017) in the sense that the marginal distribution of each point is equal to $\varrho_m d\mu$ with the Christoffel density $\varrho_m$. Hence, the individual point distributions of this so-called projection DPP and the i.i.d. Christoffel sampling are equal. But in contrast to the i.i.d. sampling, the volume sampling avoids clustering (e.g., the density is zero for all samples $(x_1, \dots, x_n)$ with $x_i = x_j$ for some $i \neq j$). The hope, which is partially supported by numerical experiments, is therefore that (some variant of) the DPP sampling can *avoid the coupon collector's problem* and fill the 'holes' in the data set more quickly than i.i.d. sampling, which may give better results. However, it is still unclear whether the hope is justified: The best error guarantees proven so far rather seem to be worse than those of i.i.d. Christoffel sampling, see (Belhadji, Bardenet and Chainais 2025, Table 2). We also refer to Nouy and Michel (2025) for related sampling strategies.

### 7.4. *Multilevel Monte Carlo methods*

We now discuss another Monte Carlo method, that can get close to the best $L_2$-approximation in a finite-dimensional space with a constant oversampling factor. It requires stronger assumptions than the least squares algorithms discussed in Section 7.3, namely, a hierarchy of good approximation spaces $(V_m)_{m\in\mathbb{N}}$ instead of a single approximation space $V_m$. On the other hand, it has a smaller runtime.

The method is a *multilevel Monte Carlo method* based on hierarchical Christoffel sampling. Foundational work on multilevel Monte Carlo methods in general has been done by Heinrich (2001), see also the article by Giles (2015) for a comprehensive study of this approach. Here, we use an idea of Wasilkowski and Woźniakowski (2007) that has been refined by Krieg (2019).

As before, $(D,\mu)$ is any measure space. We assume that we have an orthonormal system $(b_j)_{j\in\mathbb{N}}$ in $L_2$. Again, we use random samples according to the Christoffel densities

$$\varrho_m(x) \;=\; \frac{1}{m}\sum_{j=1}^{m} |b_j(x)|^2.$$

The algorithm $M_k^r$, introduced below, is a multilevel approximation, defined by an oversampling factor $r\in\mathbb{N}$ and an approximation level $k\in\mathbb{N}_0$. That is, an approximation $M_k^r f$ of $f$ will be obtained iteratively from a series of approximations $M_\ell^r f$ with $\ell < k$.

In fact, starting with $M_{-1}^r f := 0$, the approximation of level $k$ is defined by setting $m_k = 2^k$, $n_k = r\cdot m_k$, and

$$M_k^r f \;:=\; M_{k-1}^r f + \sum_{j=1}^{m_k}\left[\frac{1}{n_k}\sum_{i=1}^{n_k}\frac{\left(f - M_{k-1}^r f\right)\left(x_i^{(k)}\right)\cdot\overline{b_j}\left(x_i^{(k)}\right)}{\varrho_{m_k}\left(x_i^{(k)}\right)}\right] b_j, \qquad (7.19)$$

where the $x_i^{(k)}$ are independent random variables with distribution $\varrho_{m_k}\,d\mu$. In particular, the expectation of each term in the inner sum is

$$a_{k,j} \;:=\; \langle f - M_{k-1}^r f, b_j\rangle_2.$$

Hence, the term in the square brackets is nothing but a *standard Monte Carlo estimator* for $a_{k,j}$, i.e., an average over $n_k$ i.i.d. copies of random variables $Y_{k,j}$ with expected value $a_{k,j}$. Also note that $M_k^r f$ only requires function values of $f$ as information, assuming that the basis functions $b_j$ and the density $\varrho_m$ are known, because $M_{k-1}^r f$ is fully specified by the function values of $f$ at the sampling points from the previous levels.

The algorithm $M_k^r$ hence approximates $f$ in $k+1$ steps. In the first step, $n_0$ function values of $f$ are used for a standard Monte Carlo estimation of its $m_0$ leading coefficients with respect to the orthonormal system. In the second step, $n_1$ values of the residue are used for a standard Monte Carlo estimation of its $m_1$ leading coefficients and so on. In total, the algorithm estimates the first $m_k = 2^k$ coefficients of $f$ by using

$$\sum_{j=0}^{k} n_j \;\le\; 2r \cdot m_k$$

function evaluations of $f$. A pseudo-code for $M_k^r$ is given in Algorithm 4. The following error bound is essentially from Krieg (2019).

**Theorem 7.10.** Let $(b_j)_{j\in\mathbb{N}}$ be an orthonormal system in $L_2$ and consider the multilevel algorithm $M_k^r$ as in Algorithm 4 (which uses $n \le r2^{k+1}$ function evaluations). For any $f \in L_2$, it holds that

$$\mathbb{E}\,\|f - M_k^r f\|_2^2 \;\le\; r^{-k} \sum_{\nu=-1}^{k} r^{\nu}\,\|f - P_{2^\nu} f\|_2^2,$$

where $P_m$ denotes the $L_2$-orthogonal projection onto $V_m := \operatorname{span}\{b_j \colon 1 \le j \le m\}$.

**Remark 7.11 (Choice of parameters).** First note that, if the budget $n$ and the oversampling factor $r$ are given, then we obtain an approximation from $V_m$ with $m = n/(2r)$ based on $n$ samples by choosing the level $k = \log_2(m)$, where we assume without loss of generality that $m$ is a power of two. We write $A_n(f) := M_k^r(f)$ for the corresponding approximation. In order to choose the oversampling factor, assume first that we know a monotone positive sequence $(\sigma_m)_{m\in\mathbb{N}}$ such that for all powers $m$ of two, we have

$$d(f, V_m) := \inf_{g\in V_m} \|f - g\|_2 \;\le\; \sigma_m \qquad \text{and} \qquad \gamma := \sup_{k\in\mathbb{N}_0} \frac{\sigma_{2^k}}{\sigma_{2^{k+1}}} < \infty.$$

Then we can choose any $r \ge 2\gamma^2$ and obtain

$$e_n(f) \;:=\; \sqrt{\mathbb{E}\,\|f - A_n(f)\|_2^2} \;\le\; \sqrt{2}\cdot\sigma_{n/2r}.$$

Hence, in the above sense, $2rm$ samples suffice for best approximation from $V_m$ up to a factor $\sqrt{2}$. If no a priori bound on $d(f, V_m)$ is known, we may choose the oversampling factor $r \in \mathbb{N}$ arbitrary and obtain for all $\alpha \le \log_2(r) - 1$ and $f \in L_2$ that $d(f, V_m) \in O(m^{-\alpha})$ implies $e_n(f) \in O(n^{-\alpha})$.

*Proof.* We proceed by induction over $k$. In the case $k = -1$, both sides of the inequality are equal to $\|f\|_2^2$ and the statement is correct. For $k \ge 0$, we have by orthogonality,

$$\left\|f - M_k^r(f)\right\|_2^2 = \left\|P_{m_k} f - M_k^r f\right\|_2^2 + \left\|f - P_{m_k} f\right\|_2^2.$$

The first term satisfies

$$\left\|P_{m_k} f - M_k^r f\right\|_2^2 = \sum_{j=1}^{m_k} \left|\left\langle f - M_k^r f, b_j\right\rangle_2\right|^2.$$

**Algorithm 4** Multilevel Monte Carlo algorithm for approximating $f \in L_2$

**Input:** oversampling parameter $r \in \mathbb{N}$, level $k \in \mathbb{N}_0$, $m = 2^k$
Let $\boldsymbol{c} = \boldsymbol{0} \in \mathbb{R}^m$ {initialize vector of coefficients}
**for** $s = 0, \ldots, k$ **do**
Let $m_s = 2^s$ and $n_s = r \cdot 2^s$
Draw $x_1, \ldots, x_{n_s}$ from probability density $\varrho_{m_s}(x)\, d\mu(x)$
Let $\boldsymbol{B} = \big(b_j(x_i)\big)_{i \le n_s, j \le m_s}$ and $\boldsymbol{B}_w = \big(\varrho_m(x_i)^{-1} b_j(x_i)\big)_{i \le n_s, j \le m_s}$
Let $\boldsymbol{y} = (f(x_i))_{i \le n_s}$ {function values of $f$}
Let $\tilde{\boldsymbol{y}} = \boldsymbol{B} \cdot \boldsymbol{c}[m_s]$ {function values of current approximant}
Let $\tilde{\boldsymbol{c}} = \frac{1}{n_s} \boldsymbol{B}_w^* (\boldsymbol{y} - \tilde{\boldsymbol{y}})$ {estimated coefficiens of current remainder}
Update $\boldsymbol{c}[m_s] \leftarrow \boldsymbol{c}[m_s] + \tilde{\boldsymbol{c}}$
**end for**
**Output:** Coefficient vector $\boldsymbol{c} \in \mathbb{R}^m$. The obtained approximation of $f$ is $M_k^r f = \sum_{k=1}^m c_k b_k$, where we used $n \le 2rm$ samples.

Let us fix a realization of $M_{k-1}^r f$ and only take the expected value over the random variables from the $k$th iteration. For $j \le m_k$, we use the abbreviation

$$g_j := \frac{1}{\varrho_{m_k}} \big(f - M_{k-1}^r f\big)\, \overline{b_j}.$$

Note that $\varrho_{m_k} = 0$ implies $b_j = 0$ and we set $g_j = 0$ in this case. If $x^{(k)}$ is a random variable with distribution $\varrho_{m_k} d\mu$, we have

$$\mathbb{E}\, g_j\big(x^{(k)}\big) \;=\; \big\langle f - M_{k-1}^r f, b_j \big\rangle_2 .$$

With this and the definition of $M_k^r$, we obtain

$$\big|\big\langle f - M_k^r f, b_j \big\rangle_2\big|^2 = \left| \mathbb{E}\, g_j\big(x^{(k)}\big) - \frac{1}{n_k} \sum_{i=1}^{n_k} g_j\big(x_i^{(k)}\big) \right|^2 .$$

Hence, for any fixed realization of the random variables from the levels $j < k$ (i.e., we compute expected value and variance only with respect to the random variables from level $k$), we obtain

$$\begin{aligned} \mathbb{E}\, \big|\big\langle f - M_k^r f, b_j \big\rangle_2\big|^2 \;&=\; \mathbb{V}\left( \frac{1}{n_k} \sum_{i=1}^{n_k} g_j\big(x_i^{(k)}\big) \right) \\ &= \frac{1}{n_k} \mathbb{V}\left(g_j\big(x^{(k)}\big)\right) \;\le\; \frac{1}{n_k}\, \mathbb{E}\big|g_j\left(x^{(k)}\right)\big|^2. \end{aligned}$$

Hence, using that

$$\sum_{j=1}^{m_k} |g_j|^2 = \frac{m_k}{\varrho_{m_k}} |f - M_{k-1}^r f|^2,$$

we obtain

$$\mathbb{E}\left(\sum_{j=1}^{m_k} \left|\left\langle f - M_k^r f, b_j \right\rangle_2\right|^2\right) \le \frac{1}{n_k} \mathbb{E}\left(\sum_{j=1}^{m_k} \left|g_j(x^{(k)})\right|^2\right)$$

$$= \frac{m_k}{n_k} \mathbb{E}\left(\frac{\left|\left(f - M_{k-1}^r f\right)\left(x^{(k)}\right)\right|^2}{\varrho_{m_k}(x^{(k)})}\right) = \frac{1}{r} \left\|f - M_{k-1}^r f\right\|_2^2.$$

Altogether, this gives

$$\mathbb{E}\left\|f - M_k^r f\right\|_2^2 \le \frac{1}{r}\left\|f - M_{k-1}^r f\right\|_2^2 + \left\|f - P_{m_k}(f)\right\|_2^2$$

for every realization of $M_{k-1}^r f$. Hence, taking the expected value over the random variables from all levels, by Fubini's theorem, we get

$$\mathbb{E}\left\|f - M_k^r f\right\|_2^2 \le \frac{1}{r}\mathbb{E}\left\|f - M_{k-1}^r f\right\|_2^2 + \left\|f - P_{m_k}(f)\right\|_2^2.$$

This concludes the induction step.

□

Let us compare the multilevel method with the weighted least squares method. As we have seen, both methods give optimal error rates, but the multilevel method needs stronger model assumptions. Clearly, the multilevel approach only makes sense if there is a natural hierarchy of approximation spaces. On the other hand, the running time of the multilevel method can be much smaller. Indeed, assuming a constant oversampling factor $r$, the running time of Algorithm 4 is $\mathcal{O}(mt + m^2)$ if $t$ is the (maximal) time needed to generate a sample from one of the Christoffel densities $\varrho_{2^k}$, $k \le \log_2 m$. In comparison, the running time of the optimal weighted least squares algorithm by Dolbeault and Chkifa (2024) with points and weights as constructed in Algorithm 3 is of order $\mathcal{O}(m^4 t)$, where $t$ now is the time needed to generate a sample from $\varrho_m$. Therefore, the multilevel method seems especially advantageous if the cost $t$ of generating a Christoffel sample is high.

As a closing remark, let us mention that, just as in the deterministic setting, all the error bounds for $L_2$-approximation from Sections 7 can be lifted according to formula (3.60) to obtain error bounds for $L_p$-approximation with $p \neq 2$ (and other error norms), since both the weighted least squares method and the multilevel method return an approximation from the space $V_m$. We omit the details.

## 8. The curse of dimensionality

So far, we mostly discussed the approximation of functions $f\colon D \to \mathbb{C}$ that are defined on arbitrary sets $D$, without paying too much attention to the role of the domain $D$. In many applications, the function $f$ depends on a large number of variables, which means that the domain is a subset of $\mathbb{R}^d$ or some $d$-dimensional manifold, where the dimension $d \in \mathbb{N}$ can be quite large. It is important to study how the amount of required samples depends on the input dimension $d$.

For this purpose, it is not sufficient to study the asymptotic behavior of minimal errors (like sampling numbers), which was often the main focus in the previous sections. It is easily possible that, for a given model class $F_d$ of $d$-variate functions, we can achieve a very good and even dimension-independent rate of decay of the worst case error, while at the same time, the error remains larger than, say, 1/2 with less than $2^d$ samples. In such situations, we commonly talk about the *curse of dimensionality*. See below for a precise definition.

The term *curse of dimensionality* has been introduced by Bellman (1957) for the phenomenon that the number of needed samples increases exponentially with the (input) dimension. Informally, the opposite situation is described by the term *tractability*, referring to situations when the number of required samples shows a decent behavior both with respect to the error demand and with respect to the dimension. Depending on what one understands to be a 'decent' behavior, there are various different concepts of tractability.

These concepts of tractability are very much inspired by classifications in discrete complexity theory and we use similar concepts here. The difference is that we consider problems where a solution is only possible up to some error $\varepsilon > 0$ and we measure the *complexity* solely in terms of the amount of required samples that is needed to achieve an error less than $\varepsilon$. Since the measure of complexity is the amount of required information, we talk about the *information complexity* of the problem. It is clear that the information complexity can serve as a lower bound for the running time of numerical algorithms (in the worst case), but it can sometimes also as serve as an upper bound (e.g., for quadrature formulas). A more general treatment of the subject is given in Section 9, see also Remark 9.8.

Tractability and intractability in the setting of information-based complexity is discussed in much detail in the series of monographs Novak and Woźniakowski (2008, 2010, 2012); see also the various references in Section 8.3. In this section, we only want to convey some basic concepts, discuss a few recent developments, and relate our results from the previous sections to the question for tractability.

For a formal definition of tractability and the curse of dimensionality, it is necessary to consider a whole family of model classes $F_d \subset \mathbb{C}^{D_d}$, one for each 'dimension' $d \in \mathbb{N}$. The parameter $d \in \mathbb{N}$ does not necessarily have to represent a dimension, but it often does in the sense that $F_d$ is a class of functions on a domain $D_d \subset \mathbb{R}^d$. For each $d \in \mathbb{N}$, the error is measured in some seminormed space $Y_d \subset \mathbb{C}^{D_d}$. As in the previous sections, the minimal worst-case error for

$Y_d$-approximation on $F_d$ is described by the *sampling numbers*

$$g_n(F_d, Y_d) := \inf_{\substack{x_1,\ldots,x_n \in D_d \\ \Phi\colon \mathbb{C}^n \to Y_d}} \sup_{f \in F_d} \left\| f - \Phi(f(x_1), \ldots, f(x_n)) \right\|_{Y_d}.$$

For the discussion of tractability, it is often more convenient to work with the inverse quantity

$$n(\varepsilon, F_d, Y_d) := \min\{n \in \mathbb{N}\colon\ g_n(F_d, Y_d) \le \varepsilon\} \tag{8.1}$$

for $\varepsilon > 0$, which we shall call the **sampling complexity** of $Y_d$-approximation on $F_d$. (More generally, it is called *information complexity*, see Section 9.) Hence, the sampling complexity is the minimal amount of samples that is needed to guarantee an error smaller than $\varepsilon$.

We say that the problem of $Y_d$-approximation on $F_d$ suffers from the **curse of dimensionality** if, for some $\varepsilon_0 > 0$, there are constants $C, \delta > 0$ such that

$$n(\varepsilon_0, F_d, Y_d) \ge C\,(1+\delta)^d$$

for all $d \in \mathbb{N}$. That is, if the sampling complexity grows exponentially with the dimension $d$. On the contrary, we say that the problem is **polynomially tractable** if there are constants $C, r, s > 0$ such that

$$n(\varepsilon, F_d, Y_d) \le C\, d^s\, \varepsilon^{-r}$$

for all $d \in \mathbb{N}$ and $\varepsilon > 0$. If this is even true with $s = 0$, we call the problem **strongly polynomially tractable**.

It is clear that the complexity of a problem is not always described well by one of these three notions (e.g., a complexity of order $\varepsilon^{-2}\exp(\sqrt{d})$ is in none of the three complexity classes), and so there are various other notions of (in-)tractability in the literature. But these three notions are probably the most simple and the most widely used, and so we refrain from introducing further notions. To some extent, we regard (strong) polynomial tractability as the 'ideal' situation and the curse of dimensionality as the 'worst' scenario for high-dimensional approximation.

In this section, we discuss the following.

- We start in Section 8.1 by discussing again the introductory example of $L_p$-approximation of Lipschitz functions, this time with respect to tractability. It turns out that the problem suffers from the curse of dimensionality. We discuss whether the curse is avoided if the functions have a higher smoothness and/or if one uses randomized algorithms with the weaker Monte Carlo error criterion.
- The comparison results from the previous sections (where we estimated the sampling numbers in terms of other approximation benchmarks) do not show any dimension-dependence. We shortly discuss some consequences for tractability in Section 8.2 and present different sufficient conditions for polynomial tractability.

- In Section 8.3, we survey some concrete model classes for which the approximation problem is tractable: partially separable functions, weighted Korobov classes, functions with a bounded absolute sum of Fourier coefficients, ridge functions, symmetric functions, and rank one tensors.

### *8.1. The easy example revisited*

In Section 2, we discussed as an introductory example the problem of $L_p$-approximation on a class of univariate Lipschitz functions. We also start the discussion of tractability by considering a simple model class of Lipschitz functions. Precisely, we consider the model class

$$F^d_{\mathrm{Lip}} := \left\{ f : [0,1]^d \to \mathbb{R} : |f(x) - f(y)| \le \|x-y\|_\infty \ \text{ for } x,y \in [0,1]^d \right\}. \quad (8.2)$$

For a first insight into the problem, we will show that $L_p$-approximation on $F^d_{\mathrm{Lip}}$ suffers from the curse of dimensionality for all $1 \le p \le \infty$. It is not entirely arbitrary that we define Lipschitz continuity with respect to the maximum norm; if we modify the class and replace $\|x-y\|_\infty$ with $\|x-y\|_q$ for some $q < \infty$, then the model class becomes larger, which means that any lower bound, and in particular the curse of dimensionality, also holds for the modified class.

In fact, the curse of dimensionality already holds true for a simpler problem than $L_p$-approximation, i.e., the problem of **numerical integration** where we only want to estimate the integral

$$S_d(f) := \int_{[0,1]^d} f(x)dx,$$

by using a *sampling algorithm*

$$Q_n(f) = \varphi(f(x_1), \dots, f(x_n)) \quad (8.3)$$

with nodes $x_1, \dots, x_n \in [0,1]^d$ and $\varphi\colon \mathbb{R}^n \to \mathbb{R}$. The worst case error of the algorithm $Q_n$ on a class $F_d \subset L_1([0,1]^d)$ is given by

$$e(Q_n, F_d) := \sup_{f \in F_d} |S_d(f) - Q_n(f)|.$$

The $n$-th minimal worst case error for numerical integration is given by

$$e_n^{\mathrm{int}}(F_d) := \inf_{Q_n} e(Q_n, F_d),$$

with the infimum taken over all mappings of the form (8.3). In fact, from the general theory, it is known that it is enough to restrict the infimum to **quadrature**

**rules** of the form

$$Q_n(f) = \sum_{j=1}^{n} a_j f(x_j)$$

with weights $a_j \in \mathbb{R}$ and nodes $x_j \in [0,1]^d$ if the class $F_d$ is convex and symmetric, see Proposition 9.4, which is the case for the model class $F_d = F_{\mathrm{Lip}}^d$. The information complexity of the integration problem on $F_d$ for error demand $\varepsilon > 0$ is given by

$$n^{\mathrm{int}}(\varepsilon, F_d) \;:=\; \min\left\{n \in \mathbb{N} \colon e_n^{\mathrm{int}}(F_d) \le \varepsilon\right\}.$$

We now explain what we mean when we say that numerical integration is an easier problem than $L_p$-approximation. If $A_n$ is a sampling algorithm for $L_p$-approximation that uses $n$ sampling points, i.e., $A_n(f) = \Phi(f(x_1), \ldots, f(x_n))$ with $x_1, \ldots, x_n \in [0,1]^d$ and $\Phi \colon \mathbb{R}^n \to L_p$, then $Q_n := S_d \circ A_n$ is a sampling algorithm for integration with $n$ nodes and we have

$$|S_d(f) - Q_n(f)| \;\le\; \|f - A_n f\|_1 \;\le\; \|f - A_n f\|_p$$

for all $f \in F_d$. This implies

$$e_n^{\mathrm{int}}(F_d) \;\le\; g_n(F_d, L_p)$$

and hence also

$$n(\varepsilon, F_d, L_p) \;\ge\; n^{\mathrm{int}}(\varepsilon, F_d)$$

for all $p \ge 1$ and $\varepsilon > 0$.

For the integration problem on the model class $F_{\mathrm{Lip}}^d$, we obtain the following bounds. We remark that even more precise bounds are known, see Sukharev (1979), but we are content with the qualitative result here.

**Proposition 8.1.** For all $n, d \in \mathbb{N}$, we have

$$\frac{1}{8} n^{-1/d} \;\le\; e_n^{\mathrm{int}}(F_{\mathrm{Lip}}^d) \;\le\; n^{-1/d}.$$

In particular, for all $\varepsilon \in (0, 1/8)$, $d \in \mathbb{N}$, and $1 \le p \le \infty$ it holds that

$$n(\varepsilon, F_{\mathrm{Lip}}^d, L_p) \;\ge\; n^{\mathrm{int}}(\varepsilon, F_{\mathrm{Lip}}^d) \;\ge\; \left(\frac{1}{8\varepsilon}\right)^d,$$

so that integration and $L_p$-approximation on $F_{\mathrm{Lip}}^d$ suffer from the curse of dimensionality.

*Proof.* We start with the upper bound. We need to show that there exists a quadrature formula $Q_n$ with $e(Q_n, F_{\mathrm{Lip}}^d) \le n^{-1/d}$. For this, we choose $m \in \mathbb{N}$ maximal with $m^d \le n$ (then in particular $(2m)^d \ge n$) and write $M = m^d$. We divide $[0,1]^d$ in $M$ disjoint cubes $W_1, \ldots, W_M$ with side length $1/m$. Let $x_1, \ldots, x_M$ be the centers of these cubes, which form a product set. As quadrature formula we simply choose the product rule

$$Q_n(f) := \frac{1}{M}\sum_{j=1}^{M} f(x_j), \quad f \in F^d_{\mathrm{Lip}}.$$

Then for all $f \in F^d_{\mathrm{Lip}}$, we have the error estimate

$$\begin{aligned}
\left|\int_{[0,1]^d} f(x)dx - \frac{1}{M}\sum_{j=1}^{M} f(x_j)\right| &= \left|\sum_{j=1}^{M}\left(\int_{W_j} f(x)dx - \frac{1}{M} f(x_j)\right)\right| \\
&\le \sum_{j=1}^{M}\left|\int_{W_j} f(x)dx - \frac{1}{M} f(x_j)\right| = \sum_{j=1}^{M}\left|\int_{W_j} f(x)dx - \int_{W_j} f(x_j)dx\right| \\
&\le \sum_{j=1}^{M}\int_{W_j} |f(x) - f(x_j)|dx \le \sum_{j=1}^{M}\int_{W_j} \left\|x - x_j\right\|_\infty dx \\
&\le \sum_{j=1}^{M}\int_{W_j} \frac{1}{2m}dx \le \sum_{j=1}^{M}\frac{\mathrm{vol}(W_j)}{2m} = \frac{1}{2m} \le n^{-1/d}.
\end{aligned}$$

We now prove the lower bound. Let $Q_n$ be an arbitrary quadrature formula with nodes $x_1, \ldots, x_n$. We need to show that $e(Q_n, F^d_{\mathrm{Lip}}) \ge \frac{1}{8}n^{-1/d}$. We consider the function

$$f\colon [0,1]^d \to \mathbb{R}, \quad f(x) = \min_{j=1,\ldots,n} \left\|x - x_j\right\|_\infty.$$

Then $f(x_1) = \ldots = f(x_n) = 0$ and it holds that $\pm f \in F^d_{\mathrm{Lip}}$. The quadrature rule $Q_n$ cannot distinguish between $f$ and $-f$, we have $Q_n(f) = Q_n(-f)$. Therefore we call $f$ a *fooling function* for the quadrature rule. We have

$$\begin{aligned}
e(Q_n, F^d_{\mathrm{Lip}}) &\ge \max_{g\in\{\pm f\}} |S_d(g) - Q_n(g)| \\
&\ge \max\left\{|Q_n(f) - S_d(f)|, |Q_n(f) + S_d(f)|\right\} \ge S_d(f).
\end{aligned}$$

Therefore, it is enough to bound the integral of $f$ from below. Let $W$ be the union of all cubes $W_j$ with center $x_j$ and side length $a = (1/(2n))^{1/d}$. Then $W$ has a volume of at most $na^d = 1/2$. The set $D = [0,1]^d \setminus W$ therefore has a volume of at least $1/2$. Furthermore, for all $x \in D$ we have that $f(x) \ge a/2$, because $x$ has a distance of at least $a/2$ to every cube center $x_j$. Since $f$ is nonnegative, we get

$$\int_{[0,1]^d} f(x)dx \ge \int_D \frac{a}{2}dx \ge \frac{a}{4} \ge \frac{1}{8}n^{-1/d},$$

as required.

□

Hence, the problem of $L_p$-approximation (and even the problem of integration)

on the Lipschitz class $F_{\mathrm{Lip}}^d$ in the worst-case setting is practically unfeasible if the dimension $d$ is large. Before we discuss in the following sections how to avoid the curse of dimensionality under additional assumptions, we shortly discuss whether more smoothness or randomization affect the tractability of our easy example.

### *8.1.1. Does more smoothness help?*

A first guess might be that the reason for the curse of dimensionality lies in the relatively low smoothness of Lipschitz functions. So it is natural to ask whether the curse can be avoided if the model class consists of functions with a higher smoothness. The answer to this question is not completely clear since the term 'smoothness' can be interpreted in many different ways. Unfortunately, however, for classical notions of smoothness, the answer is typically negative. For instance, let us consider the model class

$$\mathcal{C}_d^s := \left\{ f \in C^s([0,1]^d) \colon \|f^{(k)}\|_\infty \le 1 \text{ for all } k \in \mathbb{N}_0^d \text{ with } |k| \le s \right\} \tag{8.4}$$

of $s$-times continuously differentiable functions, where $f^{(k)}$ denotes the $k$-th partial derivative with total order $|k| = k_1 + \ldots + k_d$, and $s, d \in \mathbb{N}$. For this class, it is proven by Hinrichs, Novak, Ullrich and Woźniakowski (2017) that the curse of dimensionality remains present, even with a super-exponential-in-$d$ lower bound.

**Proposition 8.2.** For each $s \in \mathbb{N}$, there are constants $c_s, C_s > 0$ such that for all $\varepsilon \in (0, 1/2)$ and $d \in \mathbb{N}$ it holds that

$$\left(\frac{c_s \cdot d}{\varepsilon}\right)^{d/s} \le n^{\mathrm{int}}(\varepsilon, \mathcal{C}_d^s) \le \left(\frac{C_s \cdot d}{\varepsilon}\right)^{d/s}.$$

In particular, integration (and hence $L_p$-approximation, $1 \le p \le \infty$) on $\mathcal{C}_d^s$ suffer from the curse of dimensionality.

See also Krieg (2019) for an even larger lower bound in the case of uniform approximation ($p = \infty$). The principal proof strategy to obtain the curse of dimensionality is the same as we have seen in the Lipschitz case: The upper bound is again proven by a product rule. For the lower bound for any set $\mathcal{P}_n$ of sampling points, one constructs a fooling function $f \in \mathcal{C}_d^s$ which equals zero on $\mathcal{P}_n$, but which is 'large' on the most part of the domain. This can be done by starting with the distance function to a small neighborhood of $\mathcal{P}_n$ and convoluting this function multiple times with the characteristic function of a small ball to make it more smooth (for the best results, it is important to consider all distances in the $\ell_1$-norm). Further results in this direction can be found in Hinrichs, Novak, Ullrich and Woźniakowski (2014*b*,*c*) and Hinrichs, Prochno and Ullrich (2019), see also the recent survey by Novak (2024).

In Example 3.25, we have already seen that the error rate $\mathcal{O}(n^{-1/d})$ for Lipschitz functions (or $\mathcal{O}(n^{-s/d})$ for $s$-times continuously differentiable functions) can be replaced with a dimension independent error rate if one considers instead function classes of *mixed smoothness*. In this regard, it would be natural to replace the class $\mathcal{C}_d^s$ with the class

$$\mathcal{C}_d^{s,\mathrm{mix}} := \left\{ f \in C^{s,\mathrm{mix}}([0,1]^d) \colon \|f^{(k)}\|_\infty \le 1 \text{ for all } k \in \mathbb{N}_0^d \text{ with } \|k\|_\infty \le s \right\}$$

with continuous partial derivatives up to mixed order $s \in \mathbb{N}$. Indeed, for this class, one can use the algorithm of Smolyak (1963), to show that the $n$-th minimal worst case errors (for integration and approximation) are in $\mathcal{O}(n^{-\alpha})$ for every $\alpha < s$, see Barthelmann, Novak and Ritter (2000, Theorem 8). Unfortunately, this does not mean that the problem becomes tractable. In fact, even for the classes

$$\mathcal{C}_d^\infty := \left\{ f \in C^\infty([0,1]^d) \colon \|f^{(k)}\|_\infty \le 1 \text{ for all } k \in \mathbb{N}_0^d \right\} \tag{8.5}$$

of infinitely smooth functions, the problem of uniform approximation suffers from the curse of dimensionality. This is proven by Novak and Woźniakowski (2009).

**Proposition 8.3.** For all $\varepsilon \in (0,1)$ and $d \in \mathbb{N}$ it holds that

$$n(\varepsilon, \mathcal{C}_d^\infty, L_\infty) \ \ge \ 2^{\lfloor d/2 \rfloor},$$

so that uniform approximation on $\mathcal{C}_d^\infty$ suffers from the curse of dimensionality.

Note that the class $\mathcal{C}_d^\infty$ is extremely small; it only contains analytic functions. In particular, since it is contained in the classes of mixed smoothness, this implies that the problem of uniform approximation suffers from the curse of dimensionality also on the classes $\mathcal{C}_d^{s,\mathrm{mix}}$. It is not known, and a very interesting open problem, whether the problem of numerical integration on the classes $\mathcal{C}_d^\infty$ suffers from the curse of dimensionality. It is only known from Wojtaszczyk (2003) that its complexity is not independent of the dimension, that is, it is not strongly polynomially tractable.

We sketch the proof of Proposition 8.3 to highlight the used *fooling functions*.

*Proof.* The reason for the lower bound is that $\mathcal{C}_d^\infty$ contains a 'large-dimensional ball' of radius one, which implies the lower bound. In fact, if we find some $V_d \subset C^\infty([0,1]^d)$ with $\dim(V_d) = 2^{\lfloor d/2 \rfloor}$ such that

$$\|f\|_\infty = \sup \left\{ \|f^{(k)}\|_\infty \colon k \in \mathbb{N}_0^d \right\} \qquad \text{for all} \quad f \in V_d,$$

then $F_d := \{ f \in V_d \colon \|f\|_\infty \le 1 \} \subset \mathcal{C}_d^\infty$. For $n \le 2^{\lfloor d/2 \rfloor} - 1$ sampling points, there is a function $f \in V_d$ with $\|f\|_\infty = 1$ that vanishes at all sampling points, and so $\pm f \in \mathcal{C}_d^\infty$ cannot be distinguished, which implies $g_n(F_d, L_\infty) \ge 1$. (This proof idea is generalized by the use of Bernstein numbers, see Proposition 9.23 below.)
As the space $V_d$, we set $m := \lfloor d/2 \rfloor$, and consider functions $f$ of the form

$$f(x_1, \ldots, x_d) = \sum_{i \in \{0,1\}^m} a_i \cdot (x_1 + x_2)^{i_1} \cdot (x_3 + x_4)^{i_2} \cdot \ldots \cdot (x_{2m-1} + x_{2m})^{i_m}$$

with $a_i \in \mathbb{R}$. Those functions are linear in each variable and satisfy $\|f^{(k)}\|_\infty \le \|f\|_\infty$ for all $k \in \mathbb{N}_0^d$. We refer to Novak and Woźniakowski (2009) for details. □

**Remark 8.4.** As noted in (Novak and Woźniakowski 2009, Remark 3), the same proof technique can be used to show the curse of dimensionality for $L_p$-approximation for the (larger) classes

$$F^{\infty}_{d,p} := \left\{ f \in C^{\infty}([0,1]^d) \colon \|f^{(k)}\|_p \leq 1 \text{ for all } k \in \mathbb{N}_0^d \right\}.$$

In fact, for $1 \leq p < \infty$, we have $n(\varepsilon, F^{\infty}_{d,p}, L_p) \geq 2^{\lfloor d/\ell \rfloor}$ for $\ell = \lceil 8(p+1)^{2/p} \rceil$ as shown by Weimar (2012*b*, Proposition 3).

Another very drastic example is due to Vybíral (2020). Proving a conjecture of Novak (1999), it is shown that the curse of dimensionality already holds for the integration problem on a class of trigonometric polynomials with frequencies in $\{-1,0,1\}^d$. Namely, Vybíral (2020) considers the $d$-fold tensor product of the reproducing kernel Hilbert space with orthonormal basis

$$b_1(x) = 1, \quad b_2(x) = \cos(2\pi x), \quad b_3(x) = \sin(2\pi x), \qquad x \in [0,1].$$

He proves that at least $2^{d-1}$ samples are needed for a quadrature error smaller than $1/\sqrt{2}$. Further recent results about the curse of dimensionality on Hilbert spaces can be found in Hinrichs, Krieg, Novak and Vybíral (2021*a*).

All this shows that classical smoothness assumptions are usually not suitable for tractability of high-dimensional approximation problems. We will see in Section 8.3 that *structural assumptions* on the target functions are often better suited for tractability. However, the question for tractability is often very subtle and depends heavily on the precise model assumptions. For example, Vybíral (2014) considers a different class of infinitely smooth functions, given by

$$\widetilde{C^{\infty}_d} := \left\{ f \in C^{\infty}([0,1]^d) \colon \sum_{k \in \mathbb{N}_0^d \colon |k| = s} \frac{\|f^{(k)}\|_{\infty}}{k!} \leq 1 \;\text{ for all } s \in \mathbb{N}_0 \right\}$$

and proves that the uniform approximation problem on this class does not suffer from the curse of dimensionality (it is quasi-polynomially tractable). Also Hinrichs, Novak and Ullrich (2014*a*) and Xu (2015) consider a modified class of infinitely smooth functions and prove a very weak form of tractability for integration and $L_p$-approximation with $p < \infty$. Hence, although smoothness does not seem to be a good model assumption for tractability studies, a general claim of the form "smoothness does not help" is certainly also incorrect.

### *8.1.2. Does randomization help?*

We have already seen in Section 7 that the (expected) error of Monte Carlo algorithms can be much smaller than the (guaranteed) error of deterministic algorithms. So a next natural question is whether we can also get rid of the curse of dimensionality by allowing this larger class of algorithms with their correspondingly relaxed error criterion.

In the case of the integration problem, the answer is a definite yes. More precisely, if we want to estimate an integral

$$S(f) = \int_D f(x) d\mu(x)$$

with respect to some probability measure $\mu$, already the simplest Monte Carlo method

$$M_n(f) := \frac{1}{n} \sum_{i=1}^{n} f(X_i), \qquad X_i \overset{\text{i.i.d.}}{\sim} \mu$$

satisfies the error bound

$$\sqrt{\mathbb{E}\,|S(f) - M_n(f)|^2} \leq \frac{\|f\|_2}{\sqrt{n}}, \tag{8.6}$$

for all $f \in L_2(D, \mu)$, which is completely independent of the dimension. Therefore, for any model class which is contained in the unit ball of $L_2(D, \mu)$, the complexity of the integration problem in the randomized setting can be bounded independently of the dimension and we have strong polynomial tractability.

Unfortunately, in the case of the approximation problem, i.e., the problem of full recovery of the function $f$, the answer is less clear. Randomization can reduce the dimension dependence of the complexity significantly in some situations, while it does not change a lot in other situations. A comprehensive discussion of this question can be found in the dissertation of Kunsch (2017).

For uniform approximation, the answer is often negative. Indeed, Kunsch (2017, Corollary 2.20) proves that the curse of dimensionality holds for uniform approximation even on the class $\mathcal{C}_d^\infty$ of infinitely smooth functions from (8.5). The proof relies on a lower bound for the error of randomized algorithms in terms of Bernstein numbers, see Section 9.3. The lower bound even holds for all randomized algorithms that use arbitrary (adaptive) linear measurements (instead of only function evaluations), and hence the curse of dimensionality reigns for the uniform approximation of infinitely smooth functions for a very general class of algorithms.

For $L_2$-approximation, on the other hand, randomization often helps to mitigate the curse of dimensionality. An example is the class

$$\mathcal{H}_d^{1,\text{mix}} := \left\{ f \in \mathbf{W}_2^1([0,1]^d) : \sum_{k \in \{0,1\}^d} \left\| f^{(k)} \right\|_2^2 \leq 1 \right\}, \tag{8.7}$$

the unit ball of the Sobolev Hilbert space of mixed smoothness $s = 1$, as already encountered in Example 3.25. For this example, it can be proven with the technique of decomposable kernels from Novak and Woźniakowski (2001), that integration and $L_p$-approximation suffer from the curse of dimensionality in the deterministic setting, see also Novak and Woźniakowski (2016) and Novak and Woźniakowski (2010, Section 11.6.1). However, there is no curse of dimensionality in the randomized setting. In order to obtain an expected error of at most $\varepsilon$ with any fixed $\varepsilon > 0$, it suffices to use $\mathcal{O}(d^q)$ random samples, where $q = \ln \varepsilon^{-1}$, see Krieg (2019, Example 2). This follows from bounds on the singular numbers (or Kolmogorov numbers) of tensor product operators between Hilbert spaces together with the error bounds for the Monte Carlo methods from Section 7 (least squares or multilevel Monte Carlo) in terms of these numbers. The complexity of the problem in the randomized setting is best described by the term **quasi-polynomial tractability**, see Gnewuch and Woźniakowski (2011).

### 8.2. *Sufficient criteria for tractability*

In the Sections 3 and 5, we gave several bounds on the sampling numbers in terms of other approximation benchmarks. From these bounds, we easily obtain several sufficient criteria for (strong) polynomial tractability of the approximation problem. In a certain sense, these results shift the question of tractability from one problem (sampling recovery) to another problem (e.g., a covering problem in the case of entropy numbers). But since these other approximation benchmarks are often much easier to handle than dealing directly with the sampling numbers, such transference results can be (but do not have to be) quite useful.

The criteria presented here all apply to the following situation: $D_d$ is a subset of $\mathbb{R}^d$ for each $d \in \mathbb{N}$, $\mu_d$ is a probability measure on $D_d$, and the model class $F_d$ is a separable subset of $B(D_d, \mu_d)$, the space of bounded complex-valued and $\mu_d$-measurable functions with the sup-norm. We consider the problem of $L_p$-approximation on $F_d$, where $L_p = L_p(D_d, \mu_d)$. In fact, these conditions are stronger than necessary, but it is quite convenient to have a simple unified assumption in mind. We start with the case $p = 2$, where the most results are available. The following is obtained from Corollary 3.31.

**Corollary 8.5.** There is an absolute constant $c \in \mathbb{N}$ such that the following holds. Let $D_d$, $\mu_d$ and $F_d$ satisfy Assumption A. Assume that

$$C_d \;:=\; \sum_{k=1}^{\infty} \frac{d_k(F_d, L_2)}{\sqrt{k}} \;<\; \infty. \tag{8.8}$$

The problem of sampling recovery on $F_d$ in the $L_2$-norm is polynomially tractable if $C_d \in \mathcal{O}(d^s)$ for some $s < \infty$ and strongly polynomially tractable if $C_d \in \mathcal{O}(1)$.

*Proof.* Corollary 3.31 implies the bound

$$g_{cn}^{\mathrm{lin}}(F_d, L_2) \;\leq\; \frac{C_d}{\sqrt{n}}, \tag{8.9}$$

which gives

$$n(\varepsilon, F_d, L_2) \;\leq\; c \cdot \left\lceil C_d^2 \cdot \varepsilon^{-2} \right\rceil$$

for all $\varepsilon > 0$ and $d \in \mathbb{N}$.

□

Hence, according to the formula (8.9), we have a situation similar to that of Monte Carlo integration (8.6) also for deterministic $L_2$-approximation whenever we have a uniform bound for the sum (8.9).

In the area of information-based complexity, it is not so common to work with the Kolmogorov widths, since these numbers do not characterize minimal errors of algorithms. It is often preferred to work instead with a different classical approximation benchmark, the **Gelfand widths**, which describe minimal errors of a family of algorithms. In fact, the Gelfand $n$-width of a convex and symmetric set $F \subset L_2$ satisfy

$$c_n(F, L_2) \asymp \inf_{\substack{\Phi\colon \mathbb{C}^n \to L_2 \\ N \in \mathcal{L}(F,\mathbb{C}^n)}} \sup_{f \in F} \left\| f - \Phi \circ N(f) \right\|_2,$$

with the hidden constants between $\frac{1}{2}$ and 2, so that it measures the worst case error of an optimal (possibly non-linear) algorithm using $n$ arbitrary (non-adaptive) linear measurements, see Section 9.3 for details. It is known that

$$d_n(F, L_2) = a_n(F, L_2) \leq \left(1 + \sqrt{n}\right) c_n(F, L_2),$$

where $a_n(F, L_2)$ are the approximation numbers (or linear widths). This result goes back to Pietsch (1974), see Mathé (1990) and Creutzig and Wojtaszczyk (2004) for the present formulation. We will discuss such relations further in Section 9.2, and prove the bound in Proposition 9.13. The reason we discuss this here is the observation that, for convex and symmetric classes $F_d$, the constant $C_d$ in Proposition 8.5 can be replaced by the (larger) constant

$$C'_d := \sum_{k=1}^{\infty} c_k(F_d, L_2) < \infty. \tag{8.10}$$

In particular, sampling recovery on $F_d$ in the $L_2$-norm is strongly polynomially tractable whenever we have a uniform summability of the Gelfand widths.

In the case of reproducing kernel Hilbert spaces, the condition (8.10) can be relaxed. Instead of uniform absolute summability it is enough to have uniform square-summability of the Gelfand numbers. This follows from the results of Section 3.3.1. In fact, it already follows from the previous work of Wasilkowski and Woźniakowski (2001), see also Novak and Woźniakowski (2012, Section 26.4).

**Corollary 8.6.** There is an absolute constant $c \in \mathbb{N}$ such that the following holds. For $d \in \mathbb{N}$, let $(D_d, \mu_d)$ be a measure space and let $F_d$ be the unit ball of a reproducing kernel Hilbert space $H_d$ that is (injectively) embedded into $L_2(D_d, \mu_d)$ and whose kernel $k_d$ has finite trace

$$W_d \;:=\; \int_{D_d} k_d(x,x)\, d\mu_d(x) \;<\; \infty.$$

The problem of sampling recovery on $F_d$ in the $L_2$-norm is polynomially tractable if $W_d \in \mathcal{O}(d^s)$ for some $s < \infty$ and strongly polynomially tractable if $W_d \in \mathcal{O}(1)$.

*Proof.* Recall that, for $L_2$-approximation in a unit ball of a Hilbert space, the Kolmogorov widths are equal to the Gelfand numbers and the singular numbers, and that the square-sum of the singular numbers equals the trace of the kernel, see Sections 3.3.1 and 9.2, or Pietsch (1987). So according to formula (3.56), we have

$$g_{4n}^{\mathrm{lin}}(F_d, L_2) \;\le\; 9\sqrt{\frac{1}{n}\sum_{k>n}\sigma_k^2} \;\le\; 9\sqrt{\frac{W_d}{n}}, \tag{8.11}$$

which gives

$$n(\varepsilon, F_d, L_2) \;\le\; 4\cdot\left\lceil 81\cdot W_d\cdot\varepsilon^{-2}\right\rceil$$

for all $\varepsilon > 0$ and $d \in \mathbb{N}$.

□

For the problem of $L_p$-approximation with $1 \le p \le \infty$, we obtain a sufficient condition for tractability in terms of the Kolmogorov widths in the uniform norm from Theorem 3.17, which reads as follows.

**Corollary 8.7.** For each $d \in \mathbb{N}$, let $(D_d, \mu_d)$ be a probability space and let $F_d \subset B(D_d)$ be a class of bounded and measurable functions. Assume that there exists some $\alpha > \max\{1/2 - 1/p, 0\}$ such that

$$R_d \;:=\; \sup_{n\in\mathbb{N}} n^{\alpha}\, d_n(F_d, B(D_d)) \;<\; \infty \tag{8.12}$$

for all $d \in \mathbb{N}$. The problem of sampling recovery on $F_d$ in the $L_p$-norm is polynomially tractable if $R_d \in \mathcal{O}(d^s)$ for some $s < \infty$ and strongly polynomially tractable if $R_d \in \mathcal{O}(1)$.

*Proof.* From Theorem 3.17, we obtain that

$$n(\varepsilon, F_d, L_p) \;\le\; 4\left\lceil\left(\frac{4R_d}{\varepsilon}\right)^{t}\right\rceil$$

with $t = (\alpha - \max\{1/2 - 1/p, 0\})^{-1}$.

□

Finally, in the case of uniform approximation, we obtain further sufficient criteria for tractability in terms of the (nonlinear) approximation benchmarks considered in Section 5.3. Namely, the relation (5.23) implies the following.

**Corollary 8.8.** For each $d \in \mathbb{N}$, let $D_d$ be a set and let $F_d$ be the unit ball of a Banach space that is embedded into $B(D_d)$. Moreover, let $s \in \{h, \varepsilon, \delta\}$ stand for the Hilbert numbers, entropy numbers, or manifold widths. Assume that there exists some $\alpha > 1$ such that

$$V_d := \sup_{n\in\mathbb{N}} (n+1)^\alpha s_n(F_d, B(D_d)) < \infty \tag{8.13}$$

for all $d \in \mathbb{N}$. The problem of sampling recovery on $F_d$ in the sup-norm is polynomially tractable if $V_d \in \mathcal{O}(d^q)$ for some $q < \infty$ and strongly polynomially tractable if $V_d \in \mathcal{O}(1)$.

*Proof.* From (5.23), we obtain that

$$n(\varepsilon, F_d, L_p) \le \left\lceil \left( \frac{4e^\alpha \cdot V_d}{\varepsilon} \right)^{1/(\alpha-1)} \right\rceil .$$

□

We only presented here the most simple criteria for (strong) polynomial tractability that we obtain immediately from the theory of Sections 3 and 5. A closer look easily reveals statements also about the exponents in the polynomial bounds on the complexity. We remark that similar transference results for **exponential tractability** are discussed by Krieg, Siedlecki, Ullrich and Woźniakowski (2023). The notion of exponential tractability is better suited to the situation of an exponential convergence of minimal errors, which typically occurs for model classes of *analytic functions*, see, for example, Dick, Kritzer, Pillichshammer and Woźniakowski (2014), Griebel and Oettershagen (2016), Sloan and Woźniakowski (2018).

### *8.3. A collection of tractable problems*

As we have seen, classical smoothness assumptions are not very well suited for high-dimensional approximation problems; the approximation problem on such model classes usually suffers from the curse of dimensionality. Here, we want to collect some other model assumptions from the literature that make the problem of sampling recovery tractable. Successful model assumptions mostly come in the form of structural knowledge, often in terms of some lower-dimensional structure that is hidden within the high-dimensional function. Of course, in practice, one cannot choose the model assumptions and the properties of the target function are given by the application. However, we think it is nonetheless interesting to collect properties that lead to tractability; it can provide motivation to examine the target function with respect to these properties and give an idea as to which properties should be exploited instead of smoothness in the approximation of high dimensional functions.

**Example 8.9 (Separable functions).** Sloan, Wang and Woźniakowski (2004) and Wasilkowski and Woźniakowski (2004) study the tractability of integration and approximation in classes of functions of the form

$$f(x_1, \dots, x_d) = \sum_{i_1,\dots,i_m \in \{1,\dots,d\}} f_{i_1,\dots,i_m}(x_{i_1}, \dots, x_{i_m}), \tag{8.14}$$

i.e., $d$-variate functions that can be written as a sum of functions that only depend on a smaller number $m < d$ of variables. Such functions are called partially separable. Wasilkowski and Woźniakowski (2004) obtain polynomial tractability for $L_2$-approximation on classes of such functions in the case that $m \in \mathcal{O}(\log d)$, and even strong polynomial tractability for constant $m$ under additional assumptions.

**Example 8.10 (Weighted Korobov spaces).** A classical approach to relinquish the curse of dimensionality, which was introduced by Sloan and Woźniakowski (1998), is by considering weighted spaces, where the weights model the influence of (groups of) certain variables. See, e.g., Ebert and Pillichshammer (2021), Novak, Sloan and Woźniakowski (2004), Kuo, Wasilkowski and Woźniakowski (2009*b*), Dick, Kuo and Sloan (2013) and the book Dick, Kritzer and Pillichshammer (2022, Ch. 13), where many more references can be found.

The most well-studied type of weights are product weights, we present a result in this setting as a generic example. Let a sequence of weight parameters $1 \geq \gamma_1 \geq \gamma_2 \geq \dots > 0$ be given and assumed that the target function $f\colon [0,1]^d \to \mathbb{C}$ is continuous with

$$\sum_{k \in \mathbb{Z}^d} |\hat{f}(k)|^2 \prod_{i=1}^{d} \max\left\{1, \frac{|k_i|^2}{\gamma_i}\right\} \leq 1. \tag{8.15}$$

The class $\mathcal{H}^{1,\mathrm{mix}}_{d,\gamma}$ of all such functions is a weighted version of the class $\mathcal{H}^{1,\mathrm{mix}}_{d}$ of mixed smoothness from (8.7); both classes are unit balls of the same Hilbert space but with respect to different equivalent norms. The sequence $\gamma$ models a decreasing importance of the input variables: For variables $i \leq d$ for which $\gamma_i$ is small, the relation (8.15) implies that the Fourier coefficients with frequencies $k_i \neq 0$ have to be small and so $f(x)$ can vary only little with respect to $x_i$. It is proven by Novak *et al.* (2004, Theorem 1) that $L_2$-approximation on this class is polynomially tractable if the weights decay fast enough, namely, if

$$\sum_{i=1}^{d} \gamma_i \in \mathcal{O}(\log(d)). \tag{8.16}$$

There are many more results of this nature and other types of weights exist which can be better suited to particular applications. For instance, so-called product-and-order-dependent (POD) weights, as introduced by Kuo, Schwab and Sloan (2012), arose in their application to PDEs with random coefficients. In fact, weighted Korobov spaces can be regarded as a generalization of the previous example; classes of separable functions are obtained by considering so-called finite-order weights, see Wasilkowski and Woźniakowski (2004).

We note that all the classes $\mathcal{H}^{1,\mathrm{mix}}_{d,\gamma}$ are unit balls of the Sobolev space of dominating mixed smoothness one, only with respect to different equivalent norms. We hence observe that changing between equivalent norms in a Banach space, although this does not change the order of convergence of minimal errors, can have a huge effect on tractability (even if the problem is properly normalized). This can also be seen in the unweighted case, as discussed by Kühn, Sickel and Ullrich (2015).

**Example 8.11 (Summable Fourier coefficients).** The $L_p$-approximation problem with $1 \leq p < \infty$, and hence also the integration problem, become tractable on many classical smoothness classes if it is additionally known that the Fourier coefficients have a bounded absolute sum

$$\sum_{k\in\mathbb{Z}^d} \left|\hat{f}(k)\right| \leq 1, \tag{8.17}$$

where $\hat{f}(k)$ are the standard Fourier coefficients from (5.11). For example, $L_p$-approximation is polynomially tractable on the class

$$F^{d,\mathrm{sum}}_{\mathrm{Lip}} = \left\{ f \in F^d_{\mathrm{Lip}} : f \text{ satisfies (8.17)} \right\},$$

where $F^d_{\mathrm{Lip}}$ is the class of Lipschitz continuous functions from (8.2). Such results follow from the theory outlined in Section 5.1, e.g., by using the square root Lasso algorithm from (5.3) with an appropriate choice of $V_N$. The above and several other such examples are discussed by Krieg (2024). Numerical integration on such classes is discussed by Dick (2014), Goda (2023), Chen and Jiang (2024), Dick, Goda and Suzuki (2025).

Another example that is covered via this approach is the class

$$\mathcal{H}^{1,\mathrm{mix}}_{d,\gamma,*} := \bigcup_{\pi\in\mathcal{S}_d} \mathcal{H}^{1,\mathrm{mix}}_{d,\pi(\gamma)}$$

where $\mathcal{H}^{1,\mathrm{mix}}_{d,\gamma}$ is the class from the previous example, defined by (8.15), and the union runs over all permutations of $\gamma$. This class models the situation that we have a decreasing influence of the input variables, but we do not know in advance which of the input variables are the important ones. We obtain polynomial tractability on this class under the same condition (8.16) as before.

In this sense, the approach via summable Fourier coefficients gives stronger tractability results than the approach via weighted spaces. However, it is also important to take into account the runtime of the algorithms; though both algorithms require only a polynomial amount of samples, the square root Lasso algorithm has to be performed on a space $V_N$ of dimension $N \geq 2^d$ and is hence not feasible for a very large dimension $d$. Further research is required to determine whether the tractability results in this example can actually be matched by an algorithm with a polynomial runtime.

**Example 8.12 (Ridge functions).** Mayer, Ullrich and Vybíral (2015), Kolleck and Vybíral (2015) and Malykhin, Ryutin and Zaitseva (2022) examine tractability under the assumption that the target function is given as a composition of a univariate function with a linear functional, i.e.,

$$f(x_1, \ldots, x_d) \;=\; g\,(a_1 x_1 + \ldots + a_d x_d)\,. \tag{8.18}$$

Such functions are called ridge functions. For example, Mayer *et al.* (2015) obtain weak tractability for the class of all functions of the form (8.18) with coefficients $\sum_{i=1}^{d} |a_i| \le 1$ and $g \in C^3([-1,1])$ with derivatives bounded by a constant.

**Example 8.13 (Symmetric functions).** Weimar (2012*a*), Nuyens, Suryanarayana and Weimar (2016), and Bachmayr, Dusson, Ortner and Thomas (2024) consider approximation and integration of symmetric and antisymmetric functions. A function $f\colon [0,1]^d \to \mathbb{C}$ is called *fully symmetric* if

$$f(x) = f(\pi(x)) \tag{8.19}$$

for all $x \in [0,1]^d$ and all permutations $\pi$. For example, it follows from Weimar (2012*a*) that randomized sampling recovery on the class

$$\mathcal{H}^{1,\mathrm{mix}}_{d,\mathrm{sym}} \;:=\; \Big\{ f \in \mathcal{H}^{1,\mathrm{mix}}_{d} \;:\; f \text{ is fully symmetric} \Big\}$$

in the $L_2$-norm is strongly polynomially tractable. Here, $\mathcal{H}^{1,\mathrm{mix}}_{d}$ is the class of mixed smoothness from (8.7). Recall that, without the symmetry assumption, we only have quasi-polynomial tractability, so the symmetry assumption significantly improves the situation in high dimensions. We note that Weimar (2012*a*) considers the minimal errors of algorithms that use arbitrary linear information (which are described by the Kolmogorov numbers in $L_2$) and we have to additionally invoke the results from Section 7.2 (i.e., the upper bounds for randomized sampling algorithms in terms of Kolmogorov numbers) in order to obtain the aforementioned result.

**Example 8.14 (Rank-one tensors).** Bachmayr, Dahmen, DeVore and Grasedyck (2014), Novak and Rudolf (2016) and Krieg and Rudolf (2019) study sampling recovery of functions of the form $f = \bigotimes_{i=1}^{d} f_i$, i.e.,

$$f(x) \;=\; \prod_{i=1}^{d} f_i(x_i) \qquad \text{for} \qquad x = (x_1, \ldots, x_d). \tag{8.20}$$

For example, it is proven in Krieg and Rudolf (2019) that the uniform approximation problem is polynomially tractable on the class

$$F^{\otimes d}_{\mathrm{Lip}} \;:=\; \Big\{ f\colon [0,1]^d \to [-1,1] \;:\; f \text{ satisfies (8.20) with } f_1, \ldots, f_d \in F_{\mathrm{Lip}} \Big\},$$

where $F_{\mathrm{Lip}}$ is the class of univariate Lipschitz functions $f\colon [0,1] \to [-1,1]$ with Lipschitz constant bounded by one.

**Example 8.15 (Bounded variation).** For numerical integration, another assumption that leads to polynomial tractability is the assumption of a bounded variation in the sense of Hardy and Krause, which is satisfied, for example, by completely monotone functions, see Aistleitner and Dick (2015). This follows from the Koksma–Hlawka inequality, which relates the integration error of a quadrature rule with the star discrepancy of the quadrature points, together with the tractability of the problem of star discrepancy as proven by Heinrich, Novak, Wasilkowski and Woźniakowski (2001). Note that the curse of dimensionality for the $L_p$-discrepancy has only been observed recently by Novak and Pillichshammer (2025). We mention this here only briefly since the main topic of this work is the approximation of functions.

## 9. A general framework for optimal algorithms

So far, we only discussed sampling algorithms in the sense that the available measurements consist solely of function evaluations. Now we want to put this into a broader context. That is, we want to discuss *optimal algorithms* that have access to a class of more general measurements or may only be allowed to use specific reconstruction maps. This can be used to study other aspects of *optimality* such as:

- Is there always an easy (e.g., linear) optimal algorithm?
- Are optimal measurements from one class as good as optimal measurements from another?
- Are randomized and/or adaptive algorithms better than deterministic nonadaptive algorithms? If yes, by how much?

Questions of this kind, where the power of one kind of algorithms is compared with the power of another kind of algorithms, are common in the area of *discrete complexity theory*. An example is the P-versus-NP problem. Here, we typically consider continuous problems (like function approximation) where an exact solution is impossible with finite information and in finite time. For continuous problems, such questions are studied in the area of **information-based complexity (IBC)**. This branch of computational mathematics was essentially initiated by Bakhvalov (1959), who studied the complexity of numerical integration. The theory was systematically developed by Traub and Woźniakowski (1980) and Traub *et al.* (1988). We also refer to the introductory articles of Traub and Woźniakowski (1991) (published in *The Mathematical Intelligencer*) or Woźniakowski (2019), as well as the books of Novak (1988) and Traub and Werschulz (1998) and the comprehensive treatment by Novak and Woźniakowski (2008, 2010, 2012).

Before we discuss in Section 10 some specific instances of the above questions, with an emphasis on relations to sampling (in the sense of function evaluations), we introduce a general framework of optimal algorithms in Section 9.1 that also allows interpreting some earlier used *benchmarks* as errors of optimal algorithms within a certain class. This might be seen as a unified approach to certain parts of IBC and approximation theory. Another perspective is given in Section 9.2, where we present the theory of *s-numbers*; a crucial concept from operator theory that provides us with a lot of powerful techniques. Finally, in Section 9.3, we discuss the relation between these concepts. We state some fundamental results from these areas along the way.

### *9.1. Concepts of optimal algorithms*

This section is primarily based on Chapter 4 of Traub *et al.* (1988).

We are given a set $F$ of possible *problem elements*, which we call a **model class**, and a metric space $Y$, which specifies the *error measure*. For convenience, we assume here that $Y$ is a (semi-)normed space. (The generalization of the definitions below is obvious.) In addition, we are given a **solution map**

$$S\colon F \to Y$$

that specifies the *problem* in the following sense: For every $f \in F$, we want to compute (an approximation of) the *solution* $S(f)$. This notion is quite flexible and also allows us to model the computation of, e.g., functionals, or solutions of differential equations (with *right-hand-side* $f$), or even optimization problems, see (Traub and Woźniakowski 1980, Chapter 8) or (Traub *et al.* 1988, Chapter 5). Also note that we no longer require the problem/solution elements to be functions, although this will be a typical assumption since function approximation is the main topic of this work. In the previous sections, we only considered the **inclusion map** (or **identity** or **embedding**) $S = \mathrm{id}$ with

$$\mathrm{id}(f) := f \qquad \text{for} \quad f \in F \subset Y,$$

which corresponds to the approximation of $f \in F$ itself in $Y$, and for most of the upcoming considerations, it is okay if the reader thinks of this special case. However, many of the following results hold directly for general (linear/continuous) $S$, and we think it is instructive to state them like that.

We first generalize some concepts that were already discussed in Section 2. Recall that by a (deterministic) **algorithm** $A_n\colon F \to Y$ we understand a mapping of the form $A_n = \Phi \circ N_n$, where the **measurement map** (aka *encoder*) and the **reconstruction map** (aka *decoder*)

$$N_n\colon F \to \mathbb{R}^n \qquad \text{and} \qquad \Phi\colon \mathbb{R}^n \to Y,$$

respectively, are mappings to be specified. Here, we follow the usual convention to consider only real-valued measurements, and also the problem elements are usually assumed to be elements of a real Banach space. We refer to Novak (1995*c*) and Woźniakowski (1999) for a discussion of the *real-number model*.

The **worst-case error** of a particular algorithm $A_n$ for the approximation of $S\colon F \to Y$ is defined by

$$e(A_n, S) := \sup_{f \in F} \|S(f) - A_n(f)\|_Y . \tag{9.1}$$

Observe that the error also depends on $F$ and $Y$ and that the short notation is due to *hiding* $F$ and $Y$ as domain and codomain of $S$.

We repeat that this quantity gives a reliable (guaranteed) a priori error bound in the class $F$. See again Remark 2.6 for a general comment on the worst-case setting and on the meaning of the class $F$. Also note that the *worst-case setting* (based on exact measurements) is only one choice, and also *randomized* (discussed below), as well as *average*, *probabilistic* and *noisy* (not discussed), settings are studied in IBC. We refer again to Traub *et al.* (1988) and Novak and Woźniakowski (2008).

We aim for *minimal errors*, i.e., the error of an optimal algorithm in some class of algorithms. Hence, let us turn to the specification of measurements and reconstructions. Until now, we restricted ourselves to function evaluations as measurements, i.e., to $N_n(f) = (f(x_1), \ldots, f(x_n))$. Now, we replace function evaluations with any class of **admissible measurements** $\Lambda \subset \{\lambda\colon F \to \mathbb{R}\}$ which, in principle, can be arbitrary real-valued functions on $F$.

We will use the symbol $\Lambda$ throughout for the class of admissible measurements, and we will add sub- and superscripts to indicate special classes. In particular, if the measurements are given by function evaluations, which requires the problem elements to be functions on some set $D$, then we write

$$\Lambda^{\text{std}} := \{\delta_x : x \in D\}\,, \tag{9.2}$$

where $\delta_x$ denotes the *Dirac functional* $\delta_x(f) := f(x)$, and the "std" stands for "standard"; in IBC function evaluations are usually called **standard information**. To give a broader picture, we will consider here also the class of *arbitrary measurements*

$$\Lambda^{\text{arb}} := \mathbb{R}^F := \{\lambda\colon F \to \mathbb{R}\}\,, \tag{9.3}$$

the class of *continuous measurements*

$$\Lambda^{\text{con}} := C(F,\mathbb{R}) := \big\{\lambda \in \mathbb{R}^F :\ \lambda \text{ is continuous}\big\}, \tag{9.4}$$

and the class of *continuous linear measurements*

$$\Lambda^{\text{lin}} := \mathcal{L}(F,\mathbb{R}) := \big\{\lambda \in C(F,\mathbb{R})\colon\ \lambda \text{ is linear}\big\}. \tag{9.5}$$

The latter two classes, of course, are only well-defined if $F$ has a topology and, in the last case, also a linear structure. An important special case, where both these structures are given and which is the setting we will mostly consider here, is that $F = B_X$ is the unit ball in some normed space $X$, in which case $F$ inherits the topology of $X$, and $\Lambda^{\text{lin}} = X'$ is the dual space of $X$.
We do not consider the class $\Lambda^{\text{arb}}$ to be of practical interest, it is mainly used to incorporate some standard approximation benchmarks in the framework. See also Remark 9.6 for a comment on other classes of measurements.

Before we discuss different ways to turn individual measurements into a measurement map $N_n\colon F \to \mathbb{R}^n$, let us consider classes of **admissible reconstruction maps**. Similarly to above, we consider the class of *arbitrary reconstructions*

$$\Psi^{\text{arb}} := \{\Phi\colon \mathbb{R}^n \to Y\colon\ n \in \mathbb{N}\}\,, \tag{9.6}$$

the class of *continuous reconstructions*

$$\Psi^{\text{con}} := \big\{\,\Phi\colon \mathbb{R}^n \to Y \text{ continuous}\colon\ n \in \mathbb{N}\big\}, \tag{9.7}$$

and the class of *linear reconstructions*

$$\Psi^{\text{lin}} := \big\{\,\Phi\colon \mathbb{R}^n \to Y \text{ linear}\colon\ n \in \mathbb{N}\big\}. \tag{9.8}$$

The reconstruction map clearly has to be a mapping from $\mathbb{R}^n$ to the target space $Y$,

where $n$ is the amount of information and has to fit to the measurement map $N_n$. For a given class $\Psi \subset \Psi^{\mathrm{arb}}$, we write $\Psi_n$ for the subclass of all reconstruction maps $\Phi \in \Psi$ with domain $\mathbb{R}^n$, i.e., with a fixed $n \in \mathbb{N}$.

Let us finally consider how algorithms can be derived from admissible measurements $\Lambda$ and reconstructions $\Psi$, respectively. For this, it remains to specify how an algorithm is allowed to *assemble* information. That is, how the actual procedure (like a *pseudo-code*) of choosing specific measurements looks like. Every *decision* for a measurement must be based solely on the available information about the particular problem element $f$, which might be the assumption "$f \in F$", but also some previous measurements.

A first and obvious approach, which is also typical in approximation theory and has been used by us until now, is to fix measurements from $\Lambda$ based on (assumptions on) $F$ alone, and independent of a particular $f$. In this case, we consider measurement maps $N_n\colon F \to \mathbb{R}^n$ of the form

$$N_n = (\lambda_1, \ldots, \lambda_n) \in \Lambda^n,$$

and we call such measurement maps, and corresponding algorithms $A_n = \Phi \circ N_n$, **non-adaptive**. Such algorithms use the same measurements for each input.

As before, but with slightly different notation, we define the **$n$-th minimal worst-case errors of non-adaptive algorithms** for approximating $S\colon F \to Y$ based on measurements $\Lambda$ and reconstructions $\Psi$ by

$$e_n^{\mathrm{non}}(S, \Lambda, \Psi) \;:=\; \inf_{\substack{N_n \in \Lambda^n \\ \Phi \in \Psi_n}} \; \sup_{f \in F} \; \Big\| S(f) - \Phi \circ N_n(f) \Big\|_Y .$$

However, there is often no reason to fix the measurement maps in advance. In general, it seems much cleverer to *choose* a measurement based on all available information. This includes the possibility of choosing measurements during the run-time of an algorithm, when already some previous measurements have been obtained and can be used in the decision which measurement to take next. Formally, this can be modeled by defining a measurement map recursively by

$$N_n(f) = \big(N_{n-1}(f),\, \lambda_n(f)\big), \tag{9.9}$$

where the choice of the $n$-th measurement $\lambda_n = \lambda_n(\cdot, N_{n-1}(f)) \in \Lambda$ may depend on the first $n-1$ measurements, i.e., the $n$-th mapping is of the form

$$\lambda_n\colon F \times \mathbb{R}^{n-1} \to \mathbb{R}$$

with $\lambda_n(\cdot, y) \in \Lambda$ for each $y \in \mathbb{R}^{n-1}$. Such measurement maps and the corresponding algorithms $A_n = \Phi \circ N_n$ are called **adaptive**. It is important to note that measurement maps $N_n$ that are based on the class $\Lambda^{\mathrm{lin}}$ of linear (or $\Lambda^{\mathrm{con}}$ of continuous) measurements are usually not linear (or continuous) as mappings from $F$ to $\mathbb{R}^n$, since the adaption step can be nonlinear and discontinuous. In fact, this is actually the standard case, since adaption often involves a discontinuous case distinction (like an if-else-element in the algorithm).
This is in contrast to $\Lambda^{\mathrm{arb}}$: Since an adaptive measurement map is still a (recursively defined) function, it can be considered non-adaptive and there is no difference between adaptive and non-adaptive measurement maps in this case. We also refer to Remark 9.9 for a general comment on adaption in our context.

The **$n$-th minimal worst-case errors** (of adaptive algorithms) for approximating $S\colon F \to Y$ based on measurements $\Lambda$ and reconstructions $\Psi$ is defined by

$$e_n(S, \Lambda, \Psi) \;:=\; \inf_{\substack{N_n \\ \Phi\in\Psi_n}} \; \sup_{f\in F} \left\| S(f) - \Phi \circ N_n(f) \right\|_Y, \tag{9.10}$$

where the infimum is over all adaptive measurement maps $N_n$ from $\Lambda$ as above.

Adaptive algorithms are not only a theoretical concept, but are quite common in applications. In fact, many famous numerical methods, like *bisection* or *Newton's method* for finding zeros and/or extrema, are obviously adaptive, and do allow for extreme speed-ups compared to all non-adaptive methods. We discuss examples below, but first give some final definitions.

If we are interested in the best that can be done based on a specific class of measurements or reconstructions, respectively, then it is natural to set the other class to "arbitrary". Here, and in IBC in general, the focus is on studying the minimal amount of information needed to solve a problem. That is, we often take $\Psi^{\mathrm{arb}}$ and consider the minimal errors

$$e_n^{(\mathrm{non})}(S, \Lambda) := e_n^{(\mathrm{non})}(S, \Lambda, \Psi^{\mathrm{arb}}). \tag{9.11}$$

As in Section 2, for a more refined analysis, we could also discuss the minimal errors for each given measurement map $N_n$, i.e., only take the infimum over all reconstructions $\Phi \in \Psi_n^{\mathrm{arb}}$ in (9.10), which leads to the *radius of information* of the measurement map $N_n$. For a discussion of this quantity, we refer to Traub and Woźniakowski (1980) and Novak and Woźniakowski (2008). In analogy to Section 2, it turns out that, for a given measurement map $N_n$, the worst case error is minimized up to a factor of two if we choose for $\Phi(y)$ any element from the set $\{S(f)\colon f \in F, N_n(f) = y\}$ of possible solutions that 'fit' to the information $y$ (in which case we talk about an *interpolatory algorithm*) and that the infimum is attained if we choose $\Phi(y)$ as a Chebychev center of this set if such a center exists (in which case we talk about a *central algorithm*).

The inverse of the quantity (9.10) is the **information complexity**

$$n(\varepsilon, S, \Lambda) \;:=\; \min\big\{n \in \mathbb{N}\colon\; e_n(S, \Lambda) \leq \varepsilon\big\}.$$

That is, every algorithm that achieves an error of at most $\varepsilon > 0$ for all $f \in F$ needs at least $n(\varepsilon, S, \Lambda)$ measurements from $\Lambda$, and it may therefore be considered a strict lower bound for the *total cost* of any such algorithm (that also has to handle this information), see Remark 9.8. Note that this concept does not deal with implementation issues.

We add that it is typical in IBC to consider the **normalized error criterion** which means that one studies the corresponding complexity

$$n^{\mathrm{nor}}(\varepsilon, S, \Lambda) := \min\{n \in \mathbb{N}\colon\ e_n(S,\Lambda) \le \varepsilon \cdot e_0(S,\Lambda)\}.$$

This is the minimal number of measurements needed to reduce the **initial error** $e_0(S,\Lambda)$ by a factor $\varepsilon \in (0,1)$. The initial error is the best that can be achieved without collecting information at all; it depends only on $S$ and $\Lambda$ and hence on the a priori knowledge. Often, the problems are *normalized*, i.e., $e_0(S,\Lambda) = 1$, and there is no difference. However, especially for high-dimensional applications as in Section 8 where the dependence on another parameter is analyzed, it is important to distinguish between achievements of an algorithm and effects of *bad normalization* of the problems.

We now discuss two examples.

**Example 9.1 (Finding a zero).** Consider the computation of the zero (aka root) of functions from the model class of monotone functions

$$F_{\mathrm{mon}} := \{f \in C([0,1])\colon f(1) = -f(0) = 1 \quad \text{and} \quad f(x) < f(y) \quad \text{for } x < y\},$$

which is modeled by the operator $S_0\colon F_{\mathrm{mon}} \to [0,1]$ such that $f(S_0(f)) = 0$, based on function evaluations $\Lambda^{\mathrm{std}}$. (As usual, the metric on $[0,1]$ is the difference in absolute value.) If we only allow non-adaptive methods, then it is basic, and somewhat reminiscent of the (graphical) analysis in Section 2, to show that

$$e_n^{\mathrm{non}}(S_0, \Lambda^{\mathrm{std}}) = \frac{1}{2(n+1)}.$$

In fact, placing $n$ equidistant points at $\mathcal{P}_n^* = \{\frac{k}{n+1}\colon k = 1,\dots,n\}$, there is, for each $f \in F_{\mathrm{mon}}$, a unique $k^* \in \{1,\dots,n\}$ with $f(\frac{k^*-1}{n+1}) < 0$ and $f(\frac{k^*}{n+1}) \ge 0$. The algorithm $A_n(f) := \frac{2k^*-1}{2(n+1)}$ satisfies the upper bound. For the lower bound, it suffices to note that every 'gap' in the point set of length $\delta$ implies a lower bound of $\frac{\delta}{2}$, since there are two indistinguishable functions whose zeros have distance arbitrarily close to $\delta$.
In contrast, an adaptive choice of sampling points allows for **bisection**: We compute the function value at $\frac{1}{2}$ and, depending on whether it is positive or negative, we continue with the midpoint of the left or right half-interval. This way, we can identify an interval of length $2^{-n}$ that contains the zero based on $n$ function values. An algorithm that outputs the midpoint of this interval implies the bound

$$e_n(S_0, \Lambda^{\mathrm{std}}) \le \frac{1}{2^{n+1}}.$$

Proving a fitting lower bound for adaptive algorithms is a little more involved; we refer to Kung (1976), where this is done for a similar model class.

The previous example is a very classical one that illustrates the clear advantage of adaptive numerical methods (in our sense). However, an important difference to the setting of the previous sections is that it is concerned with a *nonlinear solution operator*. Let us therefore also discuss an approximation problem, which is given by a linear solution operator.

**Example 9.2 (Approximation of Hölder functions).** For $0 < \alpha < 1$, we consider uniform approximation in the model class of monotone $\alpha$-Hölder-continuous functions

$$F_{\mathrm{mon}}^{\alpha} := \big\{ f \colon [0,1] \to [0,1] \colon 0 \le f(y) - f(x) \le |x-y|^{\alpha} \quad \text{for } x \le y \big\},$$

which is modeled by the solution operator $\mathrm{id}_\alpha \colon F_{\mathrm{mon}}^{\alpha} \to B([0,1])$ with $\mathrm{id}_\alpha(f) := f$. By a similar (graphical) analysis as in Section 2, we obtain

$$e_n^{\mathrm{non}}(\mathrm{id}_\alpha, \Lambda^{\mathrm{std}}) \asymp n^{-\alpha},$$

and that, again, the optimal point set for sampling is given by $\mathcal{P}_n^* = \{\frac{2k-1}{2k} \colon k = 1, \ldots, n\}$. The lower bound (for arbitrary point sets) is obtained if we exploit the largest gap in the point set, say, of size $\delta \ge \frac{1}{n+1}$. We then consider two functions $f$ and $g$, which are zero to the left and $(\delta/2)^{\alpha}$ to the right of this gap, respectively, but $f$ is constant on the left half of the gap and increases "as fast as possible" on the right half, while for g it is exactly the opposite. This way, $f$ and $g$ cannot be distinguished by the algorithm but have a distance at least $(2(n+1))^{-\alpha}$ in the supremum norm, which proves the lower bound.
In contrast, an adaptive choice of sampling points allows for the rate

$$e_n(\mathrm{id}_\alpha, \Lambda^{\mathrm{std}}) \asymp n^{-1} \log(n),$$

that does not depend on $\alpha$ (although the constants do). This is a result of Korneichuk (1994). Here, we present only the basic idea of the upper bound. In fact, an algorithm for finding an $\varepsilon$-approximation of $f \in F_{\mathrm{mon}}^{\alpha}$ is based on finding the zeros $x_k$ of the functions $g_k := f - k \cdot \varepsilon$, $k = 1, \ldots, \lfloor \varepsilon^{-1} \rfloor$. Then, clearly, the piecewise linear function that interpolates the points $(x_k, k \cdot \varepsilon)$ has an error of at most $\varepsilon$. Now, it is enough to know the points $x_k$ up to an error $\varepsilon^{1/\alpha}$ which, by Example 9.1, requires about $\frac{1}{\alpha} \log(1/\varepsilon)$ function evaluations. Overall, we need $n(\varepsilon, \mathrm{id}_\alpha, \Lambda^{\mathrm{std}}) \asymp \frac{1}{\alpha\varepsilon} \log(1/\varepsilon)$ function evaluations to achieve an error $\varepsilon$, i.e., we have an error of the order $n^{-1} \log(n)$ with $n$ function values.

Example 9.2 shows that one may *gain* almost a factor $n^{-1}$ by using adaption. Although this does not seem as extreme as Example 9.1, we will see in Section 10 that a much larger gain is not possible in the case of sampling recovery in the uniform norm for convex $F$. (Even randomization and allowing arbitrary linear measurements would not help more.)

Still, the theory of adaptive numerical methods, especially for approximation problems, is far from being complete, and it is often not clear at all how much can be gained by the use of adaption. In any case, this heavily depends on the problem $S\colon F \to Y$ and the admissible information, see, for instance, Novak (1995*a*, 1996).

It is our goal to compare the quantities above for different classes of measurements and reconstructions, ideally for large classes of $S$ at once. But let us first present some earlier benchmarks in this *language* before we state some basic results. In fact, taking into account the above terminology, we obtain the (linear) *sampling numbers* from Sections 3 and 5 by

$$g_n(F,Y) = e_n^{\mathrm{non}}(\mathrm{id}, \Lambda^{\mathrm{std}}) \quad \text{and} \quad g_n^{\mathrm{lin}}(F,Y) = e_n^{\mathrm{non}}(\mathrm{id}, \Lambda^{\mathrm{std}}, \Psi^{\mathrm{lin}}) \tag{9.12}$$

for the inclusion $\mathrm{id}\colon F \to Y$. For the same solution map, we also obtain *best-approximation* on linear subspaces, i.e., the *Kolmogorov widths* from (3.16), by

$$d_n(F,Y) = e_n(\mathrm{id}, \Lambda^{\mathrm{arb}}, \Psi^{\mathrm{lin}}) \tag{9.13}$$

(The corresponding optimal measurements are basis coefficients for the best approximation within the *optimal* subspace) and even the *manifold widths* from (5.22) by

$$\delta_n(F,Y) = e_n^{\mathrm{non}}(\mathrm{id}, \Lambda^{\mathrm{con}}, \Psi^{\mathrm{con}}).$$

Similar to $d_n$, we could define *best $m$-term approximation* by choosing $\Psi$ the class of functions with image $V_N^{(m)}$, see (5.2), and the entropy numbers from Section 5.2 can be written as $\varepsilon_n(F,Y) = e_n(\mathrm{id}, \Lambda^{\mathrm{bin}}, \Psi^{\mathrm{arb}})$, where $\Lambda^{\mathrm{bin}} := \{\lambda\colon F \to \{0,1\}\}$ are all *binary measurements*. We also stress that $e_n(S, \Lambda^{\mathrm{arb}}, \Psi^{\mathrm{arb}})$ and $e_n(S, \Lambda^{\mathrm{arb}}, \Psi^{\mathrm{con}})$ are no reasonable benchmarks, since they equal zero for all $n \geq 1$ whenever $S(F)$ is separable. We leave the details to the reader.

In the previous sections, we compared the sampling numbers with numbers like $d_n$, $\delta_n$ and $\varepsilon_n$, as they came in quite naturally from the proofs and hence this was the strongest way of writing down the error bounds. From the point of view of IBC, it is more natural to compare different ways of obtaining the information (different classes of measurements, adaptive or non-adaptive, random or deterministic) with each other. In Section 10, we will therefore present a comparison of the classes $\Lambda^{\mathrm{std}}$ and $\Lambda^{\mathrm{lin}}$ (Section 10.1), and discuss the power of adaption and randomization for $\Lambda^{\mathrm{lin}}$ (Section 10.2) and $\Lambda^{\mathrm{con}}$ (Section 10.3).

Before we introduce *randomized algorithms* in this general setting, we present some fundamental results in the field. First, we clearly have for $\Lambda_1 \subset \Lambda_2$ and $\Psi_1 \subset \Psi_2$ that

$$e_n(S, \Lambda_2, \Psi_2) \leq e_n(S, \Lambda_1, \Psi_1) \tag{9.14}$$

for every $S$ and $n \in \mathbb{N}$, since a larger class of allowed algorithms leads to a smaller error. The same holds for $e_n^{\mathrm{non}}$. The question is whether there is indeed a difference for different classes, and between $e_n$ and $e_n^{\mathrm{non}}$.

In many cases of interest, one has a linear solution map $\widetilde{S}\colon X \to Y$, like $\widetilde{S} = \mathrm{id}$, for Banach spaces $X$ and $Y$, and one considers the problem $S := \widetilde{S}|_F$ where $F$ is often assumed to be convex and symmetric, such as the unit ball $F := B_X$. In this case, and if only linear measurements, like function values, are available, then adaption as introduced above is basically useless. This is a special case of a classical result in information-based complexity, that was given (together with historical remarks) by Traub and Woźniakowski (1980, Section 2.7), see also (Traub *et al.* 1988, Section 4.5) where the necessity of the assumptions is discussed.

**Proposition 9.3.** Let $X, Y$ be normed spaces, $\widetilde{S}\colon X \to Y$ be linear, $F \subset X$ be convex and symmetric, $\Lambda \subset \mathbb{R}^F$ be a class of linear functionals, and $S := \widetilde{S}|_F$. Then,

$$e_n(S, \Lambda) \le e_n^{\mathrm{non}}(S, \Lambda) \le 2 \cdot e_n(S, \Lambda).$$

See Proposition 9.22 for a related result for non-symmetric sets. Moreover, we will see in Section 10.3 that a result like Proposition 9.3 is not true for continuous measurements $\Lambda^{\mathrm{con}}$. We stress that the word 'adaption' is used with different meanings in the literature, and we refer to (Novak and Woźniakowski 2008, Remark 4.6) for a discussion of the last result and several references. In particular, the applicability to the solution of boundary value problems for elliptic partial differential equations is discussed by Werschulz (1991).

Problems as considered in Proposition 9.3 are often called *linear problems*. The result shows that the (worst-case) error of adaptive methods can be smaller by a factor of at most two for linear problems. This also explains to some extent why we did not treat adaption earlier: The solution map was the identity, which is linear, and the class $\Lambda = \Lambda^{\mathrm{std}}$ was a class of linear measurements, so that adaption would not have helped much for all convex and symmetric model classes $F$. In many cases, even the factor two can be removed, and the inequality becomes equality. This is the case, e.g., if $S\colon F \to Y$ with $Y \in \{\mathbb{R}, B(D)\}$ or if $X$ is a Hilbert space, see (Traub *et al.* 1988, Section 4.5.2).

We do not prove Proposition 9.3 here, but note that the basic ideas have been outlined already in Section 2: Under the given assumptions, one can show that the zero information is the least informative; see (2.9). Hence, we may consider the (non-adaptive) measurement map that is obtained from (9.9) for zero info. (That is, we use zero there instead of $N_{n-1}(f)$ for finding the $n$-th measurement.) It can be shown that a corresponding non-adaptive algorithm that uses these fixed measurements has a worst-case error at most twice as large as the worst-case error of the adaptive algorithm we started with.

Another important general result in the same spirit is that nonlinear reconstruction can often be avoided (at a price) for non-adaptive measurements. As discussed above, linear methods may have several advantages. This is actually a special case of a quite deep result, which also holds for individual measurement maps. We refer to Mathé (1990) or Creutzig and Wojtaszczyk (2004) and references therein.

**Proposition 9.4.** Let $\widetilde{S} \in \mathcal{L}(X,Y)$, $F \subset X$ be convex and symmetric, $S := \widetilde{S}|_F$, and let $\Lambda \subset \Lambda^{\text{lin}}$ consist only of bounded functionals. Then,

$$e_n^{\text{non}}(S,\Lambda,\Psi^{\text{lin}}) \;\leq\; (1+\sqrt{n})\cdot e_n^{\text{non}}(S,\Lambda,\Psi^{\text{arb}}).$$

Moreover, if either

- $Y \in \{\mathbb{R}, L_\infty(\mu), B(D)\}$ <u>or</u>
- $X$ is a Hilbert space <u>or</u>
- $S \in \mathcal{L}(X, C(D))$ is compact,

then,

$$e_n^{\text{non}}(S,\Lambda,\Psi^{\text{lin}}) \;=\; e_n^{\text{non}}(S,\Lambda,\Psi^{\text{arb}}).$$

Propositions 9.3 and 9.4 may be seen in two ways: First, it shows that for some problems that are already solved by nonlinear and/or adaptive methods, there must exist a 'good' non-adaptive linear algorithm. This can be a motivation to search for such simpler and equally good algorithms for particular problems. On the other hand, for applications where it is clear from practice that nonlinear and/or adaptive algorithms give better results, this shows that the model classes and measurement classes mentioned in Propositions 9.3 and 9.4 are just not suitable to describe this superiority and hence that these are not good model assumptions for the corresponding application. We believe that it is of interest to find, for given (classes of) algorithms, even very specific ones, the classes $F$ that indicate particular features and bugs that are seen in the application.

Similarly to above, we can also restrict to continuous reconstructions in case of continuous measurements, at least for compact $F$ and up to a factor of 2. The following is from (Krieg, Novak and Ullrich 2026*a*, Lemma 9).

**Proposition 9.5.** For any continuous mapping $S\colon F \to Y$ from a compact metric space $F$ to a normed space $Y$, any $\Lambda \subset \Lambda^{\text{con}}$ and any $n \in \mathbb{N}_0$, it holds that

$$e_n^{\text{non}}(S,\Lambda,\Psi^{\text{con}}) \;\leq\; 2\cdot e_n^{\text{non}}(S,\Lambda,\Psi^{\text{arb}}).$$

We do not know whether the compactness (and the factor 2) are needed in general. By means of Propositions 9.12 and 9.21, we will see that the 2 is not needed for $S \in \mathcal{L}(H,K)$ for Hilbert spaces $H$ and $K$, and $F = B_H$ (which is in general not compact). We do not know such a bound for unit balls of Banach spaces.

Hence, we saw some sufficient conditions under which not much is lost if we restrict to linear or continuous reconstruction mappings. Similar questions can be asked for other desirable properties of the reconstruction map (especially if the above conditions are not satisfied). For instance, (positive) homogeneity of the reconstruction map is discussed by DeVore, Petrova and Wojtaszczyk (2017, Theorem 3.3) and Krieg and Kritzer (2024, Proposition 1) and measurability of the reconstruction map is discussed by Mathé (1990, Theorem 11), Novak and Ritter (1989) and Novak and Woźniakowski (2008, Section 4.3.3).

We now introduce **randomized algorithms**, as in Section 7, that allow for a random (and adaptive) choice of measurements. Formally, a randomized algorithm is a family of (adaptive) deterministic algorithms $A_n = (A_n^\omega)_{\omega\in\Omega}$ which is indexed by a probability space $(\Omega, \mathcal{A}, \mathbb{P})$. For technical reasons, we assume that the mapping $(f, \omega) \mapsto \|S(f) - A_n^\omega(f)\|_Y$ is $(\mathcal{B}_X \otimes \mathcal{A}, \mathcal{B}_Y)$-measurable, where $\mathcal{B}_X$ and $\mathcal{B}_Y$ denote the Borel $\sigma$-algebras on $X$ and $Y$, respectively. Then, formally, the desirable statement that every deterministic algorithm can be considered as randomized algorithm is not correct since we do not assume that a deterministic algorithm has to be measurable. We refer to (Novak and Woźniakowski 2008, Section 4.3.3) and (Krieg *et al.* 2025*a*, Section 6) for a discussion of this technicality.

We define the **$n$-th minimal worst-case expected error** of a randomized algorithm $A_n = (A_n^\omega)_{\omega\in\Omega}$ for approximating $S\colon F \to Y$ by

$$e_n^{\mathrm{ran}}(S, \Lambda) \ := \ \inf_{A_n} \sup_{f\in F} \mathbb{E}\, \|S(f) - A_n(f)\|_Y, \tag{9.15}$$

where the infimum is over all (adaptive) randomized algorithms $A_n = (A_n^\omega)_{\omega\in\Omega}$ that use measurements from $\Lambda$ with the above measurability assumption, and all probability spaces $(\Omega, \mathcal{A}, \mathbb{P})$.

Especially when $Y$ is a Hilbert space, it may be more common to consider the (larger) *root-mean-square error* $e_n^{\mathrm{rms}}$, as we also did in Section 7. (This often allows for additional techniques.) In the present context, however, we aim for a general comparison between randomized and deterministic methods, see Section 10.2, particularly, we want to bound from above the maximal gain of randomization. Hence, our results will be stronger when formulated with the expected error, which is smaller than the root-mean-square error. For the same reason, we do not consider additional restrictions of the reconstructions, or non-adaptive algorithms. For a discussion and further references regarding other modifications of the randomized setting (like changing the error criterion to a small error with a prescribed probability or allowing a random and adaptive cardinality $n$), we refer to Novak and Woźniakowski (2008, Section 4.3).

We close with a few remarks before discussing the related concept of s-numbers.

**Remark 9.6 (Other measurements).** The above classification of measurements and reconstructions is quite rough, and it is of great interest to study more specific classes. Some work on the case where the algorithms are Lipschitz continuous, which is related to stability and robustness against noise, has been done, e.g., by Plaskota (1996), Cohen *et al.* (2022) and Petrova and Wojtaszczyk (2023). Another interesting class of nonlinear measurements are absolute values of (linear) measurements which play in important role in imaging and *phase retrieval*. We refer to Cahill, Casazza and Daubechies (2016), Candès, Eldar, Strohmer and Voroninski (2015) and Fannjiang and Strohmer (2020) for more on this, and to Plaskota, Siedlecki and Woźniakowski (2020) for absolute value information in the context of IBC. However, for many other important classes of measurements, like values of particular averages or transforms (Radon, X-ray), the knowledge is still quite limited. A slight generalization of the results on sampling discussed in earlier sections can be found, for example, in (Sonnleitner and Ullrich 2023) and (Bartel and Dũng 2024). Also note that the assumption of continuity for linear measurements, as in $\Lambda^{\mathrm{lin}}$, is often not needed, and that this is conjectured to hold in general; see (Hinrichs, Novak and Woźniakowski 2013).

**Remark 9.7 (Optimization).** Besides comparing the power of different classes of algorithms for the same problem, it natural to compare the difficulty of different problems. We have already used this approach in Section 8, since lower bounds for numerical integration are stronger, but sometimes easier to obtain than for approximation directly.
There is another interesting relation in the context of this survey which shows that sampling numbers $g_n(F, B(D))$ are closely related to minimal errors for (global) optimization. In fact, if we consider the solution map $S_{\min}(f) := \min_{x\in D} f(x)$, defined on a convex and symmetric subset $F$ of the continuous functions and compact $D$, then we have the chain of inequalities

$$g_n(F, B(D)) \;\leq\; 2\, e_n(S_{\min}, \Lambda^{\mathrm{std}}) \;\leq\; 2\, g_n(F, B(D)).$$

A similar result holds for the problem of computing $S_{\arg}(f) := \arg\min_{x\in D} f(x)$ with the residual error criterion. See Wasilkowski (1984), Novak (1988), and Novak and Woźniakowski (2010, Section 18.4) for this and similar results. So, up to constant factors, all our considerations for sampling numbers in the uniform norm also imply results for the problem of global optimization based on function evaluations whenever $F$ is convex and symmetric.

**Remark 9.8 (Cost).** The *total cost* of solving a problem is often not proportional to the *information cost* $n$. The information complexity can be considered as a lower bound for the total cost (if this cost is also considered in a worst case setting), since any information has to be processed by the algorithm and therefore leads to additional (arithmetic) cost, see e.g. (Novak and Woźniakowski 2008, Section 4.1.2) for a discussion. Whether there exists a practical algorithm whose total cost is also bounded from above by the information complexity in some way (and whether this algorithm can be found with reasonable computational effort), depends very much on the precise problem. Some candidates for such problems are those problems where linear algorithms are optimal, in particular the integration problem, where an optimal algorithm is often just a linear quadrature rule, but also the problems mentioned in Proposition 9.4. Note that also different *cost models* can be analyzed, in which every measurement has a prescribed cost. In the context of integration over infinite-dimensional domains, this was discussed, for example, by Kuo, Sloan,

Wasilkowski and Woźniakowski (2010), Niu, Hickernell, Müller-Gronbach and Ritter (2011) and Plaskota and Wasilkowski (2011).

**Remark 9.9 (Adaption).** In our setting, only the measurement maps can be adaptive and non-adaptive. A reconstruction map, as we understand it, is always adaptive in the sense that it chooses its output (and even the set of possible outputs) depending on the available information. In other fields of numerical analysis, however, the meaning of the word 'adaptive' may be different, and the adaptivity of reconstructions plays an important role. Often, some information is used for a reconstruction which is then *tested* and, if one is not satisfied, *refined* until some *error* is reached. This can be seen as an adaptive reconstruction. However, in our view, all the refinement etc. is done during the acquisition of information, and only the *final* reconstruction is considered by us as the reconstruction map. This simplification is possible because we care primarily about the *information cost*. When the *arithmetic cost* of the algorithms is considered, too, then a more detailed analysis may be required. We refer to Plaskota (1996) and Bonito, Canuto, Nochetto and Veeser (2024) for comprehensive treatments.

### *9.2. s-numbers and widths*

Many *benchmarks* from previous sections, such as the *Kolmogorov widths* from (3.16), have an obvious geometric interpretation. For example, $d_0(F,Y)$ is the radius of $F$ (in $Y$), and $d_1(F,Y)$ is the radius of the smallest cylinder that contains $F$, etc. Taking inspiration from geometry, there are many other natural notions for the *width* of the set $F$, like the (inner or outer) radius of intersections with or projections onto finite (co-)dimensional subspaces. Such widths have been studied in various contexts, also for infinite-dimensional spaces as discussed here, and the subject has a long tradition in approximation theory. We refer to early treatments of Tikhomirov (1960) and Ismagilov (1974), as well as the foundational books of Pinkus (1985), Tikhomirov (1986), and Lorentz, von Golitschek and Makovoz (1996), and the more recent and specialized *Acta Numerica* articles by DeVore (1998), Cohen and DeVore (2015) and DeVore, Hanin and Petrova (2021).

A somewhat competing concept are *s-numbers* of operators which play an important role in operator theory and geometry in Banach spaces. Roughly speaking, they are a generalization of singular numbers (aka singular values) of matrices to linear mappings between Banach spaces (which also explains the name). They are hence a very natural concept that leads to useful tools for our purposes. Noting that singular numbers of a matrix are certain widths of the ellipsoid obtained as the image of the unit ball, the relation to the geometric concepts above is apparent.

**Singular numbers** of operators, such as matrices, on Hilbert spaces have become fundamental tools in (applied) mathematics since their introduction by Schmidt (1907) in his significant work on integral equations. For compact $S \in \mathcal{L}(H,K)$ between complex Hilbert spaces $H, K$, the *singular numbers* are defined by $\sigma_k(S) := \sqrt{\lambda_k(SS^*)}$, where the *eigenvalues* $\lambda_k(T)$ of a linear operator $T \in \mathcal{L}(X, X)$ are characterized by $Te_k = \lambda_k(T) \cdot e_k$ for some $e_k \in X \setminus \{0\}$, and ordered decreasingly. Applications, especially in computational mathematics, of these numbers and the corresponding *singular value decomposition* (traditionally referred to as *Schmidt representation*) are so manifold that we do not even try to list them here.

Building upon this theory and the *widths* from the early works mentioned above, an axiomatic theory of s-numbers was developed by Pietsch (1974, 1978, 1987, 2007) in the 1970s. The basic idea was to identify the properties that uniquely specify the function that maps an operator between Hilbert spaces to the sequence of its singular numbers, and then use the same definition for operators between Banach spaces. This has led to the following elegant definition.

In what follows, let $X, Y, Z$ and $W$ be real Banach spaces. The (closed) unit ball of $X$ is denoted by $B_X$, and the identity map on $X$ is denoted by $\mathrm{id}_X$. The dimension of a subspace $M \subset X$ is denoted by $\dim(M)$, and by $\mathrm{codim}(M) := \dim(X/M)$ we denote its codimension. The class of all bounded linear operators between Banach spaces is denoted by $\mathcal{L}$, and by $\mathcal{L}(X, Y)$ we denote those operators $S$ from $X$ to $Y$, equipped with the operator norm $\|S\| := \sup_{f \in B_X} \|S(f)\|_Y$. The rank of an operator $S \in \mathcal{L}(X, Y)$ is defined by $\mathrm{rank}(S) := \dim(S(X))$. In the context of this section, it is customary to denote elements of $X$ by $x$ etc., and to omit the specification of the norm of an element, as it is clear from the corresponding space. In particular, for $x \in X$ and $S \in \mathcal{L}(X, Y)$, we write $\|x\| = \|x\|_X$ and $\|Sx\| = \|Sx\|_Y$.

A map $S \mapsto (s_n(S))_{n \in \mathbb{N}_0}$ assigning to every operator $S \in \mathcal{L}$ a non-negative scalar sequence $(s_n(S))_{n \in \mathbb{N}_0}$ is called an **s-number sequence** if, for all $n \in \mathbb{N}_0$, the following conditions are satisfied

(S1) $\|S\| = s_0(S) \geq s_1(S) \geq \cdots \geq 0 \qquad$ for all $\quad S \in \mathcal{L}$,

(S2) $s_n(S + T) \;\leq\; s_n(S) + \|T\| \qquad$ for all $\quad S, T \in \mathcal{L}(X, Y)$,

(S3) $s_n(BSA) \;\leq\; \|B\|\, s_n(S)\, \|A\| \qquad$ where $\; W \xrightarrow{A} X \xrightarrow{S} Y \xrightarrow{B} Z$,

(S4) $s_n(S) = 0 \quad$ whenever $\mathrm{rank}(S) \leq n$,

(S5) $s_n(\mathrm{id}_{\ell_2^{n+1}}) = 1$.

We call $s_n(S)$ the $n$-th **s-number** of the operator $S$. To indicate the underlying Banach spaces, we may write $s_n(S\colon X \to Y)$.

**Remark 9.10 (Index shift).** Here, we use an index shift compared to the original definition of s-numbers, and work with $s_n$ for $n \in \mathbb{N}_0$ instead of $n \in \mathbb{N}$. The reason is that we want to use both s-numbers and widths, and even a common generalization, see Section 9.3. Since for widths, it seems more common to start counting from $n = 0$, while s-numbers typically start from $n = 1$, a clash of notation with some of the relevant literature cannot be avoided. The choice to start counting with $n = 0$ comes from our understanding of $n$ being the number of measurements (and hence the dimension or codimension of the used subspace).

For operators between Hilbert spaces, there is exactly one s-number sequence: the singular numbers, see Proposition 9.12. However, when considering operators between Banach spaces, these axioms do not lead to a unique concept. Some of these s-numbers will be more like *width*, while others correspond to minimal errors, or *quantify compactness* of an operator. As written by Pietsch (2007, 6.2.2.1),

"we have a large variety of s-numbers that make our life more interesting."

Again, there are too many applications to be summarized here. However, let us add that the research done in the classification of *operator ideals* (Pietsch 2007, 6.3), or *eigenvalue distributions* of operators, see (König 1986) or (Pietsch 2007, 6.4), has quite some similarities to the comparison of different classes of algorithms, as we are going to discuss below. For an even broader picture, we highly recommend the monumental monograph of Pietsch (2007) who gives a detailed account of the developments of the last century in –what he calls– modern Banach space theory.

There are some especially important examples of s-numbers:

- **approximation numbers**:
$$a_n(S) := \inf \left\{ \|S - L\| \colon L \in \mathcal{L}(X,Y),\ \operatorname{rank}(L) \le n \right\}$$
- **Bernstein numbers**:
$$b_n(S) := \sup \left\{ \inf_{x \in M \setminus \{0\}} \frac{\|Sx\|}{\|x\|} \colon M \subset X,\ \dim(M) = n+1 \right\}$$
- **Gelfand numbers**:
$$c_n(S) := \inf \left\{ \|S|_M\| \colon M \subset X \text{ closed},\ \operatorname{codim}(M) \le n \right\}$$
- **Kolmogorov numbers**:
$$d_n(S) := \inf \left\{ \sup_{x \in B_X} \inf_{y \in V} \|Sx - y\| \colon V \subset Y,\ \dim(V) \le n \right\}$$
- **Manifold numbers**:
$$\delta_n(S) := \inf_{\substack{N_n \in C(B_X, \mathbb{R}^n) \\ \Phi \in C(\mathbb{R}^n, Y)}} \sup_{x \in B_X} \left\| Sx - \Phi\left(N_n(x)\right) \right\|$$
- **Hilbert numbers**:
$$h_n(S) \ : = \sup \left\{ \frac{a_n(BSA)}{\|B\| \|A\|} \colon A \in \mathcal{L}(\ell_2, X),\ B \in \mathcal{L}(Y, \ell_2),\ A, B \neq 0 \right\}.$$

The numbers $a_n$, $b_n$, $c_n$ and $d_n$ were already known in the 1960s, or much earlier, see (Pietsch 2007, 6.2.3). The invention of the Hilbert numbers by Bauhardt (1977), however, to some extent demonstrates the advantage of an axiomatic approach: They are the smallest s-numbers, see Proposition 9.12, and it turned out that they are a useful tool. The manifold numbers did not appear under this name in the literature, but we think it is natural taking into account their 'width counterpart'. It is discussed in Section 9.3 that they are indeed s-numbers.

Note that these definitions lead to the corresponding *widths* considered earlier by

$$s_n(\mathrm{id}\colon X \to Y) = s_n(B_X, Y)$$

with $s \in \{d, \delta, h\}$, see (3.16), (3.54), (5.22) and (5.20). For more general $F \neq B_X$, we refer to Section 9.3 or Krieg *et al.* (2025*a*) for possible definitions of $b_n$ and $c_n$, and more references.

**Remark 9.11 (Widths and numbers).** Some of the numbers and widths have been defined, even under the same name (like Gelfand), in different areas. The differences between the concepts of s-numbers and widths can sometimes be rather subtle. For instance, if $S = \mathrm{id}$, then continuity of the measurements (and closedness of the subspaces) is postulated sometimes in $X$ and sometimes in $Y$. See e.g. the short (historical) account of Pietsch (2007, 6.2.6). In some cases, however, the differences are quite significant as has been discussed, for example, by Edmunds and Lang (2013) and Heinrich (1989). Hence, some caution is recommended when working with these concepts.

The definition of s-numbers has been changed several times over the years. Other properties of s-numbers that were used in earlier definitions are

(S2') $s_{m+n}(S + T) \leq s_m(S) + s_n(T)$ for all $S, T \in \mathcal{L}(X, Y)$,

(S3') $s_{m+n}(ST) \leq s_m(S) \cdot s_n(T)$ where $X \xrightarrow{T} Y \xrightarrow{S} Z$,

(S5') $s_n(\mathrm{id}_X) = 1$ for all $X$ with $\dim(X) > n$.

These properties are stronger than their counterparts from above, and s-numbers satisfying them are called *additive*, *multiplicative* and *strict*, respectively. Note that $a_n$, $c_n$ and $d_n$ satisfy all three properties, while $h_n$ is only additive, and $b_n$ is only strict. A detailed discussion (and proofs) of these properties can be found in Section 5.1 by Lang and Edmunds (2011), or in (Pietsch 2007, 2008).

For us, the most important aspect of the axiomatic approach is that it leads to a characterization of the smallest and largest s-number (with certain properties). See also Remark 9.18 for a characterization of $b_n$, $c_n$ and $d_n$ in this way. We will see that these s-numbers can be bounded by each other, showing that different s-numbers cannot be too different. Eventually, by placing the different minimal errors into this scale, we will also obtain inequalities between them, see Section 9.3.

We begin with the following result which is essentially due to Pietsch (1974), see also (Pietsch 1987, 2.3.4 & 2.6.3 & 2.11.9).

**Proposition 9.12.** For every s-number sequence $(s_n)$, $S \in \mathcal{L}$ and $n \in \mathbb{N}_0$, we have

$$h_n(S) \;\leq\; s_n(S) \;\leq\; a_n(S).$$

Equalities hold if $S \in \mathcal{L}(H, K)$ for Hilbert spaces $H$ and $K$.

*Proof.* We refer to (Pietsch 1987, 2.11.9) for the detailed proof that $s_n(S) = a_n(S)$ for any s-number sequence $(s_n)$, and any $S \in \mathcal{L}(H, K)$ for Hilbert spaces $H$ and $K$. Here, we only present the proof of this equality for compact $S$ in the hope to provide some insight into the matter. For this, consider the singular value decomposition $S = \sum_{k=0}^{\infty} \sigma_k \langle \cdot, u_k \rangle v_k$ for some orthonormal $u_k \in H$ and $v_k \in K$, where the $(\sigma_k)$ are the (non-increasing) singular numbers of $S$. For $n < \operatorname{rank}(S)$, we denote $S_n := \sum_{k=0}^{n-1} \sigma_k \langle \cdot, u_k \rangle v_k$. From (S2) and (S5), we see that $s_n(S) \leq s_n(S_n) + \|S - S_n\| = \sigma_n$. Now let $P := \sum_{k=0}^{n} \langle \cdot, v_k \rangle v_k$ such that $S_{n+1} = PS$ and $\|P\| = 1$, and note that the inverse of $S_{n+1}$, considered as mapping from the finite-dimensional Hilbert spaces $U := \operatorname{span}\{u_k\}$ to $V := \operatorname{span}\{v_k\}$, is given by $S_{n+1}^{-1} := \sum_{k=0}^{n} \sigma_k^{-1} \langle \cdot, v_k \rangle u_k$. Using isometries $A\colon \ell_2^{n+1} \to V$ and $B\colon V \to \ell_2^{n+1}$, we obtain from (S3) and (S4) that

$$\begin{aligned} 1 = s_n(\mathrm{id}_{\ell_2^{n+1}}) = s_n(B S_{n+1} S_{n+1}^{-1} A) &\leq \|B\| \cdot s_n(PS) \cdot \|S_{n+1}^{-1}\| \cdot \|A\| \\ &\leq \|P\| \cdot s_n(S) \cdot \sigma_n^{-1}, \end{aligned}$$

proving the equality $s_n(S) = \sigma_n(S)$ for any compact $S \in \mathcal{L}(H, K)$.

For the inequalities, we use that the uniqueness allows us to use any $(s_n)$ in the definition of $h_n$. Hence, by (S3), we obtain

$$h_n(S) \;=\; \sup\left\{ \frac{s_n(BSA)}{\|B\|\|A\|} \colon A \in \mathcal{L}(\ell_2, X),\; B \in \mathcal{L}(Y, \ell_2),\; A, B \neq 0 \right\} \;\leq\; s_n(S)$$

for any s-number sequence $(s_n)$. In addition, by (S2) and (S5), we obtain for any $L$ with $\operatorname{rank}(L) \leq n$ that $s_n(S) \;\leq\; s_n(L) + \|S - L\| \;=\; \|S - L\|$. By taking the infimum over all such $L$, we also see that $s_n(S) \leq a_n(S)$.

□

Now that we know the smallest and largest s-numbers, it is natural to ask about their maximal differences (where we mean 'difference' not in the mathematical way; the ratio is usually more interesting than the difference). Several inequalities and relations between the individual s-numbers have already been collected by Pietsch (1974), see also (Pietsch 1987, 2.10) and (Pietsch 2007, 6.2.3.14). A particularly important bound relates the approximation numbers to the Gelfand and Kolmogorov numbers, which is very much related to the treatment of linear algorithms in Proposition 9.4 (the following corresponds to the special case $F = B_X$ and $\Lambda = X'$), but also to the bounds on linear sampling numbers from Section 3, see Remark 9.20.

**Proposition 9.13.** For all $S \in \mathcal{L}$ and $n \in \mathbb{N}_0$, we have

$$a_n(S) \;\leq\; \left(1 + \sqrt{n}\right) \cdot \min\left\{c_n(S),\, d_n(S)\right\}.$$

The proof is based on important results from Banach space theory on the existence of certain projections. Moreover, the factor $\sqrt{n}$ is needed, as can be seen by the embeddings $W_1^s \hookrightarrow L_2$ and $W_2^s \hookrightarrow L_\infty$ considered in Appendix A.2.

*Proof.* Fix $\varepsilon > 0$. First, let $M \subset X$ be closed with $\operatorname{codim}(M) \leq n$ such that $\|S|_M\| \leq (1+\varepsilon) \cdot c_n(S)$. By a theorem of Garling and Gordon (1971), see also (Pietsch 1987, 1.7.17), there is a projection $P \in \mathcal{L}(X,X)$ (i.e., $P^2 = P$) such that $\{x\colon Px = 0\} = M$ and $\|P\| \leq (1+\varepsilon)\sqrt{n}$. Defining $L := SP$ with $\operatorname{rank}(L) \leq \operatorname{rank}(P) = \operatorname{codim}(M) \leq n$, we obtain

$$\begin{aligned} a_n(S) &\leq \|S - L\| = \|S(\mathrm{id}_X - P)\| \leq \|S|_M\| \cdot \|\mathrm{id}_X - P\| \\ &\leq (1+\varepsilon) \cdot c_n(S) \cdot \left(1 + (1+\varepsilon)\sqrt{n}\right). \end{aligned}$$

With $\varepsilon \to 0$, this implies the first bound.
For the second, let $V \subset Y$ with $\dim(V) \leq n$ and $\sup_{x \in B_X} \inf_{y \in V} \|Sx - y\| \leq (1+\varepsilon) \cdot d_n(S)$. By the theorem of Kadets and Snobar (1971), see also (Pietsch 1987, 1.5.5), there is a projection $P \in \mathcal{L}(X,X)$ such that $P(X) = V$ and $\|P\| \leq \sqrt{n}$. Defining $L := PS$ with $\operatorname{rank}(L) \leq \operatorname{rank}(P) \leq n$, we obtain

$$\begin{aligned} a_n(S) \leq \|S - L\| &= \|(\mathrm{id}_Y - P)S\| = \sup_{x \in B_X} \|(\mathrm{id}_Y - P)Sx\| \\ &= \sup_{x \in B_X} \inf_{y \in V} \|(\mathrm{id}_Y - P)(Sx - y)\| \leq \|\mathrm{id}_Y - P\| \cdot \sup_{x \in B_X} \inf_{y \in V} \|Sx - y\| \\ &\leq \left(1 + \sqrt{n}\right) \cdot (1+\varepsilon) \cdot d_n(S). \end{aligned}$$

With $\varepsilon \to 0$, this implies the second bound. □

We turn to bounds between all other s-numbers. This has been achieved essentially by Pietsch (1980) using a bound with the *Weyl numbers* as an intermediate step. In the following form, it is from (Pietsch 1978, 11.12.3). See also the streamlined presentation of Ullrich (2024), the extension to non-symmetric sets of Krieg *et al.* (2025*a*), as well as the history and open problems in Remark 9.19.

**Theorem 9.14.** For all $S \in \mathcal{L}$ and $n \in \mathbb{N}_0$, we have

$$\max\left\{c_n(S),\, d_n(S)\right\} \;\leq\; (n+1) \left(\prod_{k=0}^{n} h_k(S)\right)^{1/(n+1)}.$$

If $S \in \mathcal{L}(X,H)$ for a Hilbert space $H$, we can replace the factor $(n+1)$ with $\sqrt{n+1}$.

We do not give a proof of Theorem 9.14 here, because we have already seen it: The proof (for $c_n$) is one-to-one as the proof of Theorem 5.15 with function evaluations replaced by general linear functionals. We omit the details.

Bounds for individual $n$ cannot be obtained from Theorem 9.14, but combining it with the inequality $n^\alpha \leq e^\alpha\,(n!)^{\alpha/n}$ for $\alpha \geq 0$, we obtain a more handy form.

**Corollary 9.15.** For all $S \in \mathcal{L}$, $\alpha > 0$ and $n \in \mathbb{N}$, we have

$$\max\left\{c_n(S),\, d_n(S)\right\} \;\leq\; e^{\alpha}\, n^{-\alpha+1} \cdot \sup_{k\leq n}\, (k+1)^{\alpha}\, h_k(S),$$

and

$$a_n(S) \;\leq\; 2\, e^{\alpha}\, n^{-\alpha+3/2} \cdot \sup_{k\leq n}\, (k+1)^{\alpha}\, h_k(S).$$

Since $c_n(\mathrm{id}\colon \ell_1 \to \ell_\infty) \asymp d_n(\mathrm{id}\colon \ell_1 \to \ell_\infty) \asymp 1$ and $h_n(\mathrm{id}\colon \ell_1 \to \ell_\infty) \asymp n^{-1}$, see (Pietsch 2007, 6.2.3.14) and (Heinrich and Linde 1984), the first result is best possible up to constants. It is open whether the second is also.

Recalling Proposition 9.12, this gives a bound between arbitrary pairs of s-numbers. However, we may be interested in a finer comparison between particular numbers. For this, let us also mention the following bound.

**Proposition 9.16.** For all $S \in \mathcal{L}$ and $n \in \mathbb{N}_0$, we have

$$b_n(S) \;\leq\; \min\left\{c_n(S),\, \delta_n(S)\right\}.$$

Moreover, if $S \in \mathcal{L}(X,Y)$ with separable $Y$ (or $Y = L_\infty(\mu)$ with $\sigma$-finite $\mu$), then

$$\delta_n(S) \;\leq\; c_Y \cdot d_n(S)$$

for some $c_Y$ that only depends on $Y$. ($c_Y = 1$ for strictly convex $Y$)

*Proof.* For $b_n(S) \leq \delta_n(S)$, we use the Borsuk-Ulam theorem of Borsuk (1933), which implies that for every continuous (measurement) map $N_n\colon B_X \to \mathbb{R}^n$ and any $V_n \subset X$ of dimension $n+1$, there is some $x \in V_n$ with $\|x\| = 1$ and $N_n(x) = N_n(-x)$, such that any algorithm based on $N_n$ makes an error of at least $\|Sx\|$. This implies the bound.
The same argument applies to linear measurements, and can therefore be used for the inequality $b_n(S) \leq c_n(S)$. Alternatively, we can observe that $b_n$ and $c_n$ are the smallest and largest 'strict injective s-numbers', see Remark 9.18.

For the second part, we first note that, for strictly convex $Y$ and any finite dimensional subspace $V_n \subset Y$, there is a unique continuous operator $P\colon Y \to V_n$ that maps every $y \in Y$ to its best approximation from $V_n$. This proves $\delta_n(S) \leq d_n(S)$. The continuity of best-approximation remains true for any $Y$ that admits an equivalent strictly convex norm, which is implied by the conditions, see (Day 1955). □

**Remark 9.17 (Entropy numbers).** The *entropy numbers* of an operator $S \in \mathcal{L}(X,Y)$ are defined by

$$\varepsilon_n(S) \;:=\; \inf_{y_1,\dots,y_{2^n}\in Y}\; \sup_{x\in B_X}\; \min_{i=1,\dots,2^n} \|Sx - y_i\|,$$

which is in correspondence to the definition (5.16) by $\varepsilon_n(S) = \varepsilon_n(B_X, Y)$. The entropy numbers do only satisfy the conditions (S1)–(S3) and (S5), and are therefore not an s-number sequence, see (Pietsch 2007, 6.2.4). (They are also additive and multiplicative.) In fact, they satisfy $\varepsilon_n(\mathrm{id}_{\mathbb{R}}) = 2^{-n}$, and are therefore never zero. More generally, we have $\varepsilon_n(\mathrm{id}_X) = 2^{-n/m}$ for every Banach space $X$ with $\dim(X) = m$, see (Pietsch 1978, 12.1.13). The theory of entropy numbers is also well-established. We refer to Carl and Stephani (1990) for a comprehensive treatment. By means of *Carl's inequality*,

$$\varepsilon_n(S) \le C_s\, n^{-s} \cdot \sup_{k\le n} (k+1)^s \cdot \min\left\{c_k(S), d_k(S)\right\}$$

with a constant $C_s > 0$ only depending on $s > 0$ see (Carl 1981*b*) and (5.17), they are among the 'smaller' s-numbers. However, it is also known that $h_n(S) \le 2\,\varepsilon_n(S)$, see (Pietsch 1978, 12.3.1). Moreover, there is the bound

$$\max\left\{c_n(S),\ d_n(S)\right\} \le (n+1) \cdot \varepsilon_n(S)$$

from (Pietsch 1978, 12.3.2). (Recall our index shift and Remark 5.7.)

**Remark 9.18 (Further properties of s-numbers).** There are further characteristics of s-numbers that can also explain the importance of $b_n$, $c_n$ and $d_n$ without referring to their importance in approximation theory. For this, we call an s-number sequence $(s_n)$ *injective* if $s_n(S) = s_n(JS)$ for any metric injection $J \in \mathcal{L}(Y, Y_0)$, and *surjective* if $s_n(S) = s_n(SQ)$ for any metric surjection $Q \in \mathcal{L}(X_0, X)$. Here, $X_0, Y_0$ are arbitrary Banach spaces. Roughly speaking, injective s-numbers do not change if we embed $Y$ in a larger space, while surjective s-numbers do not change if we do the same with $X$. We refer to Pietsch (1974, 1987) for a detailed treatment and further properties. In turns out, that the Gelfand numbers $c_n$ are the largest injective s-numbers, and the Kolmogorov numbers $d_n$ are the largest surjective s-numbers. The Bernstein numbers are the smallest 'strict injective s-numbers'. So, many of the important s-numbers (and widths) come out of the axiomatic approach.

**Remark 9.19 (Bounds for individual $n$).** Several of the above bounds between different s-numbers bound the s-numbers from one scale in terms of the progression of s-numbers $s_0(S), \ldots, s_n(S)$ from another scale. It remains a long-standing open problem whether similar bounds exist in terms of only $s_n(S)$. In fact, it is even conjectured that $a_{cn}(S) \le c\, n\, h_n(S)$, for some constant $c \ge 1$, see (Pietsch 2009, Problem 5) or (Pietsch 1978, 11.12.4).

In the context of widths (and this survey), a particularly interesting weaker form of this conjecture is that $d_{cn}(S) \le c\, n\, b_n(S)$, see for example (Pinkus 1985, p. 24). This is an old conjecture of Mityagin and Henkin (1963), who also proved $d_n(S) \le n^2 b_n(S)$ and, apparently, the bound has not been improved since then. However, in the weaker form as in Theorem 9.14, this problem has already been solved by Pietsch (1980), see the final remarks there. It seems that this has been overlooked in the numerical analysis community until recently, see the exposition of Krieg *et al.* (2025*a*) where also the implications for numerical analysis are analyzed. We discuss that more detailed in Section 10.2.

**Remark 9.20 (Sampling projections).** The bound on the approximation numbers by Kolmogorov numbers in Proposition 9.13 is based on the existence of projections onto every $m$-dimensional space of norm at most $\sqrt{m}$. Let us stress that the results of Section 3 can partly be phrased in the same language. In fact, and as discussed in Section 3.1, the least squares algorithm from (3.3) is a projection. However, those projections are of a special form as they can only use function evaluations instead of arbitrary functionals. Hence, for $X \subset B(D)$, they are of the form $P = \sum_{i=1}^{n} \delta_{x_i} \cdot v_i$ for some $v_i \in X$ and $x_i \in D$, and may be called *sampling projections*. We refer to Krieg *et al.* (2026*c*) for more on this.

### 9.3. *Minimal errors and s-numbers*

We now connect the theory of s-numbers with the minimal errors discussed in Section 9.1. For model classes $F = B_X$, i.e., unit balls of a normed space $X$, and deterministic algorithms, the connection between the above concepts has been treated in detail by Mathé (1990), and we refer to his article for details and further references. For general convex (non-symmetric) sets, we refer to Krieg *et al.* (2025*a*) for a detailed exposition.

Here, we consider only linear solution operators $S \in \mathcal{L}(X,Y)$, because this is the setting of s-numbers. In the context of this work, this is not a major restriction, as we have been mainly concerned with the inclusion $S = \mathrm{id}$ with $\mathrm{id}(f) = f$.

Recall that the *n-th minimal worst-case errors* for approximating $S \in \mathcal{L}(X,Y)$ over $F$ based on measurements $\Lambda$ and reconstructions $\Psi$ is defined by

$$e_n^{(\mathrm{non})}(S|_F, \Lambda, \Psi) \;:=\; \inf_{\substack{N_n \\ \Phi\in\Psi}} \; \sup_{f\in F} \left\| S(f) - \Phi \circ N_n(f) \right\|_Y,$$

where the infimum is over all (non-)adaptive measurement maps $N_n$ from $\Lambda$.

Some connections to s-numbers are obvious and have already been outlined around (9.13). In fact, we have

$$a_n(S) = e_n^{\mathrm{non}}(S|_{B_X}, \Lambda^{\mathrm{lin}}, \Psi^{\mathrm{lin}}) \quad \text{and} \quad d_n(S) = e_n(S|_{B_X}, \Lambda^{\mathrm{arb}}, \Psi^{\mathrm{lin}}).$$

Moreover, it was shown by Mathé (1990, Theorem 10), who attributed it to Tikhomirov (1971), that

$$d_n(S) = e_n^{\mathrm{non}}(S|_{B_X}, \Lambda^{\mathrm{con}}, \Psi^{\mathrm{lin}})$$

for compact $S$. So, at least for compact $S$, the minimal errors with linear reconstruction and measurements from $\Lambda^{\mathrm{lin}}$, $\Lambda^{\mathrm{con}}$ or $\Lambda^{\mathrm{arb}}$ are themselves s-numbers. For arbitrary information, we also have

$$e_1^{\mathrm{non}}(S|_F, \Lambda^{\mathrm{arb}}, \Psi^{\mathrm{con}}) = e_1^{\mathrm{non}}(S|_F, \Lambda^{\mathrm{arb}}, \Psi^{\mathrm{arb}}) = 0,$$

whenever $S(F)$ is separable, see (Traub and Woźniakowski 1980, Chapter 3) or (Mathé 1990, Remark 2). They therefore do not represent s-numbers.

The remaining cases are settled by the following result of Mathé (1990).

**Proposition 9.21.** If $\Lambda \in \{\Lambda^{\rm lin}, \Lambda^{\rm con}\}$ and $\Psi \in \{\Psi^{\rm lin}, \Psi^{\rm con}, \Psi^{\rm arb}\}$, then

$$s_n(S) := e_n^{\rm non}(S|_{B_X}, \Lambda, \Psi)$$

as well as

$$\widetilde{s}_n(S) := \inf_{N_n \in \Lambda^n} \sup_{\substack{f,g \in B_X \\ N_n(f)=N_n(g)}} \frac{1}{2} \left\| S(f) - S(g) \right\|_Y$$

are s-number sequences.

Note that the second part of the proposition uses (half) the *diameter of information*; another important quantity from IBC, see (Traub and Woźniakowski 1980). It is worth noting again that, for operators $S \in \mathcal{L}(H, K)$ between Hilbert spaces, all these s-numbers coincide.

Using the last proposition and the bounds on s-numbers from Section 9.2, we could directly deduce several bounds for minimal errors. However, those are only for non-adaptive algorithms and unit balls, which is the setting for s-numbers.

For more general sets, we consider a generalization of s-numbers. In particular, the **Gelfand numbers of $S \in \mathcal{L}(X, Y)$ on $F \subset X$** are defined by

$$\begin{aligned} c_n(S, F) &:= \inf_{\lambda_1,\dots,\lambda_n \in X'} \sup_{\substack{f,g \in F: \\ \lambda_k(f)=\lambda_k(g)}} \frac{1}{2} \left\| S(f) - S(g) \right\| \\ &= \inf_{\substack{M \subset X \text{ closed} \\ \operatorname{codim}(M) \le n}} \sup_{\substack{f,g \in F: \\ f-g \in M}} \frac{1}{2} \left\| S(f) - S(g) \right\|. \end{aligned}$$

For example, $c_0(S, F) = \frac{1}{2} \operatorname{diam}\big(S(F)\big)$. This definition unifies *Gelfand numbers* ($F = B_X$) and **Gelfand widths** ($S = \mathrm{id}_X$). They characterize the minimal worst-case error of (non-adaptive, deterministic) algorithms based on arbitrary linear measurements up to a factor of two. This also implies Proposition 9.3.

**Theorem 9.22.** Let $X$ and $Y$ be Banach spaces and $S \in \mathcal{L}(X, Y)$. For every $F \subset X$ and $n \in \mathbb{N}_0$, we have

$$c_n(S, F) \le e_n^{\rm non}(S|_F, \Lambda^{\rm lin}) \le 2\, c_n(S, F).$$

If $F$ is convex and symmetric, then

$$c_n(S, F) \le e_n(S|_F, \Lambda^{\rm lin}) \le e_n^{\rm non}(S|_F, \Lambda^{\rm lin}) \le 2\, c_n(S, F).$$

We refer to Novak and Woźniakowski (2008, Sections 4.1 & 4.2) for a proof.

For adaptive or randomized algorithms on non-symmetric sets, or for more general measurements, we can only state weaker bounds. (Except possibly for operators between Hilbert spaces, where many of the numbers coincide.) In this case, we employ Bernstein numbers for a lower bound. The **Bernstein numbers of $S$ on $F$**, for general bounded sets $F \subset X$, are defined by

$$b_n(S,F) := \sup_{\substack{\dim(V)=n+1 \\ S \text{ injective on } V}} \sup \Big\{ r>0\colon\, g+U \subset F \;\text{ for some }\; g\in F \;\text{ and a ball } U := \{f\in V\colon \|S(f)\| \le r\} \Big\}.$$

They describe the radius of the largest $(n+1)$-dimensional ball in $F$ w.r.t. the norm $\|f\|_S := \|S(f)\|$. This generalization of *Bernstein numbers* ($F = B_X$) and **Bernstein widths** ($S = \mathrm{id}_X$), has been given by Krieg *et al.* (2025*a*).

The following result demonstrates the importance of the Bernstein numbers for proving lower bounds for minimal errors in different settings. The individual bounds have been proven by Novak (1995*a*), DeVore *et al.* (1989) and Kunsch (2016).

**Theorem 9.23.** Let $X$ and $Y$ be Banach spaces and $S \in \mathcal{L}(X,Y)$. For every $F \subset X$ and $n \in \mathbb{N}_0$, we have

$$b_n(S,F) \;\le\; \min\Big\{ e_n(S|_F, \Lambda^{\mathrm{lin}}),\; e_n^{\mathrm{non}}(S|_F, \Lambda^{\mathrm{con}}, \Psi^{\mathrm{con}}) \Big\}.$$

If $F$ is convex, then also

$$e_n^{\mathrm{ran}}(S|_F, \Lambda^{\mathrm{lin}}) \;\ge\; \frac{1}{30} \cdot b_{2n-1}(S,F).$$

In Section 10.2, we discuss the sharpness of these bounds. But note already here that we do not know a larger lower bound for non-adaptive randomized algorithms; it is an interesting open problem to find a better general lower bound for this class of algorithms (if possible). Moreover, we will see in Section 10.3, that a similar bound is not possible for adaptive continuous algorithms. There can be even an exponential *speed-up* of $e_n(S|_F, \Lambda^{\mathrm{con}})$ compared to all s-numbers. Unfortunately, we do not know any lower bound for $e_n(S, \Lambda^{\mathrm{con}})$ with $n \ge 2$, even if $S$ is as simple as the identity (matrix) on $\mathbb{R}^3$ (with any norm), see Section 10.3 (For $n = 1$ a lower bound follows from the Borsuk-Ulam theorem since then there can be no discontinuous adaption step and the information mapping $N_n$ is continuous). We believe that this is another interesting subject for further research.

Before we present a proof of the bounds on the deterministic errors, let us comment further on the randomized errors. In the case that $H$ and $G$ are Hilbert spaces and $S \in \mathcal{L}(H,G)$, a result of Novak (1992) implies that

$$e_n^{\mathrm{ran}}(S|_{B_H}, \Lambda^{\mathrm{lin}}) \;\ge\; \frac{1}{2} \cdot b_{2n-1}(S, B_H). \tag{9.16}$$

Note that in this case, the Bernstein numbers coincide with all the other s-numbers, and hence with many other (deterministic) minimal errors by Proposition 9.21. For a proof (of precisely this version of Novak's result), we refer to Krieg *et al.* (2025*a*, Lemma 4.3). There it was also shown that for $S \in \mathcal{L}(X, Y)$ with Banach spaces $X$, $Y$ and for general convex sets $F \subset X$, we have

$$e_n^{\mathrm{ran}}(S|_F, \Lambda^{\mathrm{lin}}) \geq \frac{1}{2} \cdot h_{2n-1}(S, F) \tag{9.17}$$

in terms of the *Hilbert numbers* that we already encountered in Section 5.3 for $S = \mathrm{id}$ being the identity. We omit the precise definition for general classes $F$ and refer to Krieg *et al.* (2025*a*) for the details. The better bound (apart from the constant) in terms of Bernstein numbers was obtained by Kunsch (2016) by a more involved method. See also Kunsch (2017) for more details and the explicit constant.

We conclude with the proof of the first part of Theorem 9.23.

*Proof.* Let $0 < \beta < b_n(S, F)$. Then there is an $(n+1)$-dimensional subspace $V \subset X$ such that $S$ is injective on $V$ as well as $m \in F$ such that, with $U := \{h \in V\colon \|Sh\|_Y \leq \beta\}$, we have $m + U \subset F$.
Now, for the bound $b_n(S, F) \leq e_n(S|_F, \Lambda^{\mathrm{lin}})$, let $A_n = \varphi_n \circ N_n$ be an algorithm based on the (adaptive) measurement map $N_n\colon F \longrightarrow \mathbb{R}^n$ with measurements from $\Lambda^{\mathrm{lin}}$. We fix the non-adaptive and linear mapping $N_n^* \in (\Lambda^{\mathrm{lin}})^n$ that is obtained if we apply $N_n$ to the 'midpoint' $m \in F$. The mapping $N_n^*$ cannot be injective on $U$ and there exists a point $\tilde{f}$ on the boundary of $U$ with $N^*(m + \tilde{f}) = N^*(m) = N^*(m - \tilde{f})$, hence $N(m + \tilde{f}) = N(m - \tilde{f})$. Then $A_n$ cannot distinguish between the two inputs and we obtain $e_n(S|_F, \Lambda^{\mathrm{lin}}) \geq \beta$.
For $b_n(S, F) \leq e_n^{\mathrm{non}}(S|_F, \Lambda^{\mathrm{con}})$, we use the Borsuk-Ulam theorem of Borsuk (1933) (as in the proof of Proposition 9.16), which implies that for every continuous (measurement) map $N_n\colon F \longrightarrow \mathbb{R}^n$, there is some $f \in V$ with $\|Sf\| = \beta$ and $N_n(m + f) = N_n(m - f)$. Since $m + f$ and $m - f$ cannot be distinguished by the information $N_n$, every nonadaptive continuous algorithm has the error at least $\beta$. This proves the first two bounds, since $\beta < b_n(S, F)$ was arbitrary.

□

## 10. Comparison of different types of measurements

In the previous sections, especially in Sections 3 and 5, we discussed many results on the minimal worst-case error of algorithms that process function values at no more than $n$ fixed sampling points, especially for the problems of $L_2$-approximation and uniform approximation. In this section, we shall see that the minimal errors do not change too much if we admit algorithms that can use the additional features:

- adaption,
- randomness,
- general linear measurements, or even
- general (Lipschitz) continuous measurements.

We shall see that any of these features individually can change the polynomial rate of convergence at most by a constant value, at least for convex model classes $F$. The precise value of this constant is known in some cases and remains an open problem in other cases. This means that, at least from the viewpoint of the information complexity (i.e., the number of required measurements), the above features can only make a big difference for examples where non-adaptive sampling algorithms have a very small or no rate of convergence.

This is even true if we combine the first three features. In contrast, for algorithms that use adaptively chosen continuous measurements, it turns out that the gain in the rate of convergence can be infinite. This last finding can be regarded only as a theoretical result, since the class of all continuous measurements is much too large to be of practical relevance. On the other hand, the result can be a motivation for future research to find and study other classes of measurements, which are more realistic, but which still show a major gain compared to function evaluations and linear measurements, see also Remark 9.6.

Most of the results in this section are corollaries of results from all over the article, which we now place into the IBC framework presented in Section 9. The structure of this section is as follows:

- In Section 10.1, we discuss the problems of $L_2$-approximation and uniform approximation. We show that algorithms using arbitrary linear measurements, possibly chosen adaptively and/or randomly, as well as non-adaptive continuous measurements, can only achieve a limited gain over linear sampling algorithms.

- In Section 10.2, we study by how much the minimal errors change with the use of adaption and/or randomization in the case of linear measurements for the approximation of general linear operators.

- In Section 10.3, we discuss the (rather unrealistic) power of algorithms that can use adaptive continuous measurements. The key

result is that an $m$-dimensional vector can, in theory, be approximated with arbitrary precision by using only $\mathcal{O}(\log m)$ adaptive and piecewise linear measurements.

We will conclude each section with some open questions.

To keep the presentation simple and focused on the main insights (and open problems), we will express the results only in terms of the polynomial rate of convergence of the minimal errors, at least in the first two sections. For this, we define the *polynomial rate of convergence* of a sequence $(z_n)_{n\in\mathbb{N}} \subset [0,\infty)$ by

$$\operatorname{rate}(z_n) := \sup\left\{\alpha > 0 \,\middle|\, \exists C \geq 0 : \forall n \in \mathbb{N} : z_n \leq Cn^{-\alpha}\right\}. \tag{10.1}$$

This quantity does not take into account logarithmic factors (which are often subject of important open problems), but it allows, e.g., to avoid the geometric means needed in some of the comparisons.

Moreover, we mostly consider convex and symmetric model classes $F$, or even the unit ball $B_X$ of $X$. In this case, we know several bounds that are sharp (up to logarithmic factors). The case of more general classes $F$, and in particular non-symmetric $F$, is widely open, see Remark 10.3, and we refer to the detailed treatment of Krieg *et al.* (2025*a*) for more on this.

**Notation.** We discuss the minimal error of different types of algorithms as defined in Section 9.1. However, since we will jump a lot between the different settings, we will expand and modify the notation a little in the hope that many of the formulas will be readable without having to scroll back to Section 9.1 too often. We write:

$$\begin{aligned} e_n^{\text{det-ada}} &:= e_n, & e_n^{\text{det-non}} &:= e_n^{\text{non}}, \\ e_n^{\text{ran-ada}} &:= e_n^{\text{ran}}, & e_n^{\text{ran-non}} &:= e_n^{\text{ran-non}} \end{aligned}$$

for the minimal error of deterministic adaptive, deterministic non-adaptive, randomized adaptive, and randomized non-adaptive algorithms, see (9.11) and (9.15). The last of these four has formally not been defined in Section 9.1, but it is of course defined as in (9.15) but with the infimum over randomized algorithms whose realizations are non-adaptive deterministic algorithms.

### *10.1. Function values vs other measurements*

We discuss the problems of approximating functions from a model class $F \subset \mathbb{C}^D$ in $Y = L_2$ or in $Y = B(D)$. For a simpler presentation, in this section, we impose the following **general assumptions**:

- $D$ is a compact Hausdorff space,
- $\mu$ is a finite Borel measure,
- $F = B_X$ for a Banach space $X$ for which $X \hookrightarrow C(D)$ is compact.

The interested reader will be able to extract more general comparison results from the referenced ingredients. The above conditions imply that $\Lambda^{\rm std} \subset \Lambda^{\rm lin} = X'$. We often consider the special case that $X$ is a Hilbert space. In this case, the above conditions imply that $X$ is, in fact, a reproducing kernel Hilbert space and that the embedding $X \hookrightarrow L_2$, where $L_2 = L_2(D,\mu)$, has a finite trace, so that the results from Section 3.3.1 apply.

In the notation of Section 9.1, we only deal with problems $\mathrm{id}\colon F \to Y$, and since the model class $F$ and error norm $Y$ bear more important information than just writing "id", we write

$$e_n^*(F,Y,\Lambda) \;:=\; e_n^*(\mathrm{id}\colon F \to Y,\Lambda)$$

with $* \in \{\text{det-non}, \text{det-ada}, \text{ran-non}, \text{ran-ada}\}$. We compare the power of non-adaptive and deterministic function evaluations with the power of more general measurements. Namely, we compare

$$e_n^{\text{det-non}}(B_X,Y,\Lambda^{\rm std}) \qquad \text{vs.} \qquad e_n^*(B_X,Y,\Lambda)$$

for $\Lambda = \Lambda^{\rm lin} = X'$ or $\Lambda = \Lambda^{\rm con}$. The main comparison results are (**10.A**)–(**10.F**).

For the use of the results from the previous sections, recall that the minimal worst-case error for non-adaptive deterministic function evaluations equals the sampling numbers, see (1.1) and (3.5), i.e.,

$$e_n^{\text{det-non}}(F,Y,\Lambda^{\rm std}) \;=\; g_n(F,Y) \;\le\; g_n^{\rm lin}(F,Y).$$

We begin with $L_2$-approximation ($Y = L_2$), as considered in Section 3.3.
For this, recall Corollary 3.31 which gives that

$$g_{cm}^{\rm lin}(B_X,L_2) \;\le\; \frac{1}{\sqrt{m}} \sum_{k>m} \frac{d_k(B_X,L_2)}{\sqrt{k}}.$$

In terms of the rate, this shows that

$$\mathrm{rate}\Big(d_n(B_X,L_2)\Big) > \frac{1}{2} \quad\Longrightarrow\quad \mathrm{rate}\Big(g_n^{\rm lin}(B_X,L_2)\Big) = \mathrm{rate}\Big(d_n(B_X,L_2)\Big). \qquad (10.2)$$

If $X = H$ is a Hilbert space, then $d_n(B_H,L_2) = d_n(H \hookrightarrow L_2)$, i.e., the Kolmogorov widths (3.54) equal the Kolmogorov numbers from Section 9.2 which equal all other s-numbers, see Proposition 9.12. Moreover, the finite trace of the embedding $H \hookrightarrow L_2$ (as implied by our assumptions) gives that $\mathrm{rate}\big(d_n(B_H,L_2)\big) \ge 1/2$. Taking into account Theorems 9.22 and 9.23, we obtain for Hilbert spaces $X = H$ that

$$\mathrm{rate}\Big(e_n^{\text{det-non}}(B_H,L_2,\Lambda^{\rm std})\Big) = \mathrm{rate}\Big(e_n^*(B_H,L_2,\Lambda^{\rm lin})\Big) \ge \frac{1}{2} \qquad (\mathbf{10.A})$$

for all $* \in \{\text{det-ada}, \text{ran-non}, \text{ran-ada}\}$. That is, non-adaptive deterministic sampling algorithms are optimal (in the sense of order) among arbitrary adaptive randomized algorithms that may use general linear information. (For rate $= 1/2$, (**10.A**) does not follow from (10.2), but directly from (3.43).) Also the rate of non-adaptive continuous measurements is equal to the rates from (**10.A**), since the corresponding minimal errors also constitute s-numbers, see Proposition 9.21. We note that our conditions on $X$, $D$ and $\mu$ are important for the result (10.2); there is a reproducing kernel Hilbert space $H$ whose unit ball satisfies

$$\operatorname{rate}\Big(e_n^{\text{det-non}}(B_H, L_2, \Lambda^{\text{std}})\Big) = 0 \quad \text{and} \quad \operatorname{rate}\Big(e_n(B_H, L_2, \Lambda^{\text{lin}})\Big) = \frac{1}{2}. \tag{10.3}$$

This was proven by Hinrichs *et al.* (2008), see also the improved bound of Hinrichs *et al.* (2022, Theorem 6). In contrast to our situation, the embedding $H \hookrightarrow L_2$ in this example has infinite trace.

For Banach spaces $X$ that are compactly embedded into $C(D)$, we still have the equality (10.2) of the rate of the sampling numbers and the Kolmogorov widths. In this more general setting, however, the Kolmogorov widths are not equivalent to the minimal errors $e_n^*(B_X, L_2, \Lambda^{\text{lin}})$. Instead, for general $X$, it follows that

$$\operatorname{rate}\Big(e_n^{\text{det-non}}(B_X, L_2, \Lambda^{\text{std}})\Big) \geq \operatorname{rate}\Big(e_n^*(B_X, L_2, \Lambda^{\text{lin}})\Big) - \frac{1}{2} \tag{\textbf{10.B}}$$

for all $* \in \{\text{det-ada}, \text{ran-non}, \text{ran-ada}\}$, whenever $\operatorname{rate}\big(e_n^{\text{ran-ada}}(B_X, L_2, \Lambda^{\text{lin}})\big) \geq 1$. That is, non-adaptive deterministic sampling algorithms are optimal within the class of all adaptive and randomized algorithms up to a loss of at most $1/2$ in the rate of convergence. Inequality (**10.B**) follows by the composition of the following results:

- $\operatorname{rate}\Big(e_n^{\text{det-non}}(B_X, L_2, \Lambda^{\text{std}})\Big) \geq \operatorname{rate}\Big(g_n^{\text{lin}}(B_X, L_2)\Big)$,
- $\operatorname{rate}\Big(g_n^{\text{lin}}(B_X, L_2)\Big) = \operatorname{rate}\Big(d_n(B_X, L_2)\Big)$, see (10.2),
- $\operatorname{rate}\Big(d_n(B_X, L_2)\Big) \geq \operatorname{rate}\Big(h_n(B_X, L_2)\Big) - 1/2$, see Theorem 9.14,
- $\operatorname{rate}\Big(h_n(B_X, L_2)\Big) \geq \operatorname{rate}\Big(b_n(B_X, L_2)\Big)$, see Proposition 9.12,
- $\operatorname{rate}\Big(b_n(B_X, L_2)\Big) \geq \operatorname{rate}\Big(e_n^*(B_X, L_2, \Lambda^{\text{lin}})\Big)$, see Theorem 9.23.

This bound is sharp, as can be seen from the unit balls $F_1^s$ of the Sobolev spaces $X = W_1^s([0,1])$ with $s \in \mathbb{N} \setminus \{1\}$, for which

$$\operatorname{rate}\Big(e_n^{\text{det-non}}(F_1^s, L_2, \Lambda^{\text{std}})\Big) = s - 1/2 \quad \text{and} \quad \operatorname{rate}\Big(e_n^{\text{det-non}}(F_1^s, L_2, \Lambda^{\text{lin}})\Big) = s,$$

see Appendix A.2 and Theorem 9.22. Since the second quantity lies between $\operatorname{rate}\big(e_n^{\text{det-non}}(F_1^s, L_2, \Lambda^{\text{std}})\big)$ and $\operatorname{rate}\big(e_n^{\text{ran-ada}}(F_1^s, L_2, \Lambda^{\text{lin}})\big)$, and the difference between these two is at most $1/2$ by (**10.B**), the difference between the latter two quantities must also be equal to $1/2$. In particular, we obtain

$$\operatorname{rate}\Big(e_n^{\text{det-non}}(F_1^s, L_2, \Lambda^{\text{lin}})\Big) = \operatorname{rate}\Big(e_n^{\text{ran-ada}}(F_1^s, L_2, \Lambda^{\text{lin}})\Big),$$

because there is no room for further gaps.

Via Proposition 9.12 and 9.21, a similar reasoning applies to the minimal errors of non-adaptive continuous algorithms. That is, since $\mathrm{rate}\big(h_n(B_X, L_2)\big) \geq \mathrm{rate}\big(e_n^{\text{det-non}}(B_X, L_2, \Lambda^{\text{con}})\big)$, we obtain from the remaining arguments above that

$$\mathrm{rate}\Big(e_n^{\text{det-non}}(B_X, L_2, \Lambda^{\text{std}})\Big) \;\geq\; \mathrm{rate}\Big(e_n^{\text{det-non}}(B_X, L_2, \Lambda^{\text{con}})\Big) - \frac{1}{2}, \qquad (\mathbf{10.C})$$

if the right hand side is at least $1/2$, showing that the gain of general continuous measurements over function values is limited in this case. Again, for the example $F = F_1^s$, we have that $\mathrm{rate}\big(e_n^{\text{det-non}}(F_1^s, L_2, \Lambda^{\text{lin}})\big) = \mathrm{rate}\big(e_n^{\text{det-non}}(F_1^s, L_2, \Lambda^{\text{con}})\big)$.

We now turn to uniform approximation ($Y = B := B(D)$).
In this case, Theorem 3.17 implies again a gap of $1/2$ in the rate. In fact, it implies

$$\mathrm{rate}\Big(g_n^{\text{lin}}(B_X, B)\Big) \;\geq\; \mathrm{rate}\Big(d_n(B_X, B)\Big) - \frac{1}{2}, \qquad (10.4)$$

Since we know that linear and non-adaptive algorithms are optimal (up to a factor of 2) for $Y = B(D)$, see Propositions 9.3 and 9.4, we also have

$$d_n(B_X, B) \leq e_n^{\text{non}}(\mathrm{id}\,|_{B_X}, \Lambda^{\text{lin}}, \Psi^{\text{lin}}) = e_n^{\text{det-non}}(B_X, B, \Lambda^{\text{lin}}) \leq 2\cdot e_n^{\text{det-ada}}(B_X, B, \Lambda^{\text{lin}}). \qquad (10.5)$$

This implies

$$\mathrm{rate}\Big(e_n^{\text{det-non}}(B_X, B, \Lambda^{\text{std}})\Big) \;\geq\; \mathrm{rate}\Big(e_n^{\text{det-ada}}(B_X, B, \Lambda^{\text{lin}})\Big) - \frac{1}{2}. \qquad (\mathbf{10.D})$$

In words, the restriction from arbitrary deterministic (adaptive) algorithms with linear measurements to non-adaptive sampling algorithms can *cost* at most an order of convergence $1/2$. Again, the sharpness is shown by the model classes $F_1^s$, $s > 1$, see Example 3.18 or Appendix A.2, for which $\mathrm{rate}\big(e_n^{\text{det-non}}(F_1^s, B, \Lambda^{\text{std}})\big) = s - 1$ and $\mathrm{rate}\big(e_n^{\text{det-ada}}(F_1^s, B, \Lambda^{\text{lin}})\big) = s - 1/2$.

For a comparison with randomized algorithms, recall Theorem 5.15 which states that

$$\mathrm{rate}\Big(g_n^{\text{lin}}(B_X, B)\Big) \;\geq\; \mathrm{rate}\Big(h_n(B_X, B)\Big) - 1 \qquad (10.6)$$

with the *Hilbert numbers* $h_n$ from (5.20). Hence, by Theorem 9.23, see also (9.17), we obtain

$$\mathrm{rate}\Big(e_n^{\text{det-non}}(B_X, B, \Lambda^{\text{std}})\Big) \;\geq\; \mathrm{rate}\Big(e_n^{\text{ran-ada}}(B_X, B, \Lambda^{\text{lin}})\Big) - 1. \qquad (\mathbf{10.E})$$

So, also for uniform approximation, the gain of randomized adaptive algorithms that use linear measurements over non-adaptive deterministic sampling algorithms is limited. We add that the Hilbert numbers also lower bound other benchmarks, like the minimal errors $e_n^{\text{det-non}}(B_X, B, \Lambda^{\text{con}})$. Hence, (**10.E**) also holds with them in place of $e_n^{\text{ran-ada}}$.

The bound (**10.E**) is also sharp. Indeed, a gain of order one can already oc-cur when we switch from non-adaptive to adaptive algorithms in the randomized setting. This has been proven only recently by Heinrich (2024*a*,*b*) for the related problem of *parametric integration* in the case of standard information. For the approximation problem, as considered in this section, an example is due to Kunsch, Novak and Wnuk (2024) and Kunsch and Wnuk (2026). The latter showed that a gain of order one occurs for approximation in $\ell_\infty = B(\mathbb{N})$ on the unit ball $F$ of a weighted $\ell_1$-space between non-adaptive and adaptive randomized algorithms in the case of linear information (and hence, also between $e_n^{\text{det-non}}(F, B, \Lambda^{\text{std}})$ and $e_n^{\text{ran-ada}}(F, B, \Lambda^{\text{lin}})$). We remark that Kunsch and Wnuk (2026) actually obtain bounds for the embedding $\ell_1^m \hookrightarrow \ell_\infty^m$ on the unweighted finite-dimensional space $\ell_1^m$, which seems a more natural problem. However, we need an infinite-dimensional example in order to talk about gaps in the rate of convergence (otherwise all the rates are infinite). We do not know a more 'natural' infinite-dimensional example where the gap of order one in (**10**.E) occurs. We come back to the adaption problem in the next section.

By the techniques used in this survey, we do not often obtain bounds on the minimal errors of (deterministic) sampling algorithms in terms of the minimal errors of more general algorithms (for the same problem) for fixed $n$. In this respect, we only have the bound

$$e_{2n}^{\text{det-non}}(B_X, B, \Lambda^{\text{std}}) \;\le\; 2 \cdot \sqrt{2n} \cdot e_n^{\text{det-ada}}(B_X, B, \Lambda^{\text{lin}}) \tag{10.7}$$

by Theorem 3.16 and (10.5). The other bounds on $e_n^{\text{det-non}}(B_X, Y, \Lambda^{\text{std}})$ are only in terms of a whole (finite or infinite) progression of minimal errors.

Another case where such a bound for individual $n$ is possible is $L_2$-approximation in reproducing kernel Hilbert spaces which satisfy an additional boundedness assumption. In this case, we obtain that linear sampling algorithms are optimal (up to a constant and oversampling) among all deterministic algorithms that use linear information. That is, we have the following theorem that has been obtained independently by Geng and Wang (2024, Theorem 2.2) and Krieg *et al.* (2025*b*, Theorem 11).

**Theorem 10.1.** There are absolute constants $b, c \in \mathbb{N}$ such that the following holds. Let $\mu$ be a finite measure on a set $D$ and $H \subset L_2(\mu)$ be a reproducing kernel Hilbert space with kernel

$$k_H(x, y) = \sum_{k=1}^{\infty} \sigma_k^2 \, b_k(x) \, \overline{b_k(y)}, \qquad x, y \in D,$$

where $(\sigma_k) \in \ell_2$ is a non-increasing sequence and $\{b_k\}_{k=1}^\infty \subset B(D)$ is an orthonormal system in $L_2$ such that there is a constant $C > 0$ with

$$\Big\| \frac{1}{n} \sum_{k=1}^{n} |b_k|^2 \Big\|_\infty \;\le\; C \tag{10.8}$$

for all $n \in \mathbb{N}$. Under these conditions we have $F := B_H \subset B = B(D)$ and

$$g_{bn}^{\text{lin}}(F, B) \le \sqrt{c \cdot C \cdot \mu(D)} \cdot c_n(F, B).$$

In particular,

$$e_{bn}^{\text{det-non}}(F, B, \Lambda^{\text{std}}) \le \sqrt{c \cdot C \cdot \mu(D)} \cdot e_n^{\text{det-ada}}(F, B, \Lambda^{\text{lin}}). \tag{10.9}$$

Recall from Theorem 9.22, that the $e_n^{\text{det-ada}}$ equal (up to a factor of 2) the Gelfand width $c_n$ in this case. As discussed in Remark 3.33, the bound (10.8) holds for several bases, such as trigonometric or Chebyshev polynomials, as well as (Haar) wavelets, with the corresponding orthogonality measure.

The proof is based on a slight modification of the lifting argument in Lemma 3.34 in reproducing kernel Hilbert spaces, see Krieg *et al.* (2025*b*, Lemma 14), and a lower bound on the Gelfand widths from Osipenko and Parfenov (1995), see also Cobos, Kühn and Sickel (2016). In particular, we obtain that

$$\text{rate}\Big(e_n^{\text{det-non}}(B_H, B, \Lambda^{\text{std}})\Big) = \text{rate}\Big(e_n^{\text{det-ada}}(B_H, B, \Lambda^{\text{lin}})\Big) \tag{\textbf{10.F}}$$

whenever $H$ is a reproducing kernel Hilbert space that satisfies the boundedness assumption of Theorem (10.1). It is an interesting open problem whether these boundedness assumptions are required, or whether (**10.F**) actually holds for more general reproducing kernel Hilbert spaces, similar to (**10.A**) for $L_2$-approximation.

**Remark 10.2 (Exponential convergence).** For many classes $F$ of smooth functions, such as the unit ball of a reproducing kernel Hilbert space with Gaussian kernel of Example 3.26, the $n$-th minimal error has a super-polynomial decay. In this case, by Theorem 9.14 or its generalization to non-symmetric sets of (Krieg *et al.* 2025*a*), also all other s-numbers have super-polynomial decay. For example, we obtain the relation

$$e_n^{\text{ran-ada}}(F, Y, \Lambda^{\text{lin}}) \in \mathcal{O}(\alpha^n) \qquad \Longrightarrow \qquad e_{cn}^{\text{det-non}}(F, Y, \Lambda^{\text{std}}) \in \mathcal{O}(\alpha^n),$$

for all $Y \in \{L_2, B(D)\}$, $F \subset C(D)$ being a compact convex set and $\alpha \in (0, 1)$, with a constant $c = c(\alpha)$. Indeed, in the case of exponential decay, the sums and geometric means in the various upper bounds simplify and additional polynomial factors are 'swallowed' by an oversampling factor $c$. In the case of $Y = L_2$, the exponential decay is discussed in more detail in Krieg *et al.* (2023).

This means that, for all such examples, there is no need for sophisticated algorithms that use randomization, adaption or general linear information, at least from the viewpoint of information complexity. In comparison to deterministic and non-adaptive algorithms that only use function evaluations, at most a factor $c$ can be gained. We do not know if this also holds for other $Y$, like $L_p$ with $p \notin \{2, \infty\}$.

*10.2. Power of adaption and randomization for linear measurements*

We now discuss linear measurements in a bit more detail. This is the most studied setting in approximation theory and IBC, possibly because its close connection to s-numbers and widths, and hence, most is known here. However, we will see that there are still fundamental open questions.

The results presented in this section hold for arbitrary linear operators $S \in \mathcal{L}(X,Y)$ with Banach spaces $X$ and $Y$. We study

$$e_n^*(S, \Lambda^{\mathrm{lin}}) := e_n^*(S : B_X \to Y, \Lambda^{\mathrm{lin}})$$

with $* \in \{\text{det-non}, \text{det-ada}, \text{ran-non}, \text{ran-ada}\}$ and $\Lambda^{\mathrm{lin}} = X'$. Again, we only discuss here the case of model classes $F = B_X$, i.e., the worst-case error is considered over the unit ball of the Banach spaces $X$, and we also focus on the polynomial rate (10.1). The case of other convex and symmetric model classes $F$ (with no essential changes) and the case of convex and non-symmetric model classes $F$ (with essential changes), is discussed in detail by Krieg *et al.* (2025*a*), see also Remark 10.3.

First, by Proposition 9.3, we have for all $S \in \mathcal{L}(X,Y)$ that

$$\mathrm{rate}\left(e_n^{\text{det-ada}}(S, \Lambda^{\mathrm{lin}})\right) = \mathrm{rate}\left(e_n^{\text{det-non}}(S, \Lambda^{\mathrm{lin}})\right).$$

If $X = H$ is a Hilbert space, and $Y = K$ is also a Hilbert space, we even have

$$\mathrm{rate}\left(e_n^{\text{det-non}}(S, \Lambda^{\mathrm{lin}})\right) = \mathrm{rate}\left(e_n^{\text{ran-ada}}(S, \Lambda^{\mathrm{lin}})\right) \quad \text{for} \quad S \in \mathcal{L}(H,K),$$

and so all types of algorithms based on linear information are of the same *power*. This result has been proven first by Novak (1992), see also (Krieg *et al.* 2025*a*, Lemma 4.3).

In more general settings, we can again rely on the fact that the deterministic minimal errors are s-numbers (Proposition 9.21) and that the randomized errors are bounded by the Bernstein numbers (Proposition 9.23), together with the bounds from Section 9.2. Since the rate of the Gelfand numbers $c_n$ differs from all s-numbers by at most 1, see Theorem 9.14, and $c_n$ is essentially $e_n^{\mathrm{non}}$, see Proposition 9.22, we obtain for every $S \in \mathcal{L}(X,Y)$ that

$$\mathrm{rate}\left(e_n^{\text{det-non}}(S, \Lambda^{\mathrm{lin}})\right) \geq \mathrm{rate}\left(e_n^{\text{ran-ada}}(S, \Lambda^{\mathrm{lin}})\right) - 1.$$

If either $X$ or $Y$ is a Hilbert space we even have

$$\mathrm{rate}\left(e_n^{\text{det-non}}(S, \Lambda^{\mathrm{lin}})\right) \geq \mathrm{rate}\left(e_n^{\text{ran-ada}}(S, \Lambda^{\mathrm{lin}})\right) - 1/2.$$

In both cases, equality can occur, as has been discussed in (Kunsch and Wnuk 2026, Corollary 4.2), (Kunsch *et al.* 2024, Remark 3.4), and (Kunsch and Wnuk 2026, Remark 4.3), respectively. In fact, for any fixed $n \in \mathbb{N}$, Kunsch *et al.* (2024) show that one may gain a factor of the order $\sqrt{n/\log n}$ for the embedding $S\colon \ell_1^m \to \ell_2^m$ and suitable (large) $m = m(n)$. Similarly, Kunsch and Wnuk (2026) prove that one may gain a factor of the order $n/\log n$ for the embedding $S\colon \ell_1^m \to \ell_\infty^m$ with appropriate $m = m(n)$. These embeddings can then be used as building blocks to construct a infinite-dimensional examples where a gain of polynomial order one (respectively one half) occurs in the rate of convergence.

This shows that the gain of adaption and randomization, even combined, is limited for the approximation of linear operators. It is natural to ask for the individual power of these *features*. For this, let us define the 'maximal gain function'

$$\mathrm{gain}(*, \square) \;:=\; \sup_{S \in \mathcal{L}} \left( \mathrm{rate}\left( e_n^{\square}(S, \Lambda^{\mathrm{lin}}) \right) - \mathrm{rate}\left( e_n^{*}(S, \Lambda^{\mathrm{lin}}) \right) \right),$$

where $*, \square \in \{\text{det-non}, \text{det-ada}, \text{ran-non}, \text{ran-ada}\}$ and the rate function is defined in (10.1). The function describes the maximal possible gain in the rate of convergence when allowing $\square$-algorithms instead of $*$-algorithms. In particular, the above discussion implies

$$\mathrm{gain}\big(\text{det-non}, \text{det-ada}\big) = 0$$

since adaption does not help in this case, and

$$\mathrm{gain}\big(\text{ran-non}, \text{ran-ada}\big) = \mathrm{gain}\big(\text{det-ada}, \text{ran-ada}\big) = \mathrm{gain}\big(\text{det-non}, \text{ran-ada}\big) = 1.$$

The very last equality is from above and for the others, note that (Kunsch *et al.* 2024, Kunsch and Wnuk 2026) actually discuss the gain of adaption for randomized algorithms, see also (Krieg *et al.* 2025*a*, Section 6) for details.

The only case where we do not know the final answer yet, and which we think is particularly interesting, is the gain of randomization for non-adaptive algorithms. So far, we only know

$$\frac{1}{2} \;\le\; \mathrm{gain}\big(\text{det-non}, \text{ran-non}\big) \;\le\; 1, \tag{10.10}$$

with the upper bound from the general statement above. The lower bound is obtained for certain Sobolev embeddings. The simplest case is from $W_2^k([0,1])$ into $L_\infty([0,1])$, where the optimal rate with deterministic algorithms is $n^{-k+1/2}$ for $k > 1/2$. This is a classical result of Solomjak and Tichomirov (1967). Using non-adaptive randomized algorithms, the upper bound $n^{-k}\log n$ has been obtained by Mathé (1991). We also refer to Byrenheid, Kunsch and Nguyen (2018), Fang and Duan (2008), Heinrich (1992) for further results and extensions.

Hence, it is not known to us whether there is some $S \in \mathcal{L}$ such that, say, $e_{2n}^{\text{ran-non}}(S, \Lambda^{\mathrm{lin}}) \le C\, n^{-\alpha}\, e_n^{\text{det-non}}(S, \Lambda^{\mathrm{lin}})$ for all $n \in \mathbb{N}$ and some $\alpha > \frac{1}{2}$. This appears to be a very interesting open problem, especially since its formulation does not require any concept of adaptivity.

**Remark 10.3 (Non-symmetric sets).** For general convex (non-symmetric) sets $F$, results that are similar to the comparisons above are given by Krieg *et al.* (2025*a*,

Section 6). The difference is that in most cases, an additional factor $\sqrt{n}$ (i.e., an extra rate of $1/2$) appears in the upper bounds on the maximal gain. In particular, the gain of adaption and/or randomization on all convex model classes is bounded above by $3/2$. However, we do not know any example with a gain larger than 1. So, it is still possible that the maximal gain of adaption and randomization is 1 for all convex model sets.

Note that in this more general case, adaption can also help in the deterministic setting. As an example, one may consider $S = \mathrm{id} \in \mathcal{L}(\ell_\infty, \ell_\infty)$, i.e., approximation in $\ell_\infty$, for inputs from

$$F = \{x \in \ell_\infty \mid x_i \geq 0, \sum x_i \leq 1, x_k \geq x_{2k}, x_k \geq x_{2k+1}\}.$$

Then one can prove a lower bound $c(\sqrt{n}\log n)^{-1}$ for non-adaptive algorithms while a simple adaptive algorithm (using "function values", i.e., values of coordinates of $x$) gives the upper bound $(n+3)^{-1}$, see (Novak 1995*b*) for details. This shows a gain of polynomial order $1/2$. On the other hand, also for the gain from non-adaptive to adaptive deterministic algorithms, we only know the upper bound $3/2$, which means that the question regarding the maximal power of adaption on convex (non-symmetric) sets is widely open even in the purely deterministic setting.

### *10.3. Power of adaption for continuous measurements*

Finally, we discuss the situation of adaptive continuous measurements. Somewhat surprisingly, it turned out that this class of measurements is (most likely) too powerful to be realistic. That is, we will show that adaption for continuous measurements allows for an exponential speed-up, even for simple problems such as approximation of a vector. This fact has only recently been observed by Krieg *et al.* (2026*a*), see Remark 10.7 for our first encounter with the problem. The solution relies on a (theoretically) simple algorithm that uses a particular 'coloring' of the model class, together with bisection. Before we turn to details, let us stress that we do not believe that this is a revolutionary approximation algorithm. On the contrary, we rather think that it demonstrates that some care is needed when working with (optimal) algorithms from such a big class. Moreover, for us, it indicates that there are further interesting phenomena to discover in the field of adaptive algorithms with nonlinear measurements, see also Remark 9.6.

The main purpose of this section is to compare

$$e_n^{\text{det-ada}}(S, \Lambda^{\text{con}}) \qquad \text{with} \qquad e_n^{\text{det-non}}(S, \Lambda^{\text{con}}).$$

That is, the $n$-th minimal errors of adaptive and non-adaptive algorithms based on continuous measurements and without restriction on the reconstruction. However, we do not have a final answer for this yet, and our results can be stated best in terms of the (larger) *manifold numbers* of a map $S\colon F \to Y$, that is, by

$$\delta_n(S) := e_n^{\text{non}}(S, \Lambda^{\text{con}}, \Psi^{\text{con}}),$$

where $F$ is a topological space (required for $\Lambda^{\mathrm{con}}$) and $Y$ is a metric space (required for the error). Note that this is again a slight generalization compared to earlier appearances with $S = \mathrm{id}$ (5.22) or with $F = B_X$ in Section 9.2.

The following result is a consequence of the finite-dimensional result (and algorithm) described below.

**Theorem 10.4.** For any mapping $S\colon F \to Y$ from a topological space $F$ to a metric space $Y$ and any $n \geq 2$, it holds that

$$e_n^{\text{det-ada}}(S, \Lambda^{\mathrm{con}}) \;\leq\; \delta_{2^{n-2}}(S).$$

Hence, optimal algorithms based on adaptive continuous measurements *always* offer an exponential speed-up compared to algorithms that use only non-adaptive continuous measurements and continuous reconstruction mappings. Note that this is not a fully satisfying comparison, since different classes of reconstructions are allowed on both sides. Unfortunately, we do not know how to resolve that in full generality, but at least in special cases.

First note that, by Proposition 9.5, we can replace $\delta_N(S)$ by $2e_N^{\text{det-non}}(S, \Lambda^{\mathrm{con}})$ whenever $S$ is a continuous map between a compact metric space $F$ and a normed space $Y$, see (Krieg *et al.* 2026*a*, Theorem 10). This, however, excludes the important case of $F$ being a unit ball of an infinite-dimensional normed space $X$. (At least if $F$ is equipped with the metric induced by $X$; if we consider continuous measurements with respect to a different metric, it might be covered.)

Alternatively, we can consider bounded linear operators $S\colon F \to Y$, where $F$ is assumed to be the unit ball of a Banach space $X$ and $Y$ is a Banach space. In this setting, we can use that the $\delta_n$ and $e_n^{\text{det-non}}(S, \Lambda^{\mathrm{con}})$ are s-numbers, see Proposition 9.21. Hence, we can employ the bounds between s-numbers from the previous sections. Especially, if $S$ maps between Hilbert spaces, so that all its s-numbers coincide (Theorem 9.12), we get the following.

**Corollary 10.5.** Let $H$ and $K$ be Hilbert spaces and $S \in \mathcal{L}(H, K)$. Then,

$$e_n^{\text{det-ada}}(S|_{B_H}, \Lambda^{\mathrm{con}}) \;\leq\; e_{2^{n-2}}^{\text{det-non}}(S|_{B_H}, \Lambda^{\mathrm{con}}).$$

This shows that adaption provides an exponential speed-up for the approximation of bounded linear operators between Hilbert spaces if continuous information is allowed. We do not know whether this remains true in more generality. Let us only mention that by the bounds between s-numbers (for linear $S$ between Banach spaces) we can obtain an exponential speed-up of $e_n^{\text{det-ada}}(S|_{B_X}, \Lambda^{\mathrm{con}})$ compared to $e_n^{\text{det-non}}(S|_{B_X}, \Lambda^{\mathrm{con}})$ for operators between Banach spaces under the assumption that the latter has a rate larger than one. We omit the details.

We now turn to the algorithm that implies the above results. In fact, the key ingredient is an algorithm for the 'elementary' problem of recovering a vector $x \in \mathbb{R}^m$, up to some error $\varepsilon > 0$ in some norm $\|\cdot\|$. Curiously, a worst-case error $\varepsilon$ over the whole space $\mathbb{R}^m$ can be achieved with only $\mathcal{O}(\log m)$ adaptively chosen (Lipschitz) continuous measurements. (We do not need to assume a bounded subset of $\mathbb{R}^m$, and also the specific norm is irrelevant.)

**Theorem 10.6.** Let $m \in \mathbb{N}$ and $\varepsilon > 0$. The algorithm $R_m^\varepsilon$ described below uses at most

$$n(m) := \lceil \log_2(m+1) \rceil + 1$$

adaptive, Lipschitz-continuous measurements and satisfies for all $x \in \mathbb{R}^m$ that

$$\|x - R_m^\varepsilon(x)\| \;\leq\; \varepsilon.$$

Note that the number of measurements required by the algorithm does not depend on the error tolerance $\varepsilon > 0$ or the norm in which the error is measured. The algorithm (which is described below) does not seem practical and the theorem is hence a purely theoretical result. Still, it demonstrates that adaptive, Lipschitz-continuous measurements are extremely more powerful than the other classes considered above: If linear measurements are used (non-adaptive or adaptive) or even if non-adaptive continuous measurements are used, then $m$ measurements are required to obtain the desired worst-case error.

We do not know whether Theorem 10.6 is optimal; in principle it is possible that a constant number of measurements suffices. We only know from the Borsuk-Ulam Theorem that one continuous measurement is not enough for $m > 1$. (There is no adaption for one measurement.)

We also stress that the *full* algorithm $R_m^\varepsilon$ above cannot be (Lipschitz-)continuous because, again, the Borsuk-Ulam theorem does not allow it. Still, from Theorem 10.6 we obtain that

$$\|R_m^\varepsilon(x) - R_m^\varepsilon(y)\| \leq \|x - y\| + 2\varepsilon,$$

which might be considered as some kind of *stability*.

We point out that it is important for Theorem 10.6 that the measurements are known exactly. If the measurements are *noisy*, only much weaker bounds can be obtained, as has been established by Krieg, Novak, Plaskota and Ullrich (2026*b*). This is already true for very small noise, which is another major difference to the other classes of measurements and algorithms, where small noise typically has small impact on the minimal worst case error. Note, however, that the upper bounds for noisy adaptive and continuous measurements, though much weaker than without noise, are sometimes still much better than everything that can be achieved with non-adaptive measurements (an example is the embedding $S\colon \ell_2^m \hookrightarrow \ell_\infty^m$, see the discussion in the aforementioned paper).

Based on Theorem 10.6, the idea for the proof of Theorem 10.4 is quite simple: Given a continuous measurement map $N\colon F \to \mathbb{R}^m$ and a continuous reconstruction map $\Phi\colon \mathbb{R}^m \to Y$, we can use $\mathcal{O}(\log m)$ adaptive continuous measurements to learn $N(f) \in \mathbb{R}^m$ up to a small error. By the continuity of $\Phi$, this suffices to learn $\Phi(N(f))$ up to a small error, and hence $\mathcal{O}(\log m)$ adaptive measurements suffice to obtain a similarly good approximation as with the $m$ non-adaptive measurements from the mapping $N$.

We now turn to a description of the algorithm $R_m^\varepsilon$ for the recovery of vectors from $\mathbb{R}^m$ up to error $\varepsilon > 0$. It is based on a suitable 'coloring' of $\mathbb{R}^m$ with the following properties: Consider a covering

$$\mathbb{R}^m = \bigcup_{i \in \mathbb{N}} D_i,$$

a (coloring) map $t\colon \mathbb{N} \to \{1, 2, \ldots, N\}$, $N \in \mathbb{N}$, and a constant $c > 0$ with:

(a) $\operatorname{diam}(D_i) \leq 1$ for each $i$;

(b) if $i \neq j$ and $t(i) = t(j)$, then $\operatorname{dist}(D_i, D_j) \geq c$,

where diameter (diam) and distance (dist) are with respect to the given norm on $\mathbb{R}^m$, and the constant can depend on $m$. The regions $D_k$, $k \in \mathbb{N}$, do not have to be pairwise disjoint.

Before we describe the algorithm and corresponding measurements, we show that such a partition exists with $N = m + 1$. Here, we only consider the maximum norm on $\mathbb{R}^m$. The general case is obtained by the equivalence of all norms on $\mathbb{R}^m$.

The idea is to split $\mathbb{R}^m$ into unit cubes and give all $r$-dimensional facets the color $r$. This does not quite work yet, since different facets of the same dimension $r$ touch. However, if we add a small neighborhood to the $(r-1)$-dimensional facets, then the $r$-dimensional facets do not touch any more, see the right hand side of Figure 10.1.

More precisely, let us denote by $F_r$ the $r$-dimensional facets, i.e., the set of vectors in $\mathbb{R}^m$ with at least $m - r$ integer components, and with $F_r^\delta$ for $\delta > 0$ the closed $\delta$-neighborhood of $F_r$ in the maximum norm. We choose $1/2 > \delta_0 > \cdots > \delta_m > 0$ and inductively define $C_0 := F_0^{\delta_0}$ and

$$C_r := F_k^{\delta_k} \setminus \bigcup_{i<k} C_i, \quad r = 1, \ldots, m.$$

The sets $C_r$ shall be colored with the color $r$. The desired covering is obtained as the family $(D_k)_{k \in \mathbb{N}}$ of the connected components of all the sets $C_r$, $0 \leq r \leq m$. It is easy to verify that this partition has the desired properties.

The above partition is certainly not the only admissible one, another partition for $m = 2$ is shown on the left hand side of Figure 10.1. However, note that it is not possible to find a partition with less than $m + 1$ colors. This is related to the Lebesgue covering dimension of $\mathbb{R}^m$, which is $m$, see Pears (1975, Chapter 2). Still, there might exist partitions with more colors that result in an algorithm that is *easier to implement*.

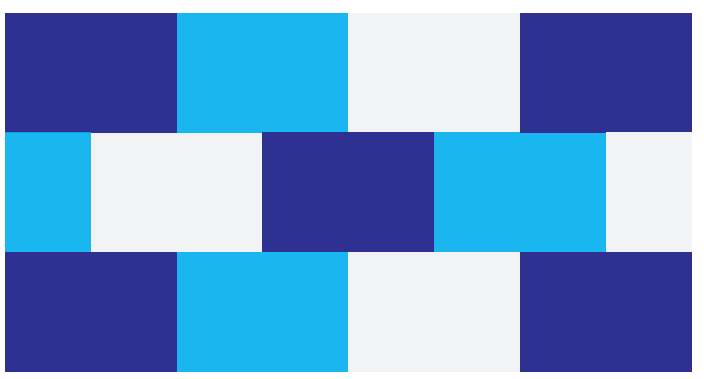
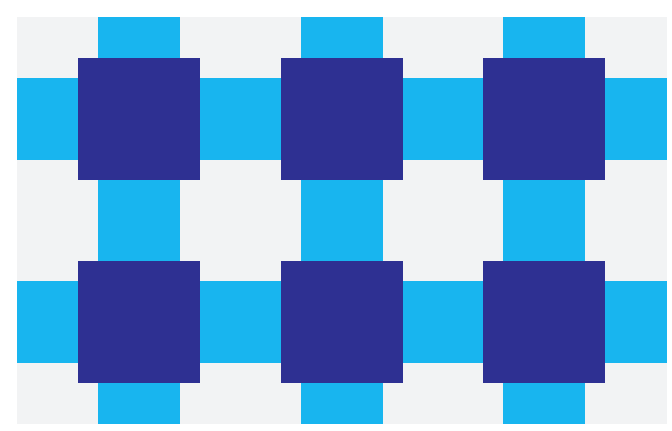

Figure 10.1. Two colorings of the plane with three colors. The second is generalized above to higher dimensions.

Now, to find an $\varepsilon$-approximation of $x \in \mathbb{R}^m$ it is clearly enough to find an index $i^*$ such that $x$ is in the closure of $\varepsilon D_{i^*}$. Define by

$$I_r := \{i \in \mathbb{N} \colon t(i) = r\} \quad \text{and} \quad E_r = \bigcup_{j \in I_r} \varepsilon D_j$$

the set of points with color $r$ in the covering $(\varepsilon D_k)_{k \in \mathbb{N}}$. A continuous measurement of the form

$$\lambda_J(x) = \operatorname{dist}\left(x, \bigcup_{j \in J} E_j\right)$$

with $J \subset \{1, \dots, m+1\}$, tells us whether or not $x$ is contained in the closure of any of the sets $E_j$, $j \in J$. We use $n = \lceil \log_2(m+1) \rceil = n(m) - 1$ such functionals and bisection to find a color $r^* = t(i^*)$ such that $x \in \overline{E}_{r^*}$. (We do not claim that $r$ and $i^*$ are unique.)

Now we can determine a correct index $i^*$ with $x \in \varepsilon\overline{D}_{i^*}$ using any continuous functional $\lambda^*$ for which the images $\lambda^*(\varepsilon\overline{D_i})$ for $i \in I_{r^*}$ are pairwise disjoint. An example is given by

$$\lambda^*(x) = \max_{i \in I_{r^*}} \left\{ \frac{c\varepsilon}{2i} - \operatorname{dist}(x, \varepsilon D_i) \right\}.$$

Note that, for each $x \in \mathbb{R}^m$, at most one of the terms in the maximum is non-negative. Since $x$ is contained in the closure of $\varepsilon D_{i^*}$ with $i^* = \frac{c\varepsilon}{2\lambda^*(x)}$, the output $R_m^\varepsilon(x)$ of the algorithm can be chosen as any element from $\varepsilon D_{i^*}$.

All the functionals above have Lipschitz-constant equal to one and, if we consider the 1- or $\infty$-norm, they are even piecewise linear.

**Remark 10.7 (The easiest case).** We first considered the "easiest case" of approximating a point from the Euclidean ball in $\mathbb{R}^3$ using only $n = 2$ measurements. We observed that adaptive continuous measurements can improve over all non-adaptive methods with worst case error 1: One can start with the functional $\ell_1(x) = x_1 - |x_3|$ and, depending on the outcome, continue with either $\ell_2(x) = x_3$ or $\ell_2(x) = x_2$. We leave it to the reader to think about the geometrical details. Observe that the used measurements are quite elementary and it would be interesting to extend this to higher dimensions.

# A. Collection of some bounds

For convenience of the reader (and ourselves), we now summarize several bounds discussed in this survey. We begin with a list of some general bounds on the sampling numbers, as discussed in Sections 3 and 5. We then present a list of the known results on various widths and the corresponding sampling numbers for our running examples, i.e., the univariate *Sobolev spaces*.

As before, we write $g_n(F,Y) \asymp n^s$ to mean that there exist positive constants $c$ and $C$, independent of $n$, although they may and generally do depend upon $F$, $Y$ and $s$, for which $cn^s \le g_n(F,Y) \le Cn^s$ for all $n$ sufficiently large, with the analogous meaning for the asymptotic inequalities $\lesssim$ and $\gtrsim$, and other benchmarks.

*A.1. General bounds on sampling numbers*

We assume here that $D$ is compact, $F \subset C(D)$ is compact and convex, and that $\mu$ is a finite Borel measure. Moreover, $1 \le p \le \infty$, $c \in \mathbb{N}$ is an absolute constant and $(a)_+ := \max\{0,a\}$. We only state the results that we know to be sharp for some $F$ and $Y$. We refer to the corresponding result and section for details.

- Theorem 3.17:
$$g_{4n}^{\mathrm{lin}}(F,L_p) \;\le\; 6 \cdot n^{(1/2-1/p)_+} \cdot d_n(F,B)$$

- Corollary 3.31:
$$g_{cn}^{\mathrm{lin}}(F,L_2) \;\le\; \frac{1}{\sqrt{n}} \sum_{k>n} \frac{d_k(F,L_2)}{\sqrt{k}}$$

- Theorem 3.35:
$$g_{cn}^{\mathrm{lin}}(F,Y) \;\lesssim\; \sum_{k>n} \frac{d(F,V_k)_2}{\sqrt{k}} \cdot \max\left\{ \frac{\Lambda_{4k}}{\sqrt{k}}, \frac{\Lambda_{4n}}{\sqrt{n}} \right\}$$
for any $Y$ that satisfies Assumption B, any sequence $(V_k)$ of nested subspaces of $L_2$ with $\dim(V_k) = k$, $d(F,V_k)_2 := \sup_{f\in F} \inf_{g\in V_k} \|f-g\|_2$, and $\Lambda_m := \sup\{\|f\|_Y \colon f \in V_m, \|f\|_2 = 1\}$ from (3.57).

- Theorem 5.8: For $s < 1/2$, we have
$$g_n(F,L_p) \;\lesssim\; n^{-s/p} \left( \sup_{k\le n} \, (k+1)^s \cdot \varepsilon_k(F,B) \right)^{1/p}$$

- Theorem 5.12:
$$g_n(F,B) \;\le\; (n+1) \cdot \varepsilon_n(F,B)$$

- Theorem 5.15: If $F$ is the unit ball of a Banach space, then, for all $s > 1$,
$$g_n(F,B) \;\lesssim\; n^{-s+1} \cdot \sup_{k\le n} (k+1)^s \cdot h_k(F,B)$$
and recall that $h_k \le \min\{2\cdot \varepsilon_k, \delta_k, d_k\}$.

*A.2. Widths of Sobolev spaces*

Let $W_q^s = \{f \in C^{s-1}([0,1]) \colon f^{(s-1)}$ abs. cont. and $f^{(s)} \in L_q\}$ and consider the approximation $(a_n)$, Bernstein $(b_n)$, Gelfand $(c_n)$, Kolmogorov $(d_n)$, manifold $(\delta_n)$, entropy $(\varepsilon_n)$ and sampling $(g_n)$ numbers, respectively. We write $s_n(X,Y) := s_n(\mathrm{id})$ for $\mathrm{id}\colon X \to Y$, see Section 9.2. In all cases below, $s_n(X, L_\infty) = s_n(X, B([0,1]))$.

**Theorem A.1.** For $s \in \mathbb{N} \setminus \{1\}$ and $1/q + 1/q' = 1$, we have:

$$a_n\big(W_q^s, L_p\big) \asymp \begin{cases} n^{-s}, & 1 \le p \le q \le \infty \\ n^{-s+1/q-1/p}, & 1 \le q \le p \le 2 \\ n^{-s+1/q-1/p}, & 2 \le q \le p \le \infty \\ n^{-s+1/q-1/2}, & 1 \le q \le 2 \le p \le \infty, q' \ge p \\ n^{-s+1/2-1/p}, & 1 \le q \le 2 \le p \le \infty, q' \le p \end{cases}$$

$$b_n\big(W_q^s, L_p\big) \asymp \begin{cases} n^{-s}, & 1 \le q \le p \le \infty \\ n^{-s}, & 1 \le p \le q \le 2 \\ n^{-s+1/q-1/p}, & 2 \le p \le q \le \infty \\ n^{-s+1/q-1/2}, & 1 \le p \le 2 \le q \le \infty \end{cases}$$

$$c_n\big(W_q^s, L_p\big) \asymp \begin{cases} n^{-s}, & 1 \le p \le q \le \infty \\ n^{-s}, & 1 \le q \le p \le 2 \\ n^{-s+1/2-1/p}, & 1 \le q \le 2 \le p \le \infty \\ n^{-s+1/q-1/p}, & 2 \le q \le p \le \infty \end{cases}$$

$$d_n\big(W_q^s, L_p\big) \asymp \begin{cases} n^{-s}, & 1 \le p \le q \le \infty \\ n^{-s}, & 2 \le q \le p \le \infty \\ n^{-s+1/q-1/2}, & 1 \le q \le 2 \le p \le \infty \\ n^{-s+1/q-1/p}, & 1 \le q \le p \le 2 \end{cases}$$

$$\delta_n\big(W_q^s, L_p\big) \asymp \varepsilon_n\big(W_q^s, L_p\big) \asymp n^{-s}$$

$$g_n\big(W_q^s, L_p\big) \asymp g_n^{\mathrm{lin}}\big(W_q^s, L_p\big) \asymp \begin{cases} n^{-s}, & 1 \le p \le q \le \infty \\ n^{-s+1/q-1/p}, & 1 \le q \le p \le \infty \end{cases}$$

We refer to Pinkus (1985, VII) for $a_n$, $c_n$, $d_n$, to Tikhomirov (1986, III.2.1) for $b_n$, to Stesin (1975) for $\delta_n$, to (Carl 1981*a*) for $\varepsilon_n$, and to Novak and Triebel (2006) for $g_n$ and $g_n^{\mathrm{lin}}$; see also (Temlyakov 2018*c*), the original references mentioned in these works, and the historical treatise of Pietsch (2007, 6.7.8.13–6.7.8.15). Results for fractional smoothness $s > 0$ can be found in Vybíral (2008) and the references therein. Note that some of the above formulas are not correct for $s = 1$, for example, it is known from Malykhin (2025) that $d_n(W_1^1, L_p) \asymp n^{-1/2} \log(n)$ for $2 < p < \infty$.

*Acknowledgments.* The authors would like to thank everyone who contributed to the content of this work, including many co-authors and friends, even if they only commented and discussed with us in the past. We believe that the progress made in this field of research in the 21st century so far is a great collective achievement. We would nevertheless like to express our thanks to Ben Adcock, Simone Brugiapaglia, Wolfgang Dahmen, Simon Foucart, Karlheinz Gröchenig, Thomas Kühn, Winfried Sickel, Mathias Sonnleitner and Jan Vybíral for their comments on a draft or discussions about the scope.

Special thanks, however, go to our former PhD advisor Erich Novak, who devoted an extraordinary amount of time to providing feedback and engaging in numerous discussions in connection with this project. More importantly, he laid the foundation for most of our results by training us both at Friedrich Schiller University Jena (starting with "Numerik 1" and "Funktionalanalysis") and beyond, and continues to inspire and motivate us to this day.

Finally, our heartfelt thanks go to Katha and Coco for their incredible patience with us during this project; we suspect that, in the end, you may have actually suffered no less than we did.